\documentclass[twoside]{article}
\usepackage[hmargin={4.3cm,4.3cm},vmargin={4.3cm}]{geometry}

\title{Theory on Structure and Coloring\\ of Maximal Planar Graphs (I):\\ Relationship between Structure and\\ Coloring}
\author{Jin Xu}
\date{}

\usepackage{amsmath,amsthm}
\usepackage{graphicx}
\usepackage{psfrag}
\usepackage{float}
\usepackage{amssymb}
\DeclareGraphicsExtensions{.pdf}

\chardef\bslash=`\\ 

 \newtheorem{theorem}{Conjecture}[section]
 \newtheorem{theorem2}{Theorem}[section]
 \newtheorem{corollary}[theorem2]{Corollary}
 
 \newtheorem{Prop}[theorem2]{Proposition}
 \newtheorem{problem}[theorem]{Problem}

\theoremstyle{definition}
\newtheorem{definition}{Definition}[section]

\theoremstyle{remark}

\newcommand{\eval}[2][\right]{\relax
  \ifx#1\right\relax \left.\fi#2#1\rvert}

\begin{document}
\maketitle \markboth{The Mathematical Proofs of Two Conjectures}
{The Mathematical Proofs of Two Conjectures}
\renewcommand{\sectionmark}[1]{}

 \begin{abstract}

Maximal planar graph refers to the planar graph with the most edges, which means no more edges can be added so that the resulting graph is still planar. The Four-Color Conjecture says that every planar graph without loops is 4-colorable. Indeed, in order to prove Four-Color Conjecture, it clearly suffices to show that all maximal planar graphs are 4-colorable. Since this conjecture was proposed in 1852, no mathematical proofs have been invented up until now. Maybe the main reasons lie in the following three aspects in terms of maximal planar graphs:  not clearing up the structures,  not figuring out the coloring types, and not straightening out the relation between structure and coloring. For this, we will write a series of articles to study the structure and coloring theory of maximal planar graphs systematically. This is our first article, which focuses mainly on the structure and coloring relations. First, we introduce a new way to construct maximal planar graphs. The advantage of this method is that it establishes an immediate relation with 4-colorings, and reveals how a given maximal planar graph is generated. Second, a special class of maximal planar graphs -- recursive maximal planar graphs is researched in depth, which lays a foundation for solving the uniquely 4-colorable planar graphs conjecture(see subsequent articles). Third, we discover an important mode for classifying 4-colorings: tree-coloring and cycle-coloring, which runs through the whole series of articles. Furthermore, this mode is applied to the research on an arbitrary 4-coloring and its corresponding structure of unions of two and three bicolored subgraphs. Finally, we introduce the concepts of black-white coloring and stamen phenomenon, and find out a necessary and sufficient condition for an even cycle to be a 2-colorable cycle.\\

 \end{abstract}

\newpage

 \begin{center}
 \textbf{Contents}
 \vspace{3mm}
 \end{center}

 1.\quad Introduction

 2.\quad Relational definitions and notations

 3.\quad Operational system to generate maximal planar graphs

 4.\quad Recursive maximal planar graphs

 5.\quad Coloring-structure of maximal planar graphs

 6.\quad Black-White coloring, and necessary and sufficient conditions for 2-colorable cycles














 Acknowledgements

 References

 Appendix.  All of the 4-colorings of maximal planar graphs with $\delta\geq 4$ and order from 6 to 11

\section{Introduction}
 The planar graph is a very important class of graphs no matter which aspect, theoretical or practical, is concerned.
In theory, there are many famous conjectures that have very significant effect on graph theory, even mathematics, such as the  Four-Color Conjecture, the Uniquely Four-Colorable planar graphs conjecture, and the Nine-Color Conjecture etc \cite{Jensen1995}. In application, planar graphs can directly be applied to the study of layout problems \cite{Josep2002}, information science \cite{Broder2000} etc.

 Particularly, maximal planar graph is one important class of planar graphs. A maximal planar graph is a simple planar graph where every face is a cycle of length 3, so it is also called triangulation. As the studying object of the well-known conjectures, i.e. the Four-Color and the Uniquely Four-Colorable planar graphs, can be confined to maximal planar graphs, many scholars have been strongly attracted to the study of this typical topic. They did research on maximal planar graphs from a number of different standpoints, such as degree sequence, construction, coloring,  traversability and generating operations, etc. The following will present some related results in brief.

 In a maximal planar graph $G$, a cycle with length $k(\geq 3)$ is called $k$-\emph{cycle} and is denoted by $C_k$. If both the interior as well as the exterior of $C_k$ contain one or more vertices, then $C_k$ is referred to as a \emph{separating} $k$-\emph{cycle}.

With respect to the connectivity of maximal planar graphs, there are obvious facts as follows.
Note that there is only one 2-connected maximal planar graph on three vertices (3-cycle), so we consider only the case $|V(G)|\geq 4$. Thus, if $\delta(G)=3$, then $G$ is 3-connected; if $\delta(G)=4$, then either $G$ is
4-connected or G has a separating 3-cycle; if $\delta(G)=5$, then either $G$ is
5-connected or $G$ has at least one of separating 3-cycles and separating 4-cycles.

In 1977-1978, Schmeichel and Hakimi \cite{Schmeichel1977,Schmeichel1978} researched degree sequences of planar graphs. They introduced the concepts of \emph{Euler sequence} and \emph{maximal Euler sequence}, and gave a necessary and sufficient condition for \emph{Euler 1-sequence} and \emph{2-sequence} to be planar graphical. They also gave a sufficient condition for a maximal planar graph that is $\delta$-connected, where $\delta$ denotes the minimum degree.

In 1982, Fanelli\cite{Fanelli1982} proved the existence of a class of $n(\geq 32)$-vertex maximal planar graphs with exactly fourteen 5-degree vertices and minimum degree $\delta$=5.

A graph is referred to as \emph{hamiltonian} if it has a cycle that contains all vertices exactly once. The problem of hamiltonian maximal planar graphs was first studied by Whitney \cite{Whitney1931} in 1931. He proved an outstanding result that each maximal planar graph without separating 3-cycles is hamiltonian.

After 1980, Nishizeki \cite{Nishizeki1980}, Dillencourt \cite{Dillencourt1991}, Harnat and Owens \cite{Harant1995}
studied hamiltonian maximal planar graphs under the toughness condition. The \emph{toughness} of a graph $G$, denoted by $t(G)$, is defined as the largest real number $t$ such that deletion of any $s$ points from $G$ results in a graph which either is connected or else has at most $s/t$ components. That is
\[
t(G)=\min\{\frac{|S|}{\omega(G-S)}: S\subset V(G),\omega(G-S)>1\}
\]
where $S$ ranges over all vertex-cutsets of $G$, and $\omega(G)$ denotes the number of components of $G$.

Since 1956, many scholars have extended Whitney's result to planar graphs. In 1977 and in 1983, Tutte \cite{Tutte1977} and Thomassen \cite{Thomassen1983} proved independently that each 4-connected planar graph is hamiltonian.
In 1990, Dillencourt \cite{Dillencourt1990} proved that Whitney's conclusion still holds if the chords satisfy a certain sparseness condition and that a Hamilton cycle through a graph satisfying
this condition can be found in linear time.
In 1997, Sanders \cite{Sanders1997} proved every 4-connected planar graph has a Hamilton cycle through any two of its edges. In 2010, G$\ddot{o}$ring and Harant \cite{Goring2010} proved that Sanders' result is the best
possibility by constructing 4-connected maximal planar graphs with three edges a large
distance apart such that any Hamilton cycle misses one of them. In 1994, Jung \cite{Jung1994} discussed infinite Hamilton paths in infinite maximal planar graphs. In 1979, Hakimi, Schmeichel and Thomassen \cite{Hakimi1979} studied the number of Hamilton cycles in a maximal planar graph, and constructed a class of $n(\geq 12$)-vertex maximal planar graphs containing exactly four Hamilton cycles. They also
proved that every 4-connected maximal planar graph on $n$ vertices
contains at least $n/(\log_2 n)$ Hamilton cycles.

In 1891, Eberhard \cite{Eberhard1891} started to consider the problem of constructing maximal planar graphs. He gave an operational system to generate all of the maximal planar graphs. We use $\langle K_4;\Phi=\{\varphi_1,\varphi_2,\varphi_3\}\rangle$ to present this system, where $K_4$, $\Phi$ and ${\varphi_1,\varphi_2,\varphi_3}$ are called
\emph{starting graph}, \emph{operation set} and \emph{generating operations}, respectively (see Figure1.1).

\begin{center}
\includegraphics [width=380pt]{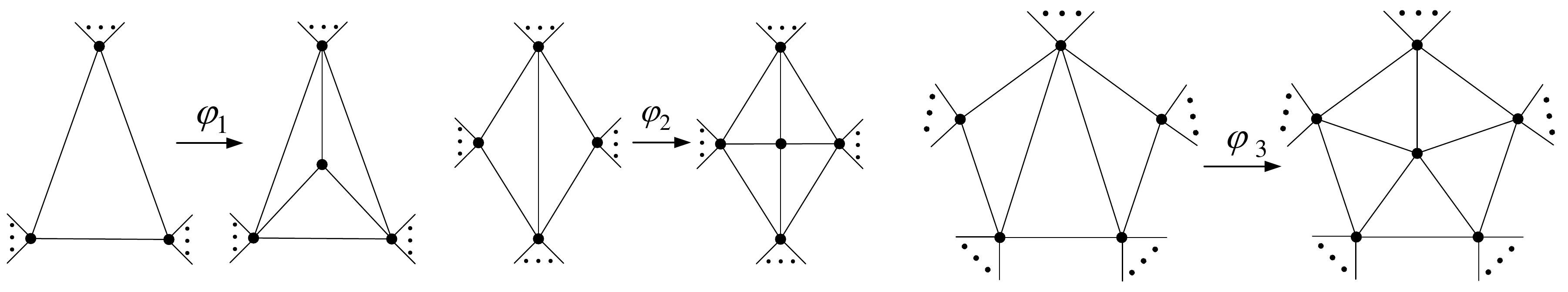}

 \textbf{Figure 1.1.} Eberhard's operations
\end{center}

For a maximal planar graph $G$ and a cycle $C$ of $G$, if the interior of $C$ contains no vertices and all the interior faces of $C$ are triangles, then we call $C$ a \emph{pure chord-cycle}. The interior edges of a pure chord-cycle are referred to as \emph{chords} of $C$. For the convenience of description, we also view triangles in maximal planar graphs as pure chord-cycles.

In fact, it can be easily seen that the implementation process of Eberhard's operations are the followings. First delete all the chords of a pure chord-cycle $C$ with length $k(=3,4,5)$; then add a new vertex inside $C$ and connect it to all
vertices of $C$ so that a wheel corresponding to $C$ is generated.

From 1999 to 2000, Wang \cite{Wang1999,Wang2000} independently proposed a method to construct maximal planar graphs. He indeed extended Eberhard's operations from pure chord-cycles with lengths $3,4,5$ to pure chord-cycles with lengths $\geq 3$.

After Eberhard's work, it was subsequently neglected for almost a century and only came to life in 1974 with the study of constructing all 5-connected maximal planar graphs by Barnette \cite{Barnette1974} and Butler \cite{Butler1974}, independently. Different from Eberhard's operational system,  Barnette and Butler's operational system is $\langle Z_{20};\Phi=\{\varphi_4,\varphi_5,\varphi_6\}\rangle$, where the starting graph $Z_{20}$ is the icosahedron. In addition, the operation set is also changed into $\{\varphi_4,\varphi_5,\varphi_6\}$ (see Figure 1.2), where the ellipses attached to the vertices in the
description of the generating rule denote any number (zero or more) of edges such that $\delta= 5$.
\begin{center}
\includegraphics [width=380pt]{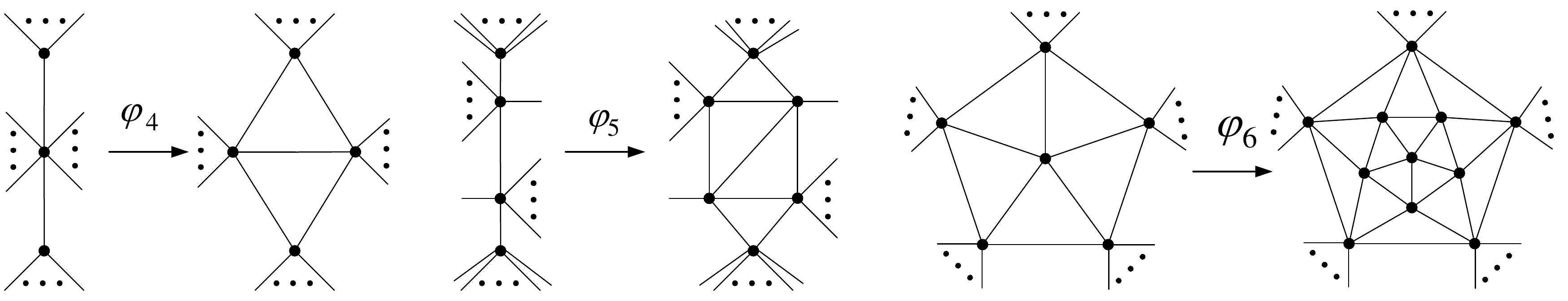}

\textbf{Figure 1.2.} Barnette and Butler's operations
\end{center}

In short, Barnette and Butler's method starts with the icosahedron graph and uses the operations
called $\varphi_4,\varphi_5$ and $\varphi_6$ given in Figure 1.2 to generate all of the 5-connected maximal planar graphs.

In 1983, Batagelj \cite{Batagelj1983} improved the method of Barnette and Butler by changing one of the generating rules. To be specific, he used a new generating operation $\varphi_7$ instead of $\varphi_6$ and kept the remainders unchanged. The new operational system is denoted by $\langle Z_{20};\Phi=\{\varphi_4,\varphi_5,\varphi_7\} \rangle$, where $\varphi_7$ is called \emph{flip} (see Figure 1.3).
\begin{center}
\includegraphics [width=160pt]{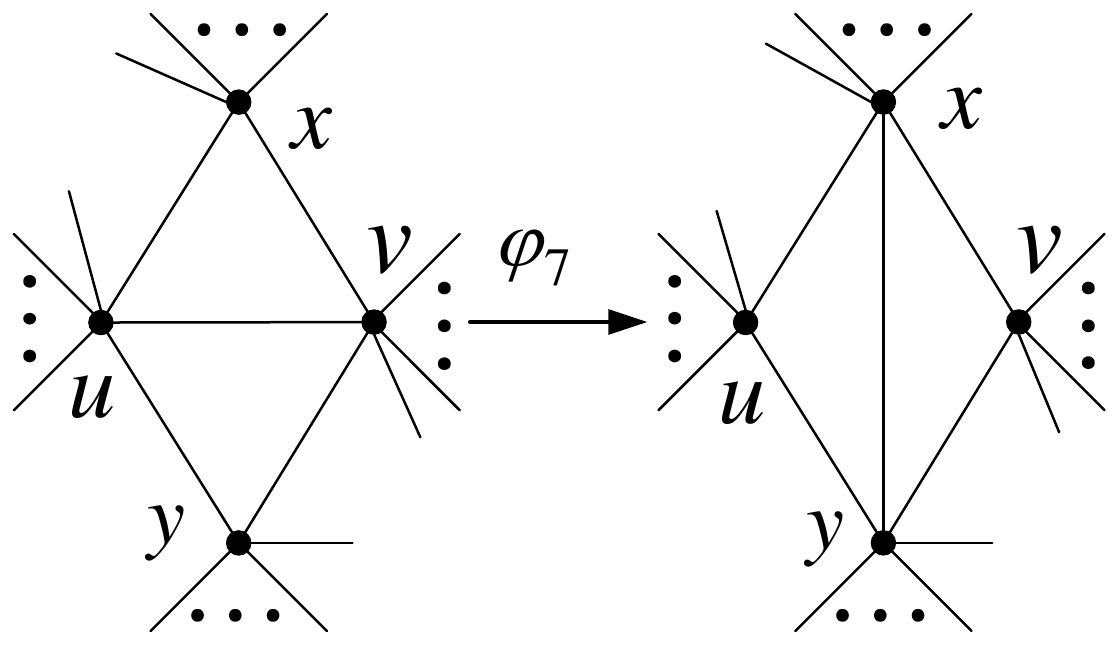}

\textbf{Figure 1.3.} Flip operation
\end{center}

In fact, the research of flip operation has been a long history. This concept was introduced first  by Wagner \cite{Wagner1936} in 1936. Up to now, the flip operation has been studied very thoroughly, so the following will give a specific discussion about it.

In 2005, further works were done by  Brinkmann and McKay \cite{Brinkmann2005} in terms of Barnette, Butlery and Batagelj's conclusions. He gave an efficient method to construct all simple maximal planar graphs of minimum degree 5. Moreover, he pointed out what condition should be satisfied for the above four generating operations, $\varphi_4,\varphi_5,\varphi_6,\varphi_7$, to construct the maximal planar graphs with minimum degree 5 that contain separating 3-cycles, 4-cycles and 5-cycles respectively. On the basis of this algorithm, he presented the results of a computer program. Particularly, he counted the number of maximal planar graphs of minimum degree 5 with orders from 12 to 40, where the numbers of 3-connected, 4-connected and 5-connected 40-vertex maximal planar graphs of minimum degree 5 are 8469193859271, 7488436558647 and 5925181102878, respectively. Note that he used the \emph{canonical construction path} method proposed by McKay \cite{McKay1998} in 1998, to avoid the generation of isomorphic copies in his computer program.

The study on algorithm to generate maximal planar graphs also stimulates many scholars' interest. In 1996, Avis \cite{Avis1996} gave a $O(r\cdot f(n,r))$-time algorithm for generating all $r$-rooted 3-connected maximal planar graphs on $n$ vertices by the reverse search technique. First, constructed a $n$-vertex canonical maximal planar graph (contains exactly two vertices of degree $n-1$); then, generated all of $r$-rooted 3-connected maximal planar graphs of order $n$ by means of the flip operation.

In 2004, Nakano \cite{Nakano2004} gave a simple algorithm to generate all of 3-connected $r$-rooted plane triangulations with at most $n$ vertices. Particularly, he showed that
all of 3-connected rooted plane triangulation having exactly $n$ vertices
including exactly $r$ vertices on the outer face can be generated in $O(r \cdot f(n, r))$ time
without duplications. Here a plane triangulation was a planar graph that each inner face has exactly three
edges on its contour; a rooted plane triangulation was a plane triangulation with
one designated vertex on the outer face; $f(n, r)$ is the number of nonisomorphic such triangulations.  The specific method of this algorithm is: first, to construct a genealogical tree $T$ so that each vertex of $T$ correspondeds to one special graph that possesses some properties; then, to generate the desired graphs based on $T$. Moreover, according to this algorithm, each one of maximal planar graphs on $n$ vertices can be generated in $O(n^3)$ time per graph.

In 2007, Brinkmann and McKay \cite{Brinkmann2007} introduced the \emph{Plantri}-operational rule depending on the canonical configuration path \cite{McKay1998}, and gave the program plantri \cite{Brinkmann}.

Let $G$ be a maximal planar graph, and $\triangle uvx$,$\triangle uvy$ be the two triangles in $G$ that have the common edge $e=uv$. An \emph{edge flip} consists of deleting the edge $e$ from $G$ and adding a new edge $e'=xy$ to the graph such that the remain is a maximal planar graph. And the edge $e$ is called \emph{flippable}(see Figure 1.3).

It is clear that edge flip converts a maximal planar graph into another one with the same number of edges. Naturally, this gives rise to a question as follow: Can any $n$-vertex maximal planar graph be transformed into any other $n$-vertex maximal planar graph through a finite sequence of flips? To our knowledge, Wagner \cite{Wagner1936} was the first to address this question directly with the positive answer. Although the number of $n$-vertex maximal planar graphs is exponential in $n$, Wagner avoided the issue of graph isomorphism by converting any given maximal planar graph into a canonical maximal planar graph, and proved any  $n$-vertex maximal planar graph could be transformed into a given $n$-vertex maximal planar graph by at most $2n^2$ edge flips. Here the canonical maximal planar graph on $n$ vertices, denoted by $\triangle_n$, is the unique one that contains exactly two vertices with degree $n-1$.

After that there are lots of scholars studying this topic, and improving the upper bound. In 1993, Negami and Nakamoto \cite{Negami1993} proved that any given $n$-vertex maximal planar graph could be converted into the canonical maximal planar graph via $n^2$ edge flips. Komuro \cite{Komuro1997} proved that any two $n$-vertex maximal planar graphs can be transformed into mutually through at most $8n-54$ (or $8n-48$) edge flips for $n\geq 13$ (or $n\geq 7$). Mori et al \cite{Mori2003} showed that any hamiltonian maximal planar graph on $n$ vertices could be transformed into $\triangle_n$ by at most $2n-10$ edge flips, preserving the existence of Hamilton cycle. He also proved that any $n$-vertex maximal planar graph could be made 4-connected by at most $n-4$ edge flips, and any two maximal planar graphs on $n$ vertices could be converted into each other through at most $6n-30$ edge flips.

In 2001, Gao et al \cite{Gao2001} proved that  every maximal planar graph on $n$ vertices contains at least $n-2$ flippable edges and that there exist some maximal planar graphs where at most $n-2$ edges are flippable. Moreover, he showed that there were at least $2n+3$ flippable edges in a maximal planar graph $G$ if $\delta(G)\geq 4$, and the bound was tight in certain cases. In addition, he generalized Wagner's theorem to the labeled graphs, and obtained a conclusion that any pair of labeled maximal planar graphs on $n$ vertices could be converted into mutually using $n\log n$ edge flips.

In 2011, Bose et al \cite{Bose2011} showed that a maximal planar graph on $n(\geq 6)$ vertices could be made to be 4-connected by at most $\lfloor(3n-6)/5\rfloor$ edge flips, and that any pair of maximal planar graphs on $n$ vertices could be transformed into each other by at most $5.2n-32.8$ edge flips.

For a 4-colorable maximal planar graph $G$, denote by $C(4)$ the color set and $C_4(G)$ the set of all 4-colorings of $G$. It is obvious that every $f\in C_4(G)$ partitions $V(G)$ into four independent sets $V_1,V_2,V_3,V_4$, where $V_i$ denotes the (possibly empty) set of vertices assigned color $i$. We refer to the subgraph of $G$ induced by $V_i\cup V_j$ $(1\leq i< j\leq 4)$, denoted by $G[i,j]$, as \emph{bicolored subgraph induced by $f$}. Clearly, any 4-coloring $f$ can yield six induced bicolored subgraphs in total. If there is a bicolored subgraph induced by $f$ containing cycles, then we call $f$ a \emph{cycle-coloring} of $G$, and $G$ a \emph{cycle-colorable graph}; otherwise, if there is no bicolored subgraph induced by $f$ containing cycles, we refer to $f$ as \emph{tree-coloring}, and $G$ as \emph{tree-colorable graph}. If $G$ contains no tree-colorings, then $G$ is called a \emph{pure cycle-coloring graph}; likewise, if $G$ contains no cycle-colorings, then $G$ is called a \emph{pure tree-coloring graph}.

According to the above classification, the elements of $C_4(G)$ can be divided into two classes: One is the kind of tree-colorings, the other is the kind of cycle-colorings. Correspondingly, all maximal planar graphs can be divided into three classes: pure tree-coloring graphs, pure cycle-coloring graphs and impure coloring graphs (contain both cycle-colorings and tree-colorings). The study of tree-coloring, cycle-coloring and the corresponding structure of a maximal planar graph will run through the whole article. Concerning the relevant study of this topic, Xu \cite{Xu2005a,Xu2005b} studied indeed the pure tree-coloring maximal planar graphs in 2005. He proclaimed that there were only two pure tree-coloring maximal planar graphs on $n$-vertices for $n\leq 43$ (see Figure 1.4). Based on this fact, he conjectured that there exist no more pure tree-coloring maximal planar graphs except for the two graphs shown in Figure 1.4. Unfortunately, we have found more pure tree-coloring maximal planar graphs with order at most $43$,  which shows the above conjecture is false. We will make a detailed research in the subsequent series of articles. Although no further researches have been done on this topic since 2005, the concept of tree-coloring proposed by Xu \cite{Xu2005a,Xu2005b} is still interesting.

\begin{center}
\includegraphics [width=170pt]{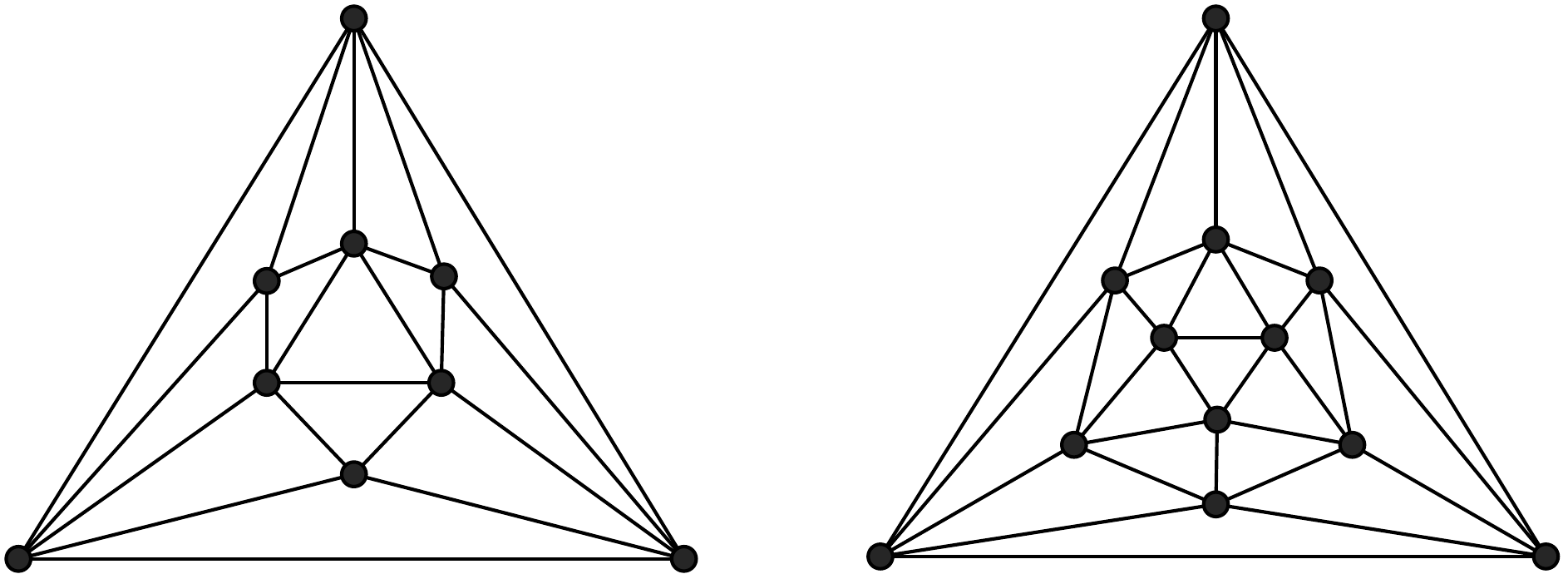}

\textbf{Figure 1.4.} Two pure tree-coloring graphs
\end{center}

 To date, much work has been done on the study of maximal planar graphs. The above is just relevant to our article.  However, there are still many topics about maximal planar graphs not being summarized, such as dominating sets, $k$-cyclable, $(m,n)$-linked, matching, triangles with restricted degrees, coloring algorithm, Chinese postman problem and enumeration, etc. Furthermore, interested readers can refer to the survey papers \cite{Zhu2013,Li2013}.

 Our research is written into five parts totally, and this is the first part that devotes mainly to the research of involving relation between structure and coloring  of maximal planar graphs. From the existing methods of generating maximal planar graphs, it is very hard to associate structure with colorings. In this paper, we introduce a new powerful tool, extending and contracting operations, to construct maximal planar graphs. Not only is this method straightforward, but also can it be related to 4-coloring, easily.

\section{Relational definitions and notations}
In this chapter we present some relational notations and terminologies of graph theory, which will be used throughout this article. Unless otherwise stated, the term \emph{graph} always presents a finite and undirected graph. For a given graph $G=(V(G), E(G))$,  $V(G)$ denotes the \emph{vertex set} of $G$  and $E(G)$ the \emph{edge set} of $G$. The cardinalities of these two sets are called the \emph{order} and \emph{size} of $G$, denoted by $n(G)$ and $m(G)$ respectively. A graph is \emph{finite} if both its vertex set and edge set are finite. The graph with no vertices is called \emph{null graph}. A graph with just one vertex is \emph{trivial}, and all other graphs are \emph{nontrivial}. We simply write $uv\in E(G)$ to denote that vertices $u,v\in V(G)$ are jointed by an edge, and say that $u, v$ are the \emph{ends} of this edge. The ends of an edge are said to be \emph{incident} with the edge, and vice versa. Two vertices that are incident with a common edge are \emph{adjacent}, as are two edges which are incident with a common vertex, and two distinct adjacent vertices are \emph{neighbors}. The set of neighbors of a vertex $v$ in  a graph $G$ is denoted by $\Gamma_G(v)$, or simply $\Gamma(v)$. The\emph{ degree} of a vertex $v$ in a graph $G$, denoted by
$d_G(v)$ (or simply $d(G)$), is the number of edges of $G$ incident with $v$.
An \emph{independent set} in a graph is a set of vertices such that no two of them are adjacent.
For a graph $H=(V(H),E(H))$, if
 $V(H) \subseteq V(G), E(H)\subseteq E(G)$, then $H$ is called a \emph{subgraph} of $G$. If $V(H)=V(G)$, then $H$ is a \emph{spanning subgraph} of $G$. And whenever $u, v \in V(H)$ are adjacent in the graph $G$,
 they are also adjacent in the graph $H$, then $H$ is called an \emph{induced subgraph} of
 $G$. An induced subgraph of $G$ under a vertex set $V'\subset V(G)$ is denoted by $G[V']$.

 A cycle with the vertices $v_1,v_2,\cdots, v_k$ and the edges $v_1v_2, v_2v_3,\cdots,v_kv_1$ is called a \emph{k-cycle}, denoted by $v_1v_2\cdots v_kv_1$. A path from $u$ to $v$ is called a \emph{P(u,v) path}. For two different vertices $u$, $v$ in a graph $G$, the \emph{distance} between $u$ and $v$ is the length of the shortest $P(u,v)$ path, denoted by $d_{G}(u,v)$.

 Two graphs $G$ and $H$ are disjoint if they have no vertex
 in common. By starting with a disjoint union of $G$ and $H$, and adding edges joining every vertex of $G$ to every vertex of $H$, one obtains the \emph{join} of $G$ and $H$, denoted by $G$+$H$. The join $C_{n}$+$K_{1}$ of a cycle and a single vertex is referred to as a \emph{wheel} with $n$ spokes, denoted by
 $W_{n}$ (the examples $W_{2},W_{3},W_{4},W_{5}$ are shown in Figure 2.1), where $C_n$, $K_1$ are called the \emph{cycle} and \emph{center} of the wheel respectively.

 A graph is \emph{k}-\emph{regular} if all of its vertices
have the same degree \emph{k}. A 3-regular graph is usually called a \emph{cubic graph}.

        \begin{center}
        \includegraphics [width=300pt]{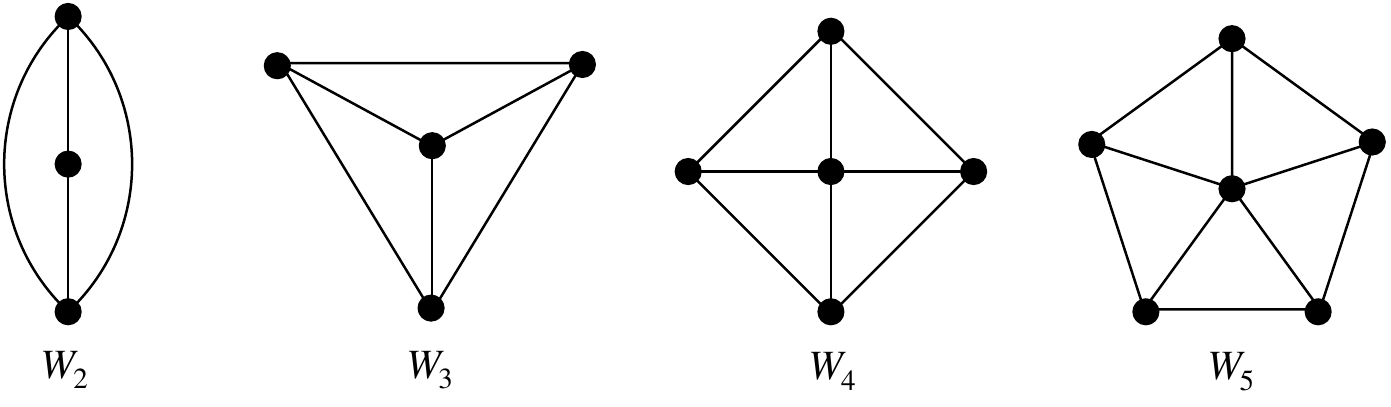}

         \textbf{Figure 2.1.} Four wheels $W_{2},W_{3},W_{4},W_{5}$.
        \end{center}

In order to \emph{identify }the nonadjacent vertices $u$ and $v$ of a graph $G$, it is necessary to replace these two vertices by a single vertex $w$, and make $w$ to be incident with all the edges which are incident to  $u$ and $v$ in $G$. We denote the resulting graph by $G \circ \{u,v\}$. To \emph{contract}
an edge $e=uv$ of a graph $G$, we delete the edge and then identify its ends.
 The resulting graph is denoted by $G \circ e$.

\subsection{Graph coloring}

    A \emph{vertex-coloring} of a graph $G$ is an assignment from color-set to its vertex-set such
    that no two adjacent vertices have the same color. A \emph{$k$-vertex-coloring}, or simply a $k$-coloring, of a
 graph \emph{G} is a mapping $f$ from  $V(G)$ to the color sets $C(k)=\{1,2,\ldots,k\}$ such
 that $f(x)\neq f(y)$ if $xy\in E(G)$.

    A graph $G$ is \emph{$k$-colorable} if it has a $k$-coloring. The minimum $k$ for which a graph $G$ is $k$-colorable
 is called its \emph{chromatic number}, denoted by $\chi(G)$. If $\chi(G)=k$, then the graph $G$ can be colored with $k$ colors, but not with $k-1$ colors, and call $G$ \emph{$k$-chromatic graph}.
Alternatively,
 each $k$-coloring $f$ of  $G$ can be viewed as a partition  $\{V_{1},V_{2},\cdots,V_{k}\}$ of $V$,
 where $V_{i}$ denotes the set of vertices assigned color $i$. So it can be written as $f=(V_{1},V_{2},\cdots,
 V_{k})$. In other words, the $k$-coloring $f$ partitions
$$V(G)=\bigcup_{i=1}^{k}V_{i},V_{i}\neq\emptyset,V_{i}\cap V_{j}= \emptyset, 1\leq i <j\leq k,
 i,j=1,2,\ldots,k \eqno{(2.1)}$$
where $V_{i}$ is an independent set
 of $G$, $i=1,2,\ldots,k$. The set of all $k$-colorings of a graph
 $G$ is denoted
 by $C_{k}(G)$. For a $k$-chromatic graph $G$, the notation $C^{0}_{k}(G)$
    denotes the set consisting of the partitions of all $k$-colorings of
    $G$, and is the \emph{partition set of $k$-color class} of $G$.
 And define
 $$\sigma_k^0(G)=|C_k^0(G)| \eqno{(2.2)}$$

 Suppose that $G$ is a $k$-chromatic graph with $k\geq 3$. Let $f\in C_k^0(G)$, and $U=\{v_1,v_2,\cdots, v_t\}(t\leq|V(G)|)$ be a subset of vertices of $G$. Now, define $f(U)=\{f(v_i)|i=1,2,\cdots, t\}$, obviously, $f(U)\in C(k)$. Particularly, when $u\in V(G)$, $f(\Gamma(u))$ denotes the set of all colors assigned to the neighbors of  $u$.

 A $k$-colorable graph $G$ is called \emph{uniquely $k$-colorable} if each $k$-coloring of $G$ induces the same
 partition of the vertex set $V$, shown in Formula (2.1).

    Similarly, an \emph{edge-coloring} of a graph $G$ \cite{Fiorini1977,Fiorini1978} is an assignment
    from color-set to its edge-set such that no two
    adjacent edges have the same color.
 A \emph{$k$-edge-coloring} of a graph $G$ is an edge-coloring with $k$ colors. A graph
 $G$ is \emph{$k$-edge-colorable} if it has a $k$-edge-coloring. The
 \emph{edge chromatic number}, $\chi^{\prime}(G)$, of $G$ is the minimum
 number $k$ for which $G$ is $k$-edge-colorable.
 A graph $G$ is \emph{uniquely $k$-edge-colorable}  if
 there is a unique $k$-edge-coloring such that any other colorings are equivalent to it. Alternatively, a graph $G$ is
 uniquely $k$-edge-colorable if there is exactly one partition of the edge-set $E(G)$ into $k$ disjoint matchings.

 In general, tow graphs  $G$ and $H$ are \emph{isomorphic}, written $G \cong H$, if there are two bijections $\theta: V(G)\rightarrow V(H)$ and $\phi:E(G)\rightarrow E(H)$ such that $\psi_G(e)=uv$ if and only if $\psi_H(\phi(e))=\theta(u)\theta(v)$; such a pair of mappings is called an \emph{isomorphism} between $G$ and $H$. A  graph $G$ is labeled if each vertex is
 assigned by a different label, traditionally represented by integer. For a labeled graph $G$, two colorings are different if
 there is at least one vertex receiving different colors. We use $f(G,t)$  to denote the number of $t$-colorings for a labeled
 graph $G$. It is called the \emph{chromatic polynomial} of a graph $G$, which is introduced first by Brikhoff to attack the
 Four-Color Problem in 1912 \cite{Birkhoff1912}. More detailed researches can be found in \cite{Beineke1978}$\sim$ \cite{Xuj1995} .

\subsection{Maximal planar graph}
    A \emph{maximal planar graph} is a planar graph to which no new edges can be added without violating planarity.
 A \emph{triangulation} is a planar graph in which every face is bounded by three edges (including its infinite face).
 It can be easily proved that a maximal planar graph is equivalent to a triangulation. Thus, we can say that each maximal
 planar graph is a triangulation. A maximal planar graph $G$ is \emph{divisible} if there exists a separating 3-cycle in $G$.

    There exists a kind of uniquely 4-colorable planar graphs, \emph{recursive maximal planar graphs}, each of which can be obtained from $K_{4}$ by embedding
 a 3-degree vertex in some triangular face continuously. In this paper, $\Lambda$ denotes the set consisting of all recursive maximal planar graphs and  $\Lambda_{n}$ the set of graphs in $\Lambda$ with order $n$.
 Let $\gamma_{n}=|\Lambda_{n}|$. Obviously, $\gamma_{4}=\gamma_{5}=\gamma_{6}=1$. The corresponding recursive maximal planar graphs are shown in Figure 2.2.
    \begin{center}
        \includegraphics [width=300pt]{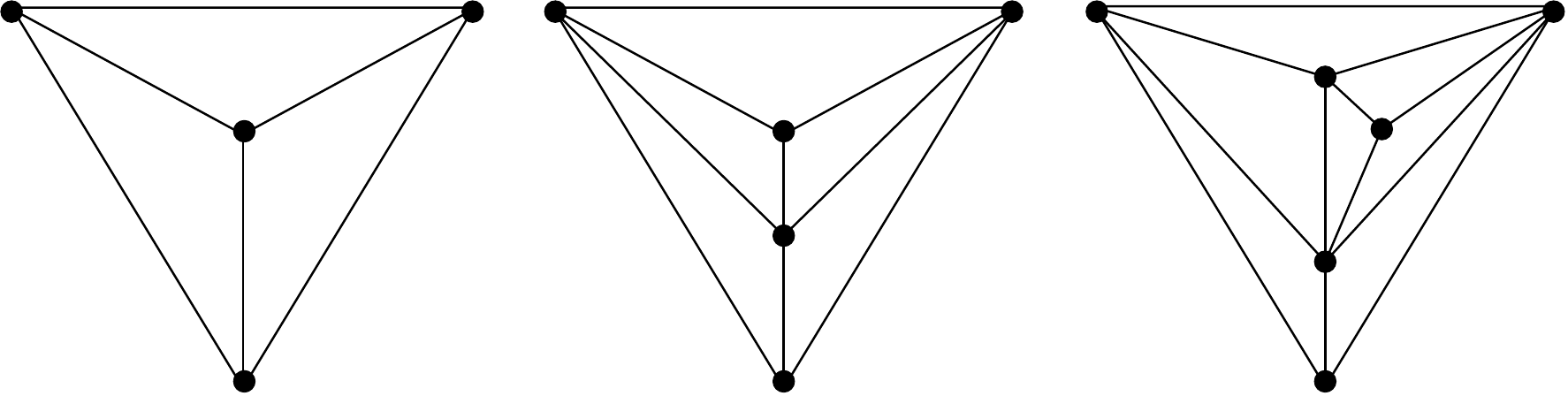}

        \textbf{Figure 2.2.} Three recursive maximal planar graphs
    \end{center}

 \noindent

In a maximal planar graph $G$, we mark a triangular face  $a-b-c$ if the vertices on its boundary are marked $a,b,c$ respectively.
A \emph{ vertex-embedding} on a triangular face $a-b-c$ is to
add a new vertex $u$ and make it to be adjacent to the vertices $a,b,c$ in this face, denoted $G+u$.
Obviously, the resulting graph is also a maximal planar graph.
We refer to a vertex-embedding as an \emph{extending 3-wheel operation}.
Another operation used in the paper is the \emph{vertex-deleting}, which is the inverse operation of a vertex-embedding. We also call a vertex-deleting a \emph{contracting 3-wheel operation}. The operations of vertex-embedding and vertex-deleting are illustrated in Figure 2.3.

The definitions and notations not mentioned can be found in \cite{Bondy2008}.

        \begin{center}
        \includegraphics [width=330pt]{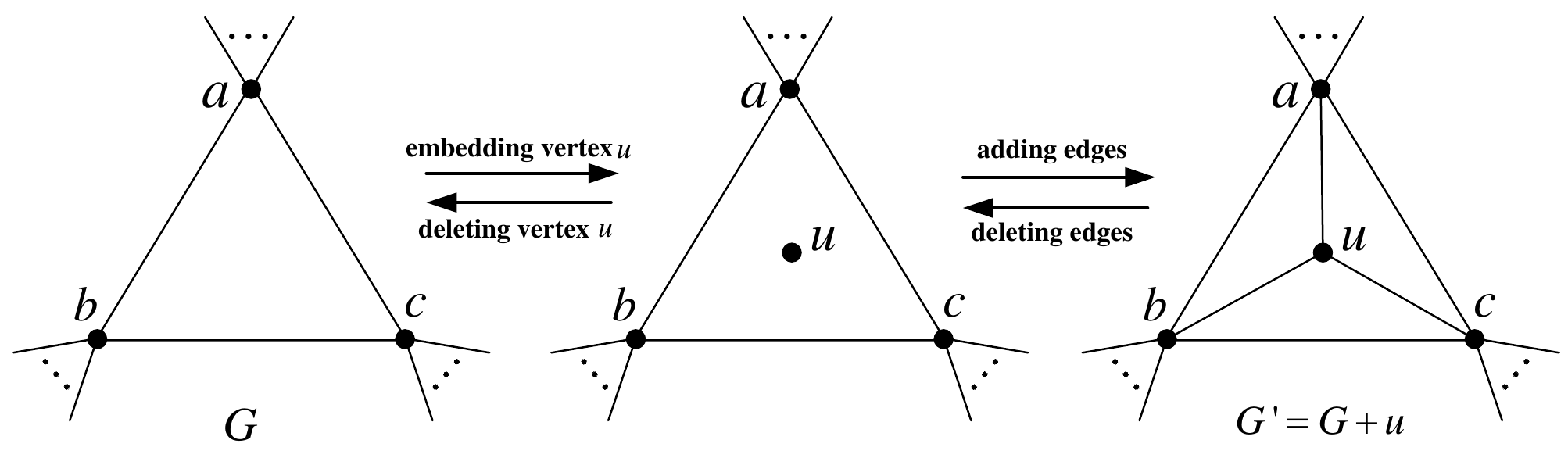}

       \textbf{Figure 2.3. } Two operations of embedding and deleting 3-degree vertex
        \end{center}


\section{Operational system to generate maximal planar graphs}

Investigating graph structure is extremely important as we study the coloring problems of maximal planar graphs. Since 1891, there have been many researchers being occupied by this topic. Correspondingly, many kinds of brilliant ideas have been created to construct maximal planar graphs, which have been introduced at length in Chapter 1. So far, lots of excellent results have been obtained with respect to the construction of maximal planar graphs, however it seems that little connection can be established with the 4-coloring. So, in this chapter, we describe a new method, the extending and contracting operation, to generate maximal planar graphs. The advantage of this method is that it can build up a direct relation with 4-coloring. In addition, with regard to the study of characteristics of maximal planar graphs, there are also many conclusions being obtained in the process of studying the computer-assisted method to attack the famous Four-Color Conjecture \cite{Appel1977a}$\sim$\cite{Ore1967}.

The operational system proposed here consists of two parts: operating objects and basic operator. Here the operating objects are maximal planar graphs; the basic operator includes four pairs of operations:
   the extending 2-wheel operation and its inverse operation, the contracting 2-wheel operation;
    the extending 3-wheel operation and its inverse operation, the contracting 3-wheel operation;
    the extending 4-wheel operation and its inverse operation, the contracting 4-wheel operation;
    the extending 5-wheel operation and its inverse operation, the contracting 5-wheel operation.
The function of this system is that starting with $K_3$ by a finite number of applications of the above four operations, it enables us to generate any given maximal planar graph. Based on this method, we construct all $(6\sim 12)$-vertex maximal planar graphs with $\delta\geq 4$. Moreover, we also discuss the extending and contracting operations under the condition of 4-coloring.


\subsection{Chromatic isomorphism}
 For the sake of convenience, we first introduce the concept of \emph{chromatic isomorphism} as follows.

 It is necessary to give an example before showing this definition clearly.
Considering the graph shown in Figure 3.1(x), it is  a 3-chromatic graph with 6 vertices, and the subgraph induced by $V'=\{v_{1},v_{2},v_{3},v_{4}\}$ is uniquely 3-colorable. We use $\{v_{1}\},\{v_{2},v_{4}\},\{v_{3}\}$ to present the unique vertex partition (see Figure 3.1(y)), which receive colors 1,2,3 respectively, where $\{1, 2, 3\}$ is the color set (see Figure 3.1(z)), that

    $$\{v_{1}\}\rightarrow 1, \{v_{2},v_{4}\}\rightarrow 2, \{v_{3}\}\rightarrow 3 \eqno{(3.1)}$$\label{3.1}

Thus, up to the chromatic isomorphism, Figure 3.1(x) has only three 3-colorings (see Figures 3.1(a), (b) and (c)), denoted by $f_{a}, f_{b}, f_{c}$, respectively.

    \[f_{a}= (\begin{array}{cccccc}
                v_{1} & v_{2} & v_{3} & v_{4} & v_{5} & v_{6} \\
                1 & 2 & 3 & 2 & 1 & 2
              \end{array}
    )\]
    \[f_{b}= (\begin{array}{cccccc}
                v_{1} & v_{2} & v_{3} & v_{4} & v_{5} & v_{6} \\
                1 & 2 & 3 & 2 & 3 & 1
              \end{array}
    )\]
    \[f_{c}= (\begin{array}{cccccc}
                v_{1} & v_{2} & v_{3} & v_{4} & v_{5} & v_{6} \\
                1 & 2 & 3 & 2 & 3 & 2
              \end{array}
    )\]

If we make (color) permutations on the Form (3.1), then  $f_{a},f_{b},f_{c}$ can correspondingly induce six 3-colorings (see Figure 3.1), respectively. In fact, Figure 3.1(x) has only three different 3-colorings, shown in Figures 3.1(a), (b) and (c). All other 3-colorings are just the color permutations of them. They indeed have the same vertex partitions with $f(a), f(b), f(c)$, respectively, in the following.

    \begin{center}

        \includegraphics [width=240pt]{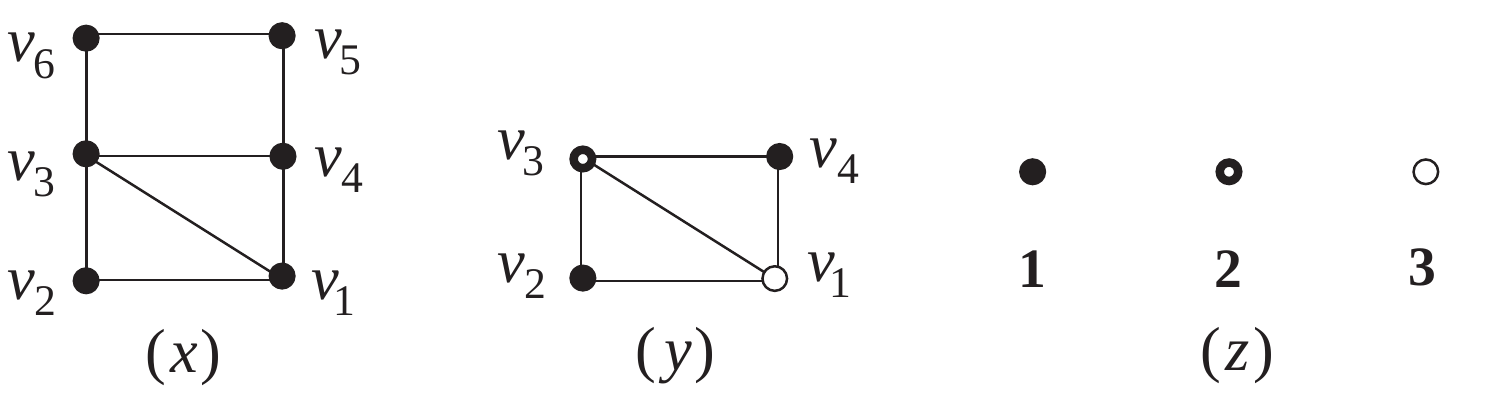}

        \vspace{5mm}
        \includegraphics [width=360pt]{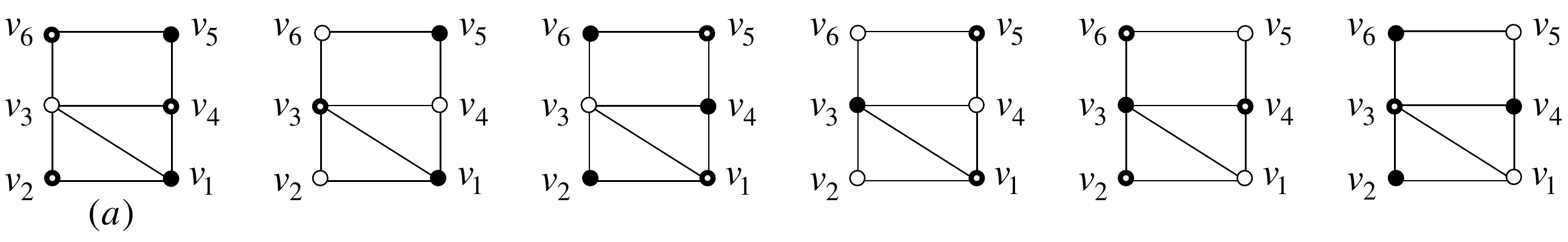}

        \includegraphics [width=360pt]{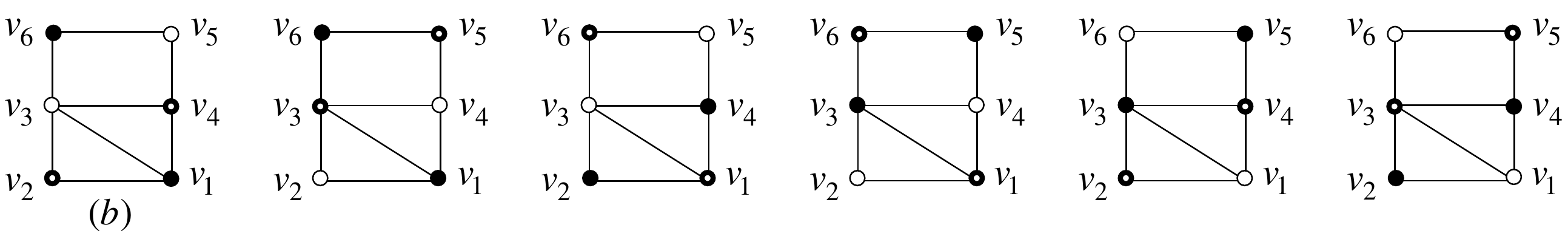}

        \includegraphics [width=360pt]{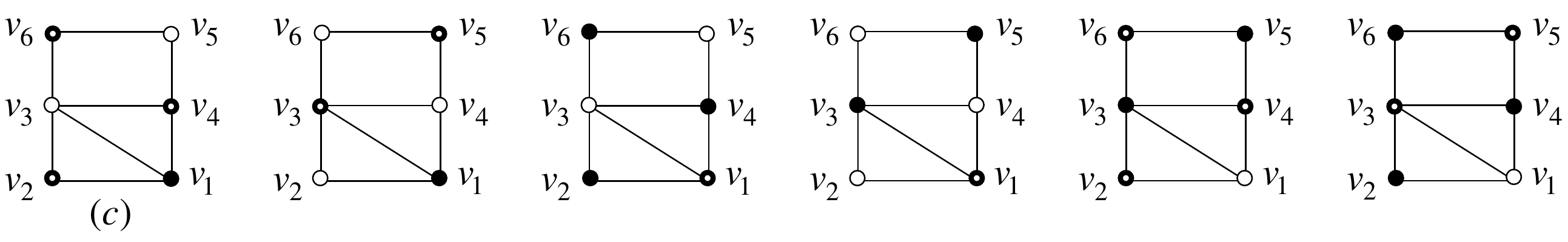}

        \vspace{5mm}
        \textbf{Figure 3.1.} The graph $G$ with 6 vertices and its 18 colorings

    \end{center}

    Recall that any $k$-coloring of a graph $G$ can partition $V(G)$ into $k$ independent subsets, which are often called the \emph{$k$-color class partition} induced by the $k$-coloring. At the same time, each such independent subset is said to be a \emph{color class} of $f$.

    \begin{definition}
        For a $k$-chromatic graph $G$ and its two colorings $f_{1},f_{2}$, we say that
    $f_{1},f_{2}$ are chromatic isomorphic if the $k$-color class partitions induced by them are identical. Obviously, chromatic isomorphic possesses the property of transitivity. A coloring $f$ chromatic isomorphic set of $G$ is the set of all colorings of $G$ that are chromatic isomorphic to $f$. If we choose one representative from each chromatic isomorphic set of $G$, then the set of all such representatives are called the chromatic isomorphic class of $G$, denoted by $C_k^0(G)$. In fact, $C_k^0(G)$ is just the set of all $k$-color class partitions of $G$, often called the $k$-color class partition set of $G$.
    \end{definition}

        If two colorings of a graph are chromatic isomorphic,
     then we can always transform one to be another one by adjusting the colors  properly. So we view such two colorings as one, and only choose one as their representative. For example, the six 3-colorings of the graph shown in Figure 3.1(x) in the same row have the identical 3-color class partition, so they are chromatic isomorphic mutually. Thus, we will say that Figure 3.1(x) has just three different 3-colorings, which are denoted by the first graph of each row, respectively.

        Clearly, we have the following theorem.

    \begin{theorem2}\label{theorem3.1}

        Let $G$ be a $k$-chromatic graph. Then for
    the chromatic isomorphic class $C^{0}_{k}(G)$ of $G$, and the
    set $C_{k}(G)$ of all colorings of $G$, we have
                             $$|C_{k}(G)|=k!|C^{0}_{k}(G)|\eqno{(3.2)}$$
    \end{theorem2}

        Based on Theorem \ref{theorem3.1}, when we analyze the properties of $k$-colorings of $G$,
   it suffices to consider the chromatic isomorphic group $C^{0}_{k}(G)$.

\subsection{Basic operational system}

    This section will be devoted to define the basic operators of the operational system for generating maximal planar graphs and some related properties under the condition of no coloring.

	The extending $2$-wheel operation means a procedure of adding a new edge between two adjacent vertices first, which will generate $2$-parallel edges (namely 2-cycle), and then, adding a new vertex in the face of the $2$-parallel edges and making the new vertex to be adjacent to the two vertices of the 2-cycle. Thus, the object of extending $2$-wheel operation is an edge of a maximal planar graph, see Figure 3.2(a).
	For a graph with 2-wheels, the contracting $2$-wheel operation means a procedure of deleting the center of a 2-wheel and the two edges incident with the center first,  and then erasing one of the parallel edges of the $2$-wheel.

	In section 2, we have introduced the extending $3$-wheel operation of maximal planar graphs as follow. First,  add a new vertex in a certain face of the maximal planar graph; second, add three edges of linking the new vertex and three vertices of the face, respectively. Thus, the object of extending $3$-wheel operation is a triangle of a maximal planar graph, see Figure 3.2(b). Correspondingly, we have also introduced the contracting $3$-wheel operation: deleting a certain $3$-degree vertex and the edges incident with it.

	For a maximal planar graph of $\delta\geq 4$, the so-called contracting $4$-wheel operation is defined in the following. First,  delete a certain $4$-degree vertex and the edges incident with it, and then identify a pair of the nonadjacent vertices in its neighbors. The extending $4$-wheel operation is the inverse operation of the contracting $4$-wheel operation. So the object of extending $4$-wheel operation is a $2$-length path of a maximal planar graph, see Figure 3.2(c). The following will show this definition in detail.

\begin{center}
        \includegraphics [width=300pt]{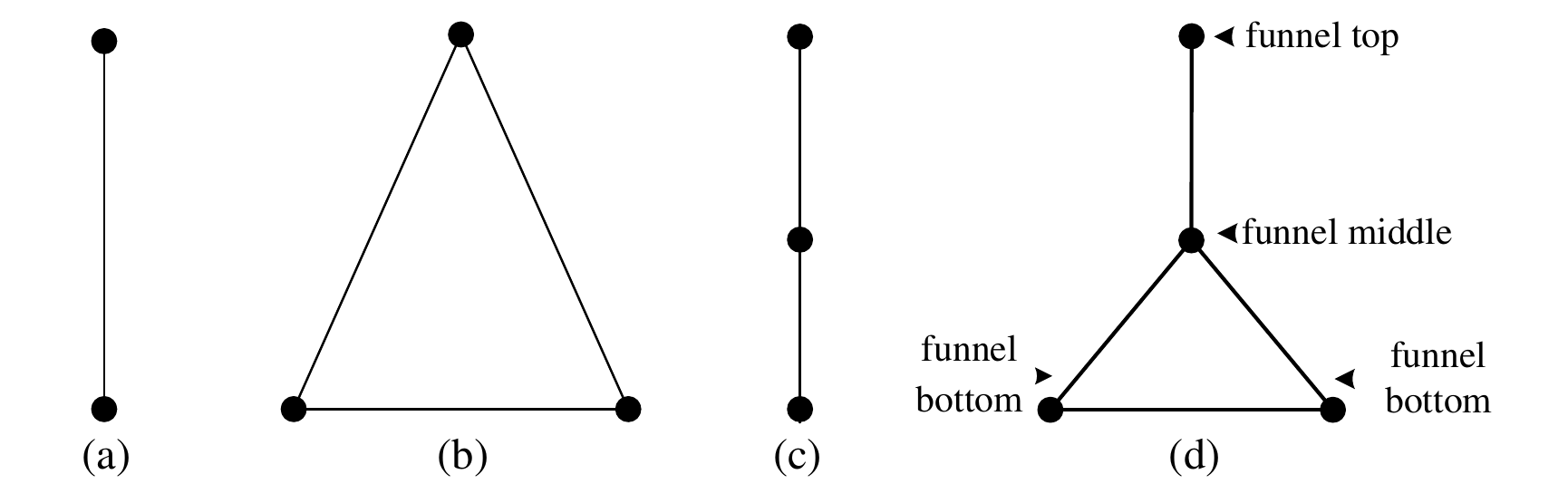}

        \textbf{Figure 3.2.} The objects of basic extending wheel operations
\end{center}

Let $G$ be a maximal planar graph, $x,u,y\in V(G)$, and $xuy$ be a path of length 2. The so-called extending $4$-wheel operation on the path $xuy$ is to replace the $xuy$ by a 4-cycle $xuu'yx$; that is, spilt the vertex $u$ into two vertices $u$ and $u'$, and split the edges $xu$  into two edges $xu,xu'$ and $uy$ into $uy,uy'$ respectively. This process is shown in Figure 3.3. Then add a new vertex $v$ in the face of the $4$-cycle $xu'yux$, and make $v$ adjacent to vertices $x,u',y,u$ respectively. The resulting graph is referred to as a \emph{generated graph} by implementing an extending $4$-wheel operation, denoted by $G*xuy$ (see the fifth graph in Figure 3.3).

The graph shown in Figure 3.2(d)  is called a \emph{funnel}, where the 1-degree vertex is the \emph{top} of the funnel, the 3-degree vertex is the \emph{middle} of the funnel and the two 2-degree vertices are the \emph{bottoms} of the funnel. As the middle and two bottoms of $L$ are vertices of a triangle, we also write $L$ by $L=v-\triangle$, where $v$ is the top of $L$.

\begin{center}
        \includegraphics [width=350pt]{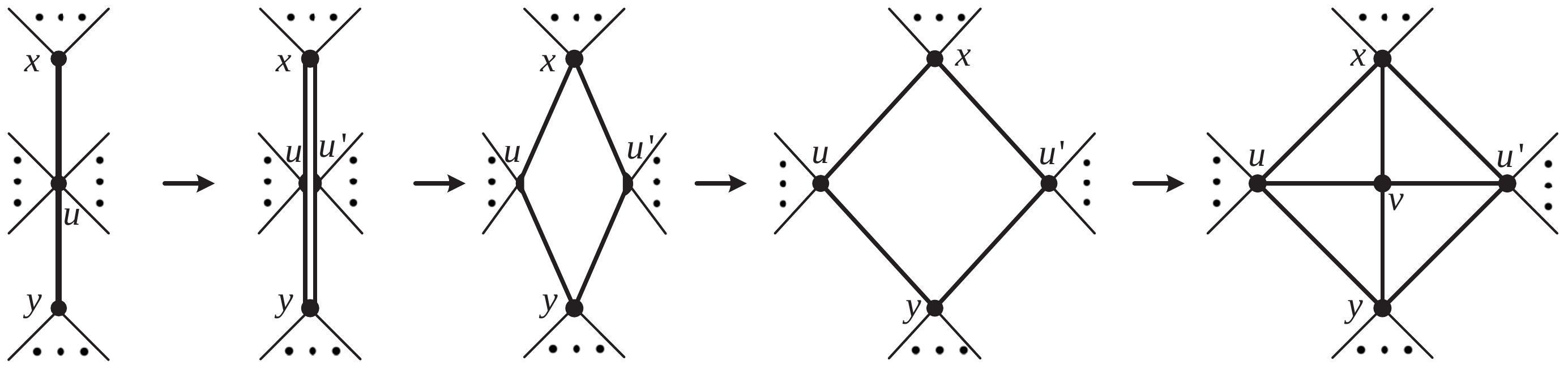}

        \textbf{Figure 3.3.} An illustration of the extending 4-wheel operation
\end{center}

	For a maximal planar graph, the contracting 5-wheel operation and the extending 5-wheel operation are similar to the contracting 4-wheel operation and the extending 4-wheel operation. The difference between them is that an extending 5-wheel operation is on a funnel, while an extending 4-wheel operation is  on a $2$-path. Here we only give a graphical illustrative definition shown in Figure 3.4, of  the contracting 5-wheel operation and the extending 5-wheel operation.
\begin{center}
        \includegraphics [width=350pt]{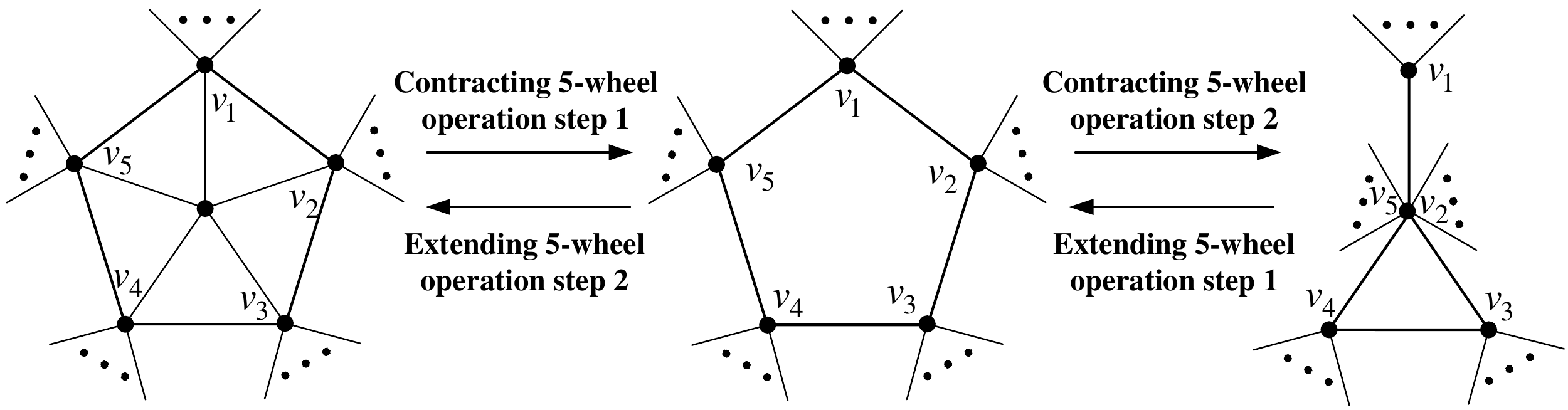}

        \textbf{Figure 3.4.} A graphical illustrative definition of an extending 5-wheel and contracting 5-wheel operations
\end{center}

	For a maximal planar graph $G$, we use $\zeta^{-}_{i}(G)$ and $\zeta^{+}_{i}(G)$ to denote the resulting graphs after implementing contracting $i$-wheel operation and extending $i$-wheel operation for $i=2,3,4,5$.

The following two propositions are easy to prove.

\begin{Prop}
	 $\zeta^{-}_{i}(G)$ and $\zeta^{+}_{i}(G)$ $(i=2,3,4,5)$ are  maximal planar graphs.
\end{Prop}

\begin{Prop}
	Let $G$ be a maximal planar graph of order $n$, then  $\zeta^{-}_{2}(G)$ and $\zeta^{-}_{3}(G)$  are  maximal planar graphs with order $n-1$;  $\zeta^{-}_{4}(G)$ and $\zeta^{-}_{5}(G)$  are maximal planar graphs with order $n-2$, that is
$$|\zeta^{-}_{2}(G)|=|\zeta^{-}_{3}(G)|=|V(G)|-1=n-1;$$
$$|\zeta^{-}_{4}(G)|=|\zeta^{-}_{5}(G)|=|V(G)|-2=n-2.$$
\end{Prop}

\begin{theorem2}\label{theorem3.4}
Suppose that $G$ is a maximal planar graph of order $n$. Then $G$ can be contracted to $K_3$ by implementing a series of contracting $i$-wheel operations for $i=2,3,4,5$.
\end{theorem2}

	\begin{proof} When $n=4$, there is only one maximal planar graph $K_4$, so the conclusion is true. Suppose that the conclusion holds for  $n\leq p~(p\geq 4)$, which means that for any maximal planar graph with order at most $p$, it can be contracted to $K_3$  by implementing contracting 2-wheel, 3-wheel, 4-wheel and 5-wheel operations, properly.

	Now we consider the case $n=p+1$. For any maximal planar graph $G$ of order $p+1$, if $G$ has a 2-degree or 3-degree vertex, then we will get a maximal planar graph with order $p$,
 $\zeta^{-}_{2}(G)$ or $\zeta^{-}_{3}(G)$, by deleting the 2-degree or 3-degree vertex and the incident edges together.
 According to the induction hypothesis, the conclusion holds. If $\delta(G)=4$ or $\delta(G)=5$, then properly implementing a contracting 4-wheel operation or a contracting 5-wheel operation for some 4-degree or 5-degree vertex, we will get a graph $\zeta^{-}_{4}(G)$ or $\zeta^{-}_{5}(G)$, which is a maximal planar graph of order $p-1$. On the basis of the induction hypothesis, they can be contracted to $K_3$ by a series of contracting $i$-wheel operations for $i=2,3,4,5$. 	
\end{proof}

	Through Theorem \ref{theorem3.4}, we clearly know that every maximal planar graph of order $n$ can be contracted to $K_3$ by implementing four basic contracting operations, properly. Of course, if we trace back to the reverses of contracting $k$-wheel operations of graph $G$, then starting with $K_3$ and doing the corresponding extending $k$-wheel operations, we can also get the original graph $G$. So,
\begin{corollary}
Any two maximal planar graphs can be transformed into each other by implementing the four pairs of contracting and extending operations.
\end{corollary}

	Here,  $\Psi=\{\zeta^{-}_{2},\zeta^{+}_{2},\zeta^{-}_{3},\zeta^{+}_{3},\zeta^{-}_{4},\zeta^{+}_{4},\zeta^{-}_{5},\zeta^{+}_{5}\}$  denotes the four basic pairs of contracting operations and extending operations, and  $S(G)=(K_3,\Psi)$  denotes the operational system of generating maximal planar graphs. So, starting with $K_3$, any maximal planar graph can be generated based on this system.

\subsection{Compound operational system}

The foregoing discussion shows that every maximal planar graph can be generated from $K_3$ through a finite sequence of extending $k$-wheel operations. In fact, the studying object of the operational system introduced in this Chapter is  the kind of maximal planar graphs with $\delta\geq 4$, so we focus mainly on such operations, including extending wheel and contracting wheel operations, by the conduction of which the resulting graph is still one with $\delta\geq 4$. Suppose $G$ is a maximal planar graph of $\delta\geq 4$, then there are three possibilities in terms of the minimum degree of $\zeta^-(G)$ when we implement contracting 4-wheel (or 5-wheel) operation to $G$: $\delta(\zeta^-(G))=2,\delta(\zeta^-(G))=3,\delta(\zeta^-(G))=4$.
When $\delta(\zeta^-(G))=2$ (or 3), obviously $\zeta^-(G)$ is not the graph we desired. So, we hope to obtain the desired graphs by conducting contracting 2-wheel and 3-wheel operations continually. That is: starting with $G$, implementing contracting 4-wheel (or 5-wheel) operation and contracting 2-wheel or 3-wheel operations repeatedly, we can always gain the graph with $\delta\geq 4$, which is also denoted by $\zeta^-(G)$.
The following proves this expectation can be achieved. We refer to the contracting 4-wheel (or 5-wheel) operation and a sequence of contracting 2-wheel and 3-wheel operations above as the \emph{compound contracting wheel operation} (or \emph{continually contracting wheel operation}, sometimes). At the same time, if we use $X$ to denote the set of vertices that are contracted in the process, $G[X\cup \Gamma(X)]$ is to be a configuration, written as $G^X$ and called a \emph{contractible subgraph} of $G$. Now we will mainly concern the structure of contractible subgraphs. Naturally, the related properties should be studied in depth. For this, we first present a straightforward result as follows.

\begin{theorem2}\label{theorem3.6}
Suppose that $G$ is a maximal planar graph of order $n$. If $G^X$ is a contractible subgraph of $G$, then the degree of vertices in $X$ is either 4 or 5, and the number of 5-degree vertices of $X$ is at most two.
\end{theorem2}

Note that there is an obvious fact as follows. For a given maximal planar graph $G$ with minimum degree $\delta\geq 4$, it is easy to find a 2-length path (or a funnel $L$) such that the resulting graph (denoted by $\zeta^+(G)$) obtained by doing an extending 4-wheel operation on $P_3$, or a pair of extending 3-wheel and 4-wheel operations, or a pair of extending 2-wheel and 4-wheel operations (or an extending 5-wheel operation on $L$, or a pair of extending 2-wheel and 5-wheel operations, or a pair of extending 3-wheel and 5-wheel operations), is still a maximal planar graph with $\delta\geq 4$. Similarly, for the resulting graph obtained by implementing some contracting wheel operations, we have

\begin{theorem2}\label{theorem3.7}
Suppose that $G$ is an $n$ $(\geq 9)$-vertex maximal planar graph with $\delta\geq4$. Then based on $G$ we can obtain an $(n-2)$-vertex or $(n-3)$-vertex maximal planar graph also with $\delta\geq 4$ through at most two contracting wheel operations.
\end{theorem2}

The above theorem shows indeed that for a maximal planar graph of minimum degree at least four, there is always a compound contracting wheel operation with $|X|\leq 2$ such that the resulting graph is also a maximal planar graph with $\delta\geq 4$, where $X$ is the set of vertices contracted in the process.

However, in the process of implementing compound contracting wheel operations, especially associated with coloring (which will be researched in the next section), an unavoidable case must be considered that the contracted vertex-set $X$ contains more than three vertices. For example, for the configuration of the join graph of a path $P_m$ ($m\geq 3$) and two isolated vertices $u,v$, it is easy to prove that the resulting graph, after implementing a sequence of contracting wheel operations, is still of this type of configuration. We refer to such configuration as the \emph{string \emph{4}-wheel graph}, which will be studied further in the later sections.

\begin{theorem2}\label{theorem3.8}
Suppose that $G$ is a maximal planar graph of $\delta\geq4$, and $G^X$ is a contractible subgraph of $G$. Then $G^X$ belongs to one of the fourteen configurations shown in Figure 3.5.
\end{theorem2}

\begin{proof}
Based on Theorem \ref{theorem3.6}, we prove this result by considering the size of $X$.

When $|X|=|\{v\}|$=1, $G^X$ is either $W_4$ or $W_5$ (see Figures 3.5(a),(b)). Tables 3.1(a) and (b) exhibit the process of these two contracting wheel operations in detail.

When $|X|=|\{u,v\}|$=2, the possible degrees of $u,v$ are in the following: 4,4; 4,5; 5,5 (see Figures 3.5(c),(d),(e)). Tables 3.1(c), (d) and (e) exhibit the process of these three contracting wheel operations in detail.

When $|X|=|\{u,v,w\}|$=3, the degree of $u,v,w$ are possible as: 4,4,4; 4,4,5; 4,5,4; 4,5,5; 5,4,5 (see Figures 3.5(f),(g),(h),(i),(j)). Tables 3.1(f),(g),(h),(i) and (j) exhibit the process of these five contracting wheel operations in detail.

When $|X|=|\{v,x_1,x_2,\cdots,x_t\}|\geq 4$, the degrees of $u,v,w$ have the following possibilities:

$\textcircled{1}$ No 5-degree vertex is in $X$. The corresponding subgraph is the string 4-wheel graph (see Figure 3.5(k)), and the process of the contracting wheel operation is shown in Table 3.1(k).

$\textcircled{2}$ One 5-degree vertex is in $X$. The corresponding subgraph and the process of the contracting wheel operations are shown in Figure 3.5(l) and Table 3.1(l), respectively.

$\textcircled{3}$ Two 5-degree vertices are in $X$. According to the adjacent relation of these two 5-degree vertices, there are two cases to be considered (see Figures 3.5(m),(n)). Tables 3.1(m) and (n) give the process of the two contracting wheel operations in detail.
\end{proof}

\begin{center}
        \includegraphics [width=350pt]{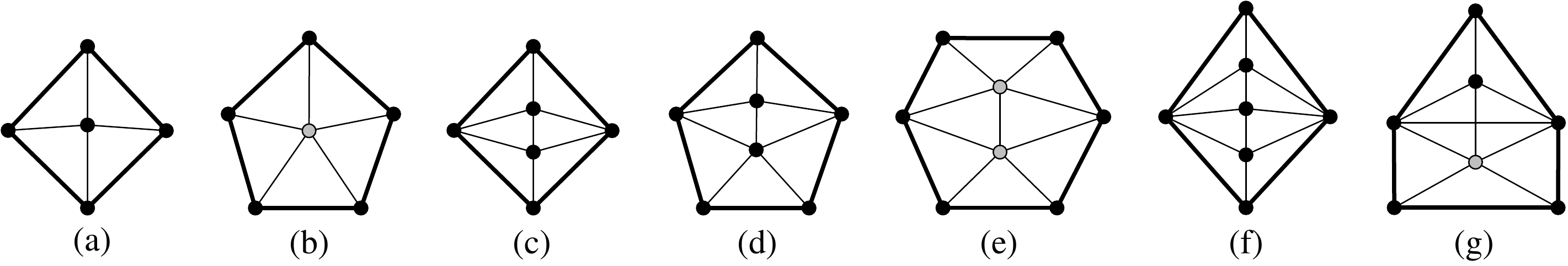}

        \includegraphics [width=350pt]{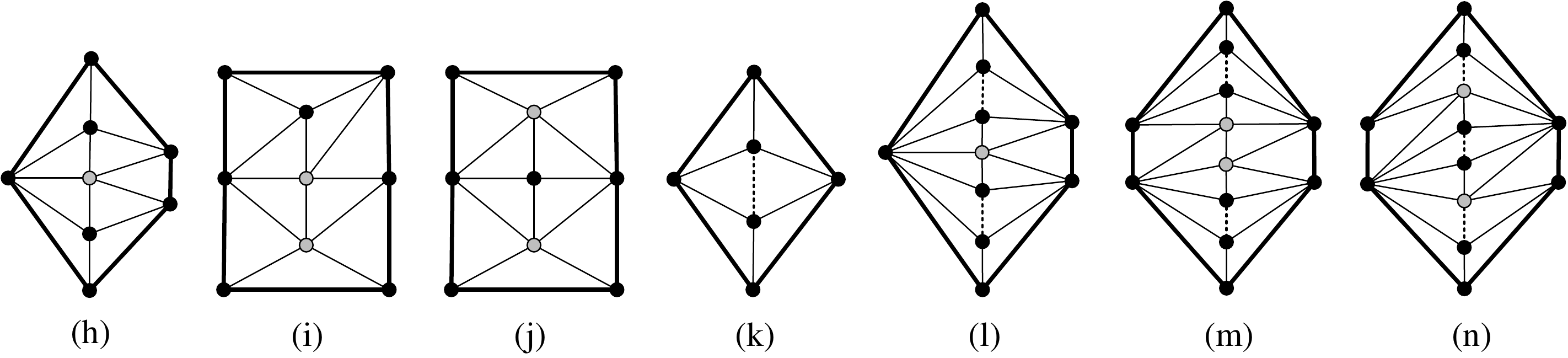}

        \textbf{Figure 3.5.}  An unavoidable-complete set of contractible subgraphs
\end{center}

\noindent\begin{tabular}{|l|l|l|}
  \hline
  ~~~~~~~~~~~~~~~\textbf{Table 3.1.} Structure of contractible subgraphs\\
  \hline
  (a) contracting 4-wheel operation, and the corresponding contracted\\
   subgraph \\
  \hline
  \includegraphics [width=360pt]{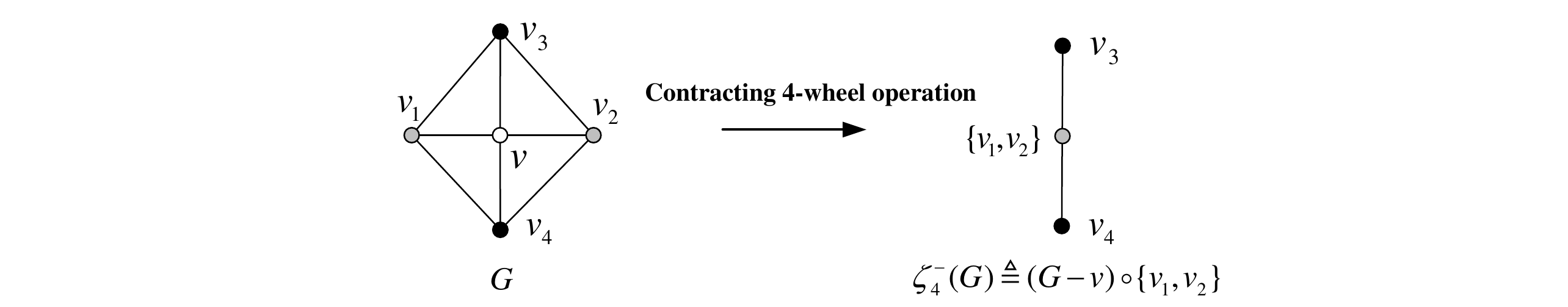} \\
  \hline
\end{tabular}
\begin{tabular}{|l|l|}
  \hline
   (b) contracting 5-wheel operation, and the corresponding contracted\\
    subgraph\\
  \hline
  \includegraphics [width=360pt]{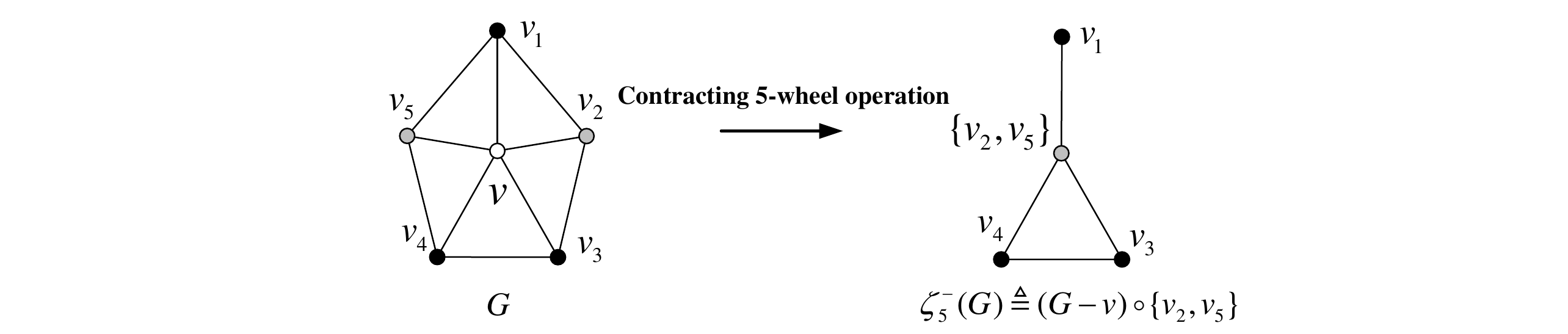} \\
  \hline
\end{tabular}
\begin{tabular}{|l|l|}
  \hline
   (c) contracting 4-wheel and 2-wheel operations, and the corresponding\\
    contracted subgraph\\
  \hline
  \includegraphics [width=360pt]{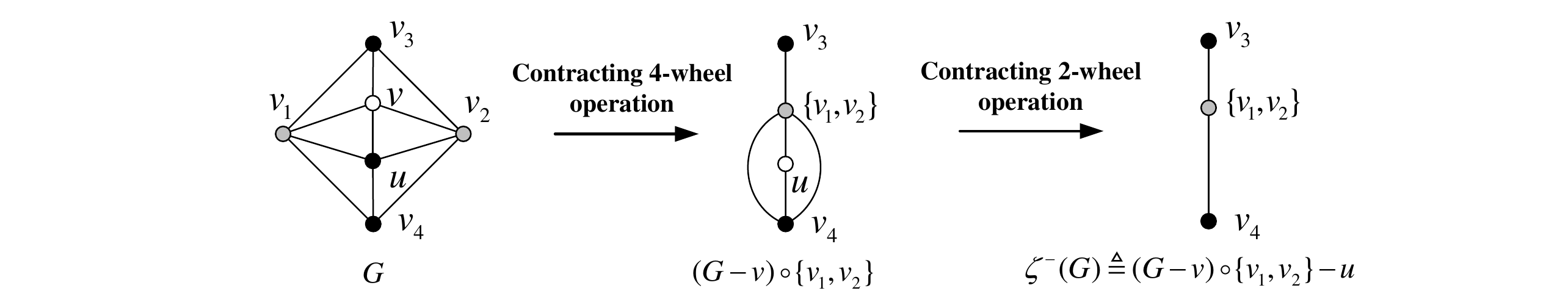} \\
  \hline
\end{tabular}
\begin{tabular}{|l|l|}
  \hline
   (d) contracting 4-wheel and 3-wheel (= 5-wheel and 2-wheel) operations, \\
 and the corresponding contracted subgraph\\
  \hline
  \includegraphics [width=360pt]{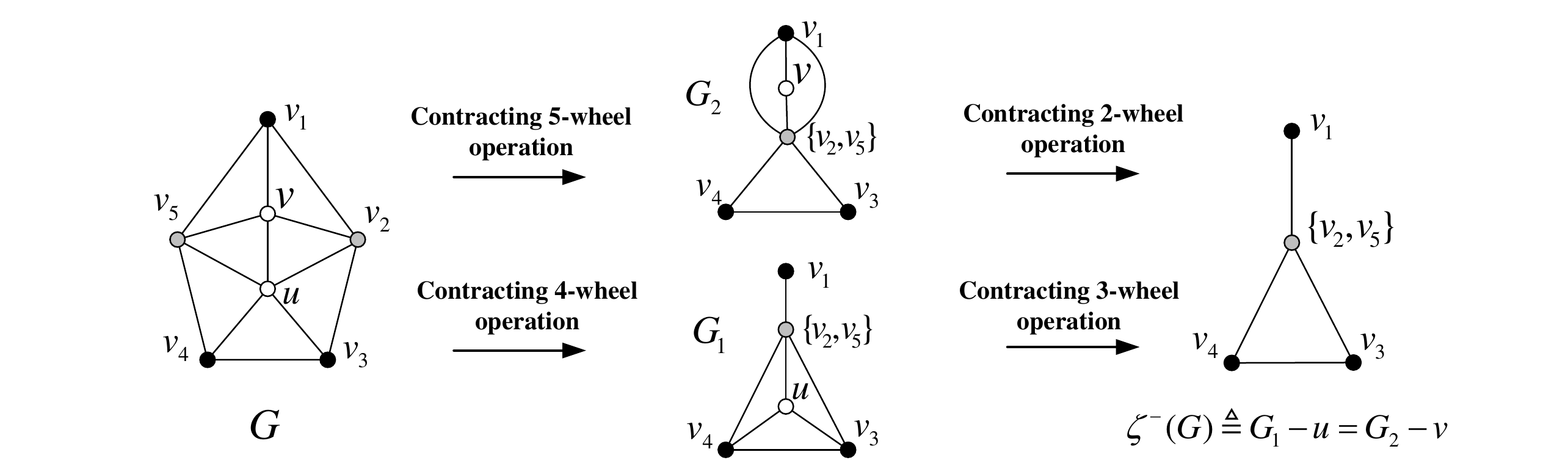} \\
  \hline
  \end{tabular}
  \begin{tabular}{|l|l|}
  \hline
  (e) contracting 5-wheel and 3-wheel operations, and the corresponding\\
   contracted subgraph\\
  \hline
  \includegraphics [width=360pt]{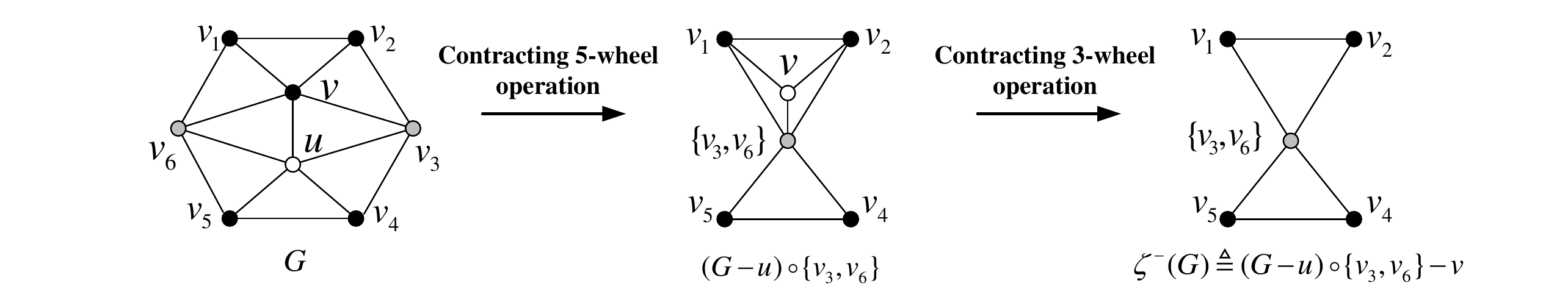} \\
  \hline
  \end{tabular} 
  \begin{tabular}{|l|l|}
  \hline
  (f) contracting 4-wheel, 2-wheel and 2-wheel operations, and the \\
  corresponding contracted subgraph\\
  \hline
  \includegraphics [width=360pt]{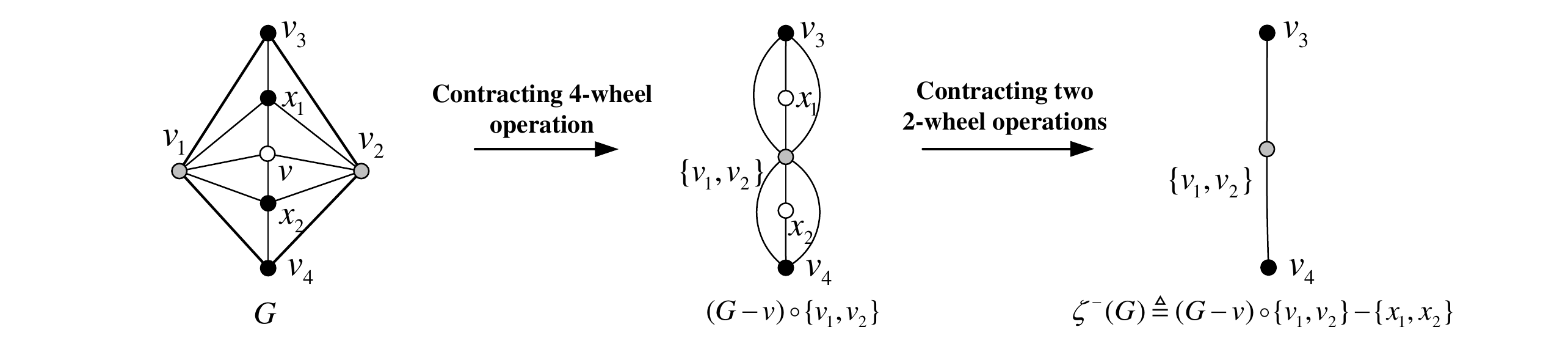} \\
  \hline
  \end{tabular} 
  \begin{tabular}{|l|l|}
  \hline
  (g) contracting 4-wheel, 2-wheel and 3-wheel (=5-wheel, 2-wheel and \\
  2-wheel) operations, and the corresponding contracted subgraph\\
  \hline
  \includegraphics [width=360pt]{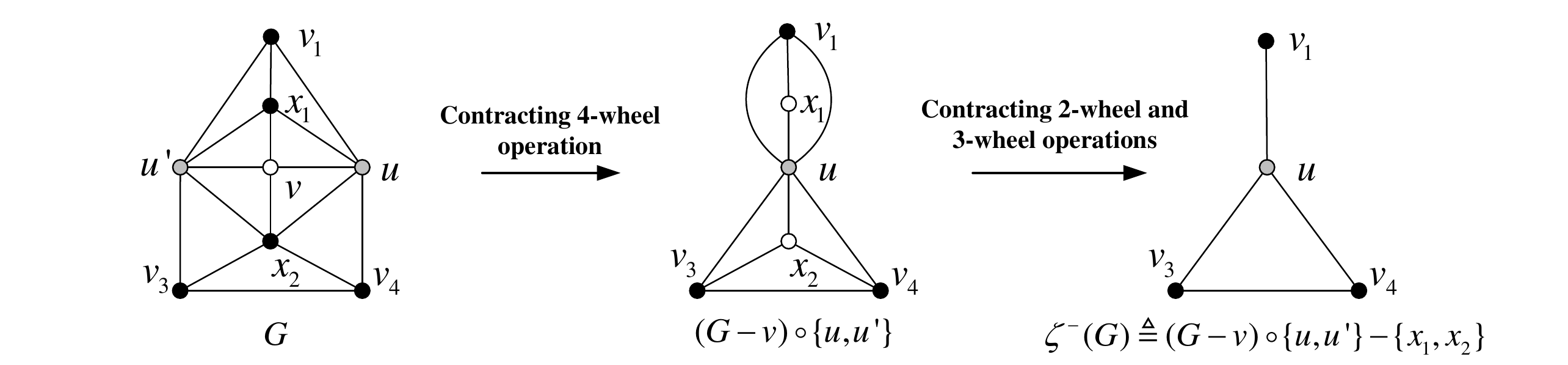} \\
  \hline
  \end{tabular} 
\begin{tabular}{|l|l|}
  \hline
   (h) contracting 5-wheel, 2-wheel and 3-wheel (= 4-wheel, 3-wheel and \\
  4-wheel)  operations, and the corresponding contracted subgraph\\
  \hline
  \includegraphics [width=360pt]{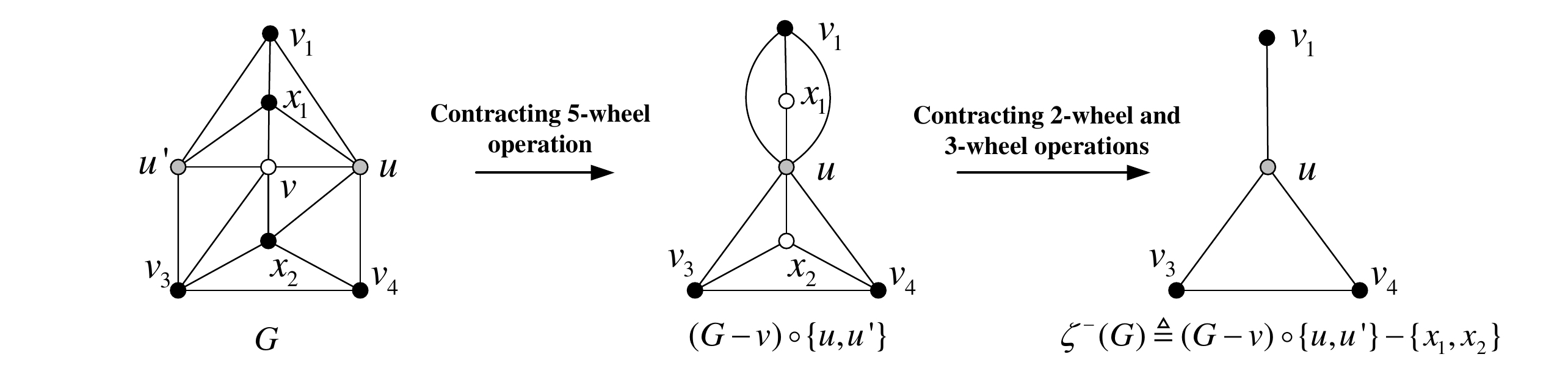} \\
  \hline
  \end{tabular} 
  \begin{tabular}{|l|l|}
   \hline
   (i) contracting 5-wheel, 3-wheel and 3-wheel operations, and the\\
    corresponding contracted subgraph\\
  \hline
  \includegraphics [width=360pt]{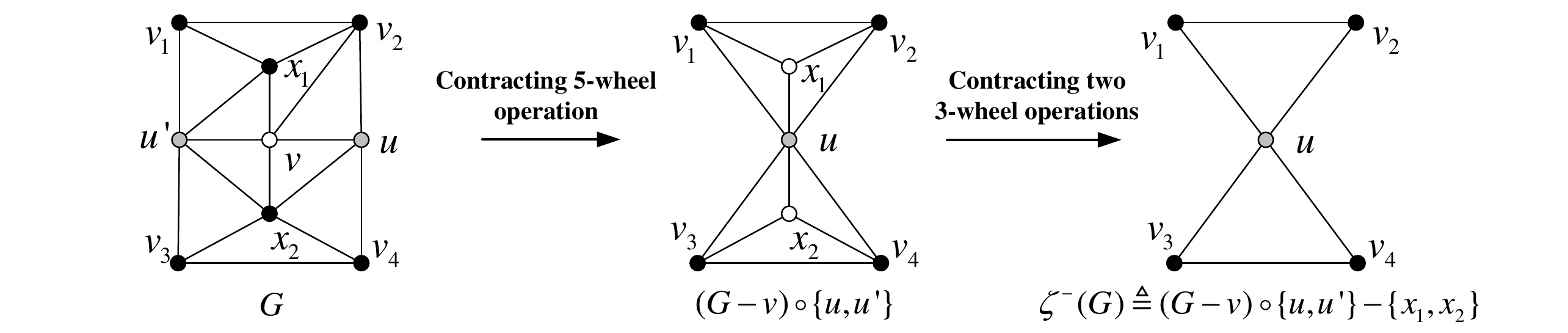} \\
  \hline
  \end{tabular}  
  \begin{tabular}{|l|l|}
  \hline
    (j) contracting 4-wheel, 3-wheel and 3-wheel (= 5-wheel, 2-wheel and\\
     3-wheel) operations, and the corresponding contracted subgraph\\
     (dumbbell transformation)\\
  \hline
  \includegraphics [width=360pt]{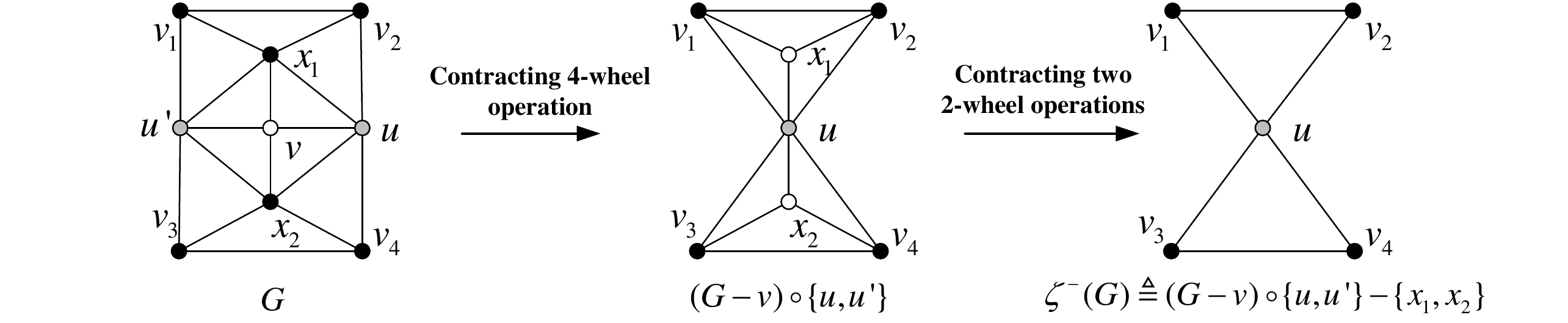} \\
  \hline
  \end{tabular} 
  \begin{tabular}{|l|l|}
  \hline
   (k) contracting 4-wheel and 2-wheel ($t$-times) operations, and the\\
    corresponding contracted subgraph\\
  \hline
  \includegraphics [width=360pt]{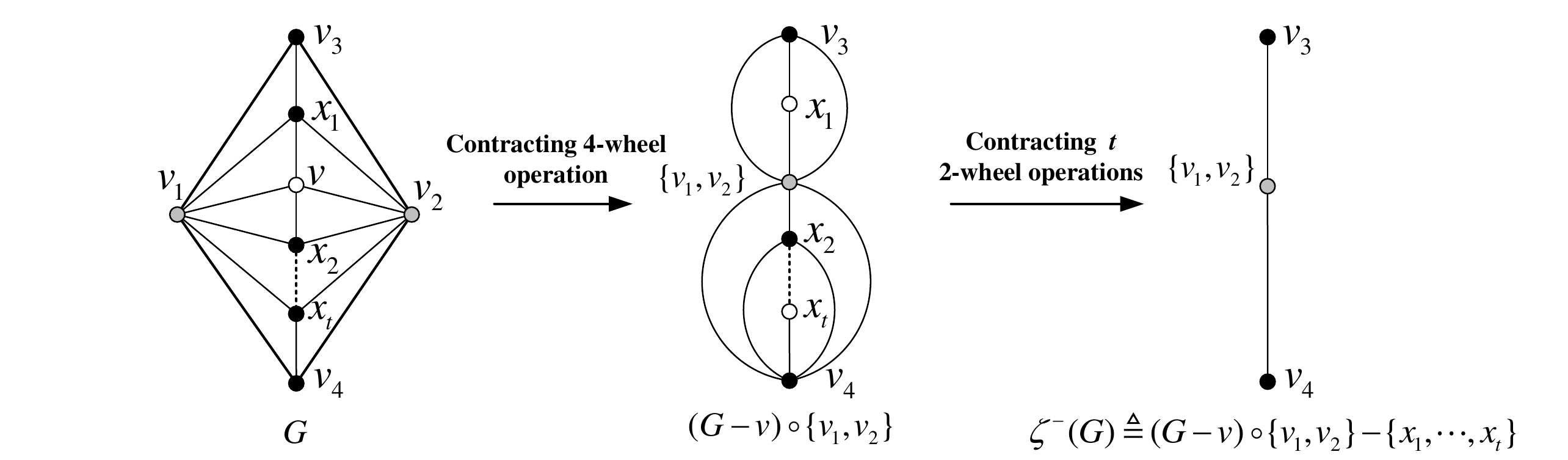} \\
  \hline
  \end{tabular}  
  \begin{tabular}{|l|l|}
  \hline
   (l) contracting 5-wheel, 2-wheel ($r$-times) and 3-wheel operations\\
 ($=$contracting 4-wheel ($r-1$  times), 2-wheel ($(t-r+1)$-times) and \\
  3-wheel), and the corresponding contracted subgraph\\
  \hline
  \includegraphics [width=360pt]{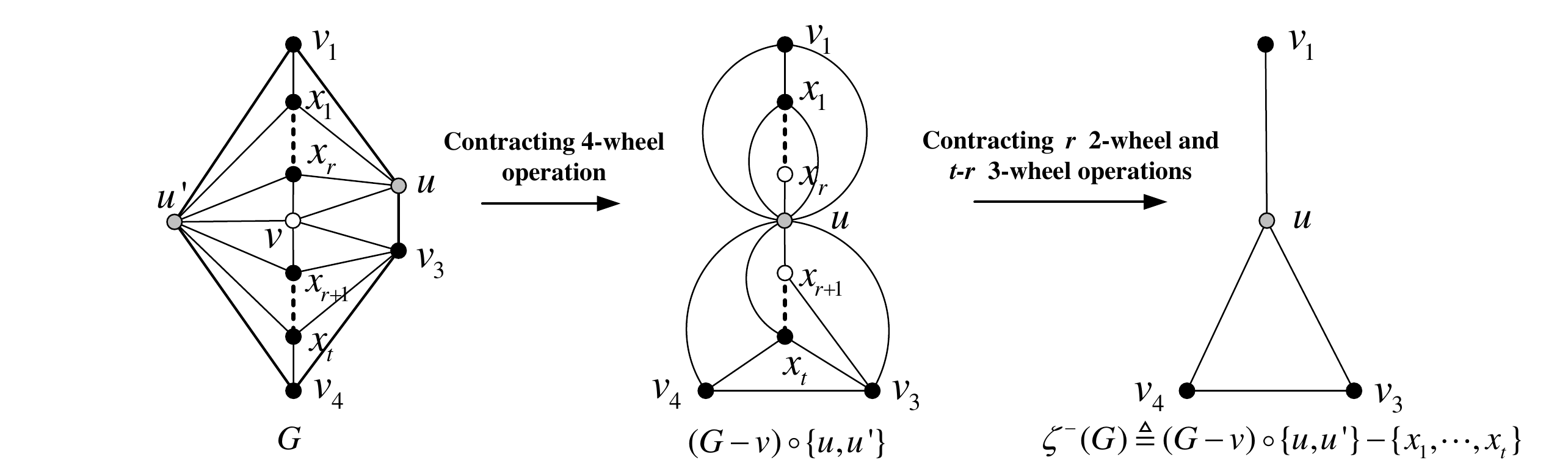} \\
  \hline
  \end{tabular} 
  \begin{tabular}{|l|l|}
  \hline
   (m) contracting 5-wheel and 3-wheel ($t$-times) operations,  and the\\
   corresponding contracted subgraph (two 5-degree vertices are adjacent)\\
  \hline
  \includegraphics [width=360pt]{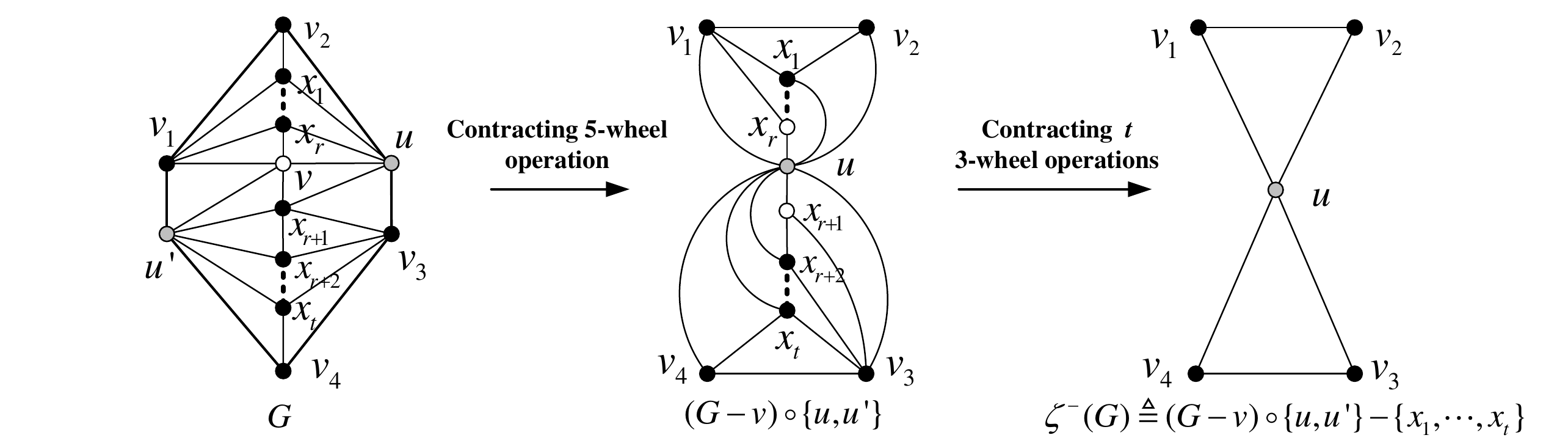} \\
  \hline
  \end{tabular}
  \begin{tabular}{|l|l|}
  \hline
   (n) contracting 5-wheel and 3-wheel ($t$-times) operations,  and the \\
  corresponding contracted subgraph (two 5-degree vertices are nonadjacent)\\
  \hline
  \includegraphics [width=360pt]{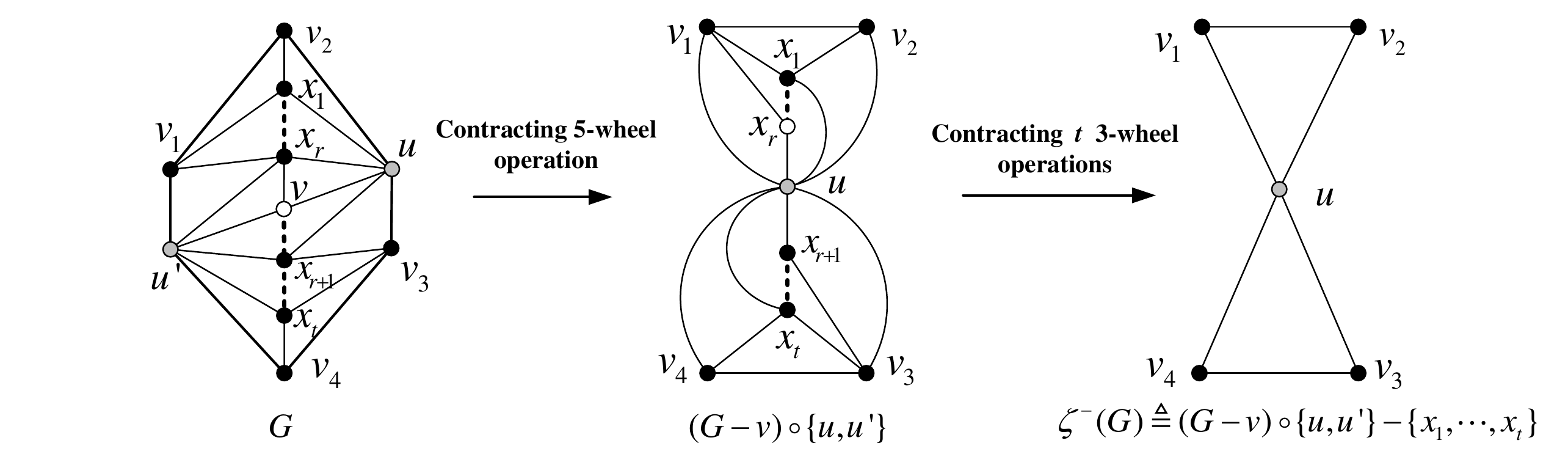} \\
  \hline
  \end{tabular}

\vspace{2mm}
According to Theorem \ref{theorem3.8},  the fourteen (class) contractible subgraphs shown in Figure 3.5 are unavoidable when we implement contracting wheel operations to maximal planar graphs with $\delta\geq 4$. In addition, they are all of the contractible subgraphs in this process, so we collect these fourteen (class) configurations together in the \emph{unavoidable-complete set of contractible subgraphs}. Similarly, as concerning extending wheel operation, there is also the concept of \emph{unavoidable-complete set of extendable subgraphs}.

Extending wheel operations mainly indicate extending 4-wheel and 5-wheel operations, where the operation objects of them are 2-length path (see Figure 3.6(a)) and the funnel (see Figure 3.6(b)), respectively. So, for a maximal planar graph $G$ with minimum degree at least four, if we conduct the extending 4-wheel operation on a 2-length path,  the resulting graph $\zeta^+(G)$ is still a maximal planar graph of $\delta=4$. However, a given $n$-vertex maximal planar graph of $\delta\geq 4$ may not be generated by just one extending 4-wheel or one 5-wheel operation from an $(n-2)$-vertex maximal planar graph. It may be generated from an $(n-3)$-vertex or an $(n-4)$-vertex maximal planar graph by implementing some extending wheel operations. This can be seen in the process of contracting wheel operation.

\begin{Prop}\label{prop3.9} Suppose that $G$ is a maximal planar graph with $\delta=2$ or 3.

\textcircled{1}  Let $\zeta_4^+(G)$ be the resulting graph from $G$ by conducting extending 4-wheel operation on a 2-length path $P_3$. Then $\delta(\zeta_4^+(G))\geq 4$ if and only if there is no 2-degree or 3-degree vertices in $G$ except for the two ends of $P_3$.

\textcircled{2}  Let $\zeta_5^+(G)$ be the resulting graph from $G$ by conducting extending 5-wheel operation on a funnel $L$. Then $\delta(\zeta_5^+(G))\geq 4$ if and only if there is no 2-degree or 3-degree vertices in $G$ except for the funnel top and one of the funnel bottom, and 2-degree vertex can only be the funnel top.
\end{Prop}

\begin{center}
        \includegraphics [width=350pt]{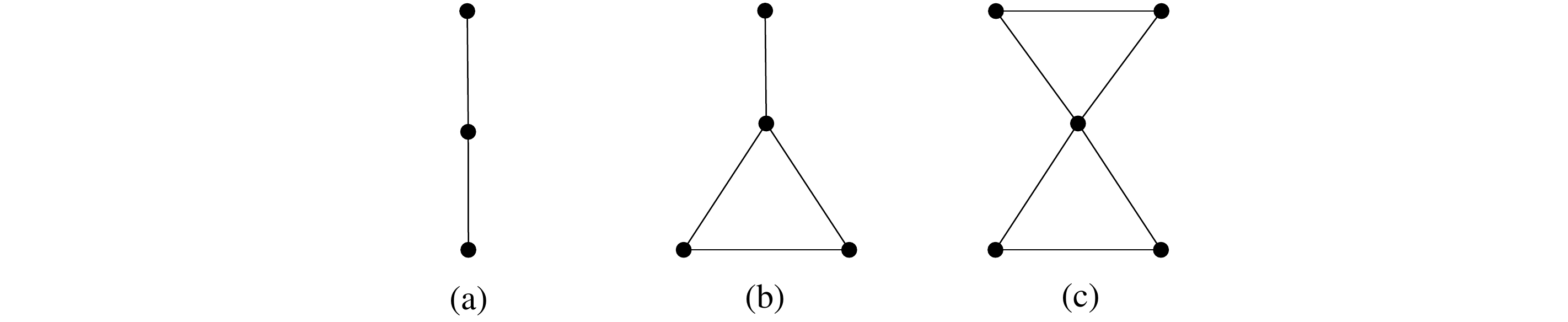}

        \textbf{Figure 3.6.} Operation objects of extending wheel operation
\end{center}

This proposition implies that for a maximal planar graph of $\delta\geq 4$, we can first do one or two extending 2-wheel or 3-wheel operations to $G$, and then implement extending 4-wheel or 5-wheel operations. When we conduct extending 4-wheel operation on a 2-length path, of which the two ends have degree 3,  we should first conduct two extending 3-wheel operations on two triangles that have a common vertex (see Figure 3.6(c)). We call the graph shown in Figure 3.6(c) the \emph{dumbbell}. It is easy to prove that there are many kinds of compound extending wheel operations based on dumbbells. For this view, we generally refer to the three subgraphs shown in Figure 3.6 as three \emph{object subgraphs} that are extended in extending wheel operation, where the first two graphs are viewed as the \emph{basic object subgraphs}.

For a maximal planar graph $G$, if there is a 2-length path $P_3=xuy$ (or a funnel $L=v_1-\triangle v_2v_3v_4$, where $v_1,v_2$ are the top and middle of $L$, and $v_3,v_4$ are the bottoms of $L$) of $G$ so that the graph $\zeta^+(G)$ generated from $G$ by conducting extending 4-wheel operation on $P_3$ (or 5-wheel operation on $L$), has $\delta(\zeta^+(G))\geq 4$, then we say $G$ is \emph{extendable}. According to Proposition \ref{prop3.9}, we can know that if $G$ is extendable and the minimum degree of $G$ is 2 or 3, then the possible 2-degree or 3-degree vertices of $G$ must belong to $\{x,y,v_1,v_3,v_4\}$. Here, we make an agreement as follow.

\textcircled{1} When $\delta(G)\geq 4$, we say both $P_3$ and $L$ are \emph{extendable} in $G$.

\textcircled{2} When $\delta(G)=2$ or 3, we refer to the union graph of $P_3$ and the 2-wheel and 3-wheel corresponding to the 2 and 3-degree vertices in $P_3$ respectively (or $L$ and the 2-wheel and 3-wheel corresponding to the 2 and 3-degree vertices in $L$ respectively) as a \emph{extendable subgraph} of $G$. For example when we conduct extending 4-wheel operation on $P_3$, if $d_G(x)=d_G(y)=2$, $P_3\cup G[x\cup \Gamma(x)] \cup G[y\cup \Gamma(y)]$ (see Figure 3.7(c)) is an extendable subgraph of $G$; if $d_G(x)=d_G(y)=3$,  Figure 3.7(j) is an extendable subgraph of $G$, etc.

We refer to the set of all possible extendable subgraphs of an extendable maximal planar graph $G$ as the \emph{unavoidable-complete set of extendable subgraphs} of $G$, denoted by $\varsigma(G)$.

\begin{theorem2}\label{theorem3.10}
Suppose that $G$ is an extendable maximal planar graph. Then  $\varsigma(G)$ consists of  eleven extendable subgraphs shown in Figure 3.7, where the 2-length paths labeled by bold are the object subgraphs of extending wheel operations.
\end{theorem2}

\begin{center}
        \includegraphics [width=350pt]{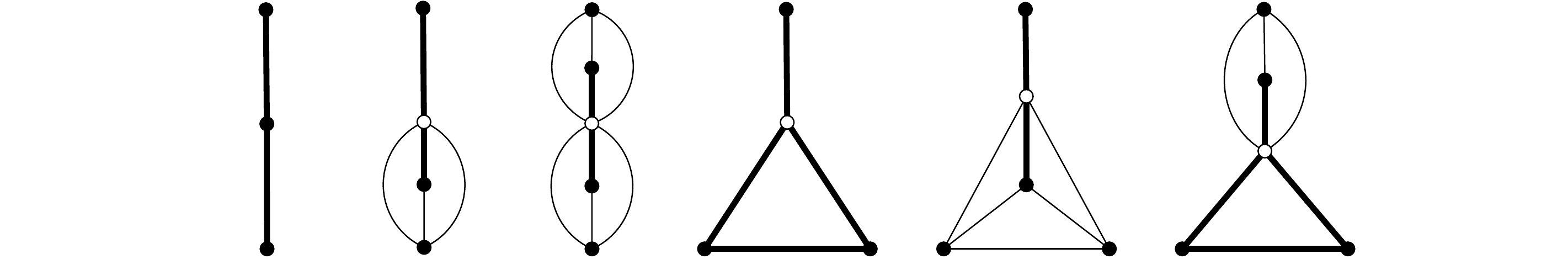}

         \vspace{2mm}
         \includegraphics [width=350pt]{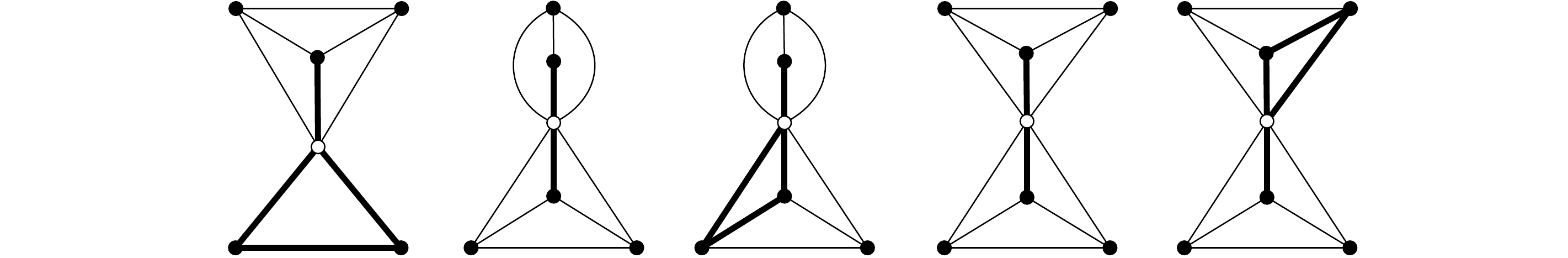}

         \vspace{2mm}
        \textbf{Figure 3.7.} The unavoidable-complete set of extendable subgraphs
\end{center}

\subsection{Constructing maximal planar graphs}

	Based on the extending wheel operation mentioned above, this section will show how to construct all $n$-vertex maximal planar graphs with $\delta(G)\geq 4$, in detail. Especially, we will present all of  $(6\sim 12)$-vertex maximal planar graphs with $\delta(G)\geq 4$.

	Let $Aut(G)$ denote the automorphism group of graph $G$, $xuy$ and $x'u'y'$ are two different paths of $G$. $xuy$ and $x'u'y'$ are called \emph{equivalent} if there is a $\sigma$  in $Aut(G)$ so that  $\sigma(x)=x'$, $\sigma(u)=u'$, $\sigma(z)=z'$, otherwise,  \emph{nonequivalent}.

\subsubsection{Constructing $n$-vertex maximal planar graphs with $\delta\geq 4$}

	\indent Step 1. Implementing extending 4-wheel or 5-wheel operations to $(n-2)$-vertex maximal planar graphs of $\delta\geq 4$.

	The detailed process is in the following. For an $(n-2)$-vertex maximal planar graph $G$, choose all of the nonequivalent 2-paths of $G$ first. For instance,  in a 7-vertex maximal planar graph, $G_7$, with $\delta(G_7)=4$ (shown in figure 3.8), there are four different 2-length paths: 444 type, 445 type, 454 type and 545 type respectively. Where 444 type means that the degree sequence of the 2-path is (444), and the other types are similar to this.  Second, conducting extending 4-wheel operation on each 2-path. For example, in $G_7$, when we
implement extending 4-wheel operations on 444 type and 454 type 2-length paths respectively, the two resulting 9-vertex maximal planar graphs are isomorphic; when we
implement extending 4-wheel operations on 545 type 2-length paths, the degree sequence of the resulting graph is  (444444477); when we implement extending 5-wheel operation on a funnel of $G_7$, the degree sequence of the resulting graph is (444455556). The processes of these extending wheel operations are shown in Figure 3.8, respectively.

	Step 2:  Implementing extending 2-wheel and 4-wheel operations (or extending 3-wheel and 5-wheel operations) to $(n-3)$-vertex maximal planar graphs of $\delta\geq 4$.

	In this step,  a 9-vertex maximal planar graph, $G_9$, can be only generated from a 6-vertex maximal planar graph, $G_6$. Because the 6-vertex maximal planar graph  with $\delta(G)=4$ is only the regular octahedron (as the first graph in Figure 3.8), it only generates a maximal planar graph $G$ of order 9 by doing extending 2-wheel and extending 4-wheel operations, and the resulting graph has degree sequence (444444666) and minimum degree $\delta=4$. Similarly, by conducting extending 3-wheel and extending 5-wheel operations, $G_6$ can also
generate only one maximal planar graph with minimum degree $\delta=4$ and degree sequence (444555555).

	We have indeed constructed all five 9-vertex maximal planar graphs with $\delta(G)=4$  in the above discussion. Note that 445 type cannot generate a 9-vertex maximal planar graphs with $\delta(G)=4$, and all other graphs generated by doing extending 5-wheel operations are isomorphic to one of these five graphs.

\begin{center}
        \includegraphics [width=230pt]{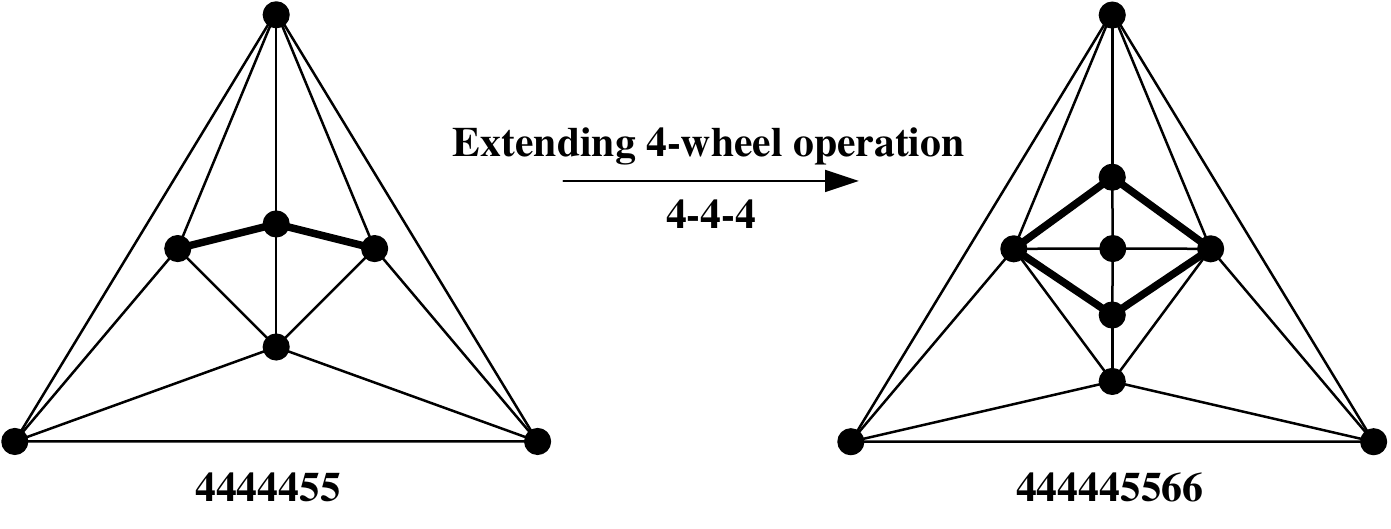}

          \vspace{2mm}
        \includegraphics [width=230pt]{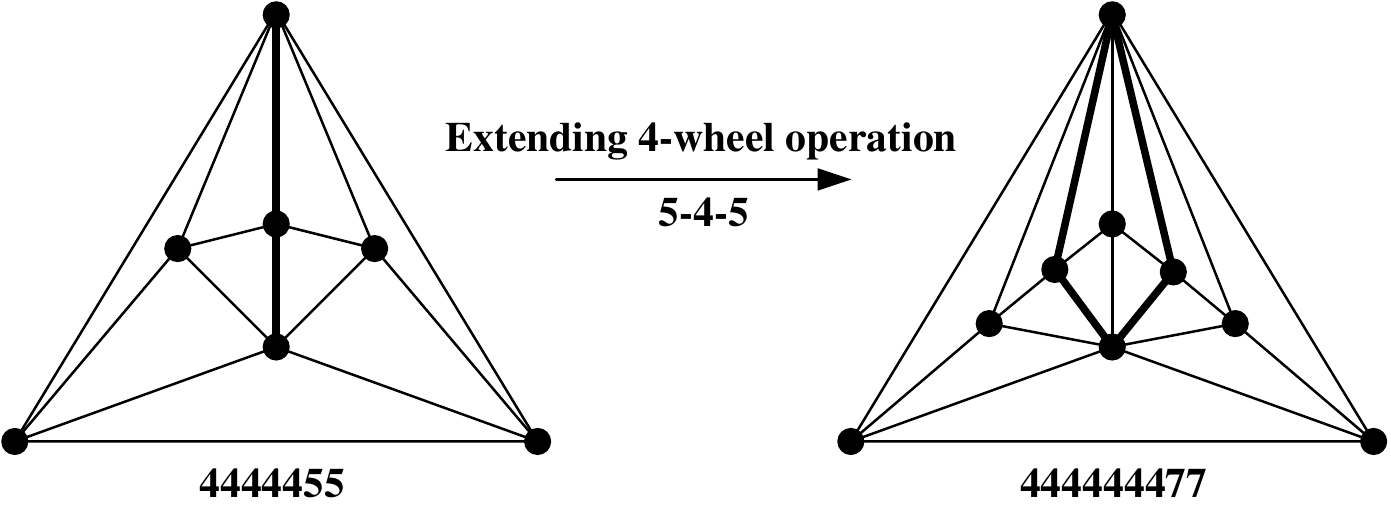}

         \vspace{2mm}
        \includegraphics [width=230pt]{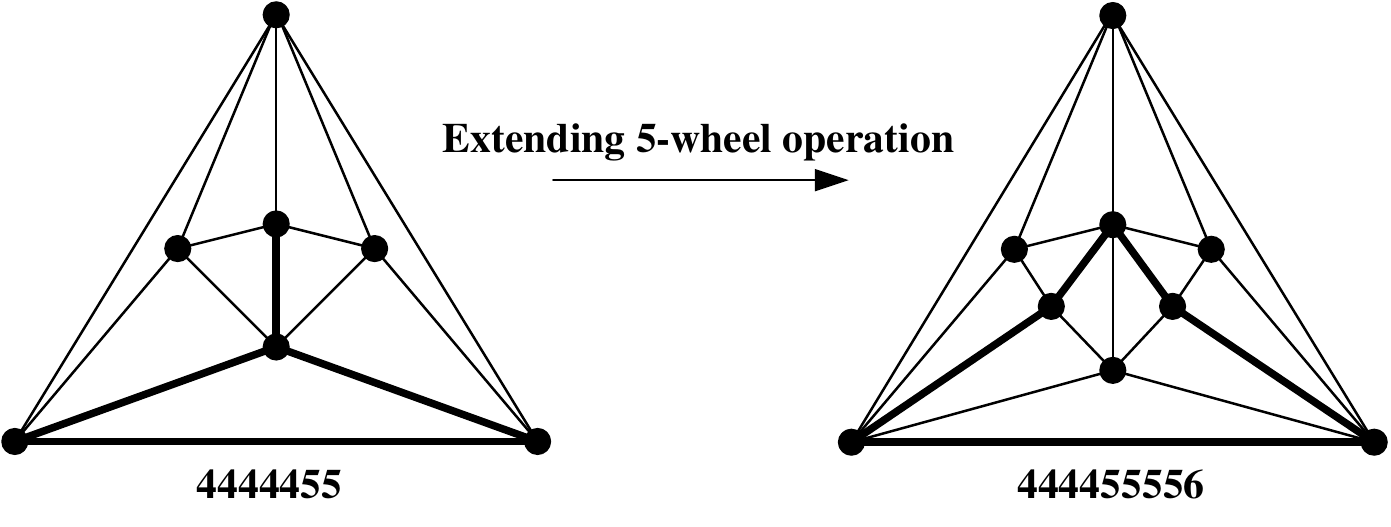}

         \vspace{2mm}
        \includegraphics [width=370pt]{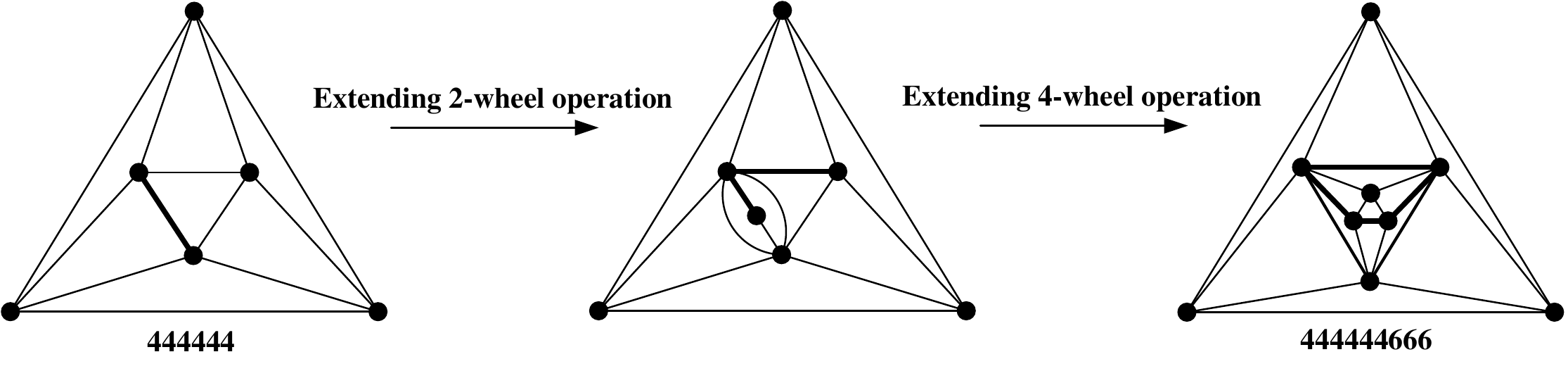}

        \vspace{2mm}
        \includegraphics [width=370pt]{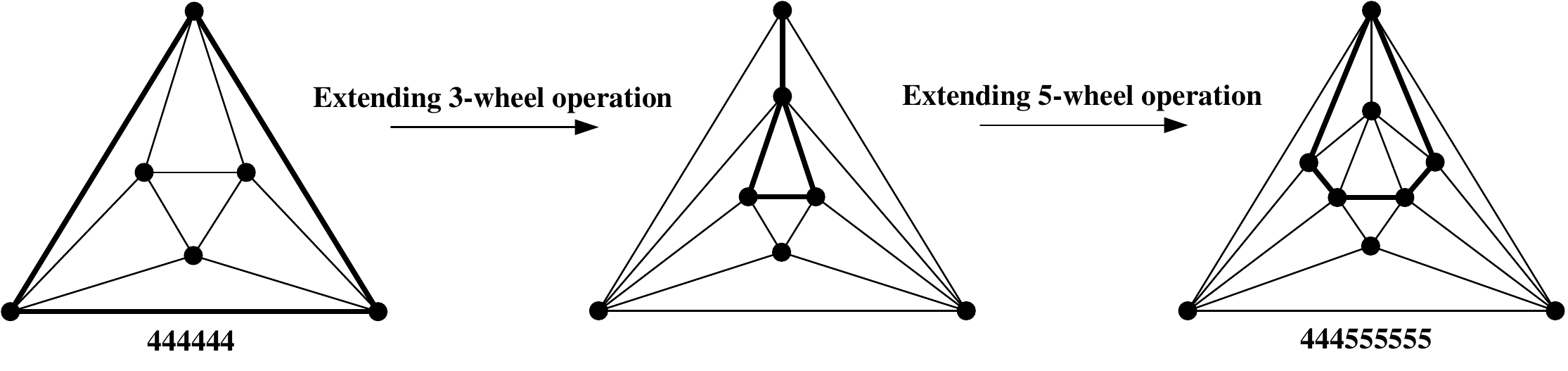}

        \vspace{2mm}
        \textbf{Figure 3.8.} Diagrams of generating 9-vertex maximal planar graphs of $\delta=4$ from 6,7-vertex maximal planar graphs
\end{center}

\subsubsection{All of the (6$\sim$12)-vertex maximal planar graphs with $\delta(G)\geq 4$ }

	In order to prove the main result  in this subsection, we need to investigate all $(6 \sim 12)$-vertex maximal planar graphs with $\delta(G)=4$. The enumeration of maximal planar graphs with minimum degree $\delta\geq 4$ has been studied, and an algorithm has been given to generate them by Brinkmann and McKay \cite{Brinkmann2007} in 2007. Here presents the number of $(6\sim 23)$-vertex maximal planar graphs of $\delta\geq 4$ (see Table 3.2).

\begin{center}
\textbf{Table 3.2.} The numbers of $(6\sim23)$-vertex maximal planar graphs with $\delta\geq 4$\label{tab1}

\hspace{0.5cm}

 \begin{tabular}{ ccccccccccccccccc }

    \hline
    Order  &6  & 7   &8  & 9 &  10\\
    \hline
    Count  &1 &  1  & 2 &  5 &  12 \\
    \hline
    Order  &  11 & 12 &13  & 14  &15 \\
    \hline
    Count  & 34  & 130  &525 &  2472  & 12400 \\
    \hline
     Order & 16 &  17&  18& 19 &20\\
    \hline
    Count &  65619 &  357504 & 1992985 & 11284042 &64719885 \\
    \hline
     Order   &21  &  22 & 23\\
    \hline
    Count   & 375126827  &  2194439398  & 12941995397\\
    \hline
  \end{tabular}
\end{center}

According to the methods of generating maximal planar graphs put forward above, we construct all $(6\sim 12)$-vertex  maximal planar graphs shown in Figure 3.9 $\sim$ Figure 3.14, respectively.

\begin{center}
        \includegraphics [width=180pt]{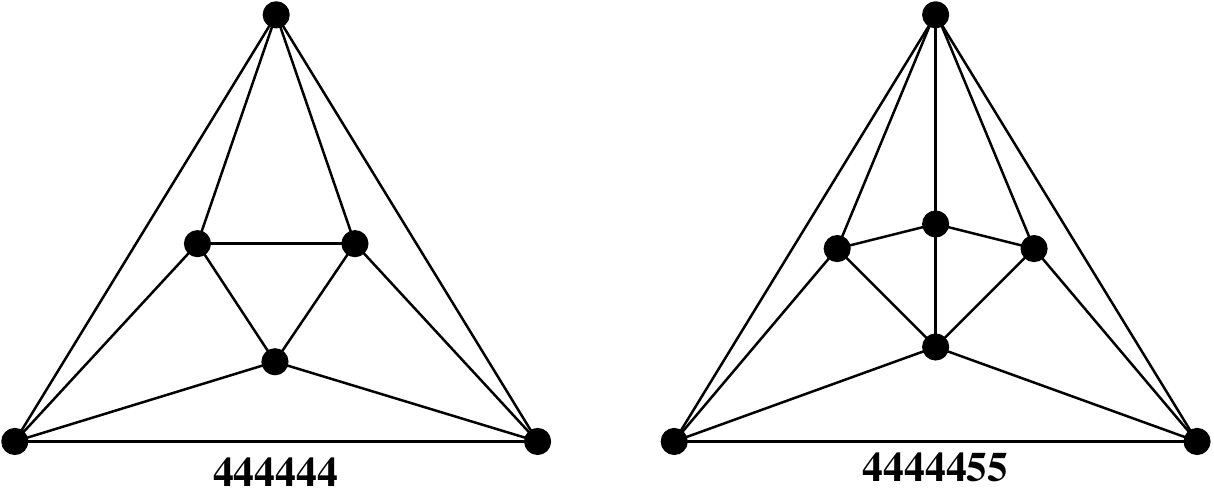}

        \vspace{2mm}
        \textbf{Figure 3.9.} The (6,7)-vertex maximal planar graphs with $\delta=4$
\end{center}

\begin{center}
        \includegraphics [width=180pt]{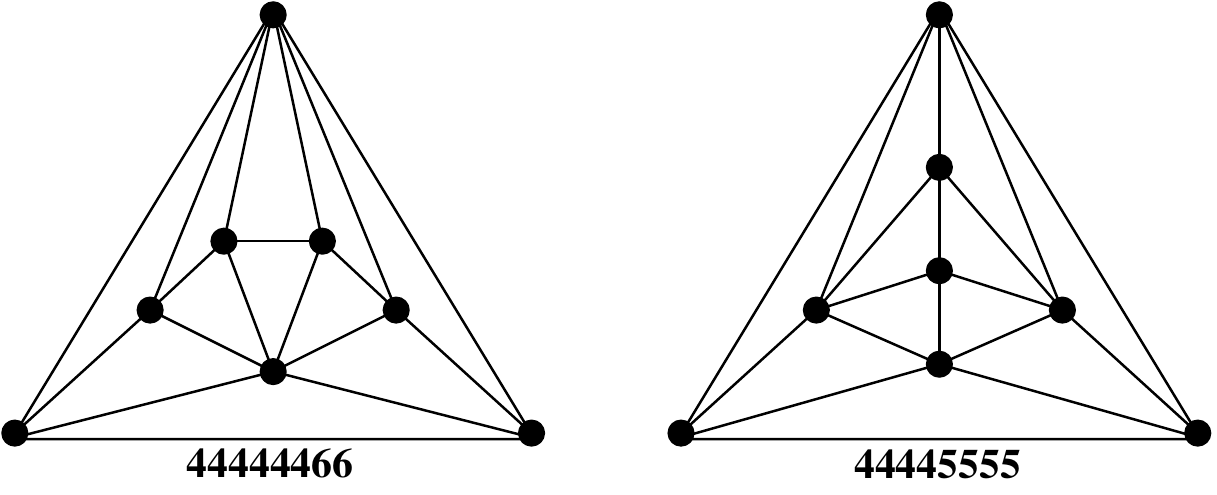}

        \vspace{2mm}
        \textbf{Figure 3.10.}  Two 8-vertex maximal planar graphs with $\delta=4$
\end{center}
\begin{center}
        \includegraphics [width=380pt]{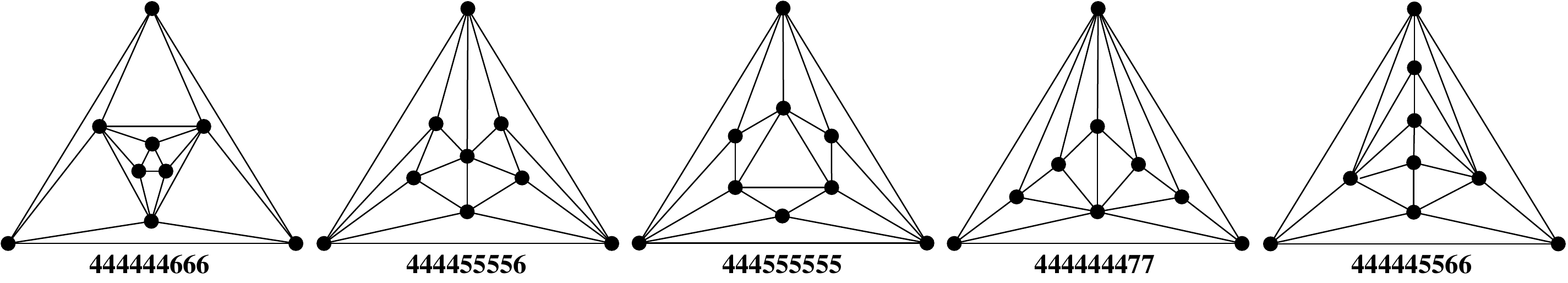}

        \textbf{Figure 3.11.}  Five 9-vertex maximal planar graphs with $\delta=4$
\end{center}
\begin{center}
        \includegraphics [width=380pt]{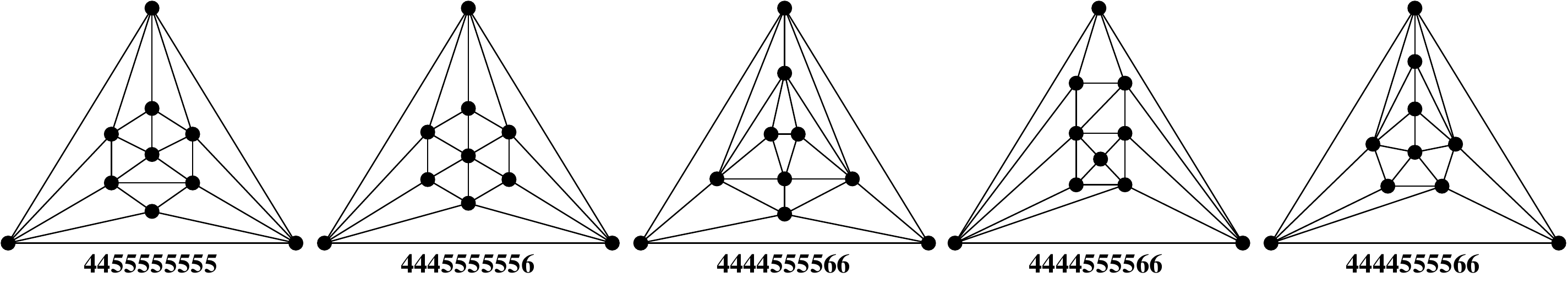}

        \vspace{2mm}
        \includegraphics [width=380pt]{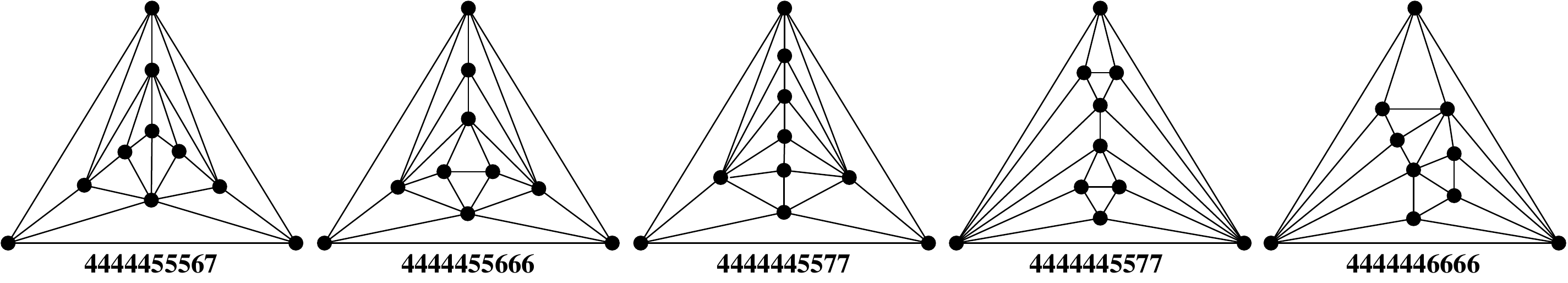}

        \vspace{2mm}
        \includegraphics [width=380pt]{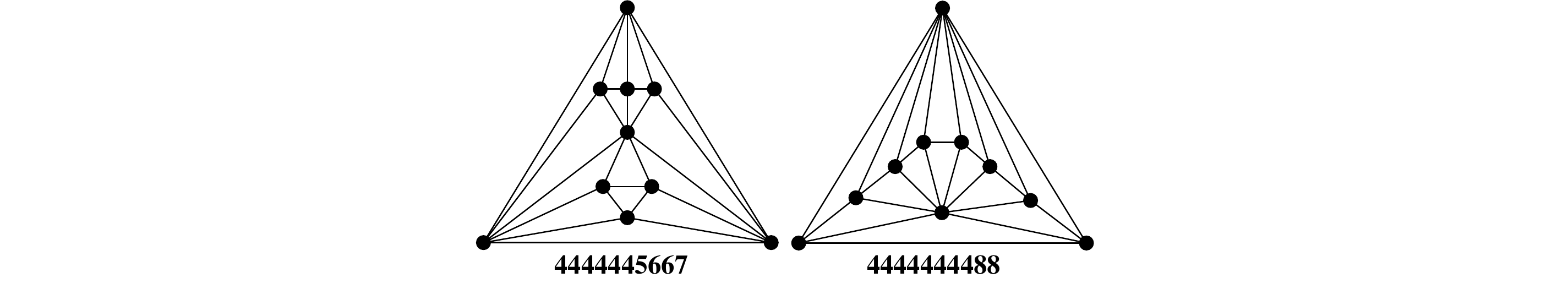}

        \vspace{2mm}
        \textbf{Figure 3.12.} Twelve 10-vertex maximal planar graphs with $\delta=4$
\end{center}
\begin{center}
        \includegraphics [width=380pt]{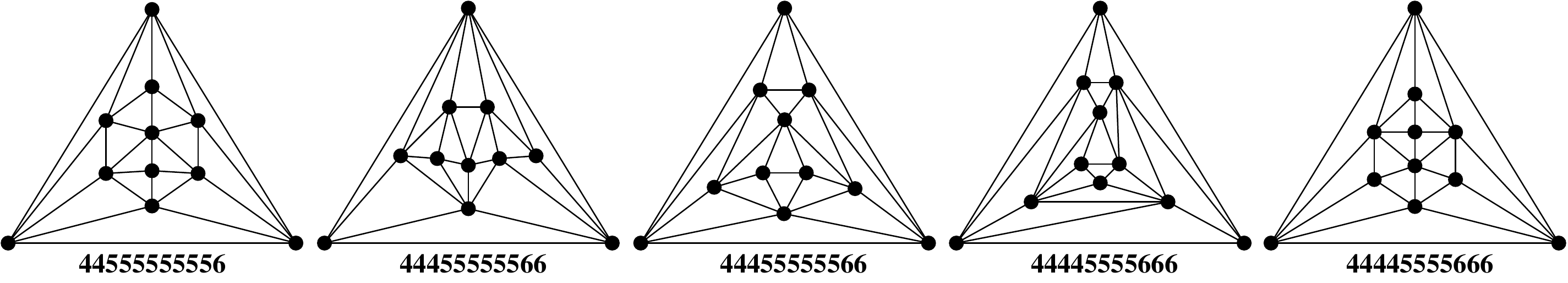}

        \vspace{2mm}
        \includegraphics [width=380pt]{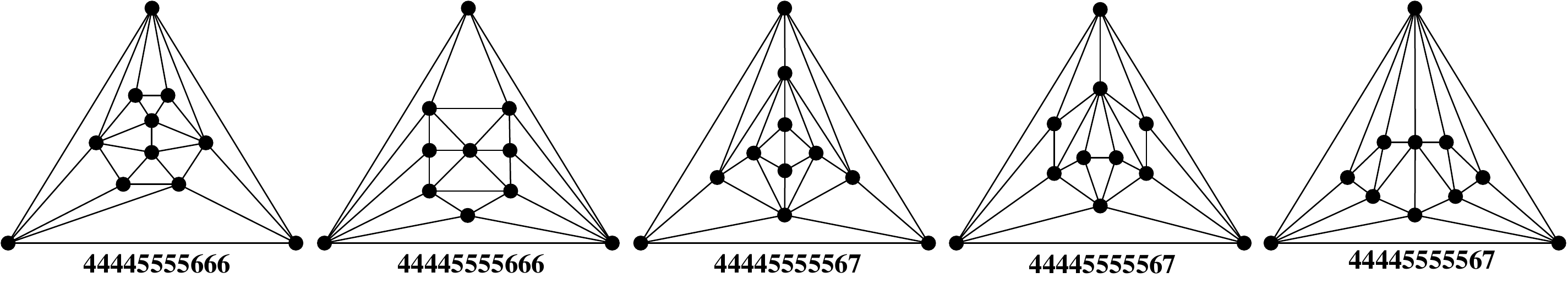}

        \vspace{2mm}
        \includegraphics [width=380pt]{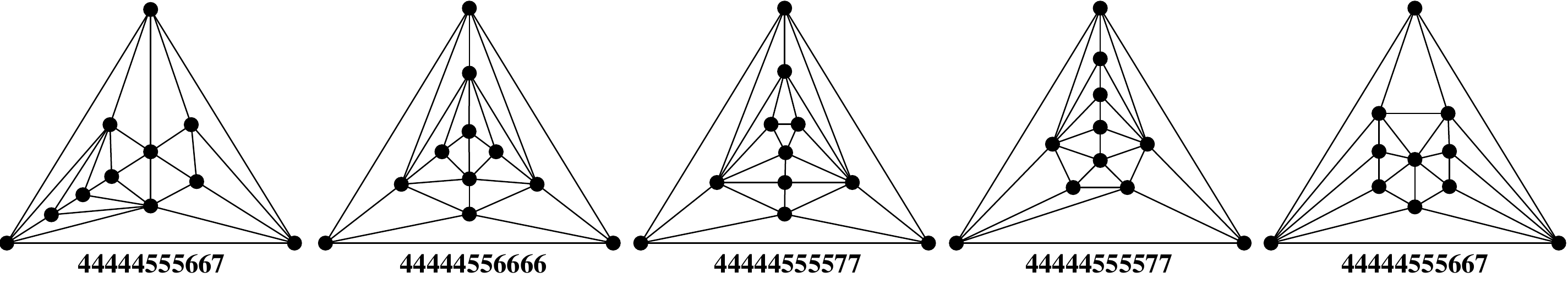}

        \vspace{2mm}
        \includegraphics [width=380pt]{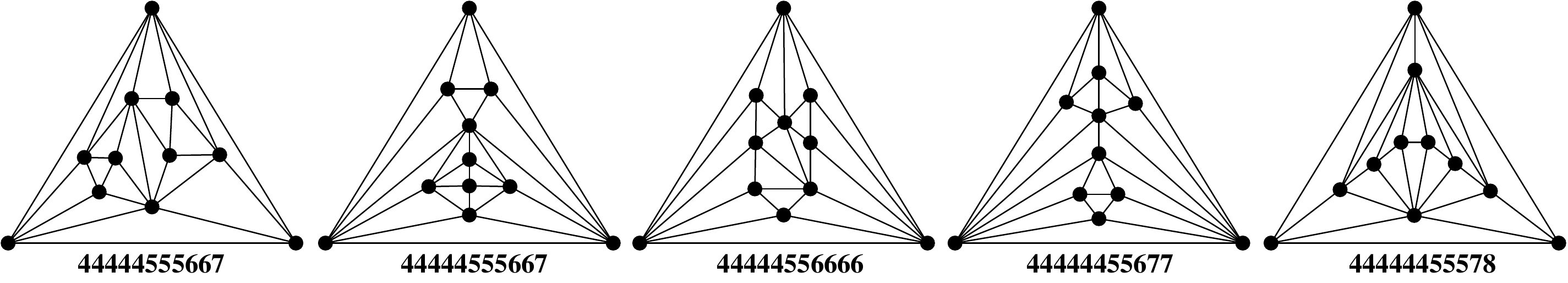}

        \vspace{2mm}
        \includegraphics [width=380pt]{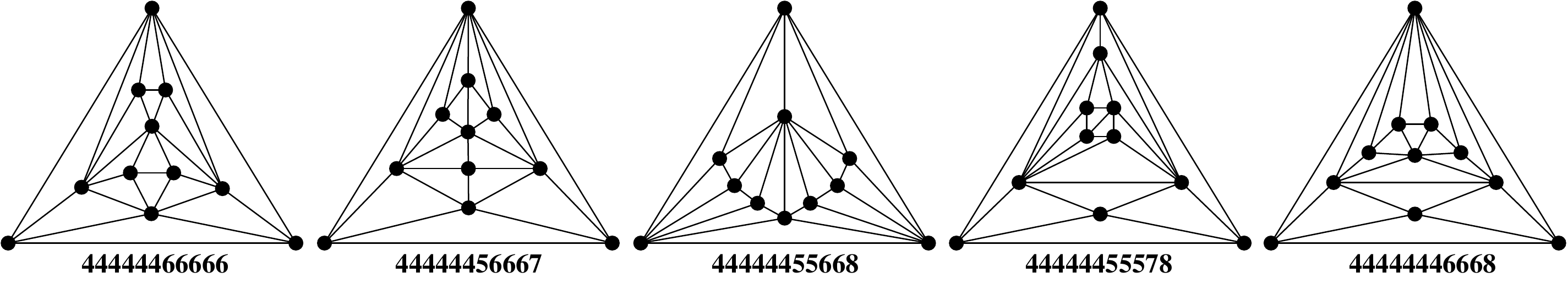}

        \vspace{2mm}
        \includegraphics [width=380pt]{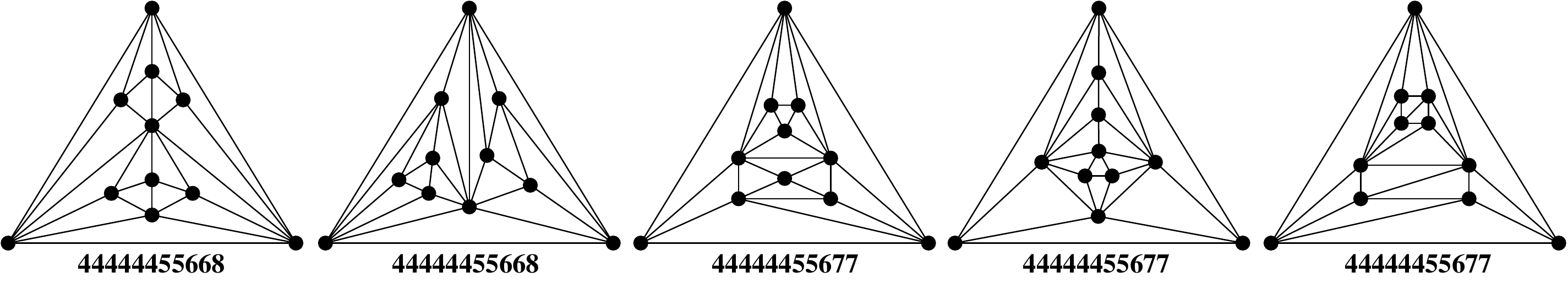}

        \vspace{2mm}
        \includegraphics [width=380pt]{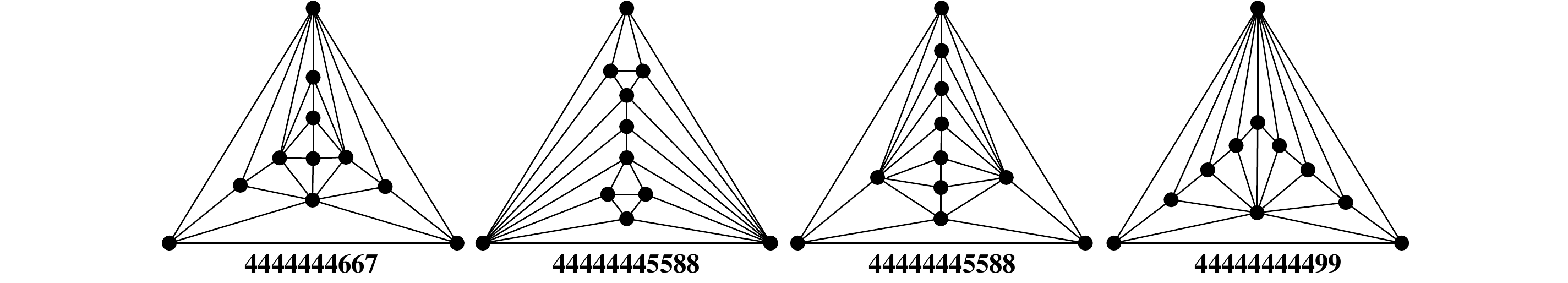}

        \vspace{2mm}
        \textbf{Figure 3.13.} Thirty-four 11-vertex maximal planar graphs with $\delta(G)=4$
\end{center}

\begin{center}
       \includegraphics [width=380pt]{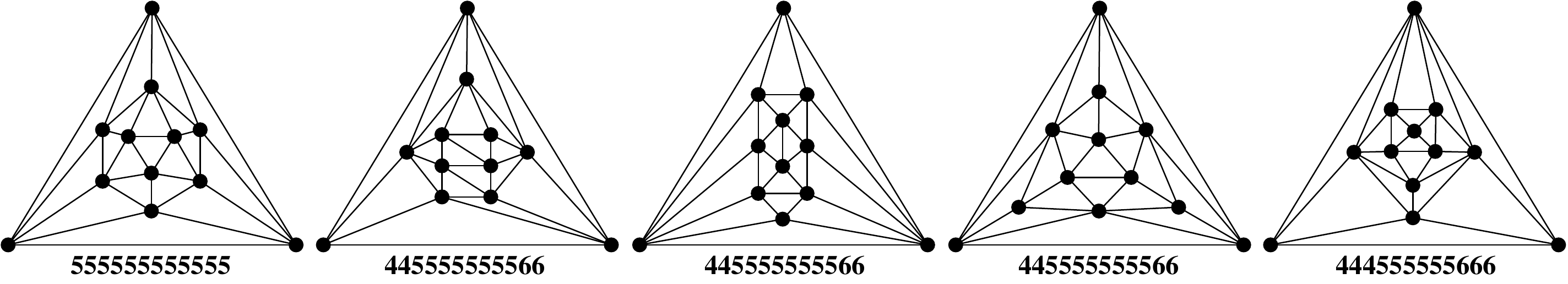}

        \vspace{2mm}
       \includegraphics [width=380pt]{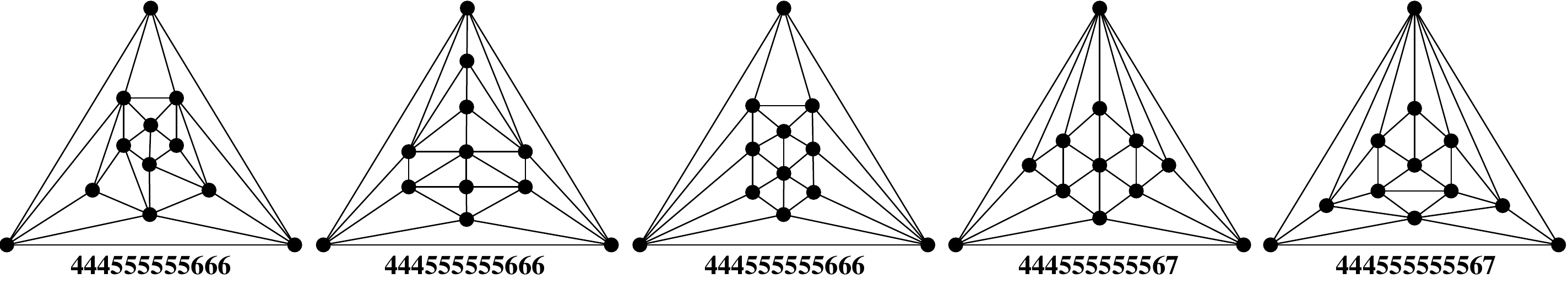}

        \vspace{2mm}
       \includegraphics [width=380pt]{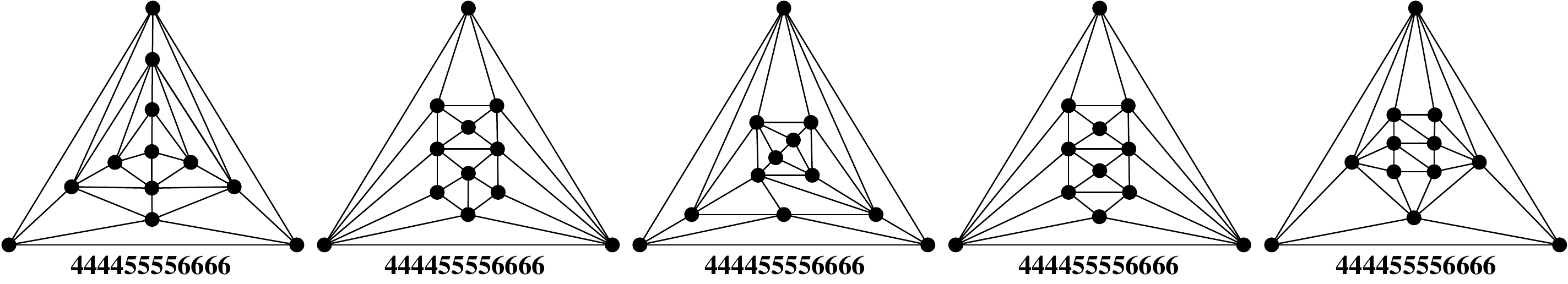}

        \vspace{2mm}
       \includegraphics [width=380pt]{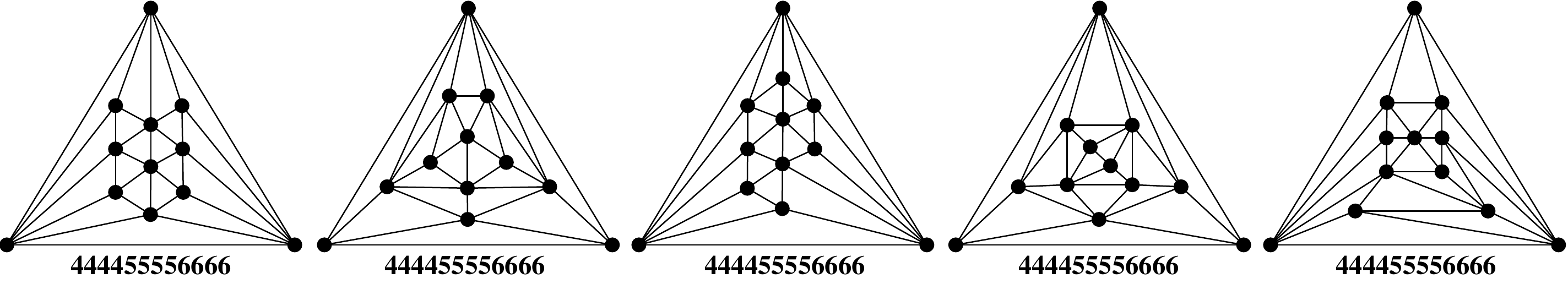}

        \vspace{2mm}
       \includegraphics [width=380pt]{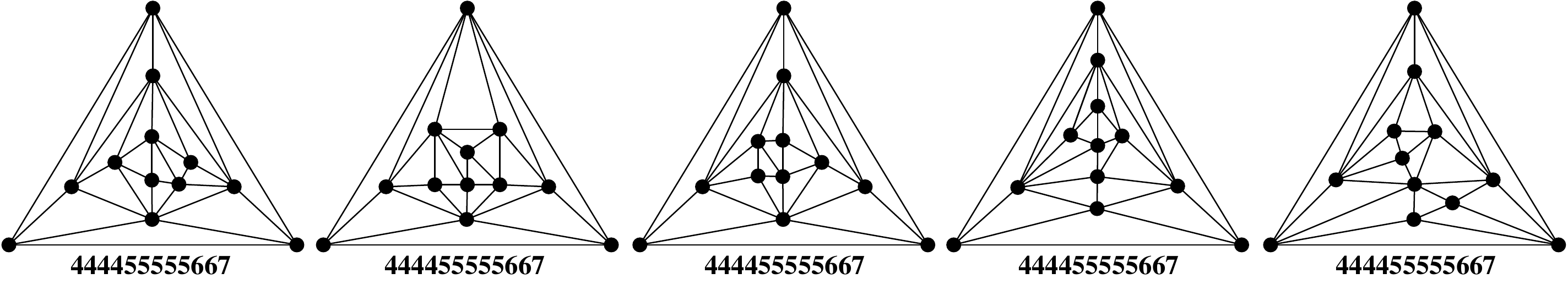}

        \vspace{2mm}
       \includegraphics [width=380pt]{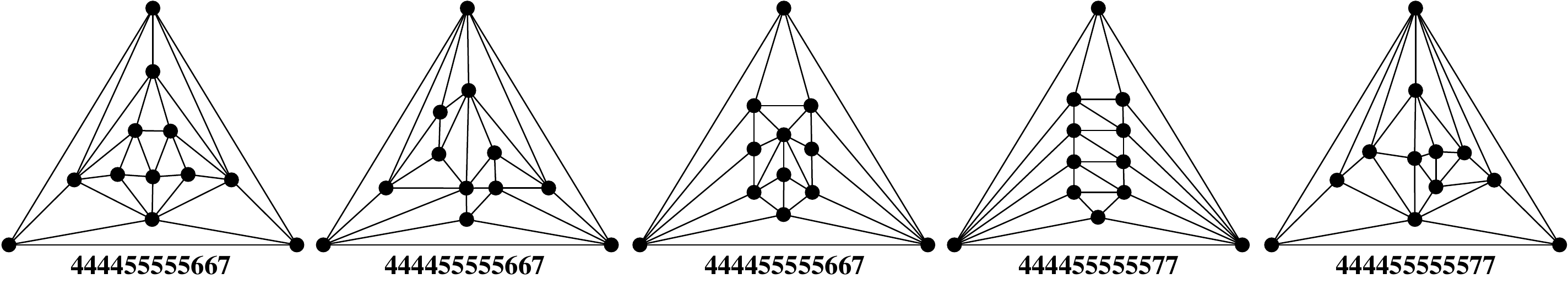}

        \vspace{2mm}
       \includegraphics [width=380pt]{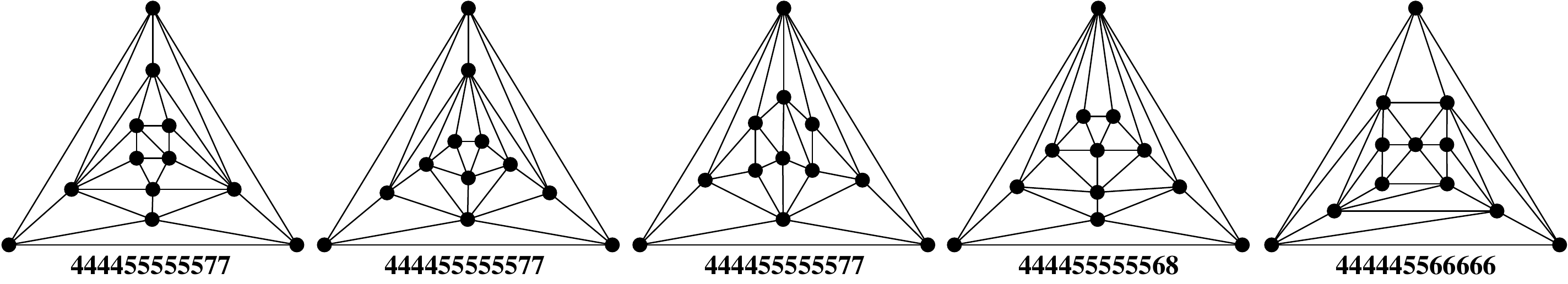}

        \vspace{2mm}
       \includegraphics [width=380pt]{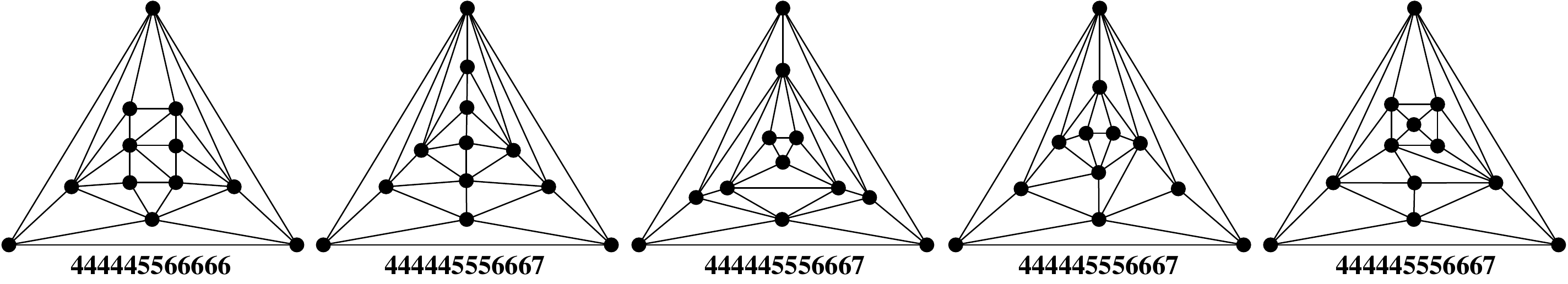}

        \vspace{2mm}
       \includegraphics [width=380pt]{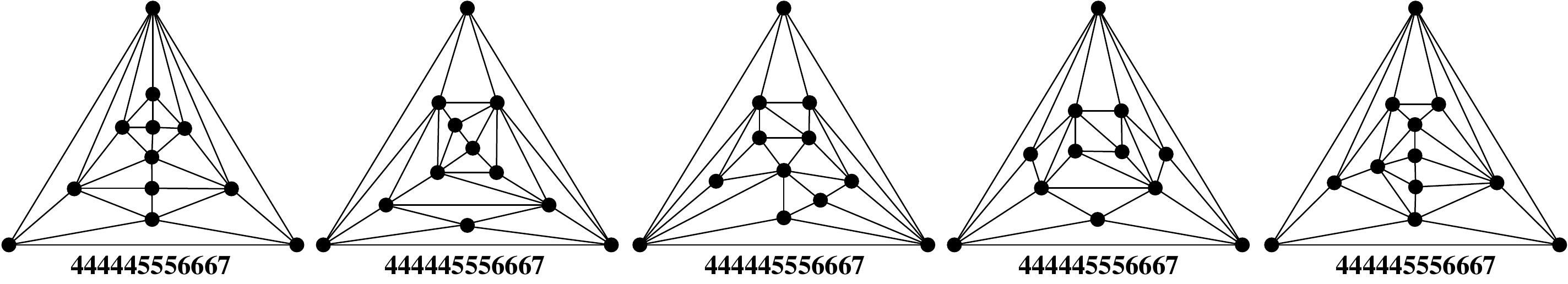}

        \vspace{2mm}
       \includegraphics [width=380pt]{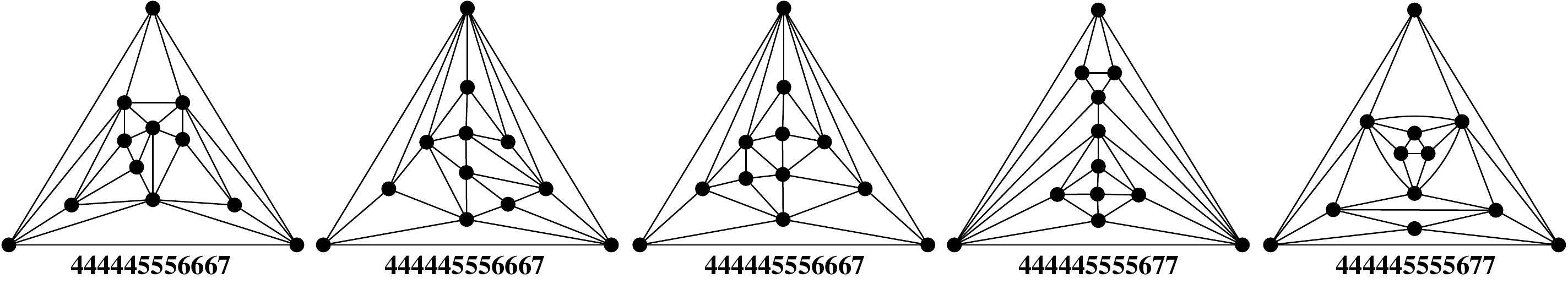}

        \vspace{2mm}
       \includegraphics [width=380pt]{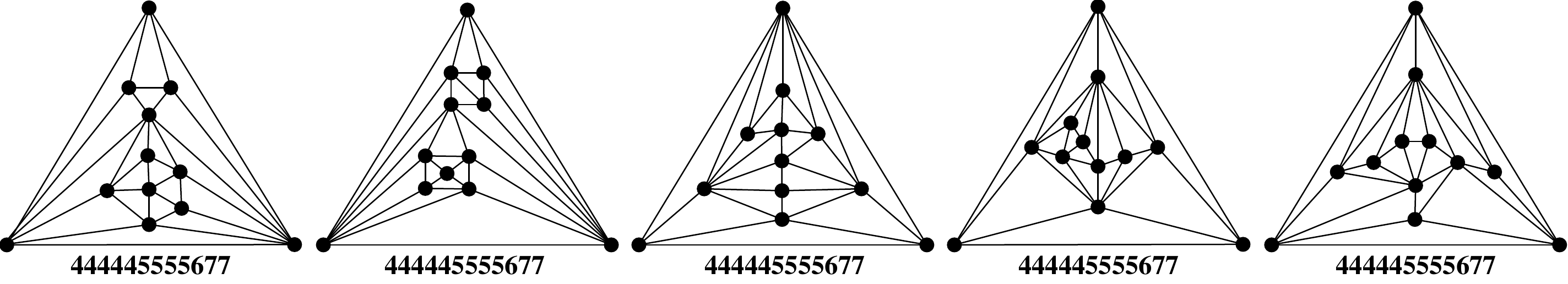}

        \vspace{2mm}
       \includegraphics [width=380pt]{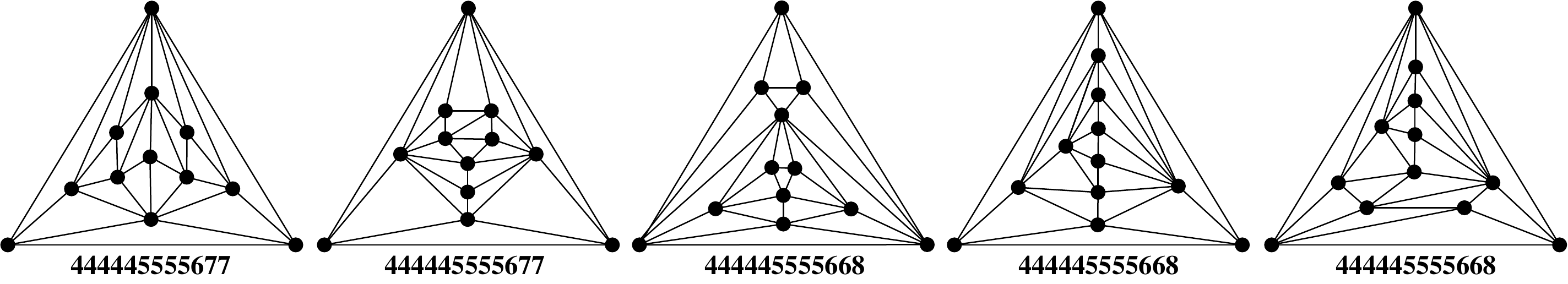}

        \vspace{2mm}
       \includegraphics [width=380pt]{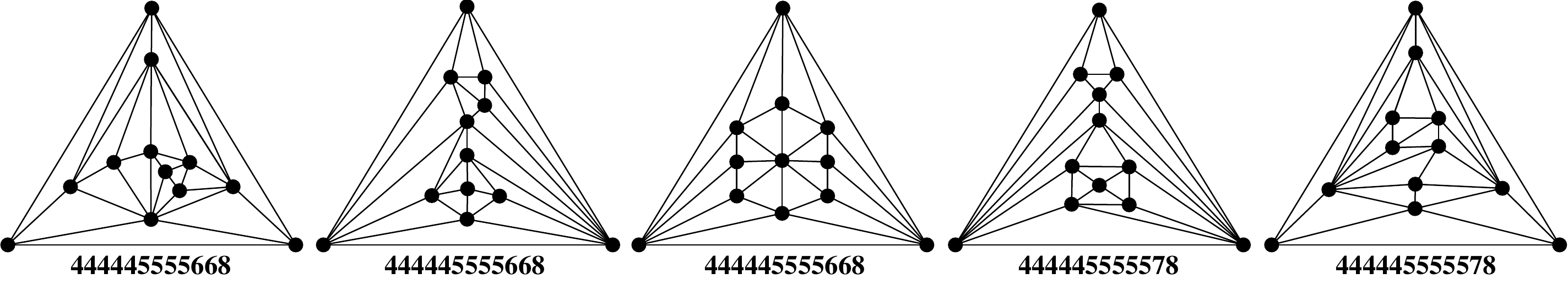}

        \vspace{2mm}
       \includegraphics [width=380pt]{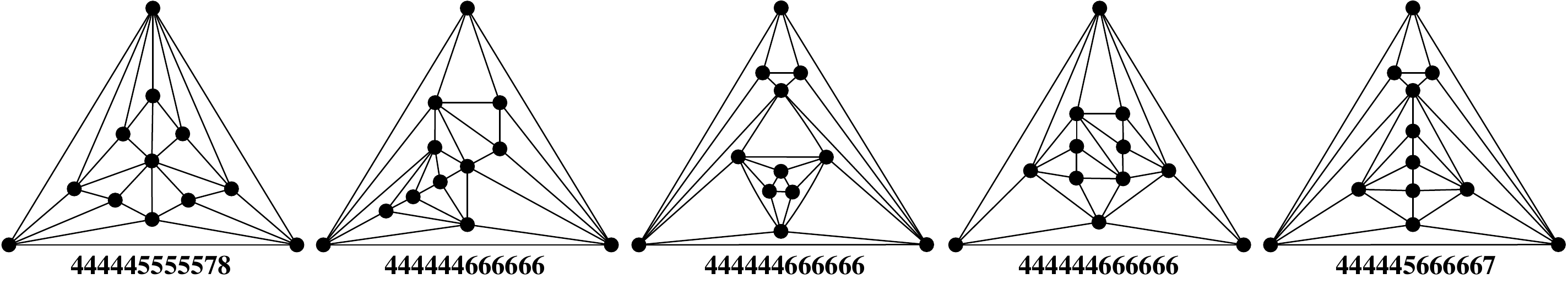}

        \vspace{2mm}
       \includegraphics [width=380pt]{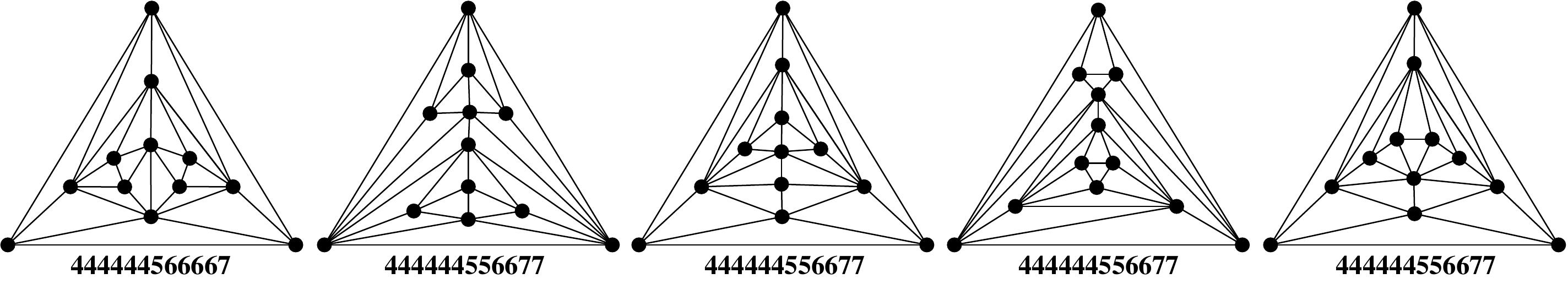}

        \vspace{2mm}
       \includegraphics [width=380pt]{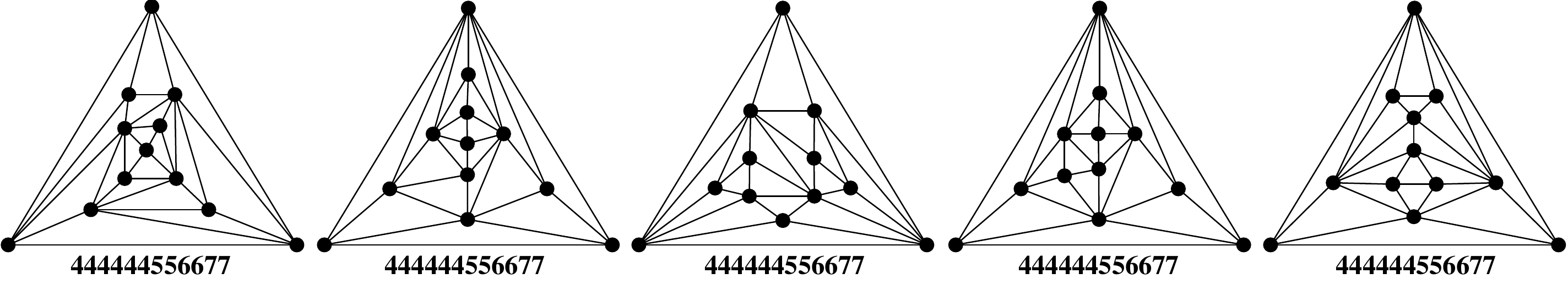}

        \vspace{2mm}
       \includegraphics [width=380pt]{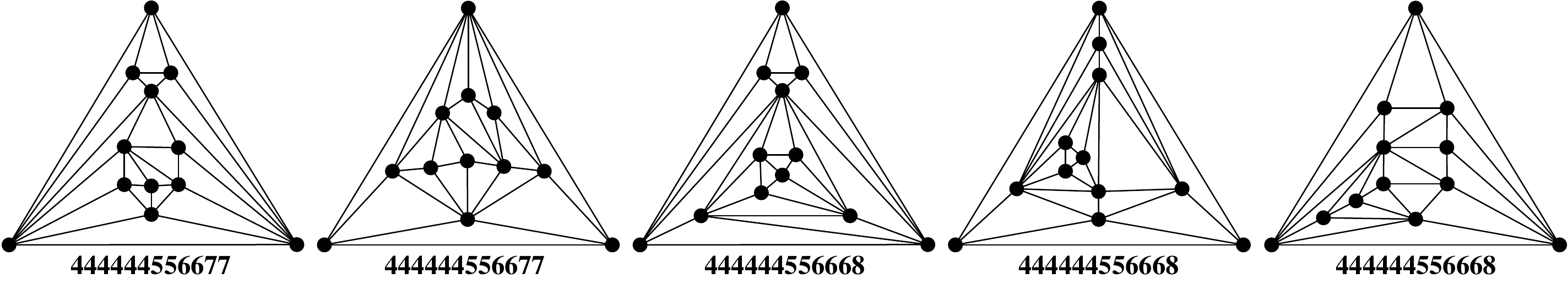}

        \vspace{2mm}
       \includegraphics [width=380pt]{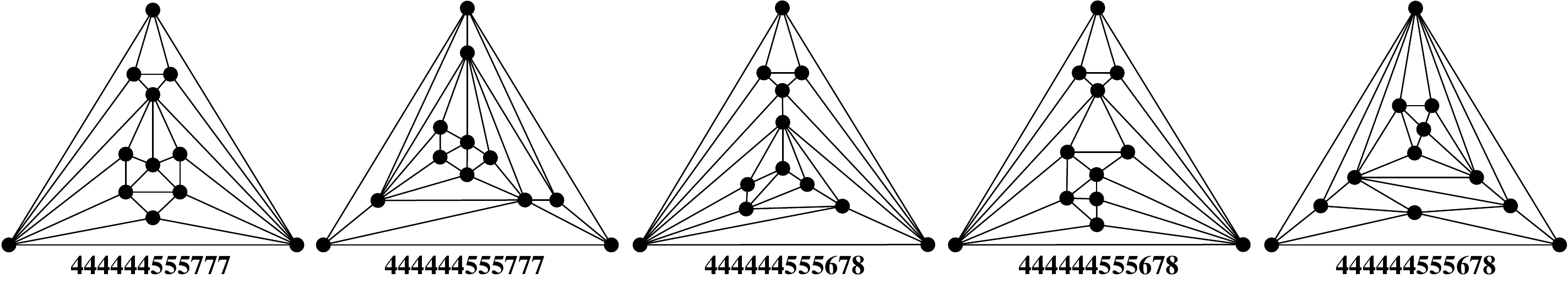}

        \vspace{2mm}
       \includegraphics [width=380pt]{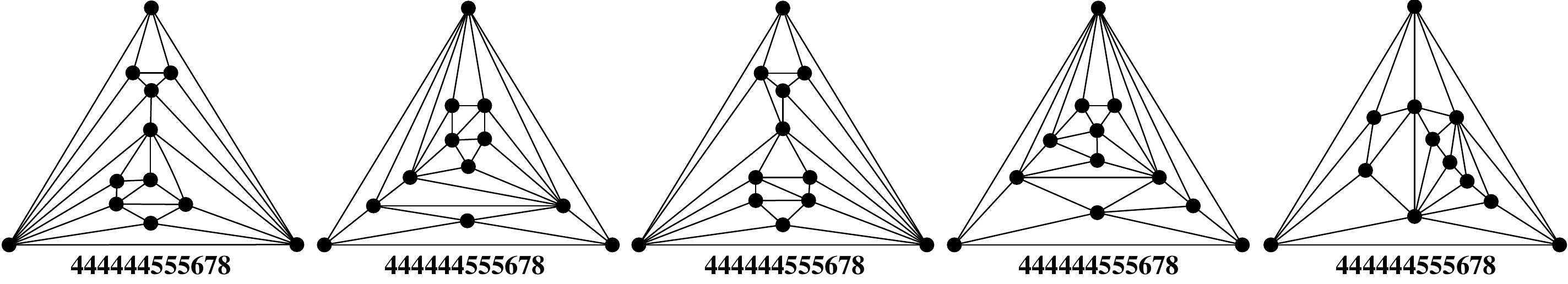}

        \vspace{2mm}
       \includegraphics [width=380pt]{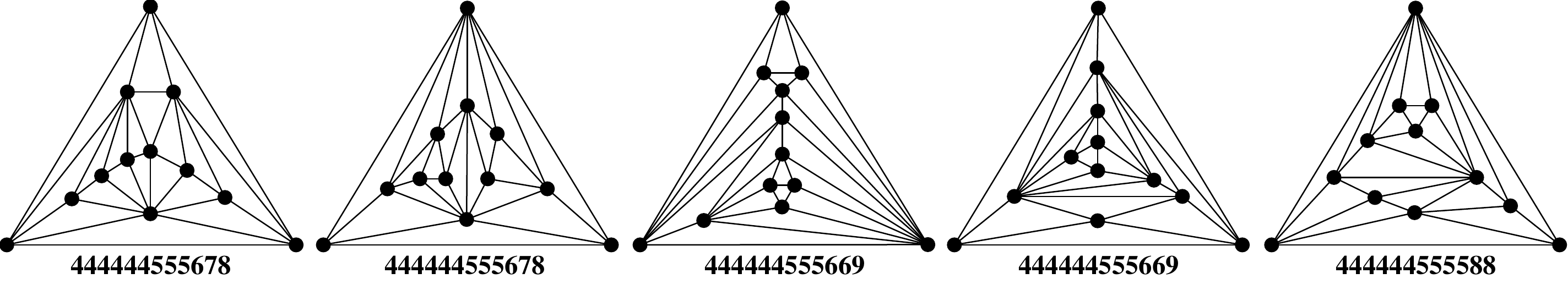}

        \vspace{2mm}
       \includegraphics [width=380pt]{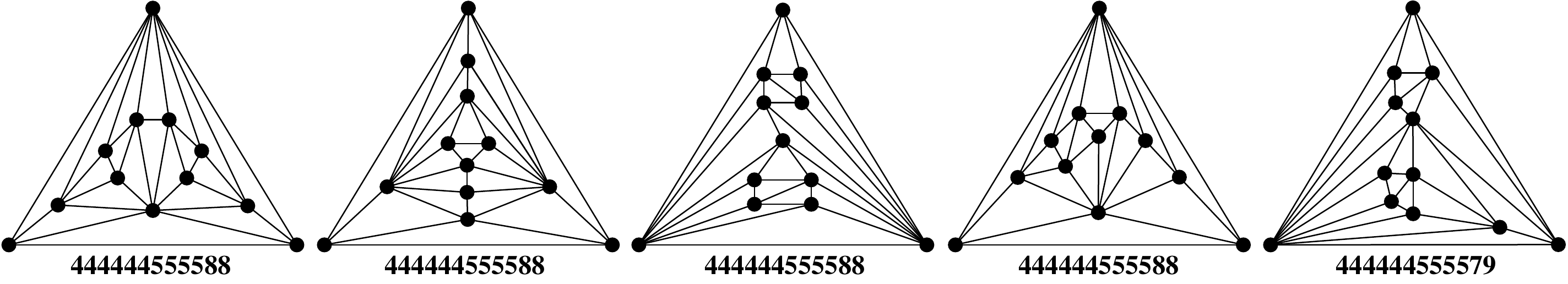}

        \vspace{2mm}
       \includegraphics [width=380pt]{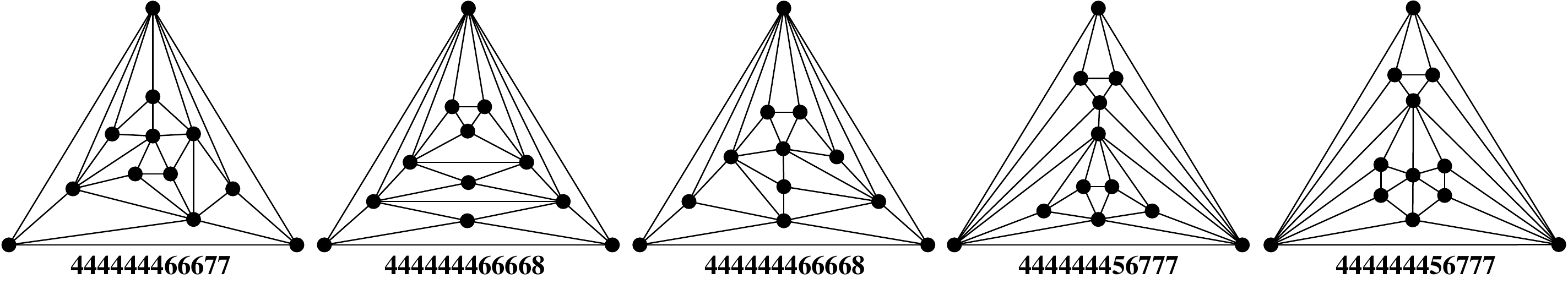}

        \vspace{2mm}
       \includegraphics [width=380pt]{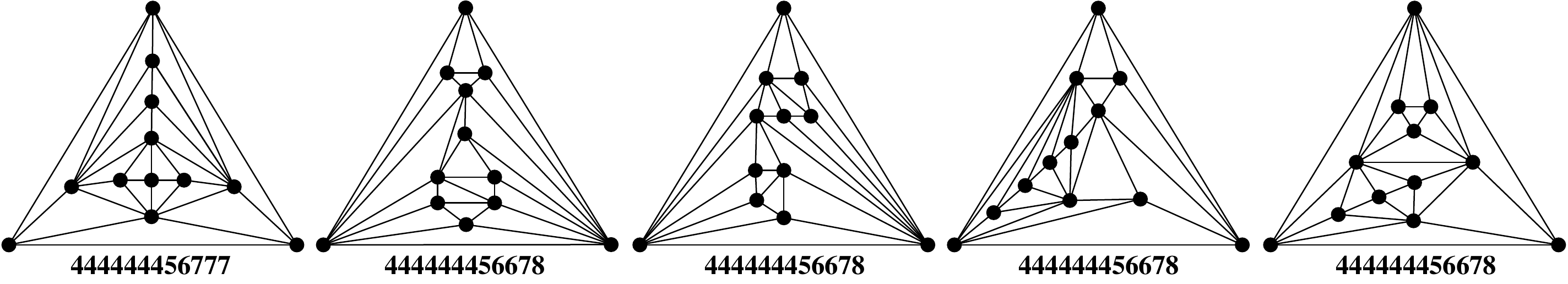}

        \vspace{2mm}
       \includegraphics [width=380pt]{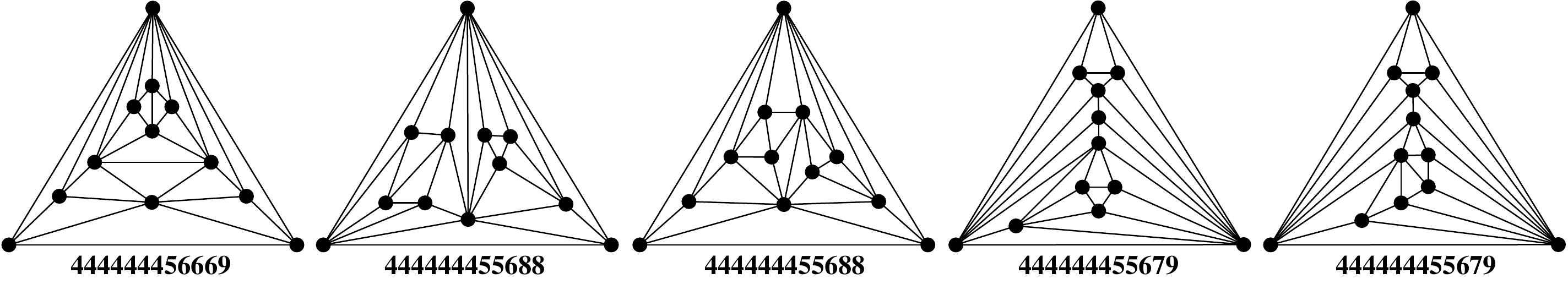}

        \vspace{2mm}
       \includegraphics [width=380pt]{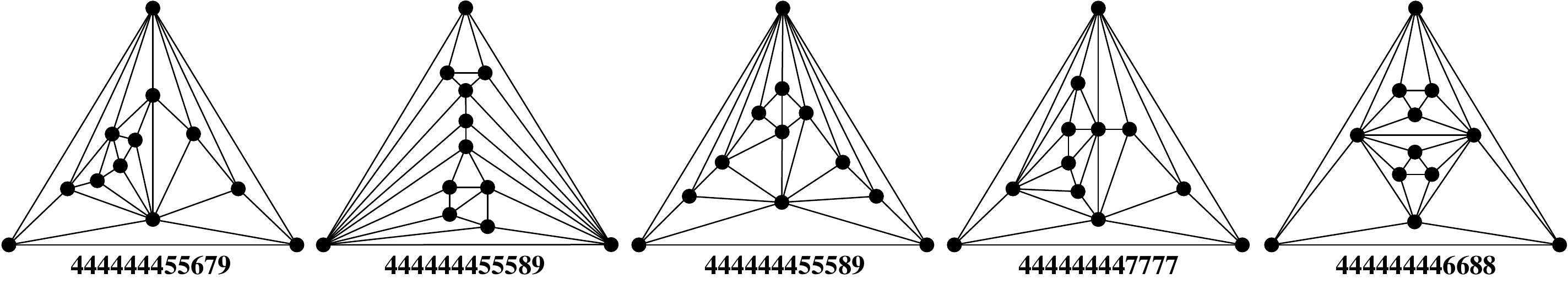}

        \vspace{2mm}
       \includegraphics [width=380pt]{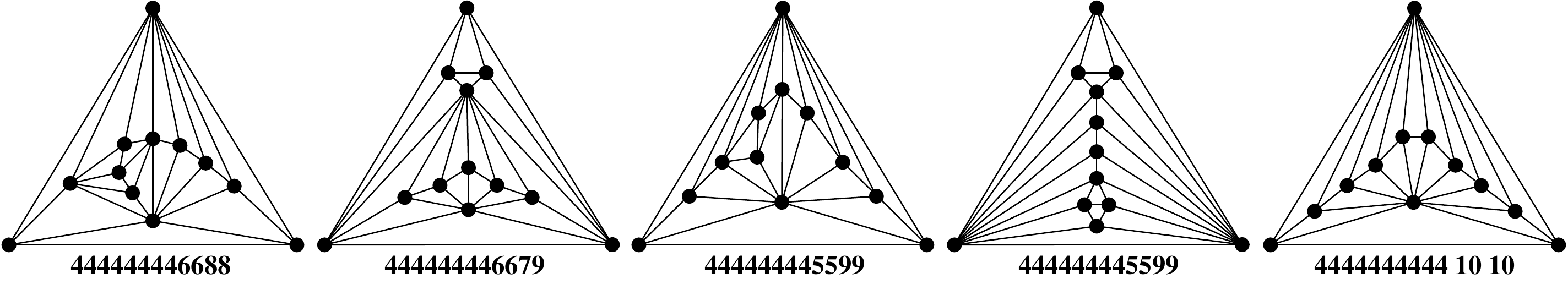}

        \vspace{2mm}
      \textbf{Figure 3.14.} A hundred and thirty 12-vertex maximal planar graphs with $\delta\geq 4$
\end{center}

\subsection{Basic extending and contracting operations based on coloring}

	On the basis of the previous two sections,  the \emph{contracting $k$-wheel operation under coloring} and its inverse operation--the \emph{extending $k$-wheel operation under coloring} to a maximal planar graph will be introduced in this section, where $k=2,3,4,5$. Meanwhile, some related properties are studied as well.

	The contracting 2-wheel operation and the extending 2-wheel operation under a coloring are almost as the same as the corresponding operations under the condition of no coloring. We need only concern how to assign a color to the center of the 2-wheel. This process is simple, so no more discussion is given here. Please see section 2.2.

	Let $G$ be a 4-colorable maximal planar graph. If $v$ is a 3-degree vertex of $G$, and $\Gamma(v)=\{v_1,v_2,v_3\}$, then  $\forall f\in C^{0}_{4}(G)$ ,\emph{the contracting 3-wheel operation under $f$} on $v$, means  deleting vertex $v$ from  $G$.
Naturally, the resulting graph $G-v$ is still a 4-colorable maximal planar graph. Meanwhile, the extending 3-wheel operation on the face $v_1-v_2-v_3$ under a coloring $f\in C^{0}_{4}(G-v)$ is in the following. First, we add a new vertex $v$ on this face, and make $v$ adjacent to $v_1,v_2,v_3$ respectively, and then assign $v$ a color differing from $f(v_1),f(v_2),f(v_3)$.

	Let $G$ be a 4-colorable maximal planar graph, $v$ be a 4-degree vertex of $G$, and  $\Gamma(v)=\{v_1,v_2,v_3,v_4\}$. It is clear that for any $f\in C^{0}_{4}(G)$, there is either  $f(v_1)=f(v_3)$ or $f(v_2)=f(v_4)$. From now on, we always assume $f(v_1)=f(v_3)$ (see Figure 3.15(a)). Then a \emph{contracting 4-wheel operation under $f$} on $v$ is to delete vertex $v$ from graph $G$, and identify vertices $\{v_1,v_3\}$ (see Figure 3.15(b)). An \emph{extending 4-wheel operation under $f$} means that for a 2-path $xuy$ of a 4-colorable maximal planar graph, first, conduct an extending 4-wheel operation considering no coloring (see Figure 3.15(b)); second, assign the new center $v$ of the wheel a color differing from $f(x),f(u),f(y)$.

\begin{center}
        \includegraphics [width=380pt]{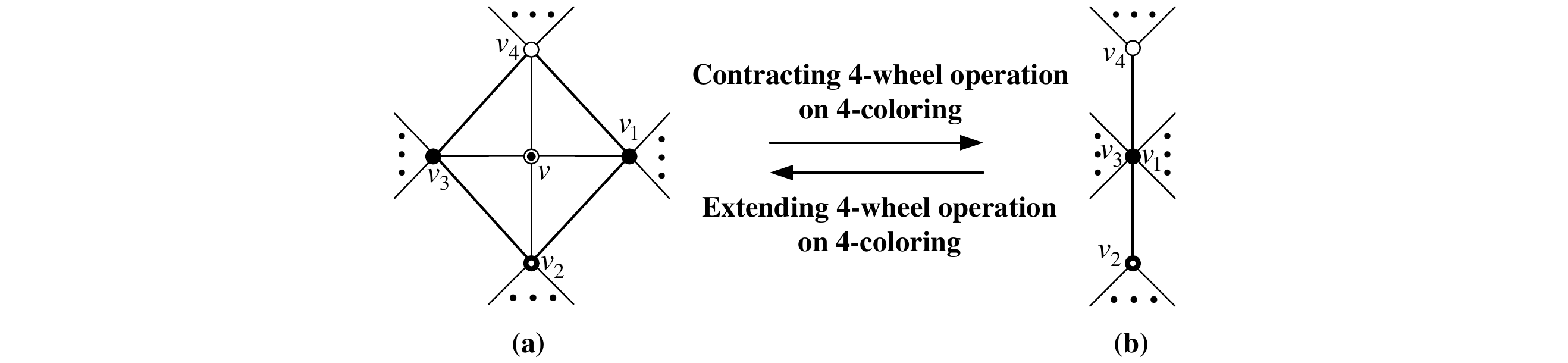}

        \textbf{Figure 3.15.} Schematic diagram of contracting 4-wheel and extending 4-wheel operations based on coloring
\end{center}

	Let $v$ be a 5-degree vertex of a 4-colorable maximal planar graph $G$, and $\Gamma(v)=\{v_1,v_2,v_3,v_4,v_5\}$. When $f\in C^{0}_{4}(G)$, without loss of generality, let $f(v_1)=f(v_3)$, $f(v_2)=f(v_5)$, shown in Figure 3.16(a). Then a \emph{contracting 5-wheel operation under $f$} is to delete vertex $v$ from graph $G$, and identify $\{v_2,v_5\}$ or $\{v_1,v_3\}$. Here we identify vertices $\{v_2,v_5\}$ to describe the process (see Figure 3.16(b)). For a 4-colorable maximal planar graph $G$, the object of an \emph{extending 5-wheel operation under $f$} is a funnel $L$ of $G$, in which both the top and one of bottoms have the same color under $f$ (here let $f(v_1)=f(v_3)$), see Figure 3.16(b). The specific steps of the \emph{extending 5-wheel operation under $f$} is shown in the following. First, conduct the extending 5-wheel operation considering no coloring (see Figures 3.4 and 3.16(b)); second, assign to the new center $v$ of the wheel a color differing from $f(v_1)=f(v_3),f(v_2),f(v_4)$ (see Figure 3.16(a)).

\begin{center}
        \includegraphics [width=380pt]{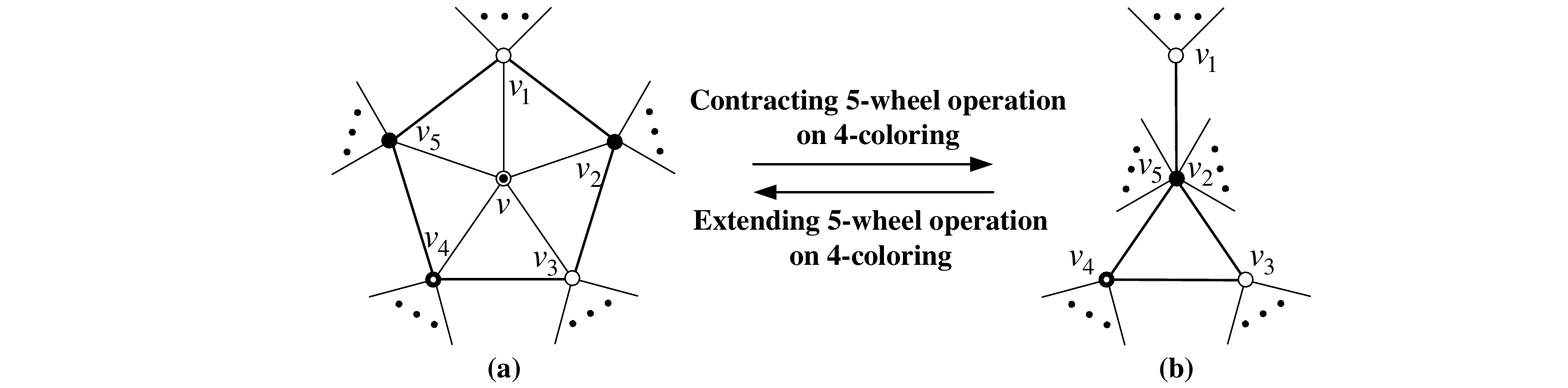}

        \textbf{Figure 3.16.} A diagram for contracting 5-wheel and extending 5-wheel operations
\end{center}

	The above discussion describes the detail of contracting $k$-wheel and  extending $k$-wheel operations of a 4-colorable maximal planar graph $G$ ($2\leq k\leq 5$) based on coloring.
 It is easy to see that, when $2\leq k \leq 5$, the extending $k$-wheel operation and the contracting $k$-wheel operation are one-one correspondence. But for $k\geq 6$, the similar extending $k$-wheel operation and contracting $k$-wheel operation are not one to one, may be one to many. In the following, we first introduce the definition of  extending 6-wheel operation and contracting 6-wheel operation of a 4-colorable maximal planar graph $G$ based on coloring, and then define the extending $k$-wheel operation and contracting $k$-wheel operation of $G$ based on coloring.

\begin{center}
        \includegraphics [width=360pt]{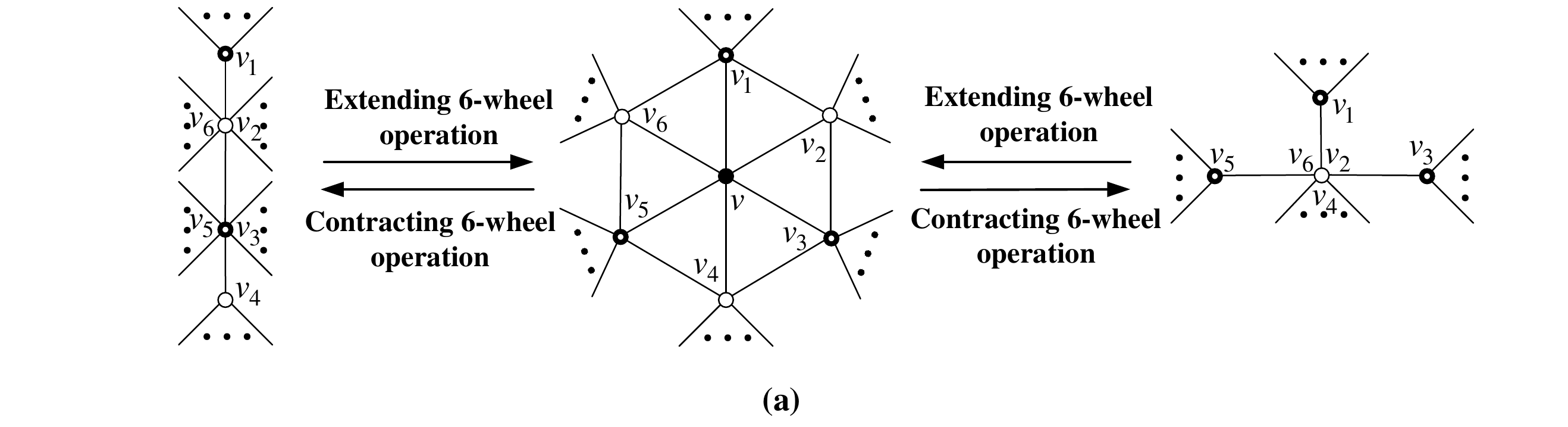}

        \includegraphics [width=360pt]{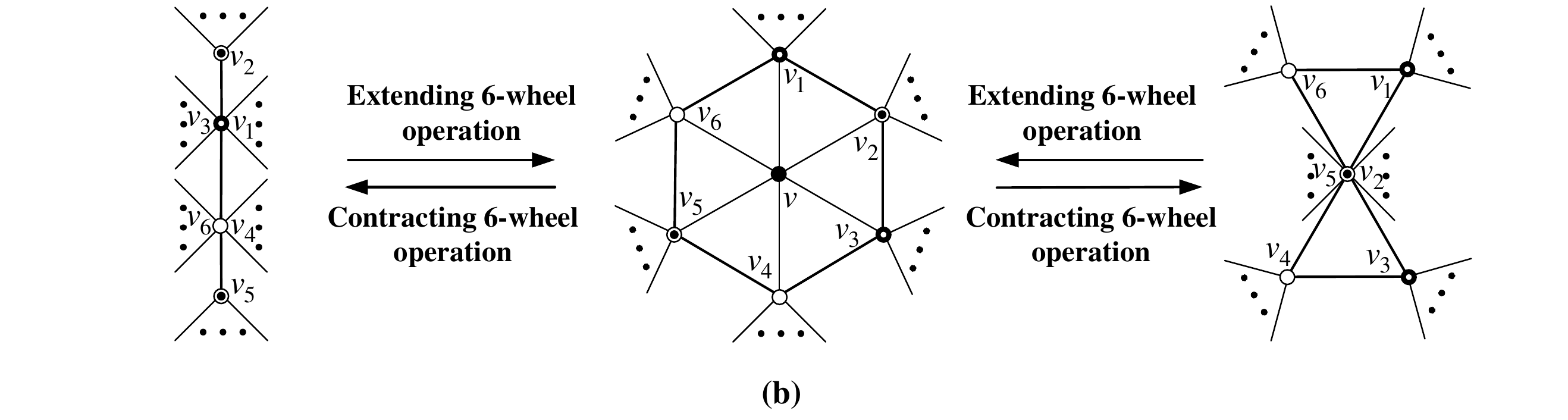}

        \includegraphics [width=360pt]{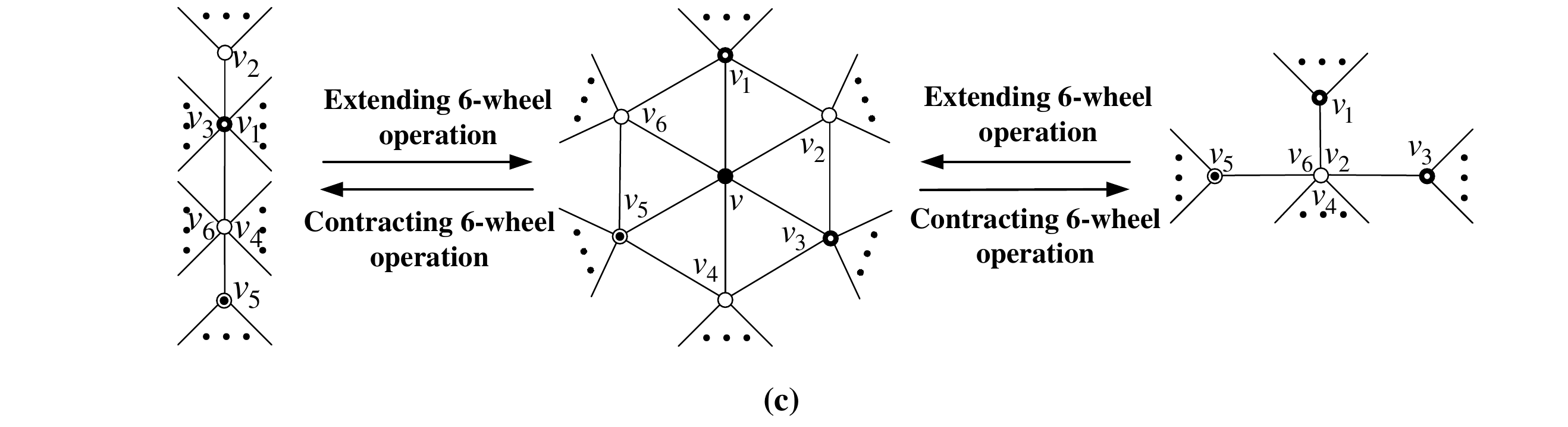}

        \includegraphics [width=360pt]{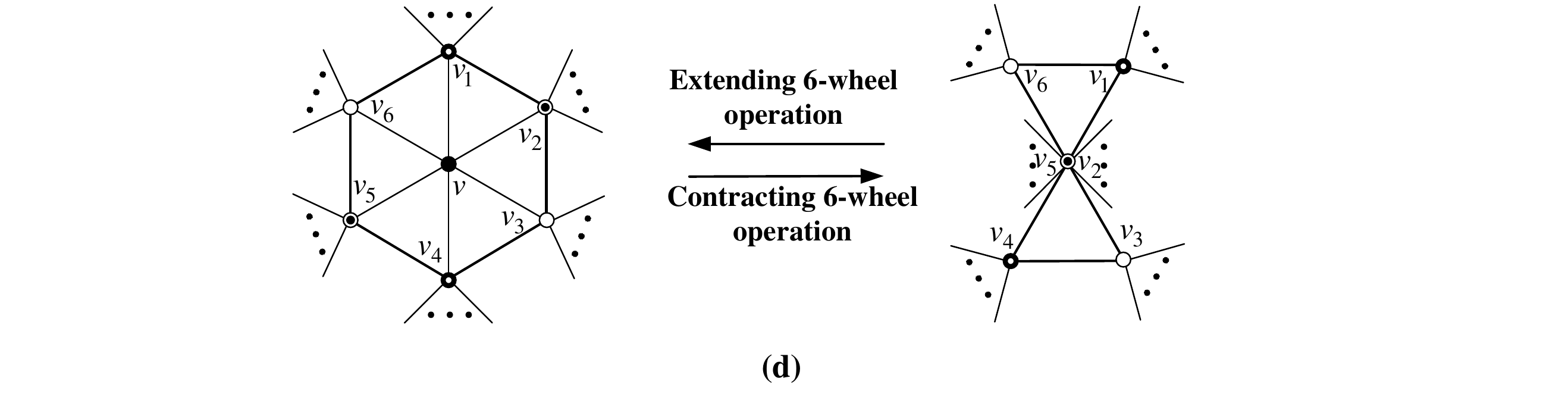}

        \textbf{Figure 3.17.} Diagrams for extending 6-wheel and contracting 6-wheel operations
\end{center}

\begin{definition}
	Let  $v$ be a 6-degree vertex of a 4-colorable maximal planar graph $G$, and $\Gamma(v)=\{v_1,v_2,v_3,v_4,$ $v_5,v_6\}$. When $f\in C^{0}_{4}(G)$, for the 6-wheel $W_6=G[\Gamma(v)\cup \{v\}]$, there might be four kinds of colorings, shown in Figure 3.17.  For every kind of colorings, the definition of the relevant extending 6-wheel operation and contracting 6-wheel operation is shown in the graphs shown in Figure 3.17.
\end{definition}

\begin{definition}
	Let  a maximal planar graph $G$ be 4-colorable, and  $f\in C_4^0(G)$. Suppose that $v$ is a vertex with degree $k\geq 3$ of $G$, and $\Gamma(v)=\{v_1,v_2,\cdots,v_k\}$. The so-called \emph{contracting $k$-wheel operation under $f$} on the wheel $W_k=G[\Gamma(v)\cup \{v\}]$ is to delete vertex $v$, and to identify its neighbors received the same colors. If the resulting graph is a maximal planar graph, then the contracting $k$-wheel operation is completed; otherwise, the new formed face has degree (the number of vertices incident to the face) at least 4, and  at least two vertices receiving the same color, then we identify these vertices again.  We should do this process repeatedly until the last  resulting graph is a maximal planar graph. Clearly, for this resulting graph, if we conduct the corresponding inverse operations of the contracting $k$-wheel operations step by step, $G$ will be obtained again.
We call the above process an \emph{extending $k$-wheel operation based on coloring}.
\end{definition}

	We  use  $\zeta^{-}_{k}(G,v)$ to denote the graphs obtained by contracting $k$-wheel operation on a $k$-degree vertex $v$ of a maximal planar graph $G$, and use $\zeta^{+}_{k}(G,v)$ to denote the graphs obtained by  extending $k$-wheel operation of a maximal planar graph $G$, where $v$ denotes the new  center of the $k$-wheel. If no confusion, they can be abbreviated as $\zeta^{-}_{v}(G)$ and $\zeta^{+}_{v}(G)$ respectively.

Finally, a foundational result will be given in terms of the operations of contracting $k$-wheel and extending $k$-wheel.

For a maximal planar graph $G$ on $n$ $(\geq 5)$ vertices,  there must be at least one vertex of degree three, or four, or five. When $G$ is a 4-colorable maximal planar graph, it is possible to obtain a maximal planar graph of order $n-1$, $n-2$ or $n-3$ from $G$ by doing the (compound) operations of contracting 3-wheel, 4-wheel, 5-wheel, 5-wheel and 3-wheel, 5-wheel and 2-wheel, or 4-wheel and 2-wheel. It means that any 4-colorable maximal planar graph can be obtained by conducting the (compound) operations of extending 3-wheel, 4-wheel, 5-wheel, 5-wheel and 3-wheel, 5-wheel and 2-wheel or 4-wheel and 2-wheel.

\begin{theorem2}\label{theorem3.11}
Any given 4-colorable maximal planar graph $G$ of order $n$ $(n\geq 8)$ can be obtained from some graphs of order $n-1$, $n-2$, or $n-3$ by conducting the $($compound$)$ operations of extending 3-wheel, 4-wheel, 5-wheel, 3-wheel and 5-wheel, 2-wheel and 5-wheel or 2-wheel and 4-wheel.
\end{theorem2}

\subsection{Compound extending and contracting operation based on coloring}

The last section in this chapter will describe the basic extending and contracting operations under the coloring constraint. Being analogous to the compound extending and contracting operations without the coloring constraint, the corresponding operations under colorings will be analyzed in this section.

\textcircled{1}  \emph{An unavoidable-complete set of contractible subgraphs under 4-coloring}. Similarly with the case of no coloring, we only give a kind of 4-colorings for the unavoidable-complete set of contractible subgraphs under no coloring (see Figure 3.5). Figure 3.18 presents  fourteen (class) contractible subgraphs with a 4-coloring.

\begin{center}
        \includegraphics [width=350pt]{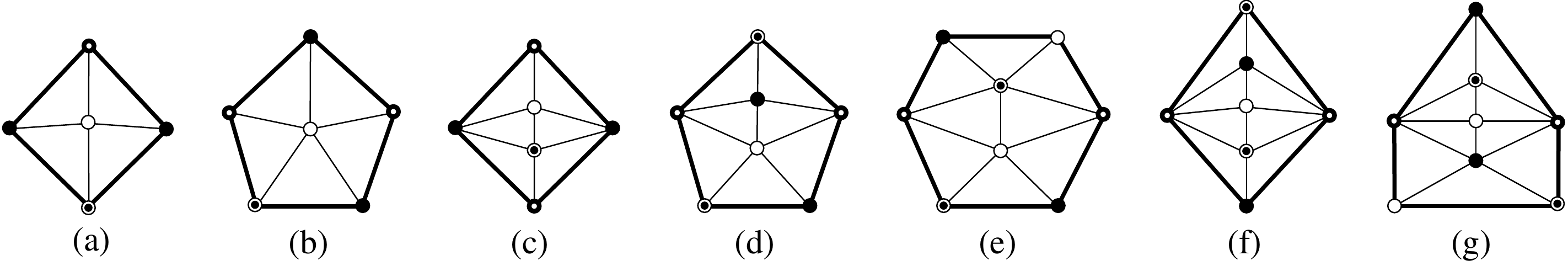}

        \includegraphics [width=350pt]{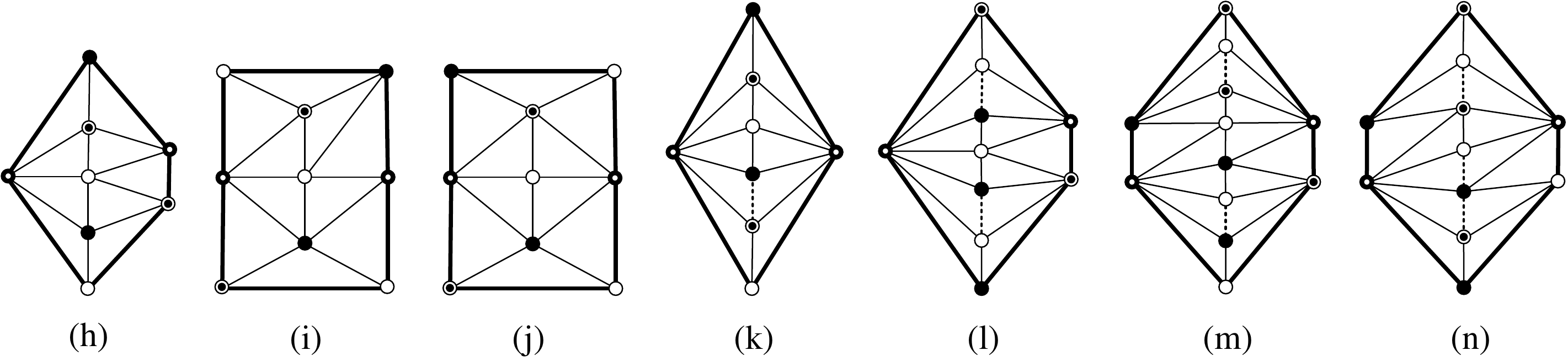}

        \textbf{Figure 3.18.} An unavoidable-complete set of contractible subgraphs under a coloring
\end{center}

\textcircled{2} \emph{The contractible subgraphs based on coloring and corresponding structure analysis}. Table 3.3 not only analyzes the structure of the contractible subgraphs, but also gives a proof of  completeness of the unavoidable set shown in Figure 3.18.\\
\noindent\begin{tabular}{|l|l|l|}
  \hline
  ~~\textbf{Table 3.3.} Structure of contractible subgraphs based on a 4-coloring\\
  \hline
  (a) contracting 4-wheel operation, and the corresponding contracted\\
   subgraph \\
  \hline
  \includegraphics [width=360pt]{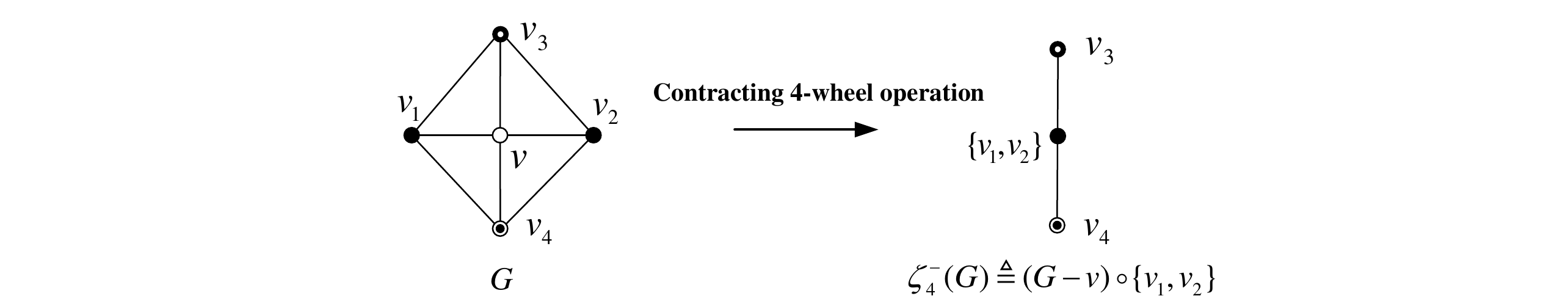} \\
  \hline
\end{tabular}
\begin{tabular}{|l|l|}
  \hline
   (b) contracting 5-wheel operation, and the corresponding contracted\\
    subgraph\\
  \hline
  \includegraphics [width=360pt]{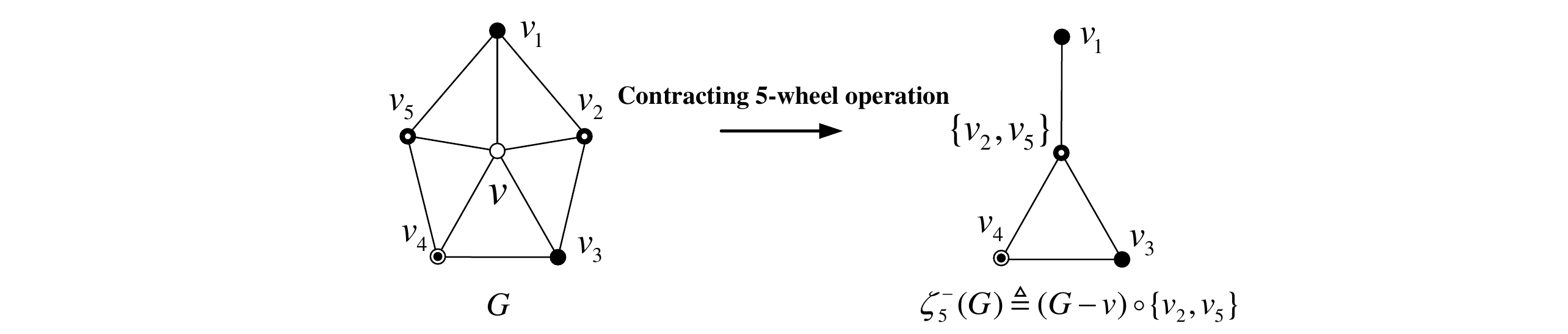} \\
  \hline
\end{tabular}
\begin{tabular}{|l|l|}
  \hline
   (c) contracting 4-wheel and 2-wheel operations, and the corresponding\\
    contracted subgraph\\
  \hline
  \includegraphics [width=360pt]{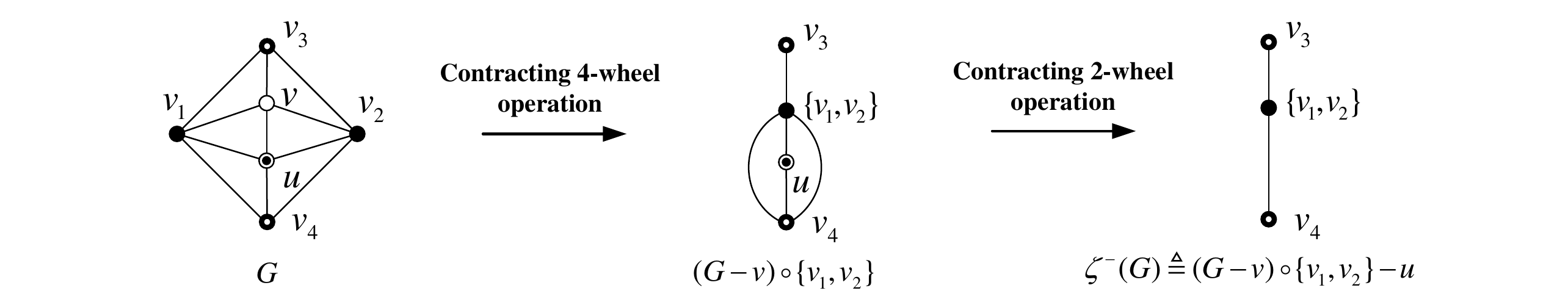} \\
  \hline
\end{tabular}
\begin{tabular}{|l|l|}
  \hline
   (d) contracting 4-wheel and 3-wheel (= 5-wheel and 2-wheel) operations, \\
 and the corresponding contracted subgraph\\
  \hline
  \includegraphics [width=360pt]{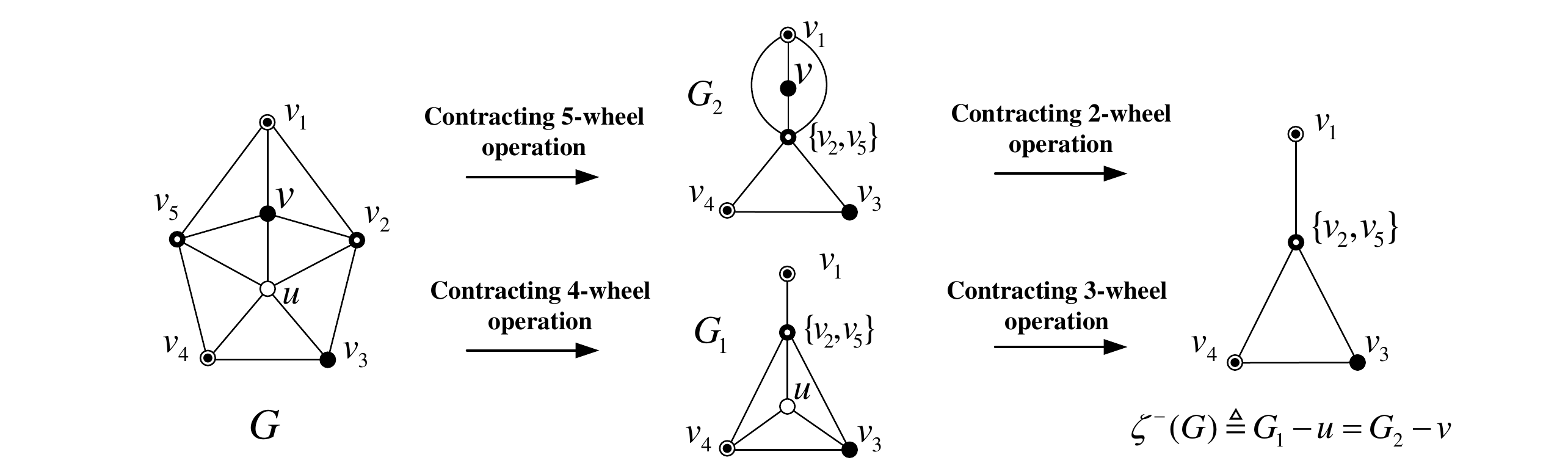} \\
  \hline
  \end{tabular}
  \begin{tabular}{|l|l|}
  \hline
  (e) contracting 5-wheel and 3-wheel operations, and the corresponding\\
   contracted subgraph\\
  \hline
  \includegraphics [width=360pt]{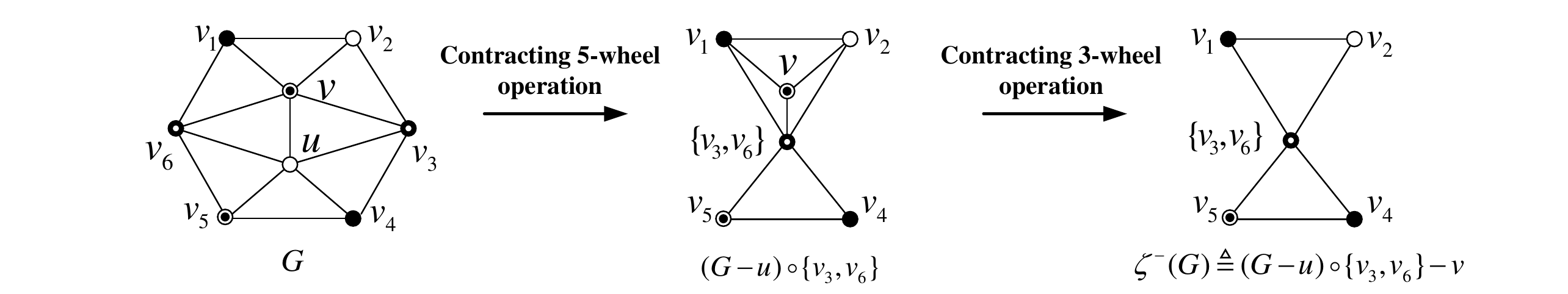} \\
  \hline
  \end{tabular}
  \begin{tabular}{|l|l|}
  \hline
  (f) contracting 4-wheel, 2-wheel and 2-wheel operations, and the \\
  corresponding contracted subgraph\\
  \hline
  \includegraphics [width=360pt]{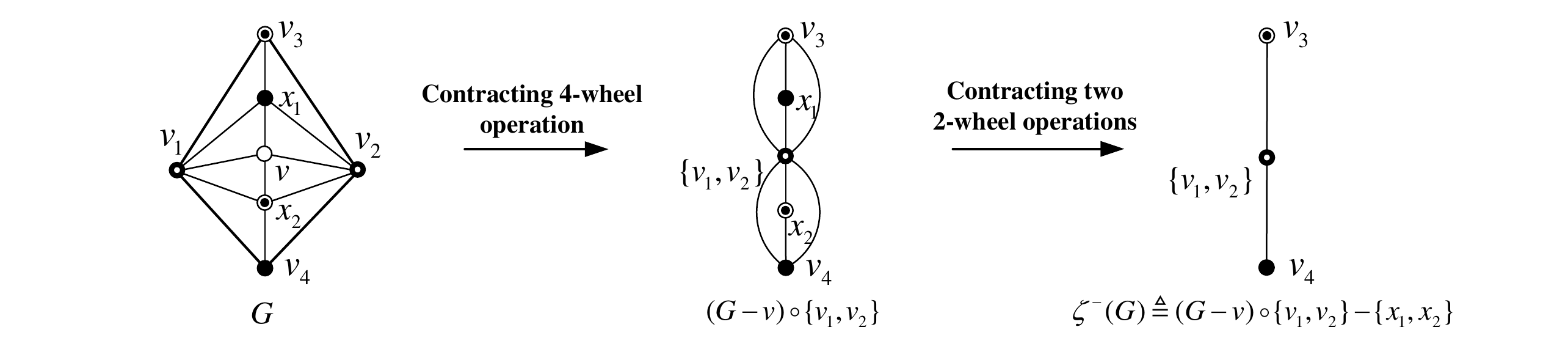} \\
  \hline
  \end{tabular}
  \begin{tabular}{|l|l|}
  \hline
  (g) contracting 4-wheel, 2-wheel and 3-wheel (=5-wheel, 2-wheel and \\
  2-wheel) operations, and the corresponding contracted subgraph\\
  \hline
  \includegraphics [width=360pt]{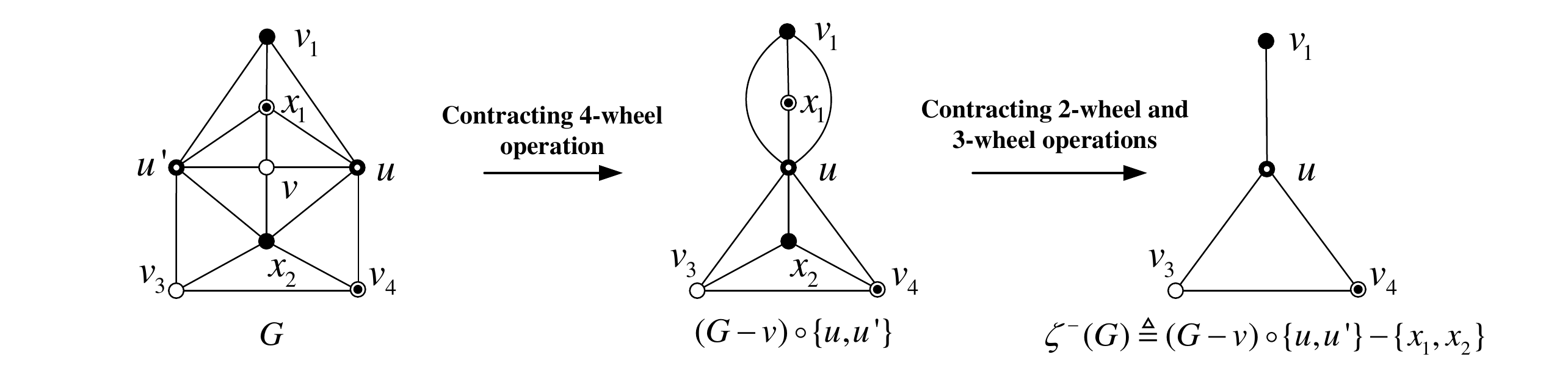} \\
  \hline
  \end{tabular}
\begin{tabular}{|l|l|}
  \hline
   (h) contracting 5-wheel, 2-wheel and 3-wheel (= 4-wheel, 3-wheel and \\
  4-wheel)  operations, and the corresponding contracted subgraph\\
  \hline
  \includegraphics [width=360pt]{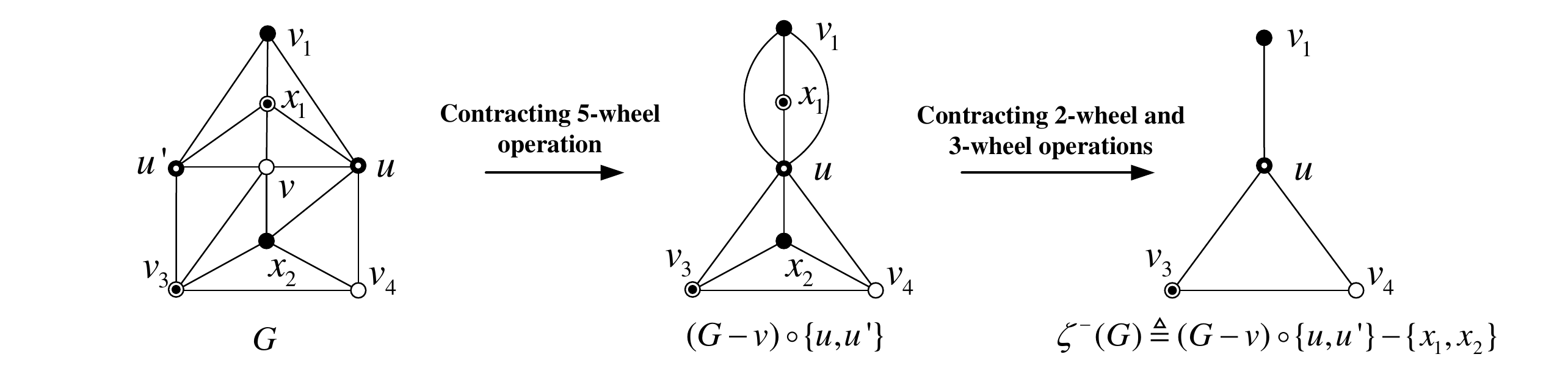} \\
  \hline
  \end{tabular}
  \begin{tabular}{|l|l|}
   \hline
   (i) contracting 5-wheel, 3-wheel and 3-wheel operations, and the\\
    corresponding contracted subgraph\\
  \hline
  \includegraphics [width=360pt]{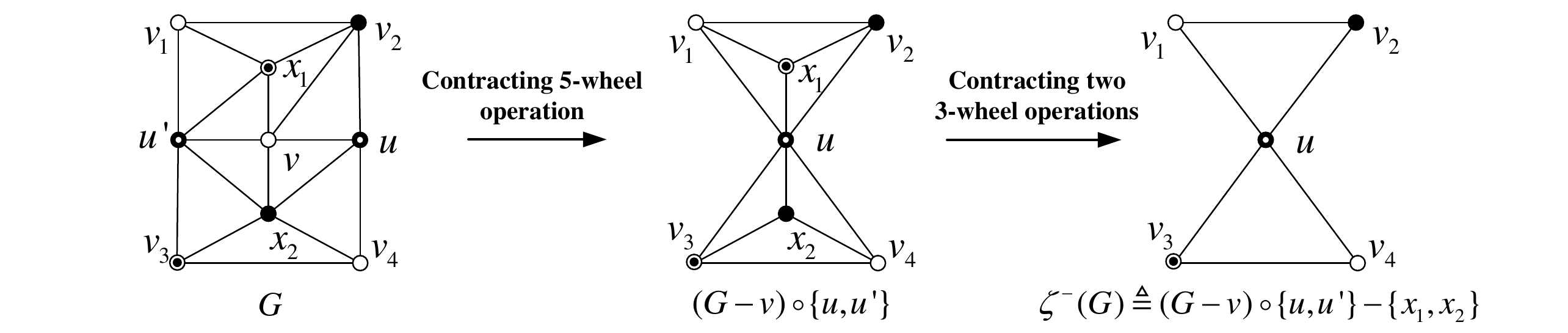} \\
  \hline
  \end{tabular}
  \begin{tabular}{|l|l|}
  \hline
    (j) contracting 4-wheel, 3-wheel and 3-wheel (= 5-wheel, 2-wheel and\\
     3-wheel) operations, and the corresponding contracted subgraph\\
     (dumbbell transformation)\\
  \hline
  \includegraphics [width=360pt]{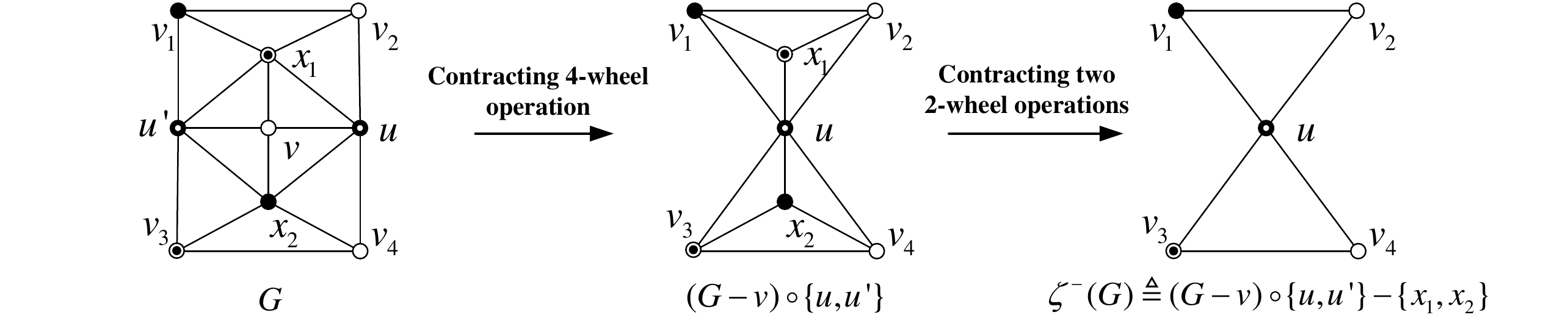} \\
  \hline
  \end{tabular}
  \begin{tabular}{|l|l|}
  \hline
   (k) contracting 4-wheel and 2-wheel ($t$-times) operations, and the\\
    corresponding contracted subgraph\\
  \hline
  \includegraphics [width=360pt]{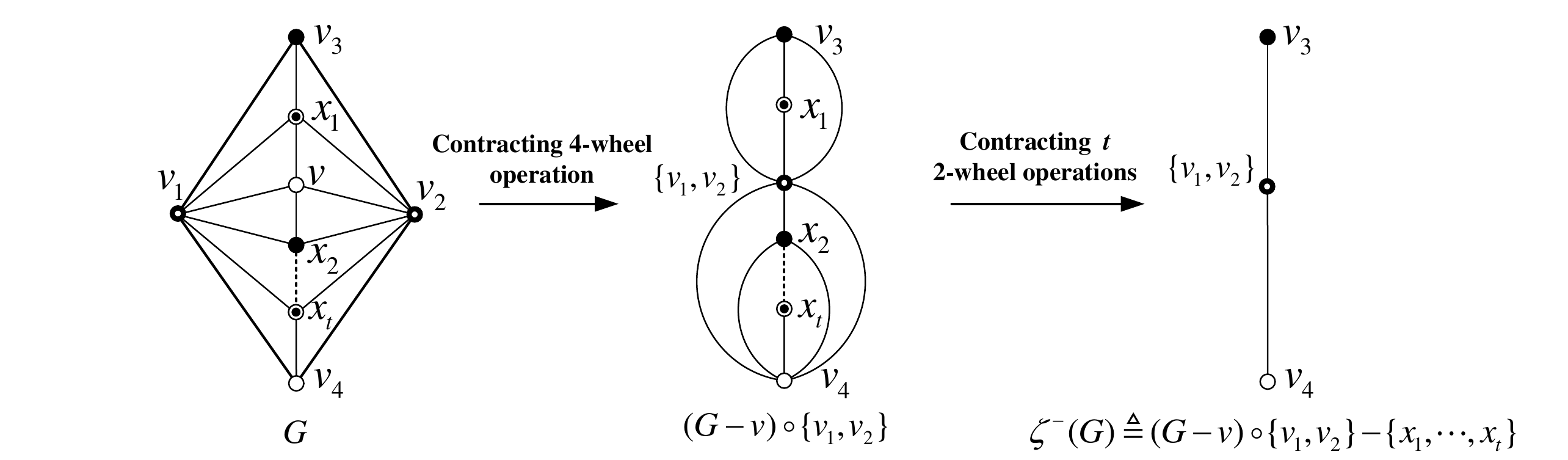} \\
  \hline
  \end{tabular}
  \begin{tabular}{|l|l|}
  \hline
   (l) contracting 5-wheel, 2-wheel ($r$-times) and 3-wheel operations\\
 ($=$contracting 4-wheel ($r-1$  times), 2-wheel ($(t-r+1)$-times) and \\
  3-wheel), and the corresponding contracted subgraph\\
  \hline
  \includegraphics [width=360pt]{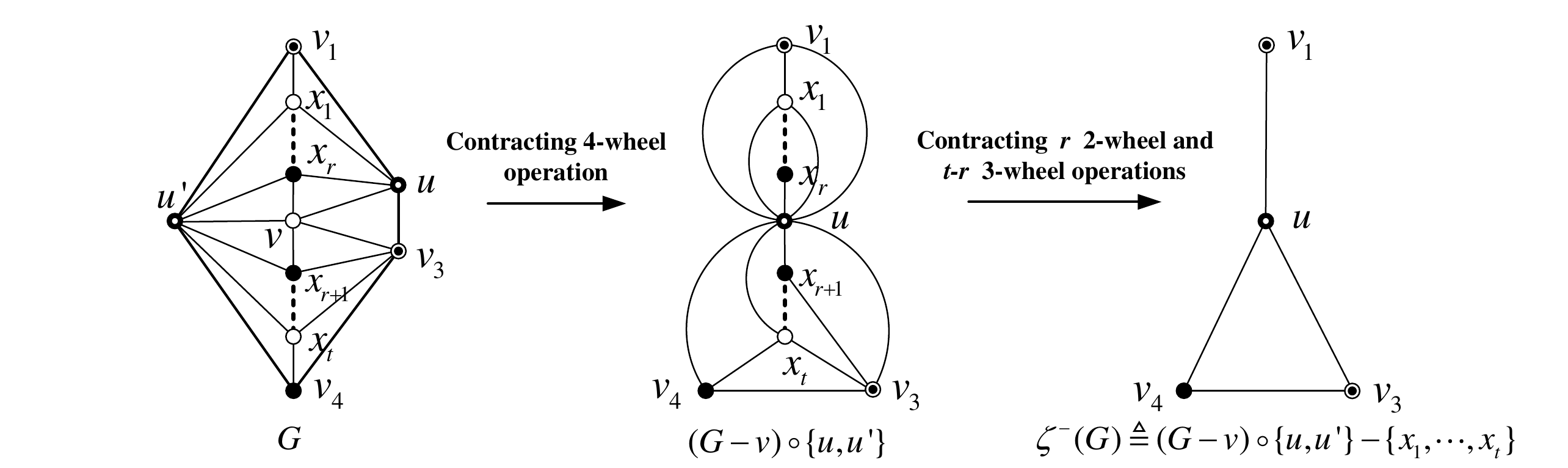} \\
  \hline
  \end{tabular}
  \begin{tabular}{|l|l|}
  \hline
   (m) contracting 5-wheel and 3-wheel ($t$-times) operations,  and the\\
   corresponding contracted subgraph (two 5-degree vertices are adjacent)\\
  \hline
  \includegraphics [width=360pt]{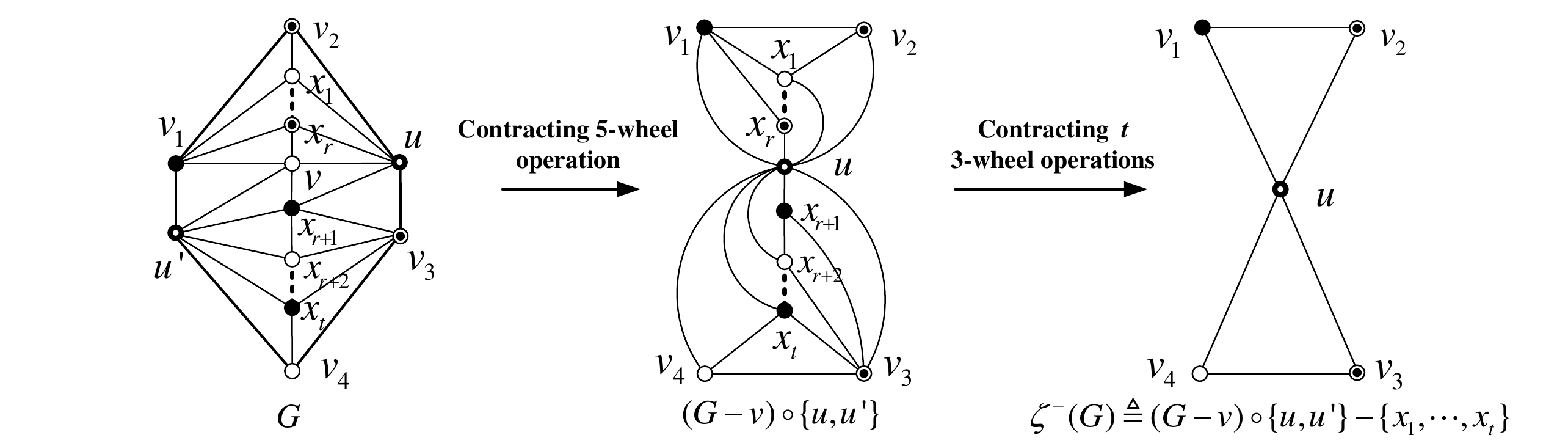} \\
  \hline
  \end{tabular}
  \begin{tabular}{|l|l|}
  \hline
   (n) contracting 5-wheel and 3-wheel ($t$-times) operations,  and the \\
  corresponding contracted subgraph (two 5-degree vertices are nonadjacent)\\
  \hline
  \includegraphics [width=360pt]{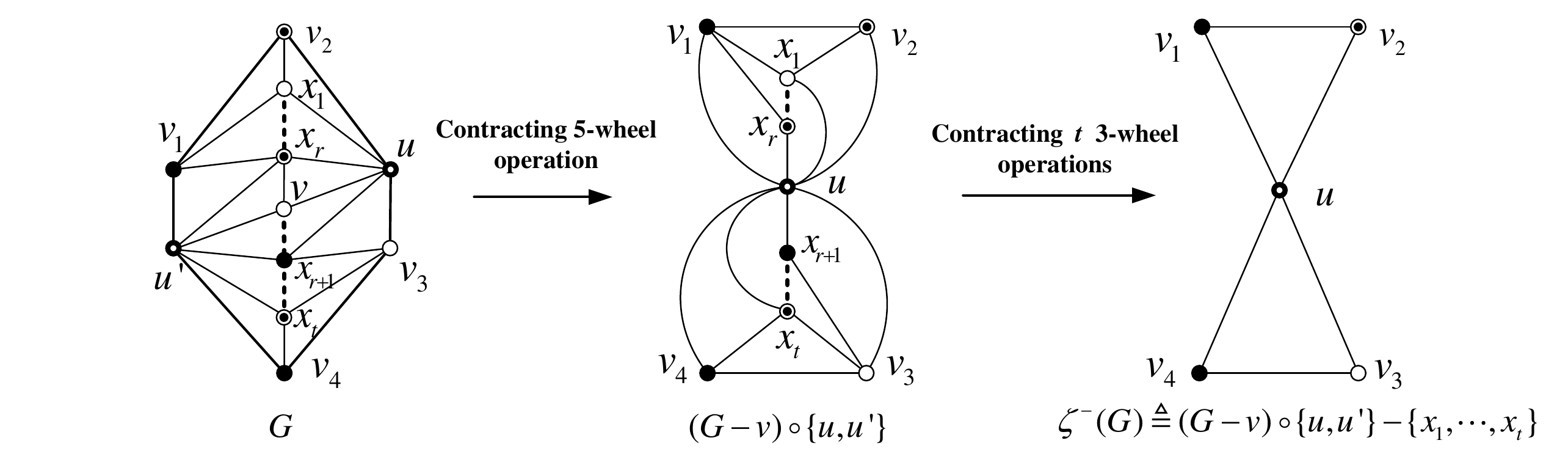} \\
  \hline
  \end{tabular}

\subsection{Parents and children of maximal planar graphs}

On the basis of the foregoing discussion, this section will mainly investigate two interesting problems of maximal planar graphs. First, for a given maximal planar graph $G$, we want to know where $G$ comes from, specifically, by starting with what graphs the graph $G$ can be generated through a sequence of extending $k$-wheel operations. Second, we desire to know how many nonisomorphic maximal planar graphs can be induced by implementing extending $k$-wheel operations to $G$. For this, we put forward the concepts of \emph{parents} and \emph{children} as follows.

For a maximal planar graph $G$ with $\delta\geq 4$, we have known that $G$ can be generated from some such maximal planar graph $\zeta^-(G)$ of lower order by extending wheel operations. In addition, if we conduct the extending wheel operations on $G$, we can also obtain some $\delta \geq 4$ maximal planar graph $\zeta^+(G)$ of higher order. We refer to $\zeta^-(G)$ as a \emph{parent} of $G$, and  $\zeta^+(G)$ as a \emph{child} of $G$. Likewise, $G$ is a parent of $\zeta^+(G)$ and a child of $\zeta^-(G)$.

If not stated otherwise, the term \emph{maximal planar graphs} will be used in the following to represent the graphs with minimum degree $\delta \geq 4$.

In general, an $n (\geq 10)$-vertex maximal planar graph $G$ may have many different parents, written $\zeta^{-1}(G),\zeta^{-2}(G),\cdots, \zeta^{-m}(G)$. We refer to the set of these $m$ parents to the \emph{parent set} of $G$, denoted by $\Upsilon^-(G)$. Namely
\[
\Upsilon^-(G)=\{\zeta^{-1}(G),\zeta^{-2}(G),\cdots, \zeta^{-m}(G)\} \eqno{(3.3)}\label{3.3}
\]
Of course, there are also many children of $G$ when we implement extending wheel operations to $G$, written $\zeta^{+1}(G),\zeta^{+2}(G),\cdots, \zeta^{+m}(G)$. Analogously, we refer to the set of these $m$ children to the \emph{child set} of $G$, denoted by $\Upsilon^+(G)$. Namely
\[
\Upsilon^+(G)=\{\zeta^{+1}(G),\zeta^{+2}(G),\cdots, \zeta^{+m}(G)\} \eqno{(3.4)}\label{3.4}
\]

It is easy to prove that icosahedron, the second graph in Figure 1.4, has only one parent (see Figure 3.19(a)), and twelve children (see Figures 3.19(b)$\sim$ (m)), where the twelve children are obtained by implementing eleven different extending wheel operations to $G$.
\begin{center}
        \includegraphics [width=360pt]{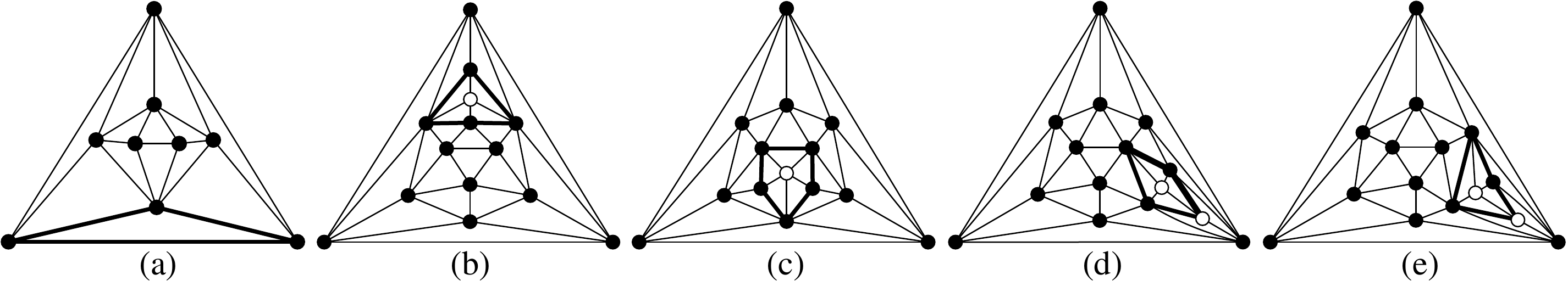}

        \vspace{2mm}
        \includegraphics [width=360pt]{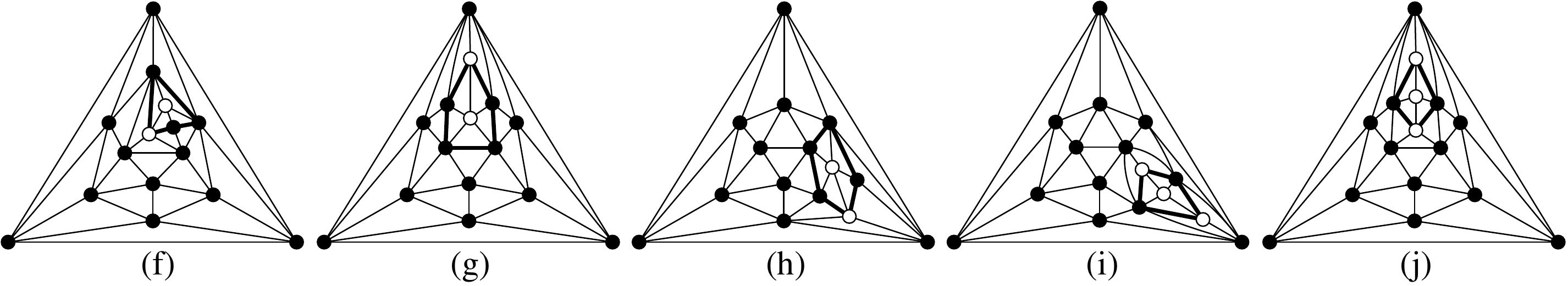}

        \vspace{2mm}
        \includegraphics [width=360pt]{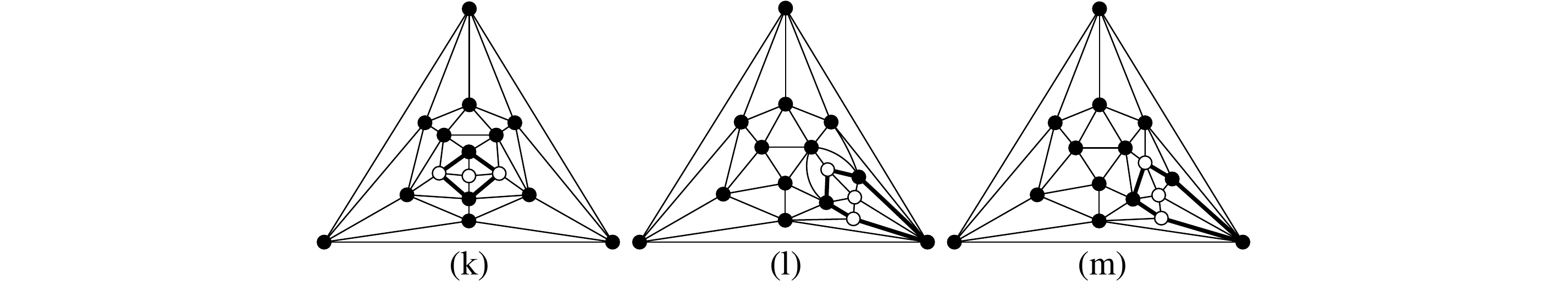}

         \vspace{2mm}
        \textbf{Figure 3.19.} A parent and twelve children of icosahedron
\end{center}

For a given maximal planar graph $G$, here we introduce the concept of equivalent subgraphs to investigate $\Upsilon^-(G)$ and $\Upsilon^+(G)$ clearly. Suppose that $Aut(G)$ is the automorphism group of $G$, and $H,H'$ are two isomorphic subgraphs of $G$. If $\exists \sigma\in Aut(G)$ such that $\sigma(H)=H'$, then $H$ and $H'$ are called to be \emph{equivalent}, otherwise, \emph{nonequivalent}.

We use $\Im^H_G$ to denote the set of all nonequivalent subgraphs $H$ of $G$. For example, $\Im^L_G$ denotes the set of all nonequivalent funnels of $G$;  $\Im^{P_3}_G$ denotes the set of all nonequivalent 2-length paths of $G$; $\Im^{P_2}_G$ denotes the set of all nonequivalent edges of $G$, and so on.

Obviously, $\Im^H_G$ has a close relation with $Aut(G)$. When the topology structure of $G$ has stronger symmetry, the size of $\Im^H_G$ is very small.  For instance, for the icosahedron $G_{20}$, we have
\[
 |\Im^{P_2}_{G_{20}}|=|\Im^{P_3}_{G_{20}}|=|\Im^{K_3}_{G_{20}}|=|\Im^{L}_{G_{20}}|=1\eqno{(3.5)}\label{3.5}
\]
However, when $Aut(G)$ is the unit group, $|\Im^H_G|$ will be larger.  The following theorems are obvious.

\begin{theorem2}\label{theorem3.12}

Suppose that $G$ is a maximal planar graph of $\delta\geq 4$ and $H,H'$ are two isomorphic subgraphs of $G$. If we implement the same extending wheel operation on $H$ and $H'$ respectively, and the resulting graphs are written by $\zeta^+_H(G), \zeta^+_{H'}(G)$, then $\zeta^{+}_H(G)\cong \zeta^+_{H'}(G)$ if and only if $H$  and $H'$ are equivalent.
\end{theorem2}

\begin{theorem2}\label{theorem3.13}
Suppose that $G$ is an $n$ $(\geq 6)$-vertex maximal planar graph with minimum degree $\delta\geq 4$. Then

\textcircled{1} The set of $(n+2)$-vertex maximal planar graphs with $\delta \geq 4$, induced by $G$, is
$$
\{\zeta^+_H(G): H\in \{P_3,L\}, P_3\in \Im^{P_3}_G, L\in \Im^{L}_G\} \eqno{(3.6)}\label{3.6}
$$

\textcircled{2} The set of $(n+3)$-vertex maximal planar graphs with $\delta \geq 4$, induced by $G$, is
$$
\{\zeta^+_H(G): H\in \{P_3,L\}, P_3\in \Im^{P_3}_{\zeta^+_2(G)}, L\in \{\Im^{L}_{\zeta^+_2(G)},\Im^{L}_{\zeta^+_3(G)}\} \} \eqno{(3.7)}\label{3.7}
$$
In Form (3.7), one of the ends of $P_3$ has degree 2 in graph $\zeta^+_2(G)$ when we implement the extending 2-wheel 4-wheel operation; the top of the funnel $L$ has degree 2 in graph $\zeta^+_2(G)$ when we implement the extending 2-wheel 5-wheel operation; the top or one of bottoms of the funnel $L$ has degree 3 in graph $\zeta^+_3(G)$ when we implement extending 3-wheel 5-wheel operation.

\textcircled{3} The set of $(n+4)$-vertex maximal planar graphs with $\delta \geq 4$, induced by $G$, is
$$
\{\zeta^+_H(G): H\in \{P_3,L\},
$$
$$
P_3\in \{\Im^{P_3}_{\zeta^+_{2-2}(G)},\Im^{P_3}_{\zeta^+_{2-3}(G)},\Im^{P_3}_{\zeta^+_{3-3}(G)}\},
L\in \{\Im^{L}_{\zeta^+_{2-3}(G)},\Im^{L}_{\zeta^+_{3-3}(G)}\}\} \eqno{(3.8)}\label{3.8}
$$
In Form (3.8), the degrees of two ends of the path $P_3$ are 2 in graph $\zeta^+_{2-2}(G)$, 2 and 3 respectively in  $\zeta^+_{2-3}(G)$, and 3 in  $\zeta^+_{3-3}(G)$;  the degrees of the top and  one of bottoms of the funnel $L$ are 2 and 3 respectively in $\zeta^+_{2-3}(G)$, and 3 in $\zeta^+_{3-3}(G)$.
\end{theorem2}

 For an $n$ $(\geq 6)$-vertex maximal planar graph $G$ with minimum degree $\delta\geq 4$, we can understand by means of Theorem 3.13 that the numbers of $(n+2)$-vertex, $(n+3)$-vertex, and $(n+4)$-vertex  maximal planar graphs in $\Upsilon^+(G)$ are defined by the number of the nonequivalent subgraphs of $G$. Particularly, the number of $(n+2)$-vertex maximal planar graphs in $\Upsilon^+(G)$ is equal to $|\Im_G^{P_3}|+|\Im_G^{L}|$.

\section{Recursive maximal planar graphs}

Recall that the concept of recursive maximal planar graphs has been introduced in section 2: they can be obtained from $K_4$, embedding a 3-degree vertex in a triangular face continuously. The set of recursive maximal planar graphs is denoted by $\Lambda$, and the set of recursive maximal planar graphs with order $n$ is denoted by $\Lambda_n$, write $\gamma_n=|\Lambda_n|$. Based on section 4, an exact definition for recursive maximal planar graphs is defined as follow: conducting extending 3-wheel operation continuously from $K_3$ or $K_4$, namely $\Lambda=(K_3,\zeta_3^+)$.

JT Conjecture states that a 4-colorable maximal planar graph is uniquely 4-colorable if and only if it is a recursive maximal planar graph. So, the foundation to attack this conjecture is to further study recursive maximal planar graphs. This kind of graphs is also
 called the FWF graphs. In the process of study, one class of graphs called the (2,2)-FWF
 graphs is actually the main class of recursive maximal planar graphs, which is indispensable in the proof of JT Conjecture.  Below, we give some related
 properties of  FWF graphs, especially for the (2,2)-FWF graphs.

 \subsection{Basic properties}
   \begin{theorem2}\label{th4.1}
   If $G$ is a FWF graph of order $n$, then it has at least two vertices of degree three.
 And when $n\geq 5$, any two vertices of degree three are not adjacent to each other.
   \end{theorem2}

   \begin{proof}
     By induction on the number of vertices. When $n=4,5,6$, $\gamma_{4}=\gamma_{5}=\gamma_{6}=1$,
 and the corresponding graphs are shown in Figure 2.2.
 So the result is true obviously.

 Assume that the theorem holds when the number of vertices is $n$. That is, for any FWF graph $G$
 with $n$ vertices, it has at least two $3$-degree vertices, and all the vertices of $3$-degree
 are not adjacent to each other.

     A graph $G \in \Lambda_{n+1}$ is constructed by adding a $3$-degree vertex $v$
 in any triangular face of a FWF graph with $n$ vertices, assuming $\Gamma_{G}(v)=\{v_{1},v_{2},v_{3}\}$. By induction, there are at least two $3$-degree vertices, and all the vertices of $3$-degree
 are not adjacent to each other in $G-v$. For $G-v$, if 3-degree vertices are included in $\{v_1,v_2,v_3\}$, then one exists at most, saying $v_1$. Obviously, there is at least another 3-degree vertex except for $v_1$ in $G-v$, and all those $3$-degree vertices
 are not adjacent to each other. Since $v$ is a 3-degree vertex of $G$ and it is not adjacent to any other vertices except for $\{v_1,v_2,v_3\}$, so  there are also at least two 3-degree vertices in $G$, and all those $3$-degree vertices
 are not adjacent to each other. Thus, the conclusion holds. For $G-v$, if 3-degree vertices are not in $\{v_1,v_2,v_3\}$, the conclusion holds by the same way above. The theorem follows by the principle of induction.
 \end{proof}

 \begin{theorem2}\label{th4.2}

 $(1)$ There exists no maximal planar graph having exactly two adjacent vertices of degree 3.

 $(2)$ There exists no maximal planar graph having exactly three vertices of degree three such that one of them is adjacent to the rest.
 \end{theorem2}
 \begin{proof}
   By contradiction.
   Assume that $G$ is a maximal planar graph with two adjacent vertices $u,v\in V(G)$ exactly, satisfying $d(u)=d(v)=3$.
 Since $u$ is also a 3-degree vertex, $\Gamma(u)=\{v,x,y\}$. Notice that $G$ is a maximal planar
 graph and  $u$ must be in a triangular face which consists of the vertices $v$, $x$ and $y$. In other words,
 $v$ is adjacent to vertices $x$ and $y$. These four vertices can form a subgraph $K_{4}$ (shown in Figure 4.1).
 Since $G$ is a maximal planar graph, if there exist any other vertices of $G$, then it can
 form a triangle with $u$ or $v$. It contradicts  $d(u)=d(v)=3$. Otherwise,
 if there exists no other vertices, then $G$ is isomorphism to $K_{4}$ with
 four vertices of 3-degree. Therefore, there exists no maximal planar graph with two adjacent vertices
 of 3-degree exactly.
 \begin{center}
 \includegraphics[width=130pt]{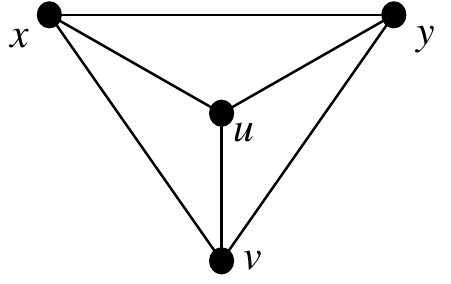}

 \textbf{Figure 4.1.} The schematic for the proof of Theorem 4.2
 \end{center}

    To consider the assertion (2), we assume that $u,v,x \in V(G)$, $d(u)=d(v)=d(x)=3$, and
 $uv, ux, vx \in E(G)$. There exist three vertices in $\Gamma(u)$, two of which are $v$ and $x$. Let $y$ denote the other vertices adjacent to $u$. So $\Gamma(u)=\{v,x,y\}$. Since a face can be constructed with three vertices $v$, $x$
 and $y$, any two of which are adjacent, so they  induce a subgraph $K_{4}$, shown in Figure 4.1. Since $G$ is a maximal
 planar graph, if there exist any other vertices of $G$, we get a triangle with vertices $u$, $v$ or $x$. It contradicts
 the fact that $d(u)=d(v)=3$. Otherwise, $G\cong K_{4}$. Obviously, $K_{4}$ contains four vertices of 3-degree.
 Thus, there exists no maximal planar graph with three vertices of 3-degree, any two of which are adjacent.
 \end{proof}

 \begin{theorem2}\label{th4.3}
     If $G$ is a maximal planar graph having only one vertex of $3$-degree, then a graph without any
 $3$-degree vertex can be obtained from $G$ by deleting $3$-degree vertices repeatedly.
 \end{theorem2}
 \begin{proof}
    Let $v$ be the unique vertex of 3-degree in graph $G$, and $\Gamma_{G}(v)=\{u_{1},u_{2},u_{3}\}$.
 Thus, these three vertices can form a triangle, any two of which are adjacent. $G_{1}$=$G-v$ is also a maximal planar graph. There may exist four cases in $G_1$ as follows:

    (1)  $\delta(G_{1})\geq 4$;

    (2)  There exists only one 3-degree vertex;

    (3)  There exactly exist  two 3-degree vertices;

    (4)  There exactly exist three 3-degree vertices.

    For case (1), the theorem holds naturally. The cases (3) and (4) do not exist by Theorem 4.2.
 So we just need to consider the case (2). In this case, there exists a 3-degree vertex in subgraph $G_{1}$,
 denoted by $v_{1}$. Let $G_{2}=G_{1}-v_{1}$. Like the method mentioned above, if $\delta(G_{2})\geq 4$, then the
 theorem holds.
 Otherwise, the graph $G_{2}$ must contain a 3-degree vertex. In this way, we can get $\delta(G_{m})\geq 4$
 within finite $m$ steps. Otherwise, $G_{m} \cong K_{4}$ when $G_{m}$ contains only four vertices.
 It means that the graph $G$ is a FWF graph. But there is only one 3-degree vertex in $G$.  It contradicts Theorem
 4.1.
 \end{proof}

 \subsection{(2,2)-FWF graphs}
    In this section, we introduce and study (2,2)-FWF graphs,  which are a special class of FWF graphs. A FWF graph is called a \emph{(2,2)-FWF graph} if it contains only two vertices of 3-degree,
 and the distance between them is 2. It is easy to prove that there exist only one $(2,2)$-FWF graph with order 5 and 6 respectively, shown in Figures 4.3(a) and (b).

 To understand the structure of a (2,2)-FWF graph, we can divide three inner faces of the complete graph $K_{4}$
into three regions, and label its vertices correspondingly. In Figure 4.2, the triangle
 is called the \emph{outside triangle} when its vertices are labeled by 1, 2, 3, and the vertex $u$ (also labeled by 4) is called
 the \emph{central vertex}. Here we define that the vertices 1,2,3,4 are colored with yellow, green, blue and red, respectively. The four
 vertices and their corresponding colorings are called the \emph{basic axes} in the \emph{color-coordinate system} of a (2,2)-FWF graph.
 Four color axes are 1 (yellow), 2 (green), 3 (blue) and $u$ (red). Obviously, there exists no (2,2)-FWF graph of order 4; and
 there is only one (2,2)-FWF graph with 5 vertices under the isomorphism of view, which can be obtained by embedding a 3-degree vertex in the region
 $I$, $II$ or $III$ of the graph $K_{4}$(shown in Figure 4.2). Without loss of generality, we make an agreement that new vertices are only added in the region $II$. Thus, the vertex is colored by yellow (Figure 4.3 (a));
  the non-isomorphic (2,2)-FWF graphs of order $6$ can be obtained by embedding a 3-degree vertex in any region of the (2,2)-FWF graph of
 5-order. It is easy to prove that this kind of graphs with $6$ vertices obtained by embedding a new vertex in any face are isomorphic.
 Therefore, the number of (2,2)-FWF graphs of order $6$  is one. In general, we make an agreement that the 6th vertex is embedded in
 the face composed of the vertices 2, 4, 5 (i.e. the sub-region $I$ of the region $II$), which is colored by blue (Figure 4.3 (b)).
 Further, for (2,2)-FWF
 graphs of higher orders, we restrict that new vertices are only added in the regions $I$ and $II$, but not in the region $III$.

 \begin{center}
    \includegraphics[width=240pt]{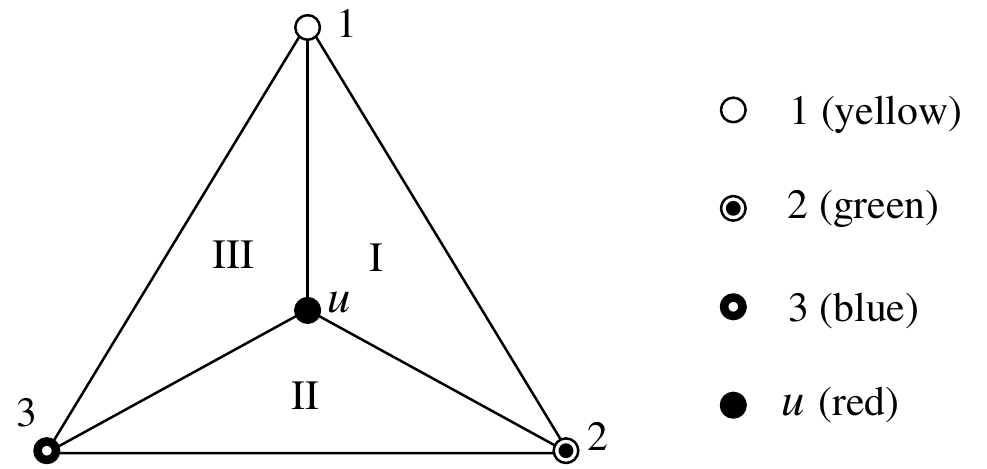}

    \textbf{Figure 4.2.}  The basic framework of the color-coordinate system
  \end{center}

 \begin{center}
   \includegraphics[width=320pt]{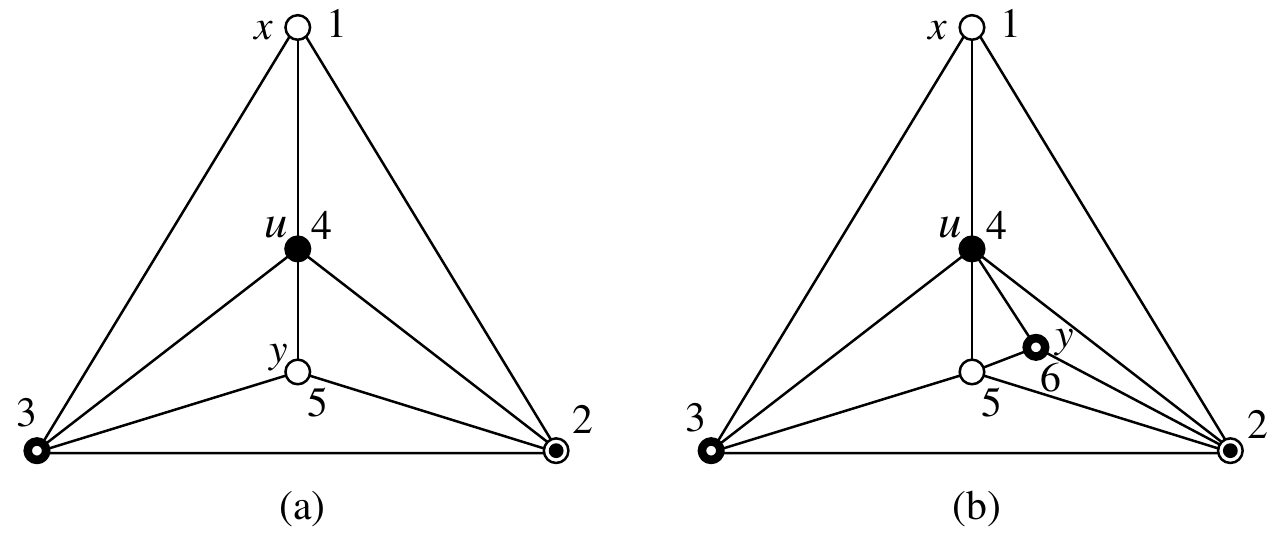}

    \textbf{Figure 4.3.}  Two (2,2)-FWF graphs \\ (a) a graph of order 5; (b) a graph of order 6
  \end{center}

    Based on the agreement above, we can discuss  the classification of (2,2)-FWF graphs.
 Two methods for this purpose will be introduced as follows.

    The first is based on the region where the 3-degree vertices are embedded: (1) A (2,2)-FWF graph is obtained by
    successively
 embedding  3-degree vertices only in the region $II$. Such graphs are shown in Figure 4.4; (2) A (2,2)-FWF graph is  obtained by
 successively and randomly embedding 3-degree vertices in the region $I$ and $II$, shown in Figure 4.5.  we have a straightforward fact as follows.
    \begin{Prop}$^{[56]}$
    Any face in a $($maximal$)$ planar graph can become the infinite outside face.
    \end{Prop}

    That is, the (2,2)-FWF graphs mentioned above are obtained by embedding 3-degree vertices randomly in the region $I$ and $II$.  We can transform
 any one 3-degree vertex in the region $I$ or $II$ to the outside triangular face by Proposition 4.4, which is equivalent to the first classification.
 It means that this kind of (2,2)-FWF graphs are obtained by successively embedding 3-degree vertices only in the region $II$. Therefore, we only consider
 this kind of graphs in the later sections.

  \begin{center}
   \includegraphics[width=380pt]{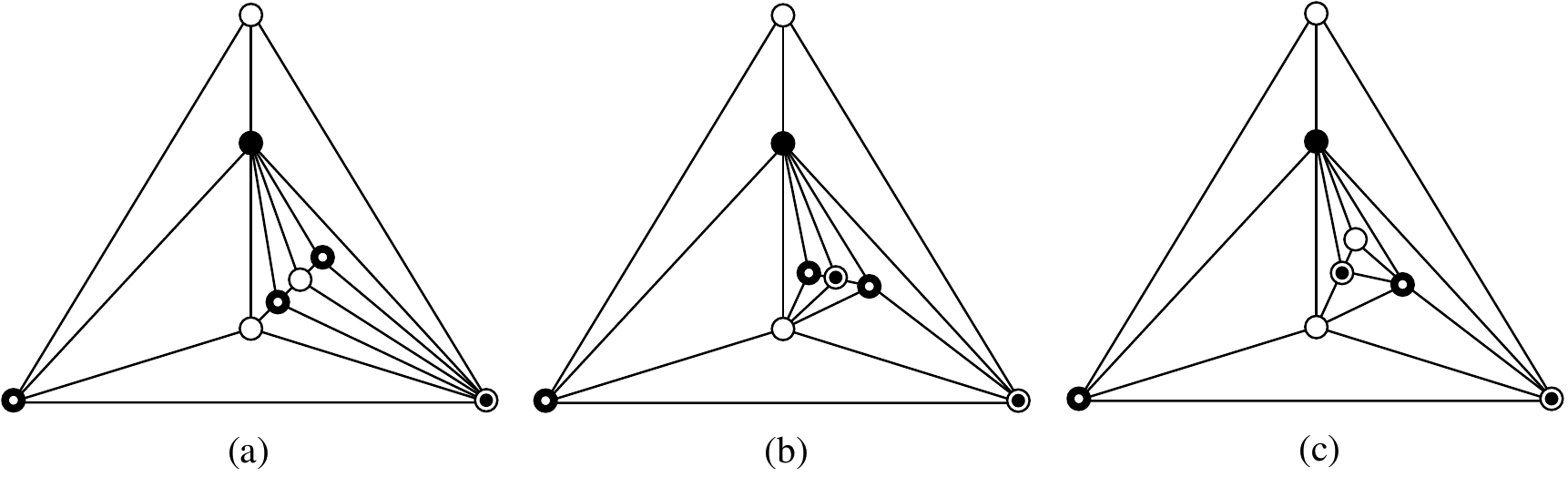}

    \textbf{Figure 4.4.}  The (2,2)-FWF graphs obtained by embedding 3-degree vertices only in the region $II$,
    \\(a) the adjacent type;(b) and (c) the non-adjacent type
   \end{center}
    \begin{center}
    \includegraphics[width=300pt]{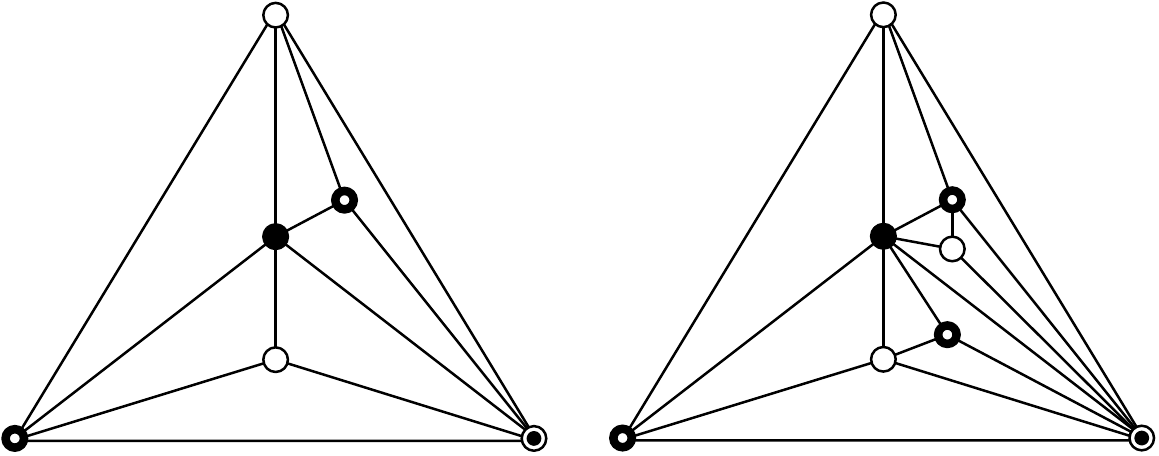}

    \textbf{Figure 4.5.} The (2,2)-FWF graphs obtained by embedding the 3-degree vertices in the region $I$ and $II$ randomly
    \end{center}

    The second method is based on whether there exists a common edge between the two triangular surfaces of two $3$-degree vertices or not.
  It is called the adjacent type if there is a common edge; otherwise, the non-adjacent type. As shown in Figure 4.4, the first graph belongs to the adjacent type, whereas the last two graphs belong to the non-adjacent type.

  From the two classification methods above, all $(2,2)$-FWF graphs can be divided into the adjacent type of region $II$ and  the nonadjacent type of region $II$.

  In the Figure 4.3(a), the $(2,2)$-FWF graph of order 5 is a double-center wheel, and the degree of vertices in the neighbor of each center is 4, where the \emph{double-center wheel} is a specific maximal planar graph constructed by the join of a cycle and two single vertices.  When the order of $(2,2)$-FWF graph is not less than 6, we have the following result.

  \begin{theorem2}\label{th31}
    $(1)$ Let $G$ be a $(2,2)$-FWF graph with order $n$ $(n\geq 6)$, then for each 3-degree vertex $v$ in $G$, there only exists one vertex with order 4 in $\Gamma(v)$;  $(2)$ Every (2,2)-FWF graph $G$ of nonadjacent type with order $n (n\geq 5)$ has one and only one $(n-1)$-degree vertex, and it is called the \emph{central vertex} of the graph $G$,
   denoted by $u$.
    Furthermore, in any partitions of color class in $G$, only the central vertex is colored with red; $(3)$ For
   the (2,2)-FWF graphs of adjacent type obtained by embedding the 3-degree vertices only in the region $II$, only its central vertex is colored with red and also only
   its color axis $2$ is colored with green.
 \end{theorem2}
 \begin{proof}
    By induction. There is only one maximal planar graph of order 5 (shown in Figure 4.3(a)), also a double-center wheel, so all triangular faces are equivalent.  Therefore, in the isomorphism of view, there exist only one FWF graph with order 6, also a $(2,2)$-FWF graph (shown in Figure 4.3(b)). Thus, the theorem holds when $n=6$.

   Assume that the theorem holds when $n$ $(n\geq 6)$. We consider a $(2,2)$-FWF graph $G$ of order $n+1$. Suppose that $v$ is a 3-degree vertex in $G$, there are two cases as follows:

   First, two or three vertices of 4-degree
are included in $\Gamma(v)$, then $G-v$ is also a FWF graph with order at least 5 which contains two or three vertices with 3-degree adjacent to each other, which contradicts with Theorem 4.1.

  Second,  the vertices of degree 4 of $G$
is not included in $\Gamma(v)$, then $G-v$ is also a FWF graph with order at least 5 which contains only one vertex of 3-degree, which contradicts with Theorem 4.1.

In conclusion, we have proved that only one vertex of degree 4 is included in $\Gamma(v)$.

Further, $G$ is a $(2,2)$-FWF graph of order $n+1$, then it exactly contains two vertices of 3-degree and the distance between them is two. Hence, there must be a vertex $u\in V(G)$ making other vertices of $G$ adjacent to $u$,
namely $d(u)=n-1$, which can be proved by the gradual construction of $(2,2)$-FWF graphs. The theorem follows the principle of induction.
   \end{proof}

According to Theorem 4.5, we now define some special triangular faces as follows: for a triangular face containing a vertex of 3-degree, if the degrees of three vertices in this triangular face are 3,4 and $n-1$ respectively, then this triangular face is called an \emph{$I$-type face}; if the degrees of three vertices in this triangular face are 3, $m$ and $n-1$ respectively, then this triangular face is called a \emph{$II$-type face}; if the degrees of three vertices in this triangular face are 3,4 and $m$ respectively, then this triangular face is called a \emph{$III$-type face}; where $5\leq m\leq n-1$.

\begin{theorem2}\label{thnew1}
$\gamma_5=\gamma_6=1,\gamma_7=2, \gamma_8=3, \gamma_9=6.$
\end{theorem2}

The corresponding $(2,2)$-FWF graphs in Theorem 4.6 are shown in Figures 4.3, 4.6, 4.7 and 4.8 respectively.

\begin{center}
    \includegraphics[width=260pt]{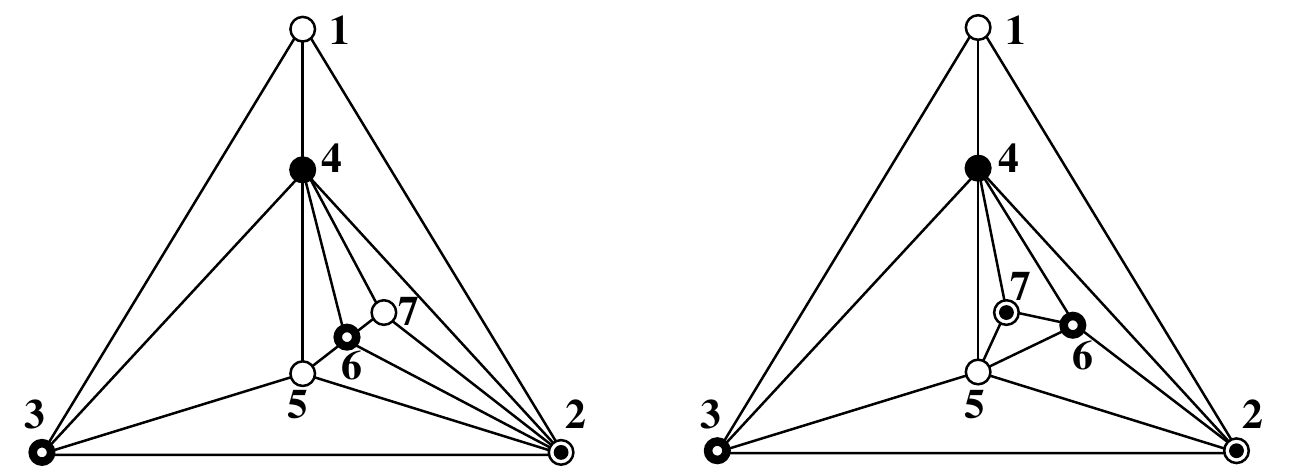}

   \textbf{Figure 4.6.} All of the two (2,2)-FWF graphs with
   order 7
\end{center}

\begin{center}
    \includegraphics[width=360pt]{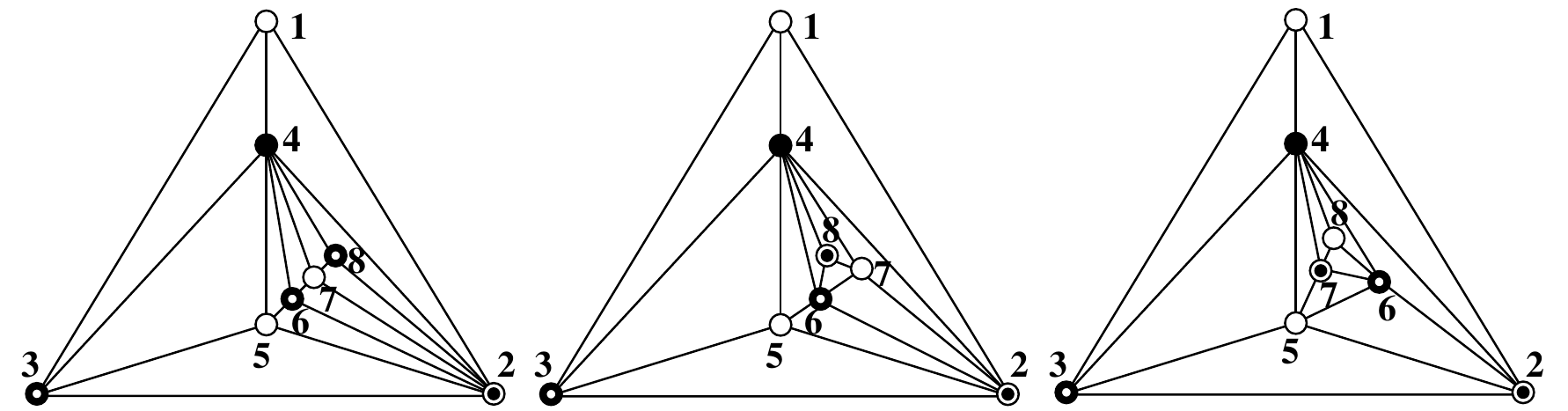}

   \textbf{Figure 4.7.} All of the three (2,2)-FWF graphs with
   order 8
\end{center}

\begin{center}
    \includegraphics[width=360pt]{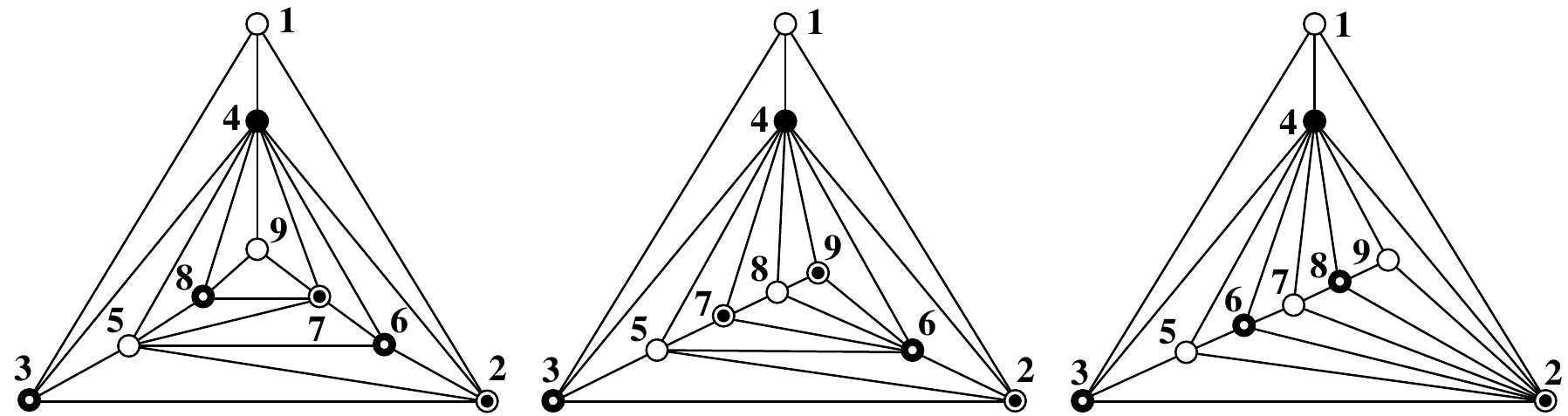}
     \includegraphics[width=360pt]{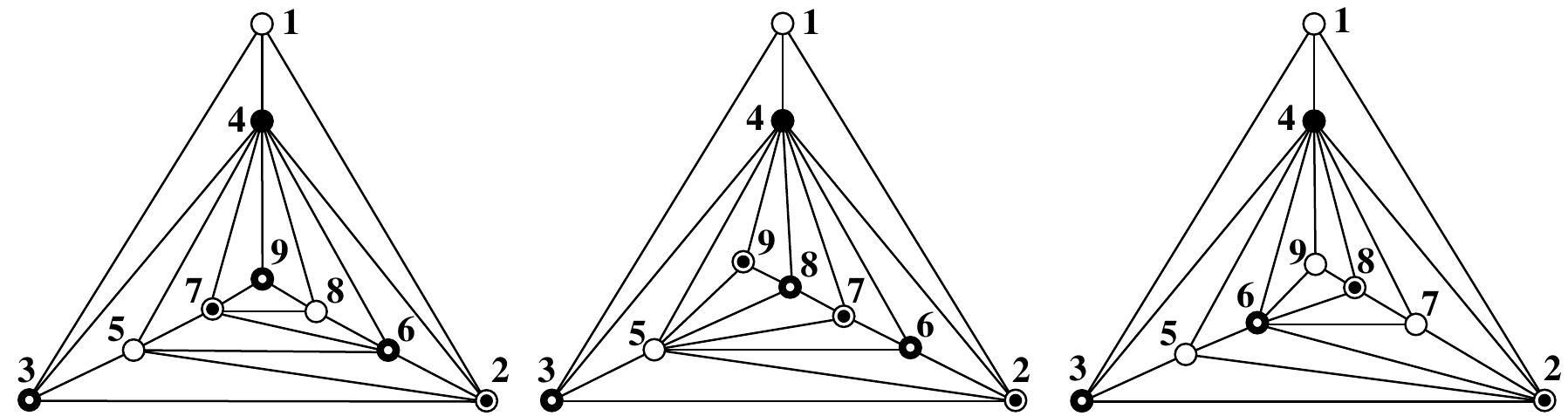}

   \textbf{Figure 4.8.} All of the six (2,2)-FWF graphs with
   order 9
\end{center}

\subsection{The color sequence of a (2,2)-FWF graph}
    Without loss of generality, we can always assume that the (2,2)-FWF graph $G$ is obtained by embedding $3$-degree vertices only in the region $II$
 in the following discussion. Thus, a (2,2)-FWF graph $G$ can be uniquely represented by its color sequence. The specific method is shown in the following.

    Let $V(G)=\{1,2,3,4=u,5,\ldots,n\}$, where vertex $1$($x$) indicates the first fixed  vertex of 3-degree, while the vertex $n(y)$ indicates the second
  vertex of 3-degree; the vertices $1(x)$, 2, 3, and $4(u)$ indicate the 1st, 2nd, 3rd, 4th color axis respectively; while the vertex $4(u)$ is the central
 vertex; the vertex $n-1$ signifies the  3-degree vertex  of the subgraph $G_{n-1}=G-n$; the vertex $n-2$ denotes the 3-degree vertex of the subgraph
 $G_{n-1}-(n-1)$; the rest can be deduced in the same way. The sequence $c_{1}c_{2},\ldots,c_{n}$ is used to indicate the corresponding color sequence of
 the sequence $(1,2,3,4=u,5,\ldots,n)$, and the parameter $c_{i}$ is the color of the vertex $i$ in the (2,2)-FWF graph $G$. So we can obtain

    $$c_{i}\in \{1=y(yellow),2=g(green),3=b(blue),4=r(red)\}.$$

    According to the definition of a (2,2)-FWF graph, we can know that this representation also determines the structure of a graph.
 This structure starts from $K_{4}$ (shown in Figure 4.2), and selects a triangular face embedded the vertices according to the coloring of each vertex.

 \textbf{Example 4.1.} For the color sequence $c_{1}c_{2}c_{3}c_{4}c_{5}c_{6}c_{7}c_{8}c_{9}=ygbrybgyg$, its corresponding (2,2)-FWF graph is shown in Figure 4.9.

    \begin{center}
    \includegraphics[width=200pt]{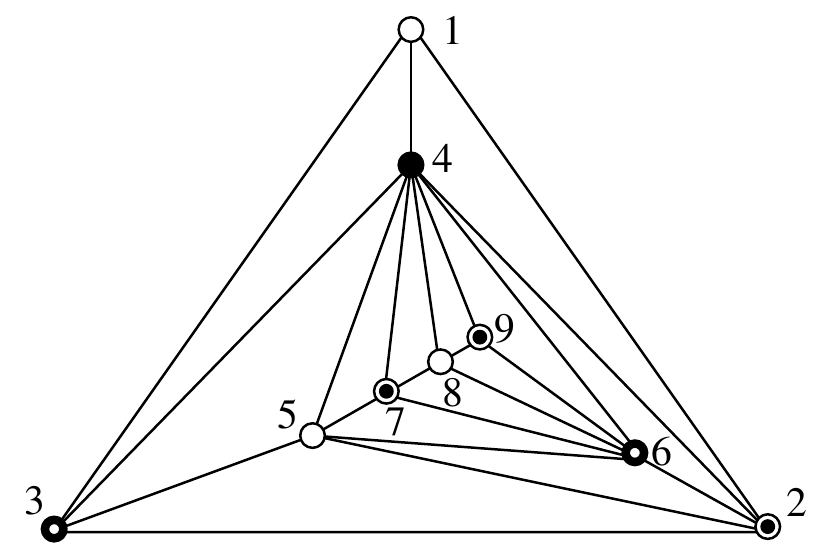}

    \textbf{Figure 4.9. } A color sequence $c_{1}c_{2}c_{3}c_{4}c_{5}c_{6}c_{7}c_{8}c_{9}=ygbrybgyg$ and its corresponding (2,2)-FWF graph
    \end{center}

    For the color sequence of a (2,2)-FWF graph, we can obtain the following theorem:
    \begin{theorem2}
        Let $c_{1}c_{2}\ldots c_{n}$ be the color sequence of a (2,2)-FWF graph. With the agreement in Section 4.2,
    the colors of the first six vertices in this sequence are determined, namely $c_{1}=y, c_{2}=g, c_{3}=b, c_{4}=r, c_{5}=y, c_{6}=b$;
     if $G$ belongs to the adjacent type, then $c_{7}=y$; otherwise, $c_{7}=g$.
    \end{theorem2}

 \subsection{ Chromaticity of graphs induced by extending 4-wheel operation}

    In this section, we are devoted to discussing the vertex coloring problem of the induced graph from
 a (2,2)-FWF graph by extending $4$-wheel operation. We know that a given (2,2)-FWF graph is uniquely $4$-colorable, and according to the definition of the
 color-coordinate system in Section $4.3$, every vertex can also be colored determinately.
    \begin{definition}
 Let $G$ be a $(2,2)$-FWF graph, $f$ be the unique $4$-coloring of $G$,  and $xuy$ be a path of length $2$ in $G$.  Obviously, there exists a coloring $f^*$ of the graph $G*xuy$ that is induced from $G$ by extending $4$-wheel operation on the path $xuy$, and
$$
f^*(x)=\left\{
\begin{array}{cc}
   f(u)& if~ x=u';  \\
   4&  if~ x=v;\\
   f(x)& otherwise.
\end{array}\right. \eqno{(4.1)}
$$
Namely, vertices $u$ and $u'$ are assigned the same color under $f^*$, and the new added vertex $v$ is assigned the different color $4$ from vertices $x,y,u$, while the colors of the rest vertices remain unchanged. We refer to $f^{*}$ as the natural $4$-coloring of the graph  $G*xuy$.
\end{definition}

    Naturally, one question is proposed about whether the induced graph $G*xuy$ obtained by extending $4$-wheel operation
 is uniquely 4-colorable or not. This question is  a key problem in this section. Definitely, the answer is negative, that is,
 $|C_{4}^0(G*xuy)|\geq1$.

 Here the definition of the color neighbor is introduced as follows:
 \begin{definition}
    Let $G$ be a $k$-chromatic graph, and $f \in C_{k}(G)$. The color neighbor of a vertex $u$ of $G$ on
 coloring $f$ is the set which consists of all colors assigned to $\Gamma(u)$ under $f$, denoted as $C(f, \Gamma(u))$.
 \end{definition}

 \begin{theorem2}\label{th4.8}
    Let $G$ be a $(2,2)$-FWF graph with order $n$ and $f$ be the unique $4$-coloring of it. The vertices $x$, $y$ are two vertices of degree $3$ and the vertex $u$ is
 the central vertex of $G$. Then, the induced graph $G*xuy$ is not uniquely $4$-colorable.
 \end{theorem2}
\begin{proof}
    Obviously, if $f(x)=f(y)$,
 the vertex $v$ in the graph $G*xuy$ has two possible colors to choose when both vertices $u$ and $u'$ are colored by red. Hence, the graph $G*xuy$ is not uniquely 4-colorable.
 So we only need to consider the case of $f(x)\neq f(y)$.

    According to the classification in Section $4.2$, all (2,2)-FWF graphs can be classified into two types:
     the adjacent type of region $II$ and  non-adjacent type of region $II$.

    \textbf{Case 1 :} The $(2,2)$-FWF graph $G$ belongs to the adjacent type of region $II$.

     Based on Theorem 4.5, we know that vertices $x$,$2$,$3$ are coordinate axes $1,2,3$, colored by yellow, green and blue respectively;
 and the central vertex $u$ is colored by red. Since all 3-degree vertices can only be embedded in the subregion $I$ of
 the region $II$, the vertex $y$ can be colored with yellow or blue, illustrated in Figure 4.10(a). But when the vertex $y$ is colored
 with yellow, which is the same with vertex $x$, this case is not needed  considering. So we only discuss the case that the
 vertex $y$ is colored with blue. With the definition of extending 4-wheel
 operation,
an extending 4-wheel operation on the path $xuy$ can be done and the graph $G*xuy$ is obtained.

  \begin{center}
   \vspace{5mm}
     \includegraphics[width=360pt]{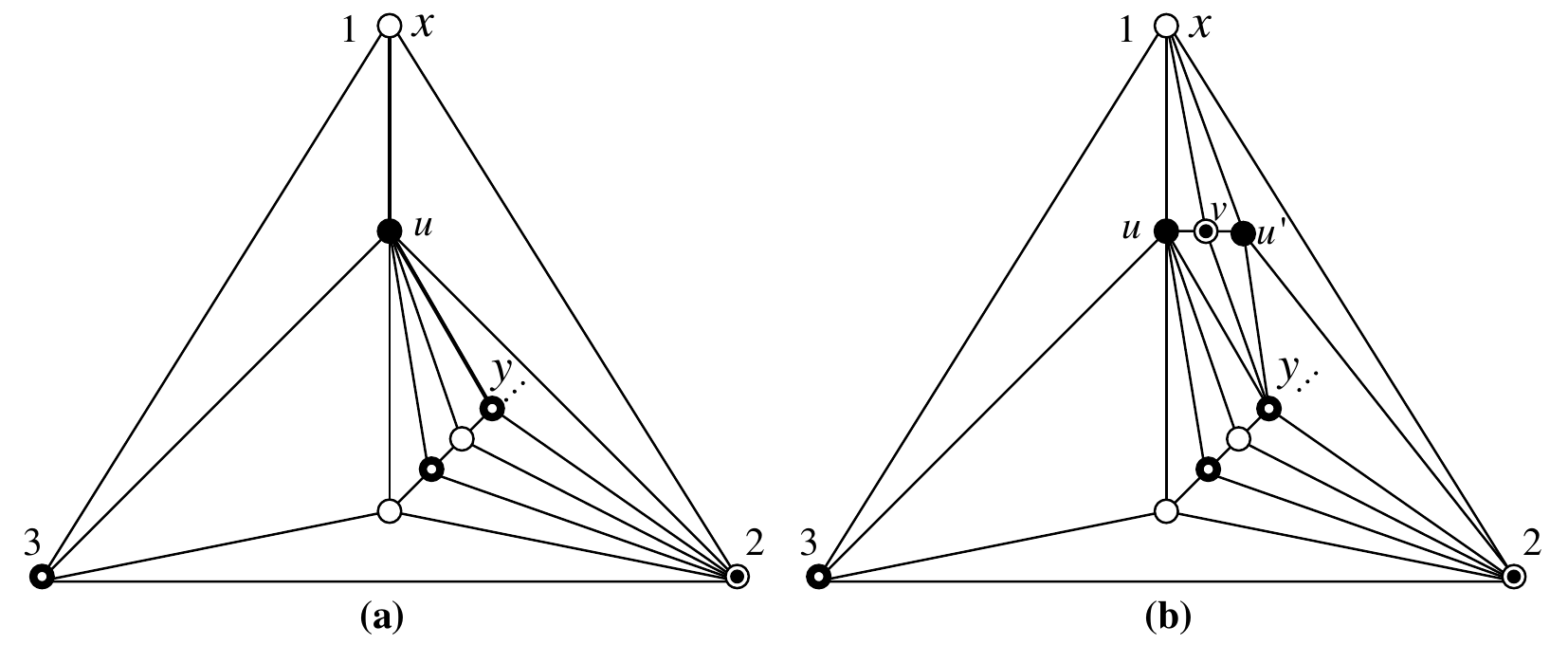}

     \includegraphics[width=180pt]{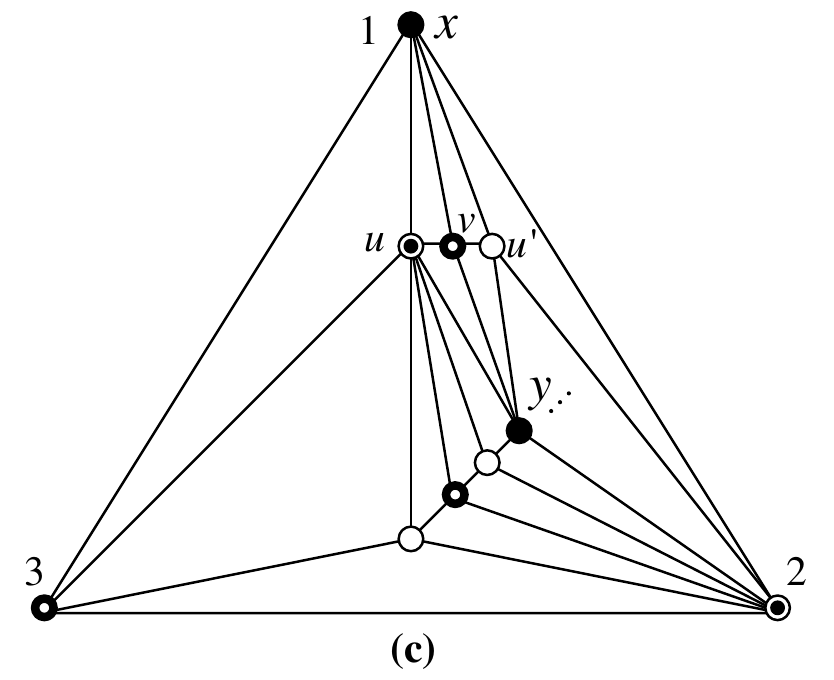}

  \textbf{Figure 4.10.} A graph of the adjacent type of region $II$ and two colorings of its induced graph by extending 4-wheel operation
  \end{center}

    It is easy to obtain two colorings of the graph $G*xuy$ as follows: one is the natural 4-coloring $f$ in which
 the vertex $u'$ is colored with red and the vertex $v$ embedded newly is colored with green.
 And the colorings of other vertices remain
 unchanged. Obviously, it is a coloring of the graph $G*xuy$, illustrated in Figure 4.10(b). Besides,
from the discussion above, there is only one vertex $v$ colored with green in $\Gamma(u)$ and only one vertex $x$ colored with yellow in $\Gamma(u')$, under $f$.
Further, for graph $G*xuy$, only two vertices $u$ and $u'$ are colored with red under its natural coloring $f$, so we can obtain a new 4-coloring $f'$ of the graph $G*xuy$ in the way that:
let the vertices
 $u$,$u'$,$x$,
 $y$ and $v$ be recolored  with green, yellow, red, red and blue respectively, other vertices remain unchanged on the basis of $f$.

 Since there is only one vertex $v$ colored with green in $\Gamma(u)$ under $f$,
 then change the color assigned to vertex $u$ from red to green, and only vertex $v$ receives the same green color.
 Similarly, since there is only
 one vertex $x$ colored with yellow $\Gamma(u')$ under $f$,
then change the color assigned to vertex $u'$ from red to yellow,
 and only vertex $x$ receives the same yellow color.
 After that, recolor nonadjacent vertices $x$ and $y$ with red, and change the color assigned to vertex $v$ from green to blue. Then we will obtain a new 4-coloring $f'$ of the graph $G*xuy$ when remain the colors of other vertices
 unchanged, see an illustration shown in Figure 4.10(c).
 These two colorings $f$ and $f'$ are different apparently. Hence, the case 1 is proved.

    \textbf{Case 2:} The $(2,2)$-FWF graph $G$ belongs to the nonadjacent type of region $II$.

    According to Theorem 4.7, the first six vertices of all the (2,2)-FWF graphs are colored in the same way, illustrated as follows:
   $$\left(
    \begin{array}{ccccccc}
        1 & 2 & 3 & 4 & 5 & 6 & \ldots \\
        y & g & b & r & y & b & \ldots \\
    \end{array}
   \right)\eqno{(4.2)}$$
 Namely, the color coordinate axes vertex 1 (or vertex $x$), 2, 3 and 4 (or vertex $u$) receive yellow, green, blue and red colors respectively.
 Vertex 1 is a vertex of degree 3 and adjacent to the central vertex 4 colored with red, the vertex 2 colored with green and the  vertex 3 colored with blue.
 Therefore, vertex 1 is a vertex of degree 5 in the graph $G*xuy$. And in the natural coloring of $G*xuy$, it is adjacent to the vertex 2 colored with green, the vertex 3 colored with blue, the vertex 4 (or vertex $u$) colored with red, the vertex $u'$ colored with red and the vertex $v$ colored with blue respectively.

    Since the graph $G$ belongs to the nonadjacent type, so the vertices $7,8,\ldots,n$ must be added in the triangular face formed by vertices 4, 5 and 6, shown in the Figures $4.6(c)$ and $4.11(a)$. According to Theorem 4.7, the 7th vertex  can only be colored with green, this case is illustrated as
 follows:
   $$\left(
    \begin{array}{cccccccc}
        1 & 2 & 3 & 4 & 5 & 6 & 7 & \ldots \\
        y & g & b & r & y & b & g & \ldots \\
    \end{array}
   \right)\eqno{(4.3)}$$

    Then, it can be known easily that the vertex 2 colored with green, which is also a color coordinate axis and a vertex of degree 5. The neighbors of vertex 2 are
 vertex 1 (yellow), vertex 3 (blue), vertex 5 (yellow), vertex $u'$ (red) and vertex 6 (blue). Hence
    $$C(f,\Gamma(1))=\{g\{2\},b\{3,v\},r\{u,u'\}\} \eqno{(4.4)}$$
    $$C(f,\Gamma(2))=\{y\{1,5\},b\{3,6\},r\{u'\}\} \eqno{(4.5)}$$

    Now we take the representative graph in Figure 4.11(a) as an example, and the detailed steps of a new $4$-coloring induced by the natural 4-coloring  of graph $G*xuy$ (Figure 4.11(b)) is given as follows:

    First, change the color assigned to vertex 1 from yellow to green. By formula (4.4), the two
 ends  of edge $\{1,2\}$ are both colored with green, and it is the unique pseudo color edge (the two ends of this edge are not colored properly). Other vertices are colored properly, illustrated
 by Figure 4.11(c).

    Second, change the color assigned to vertex 2 from green to red. Thus,  the coloring of the two ends of edge $\{1,2\}$ becomes proper, while
  $\{u',2\}$ becomes a pseudo color edge, for its two ends are both colored with red.  The coloring of other vertices is proper,
 illustrated in Figure 4.11(d).

    Third, the vertex $u'$ is recolored with yellow. Thus, the pseudo color edge $\{u',2\}$ becomes proper. There may be
 several vertices colored with yellow in the neighbor of the vertex $u'$, which can form a set $C_{4}(u',yellow)$. Therefore,
 this step generates several pseudo color edges whose number is $|C_{4}(u',yellow)|$. Obviously, other edges are all proper,
 it is illustrated in Figure 4.11(e).

    Fourth, all the vertices in $C_{4}(u',yellow)$ are colored with red. Since in the $4$-coloring of the third step, only two
 vertices $u$ and 2 are colored in red. Obviously, $C_{4}(u', yellow)\subset \Gamma(u')$. So all vertices in $C_{4}(u', yellow)$
 are not adjacent to the vertex $u$. In the set of neighbors of the  vertex 2 after the third step, the vertices 1, 3 and 6 are colored
with green, blue and blue.
Although the vertex 5 is colored with yellow, it is a vertex of degree 5 and not in $C_{4}(u', yellow)$. Therefore,
 the edges between vertex 2 and all  red vertices in $C_{4}(u', yellow)$ are proper. Moreover, the vertices in
 $C_{4}(u', yellow)$ form an independent set of the graph. So they can not generate pseudo color edges by themselves. This step is
 illustrated in Figure 4.11(f).

    Thus, based on the natural 4-coloring of the graph $G*xuy$, we can obtain a new coloring differing from the natural $4$-coloring, which means that the induced graph constructed from nonadjacent (2,2)-FWF graph by extending 4-wheel operation is not uniquely $4$-colorable.

    To sum up the cases 1 and 2, this theorem holds.
 \end{proof}

\textbf{Remark}:  Theorem 4.8 shows that the induced graph $G*xuy$ obtained from a (2,2)-FWF graph $G$ by extending 4-wheel operation on  a path $xuy$ of $G$ is not uniquely $4$-colorable. Here the degree of each end of the path $xuy$ is 3. However, it is possible that $G*xuy$ is uniquely $4$-colorable when the ends of the path $xuy$ are not the vertices of degree 3 in $G$.
 \begin{center}

  \includegraphics[width=240pt]{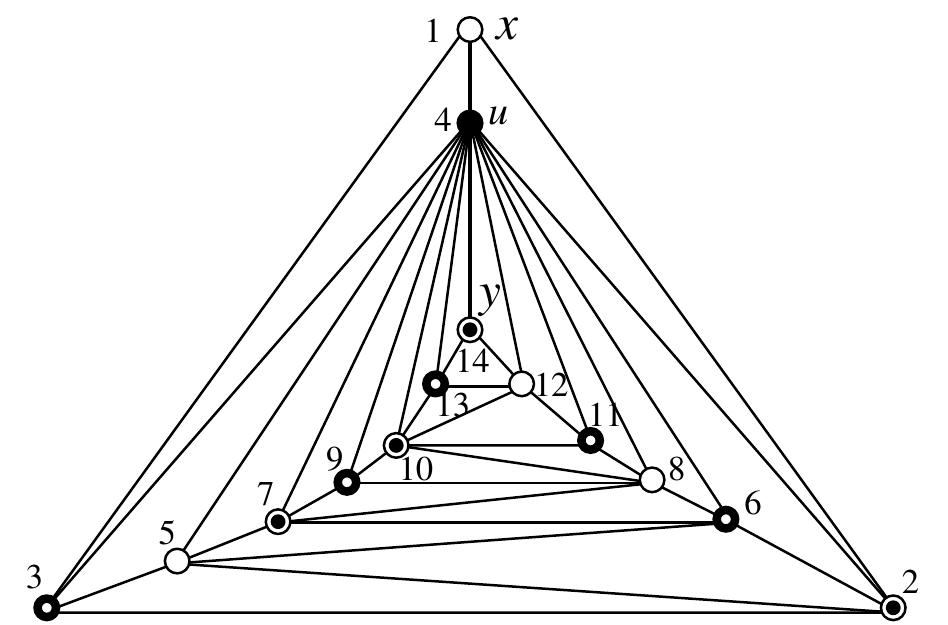}

  \textbf{(a)}  A representative (2,2)-FWF $G$

  \includegraphics[width=240pt]{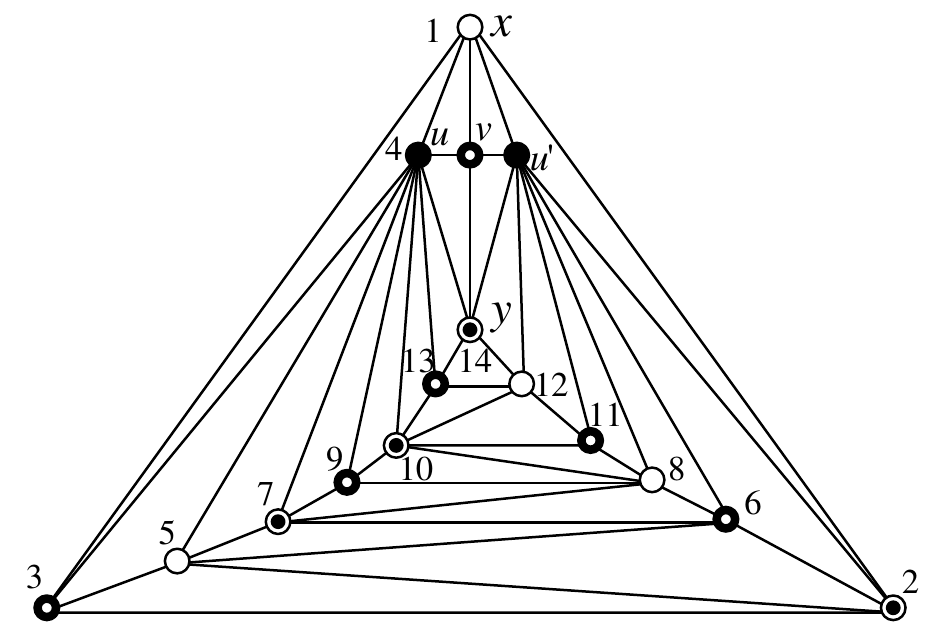}

  \textbf{(b)} The natural 4-coloring of the induced graph $G*xuy$ by extending 4-wheel operation to $G$

  \includegraphics[width=240pt]{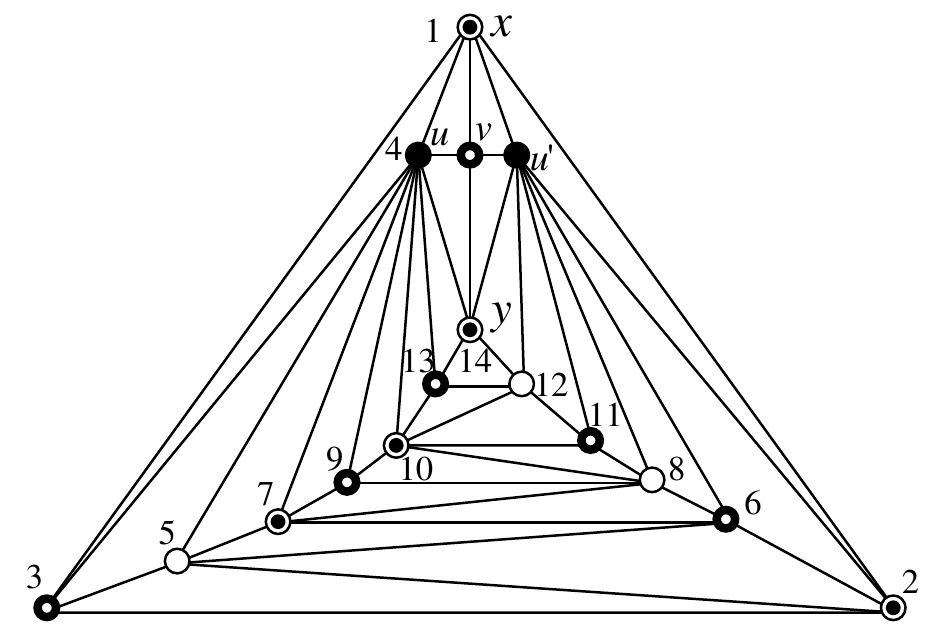}

  \textbf{(c)} The coloring of the induced graph $G*xuy$ when vertex 1 is recolored with green,
  which generates a pseudo color edge \{1,2\}

  \includegraphics[width=240pt]{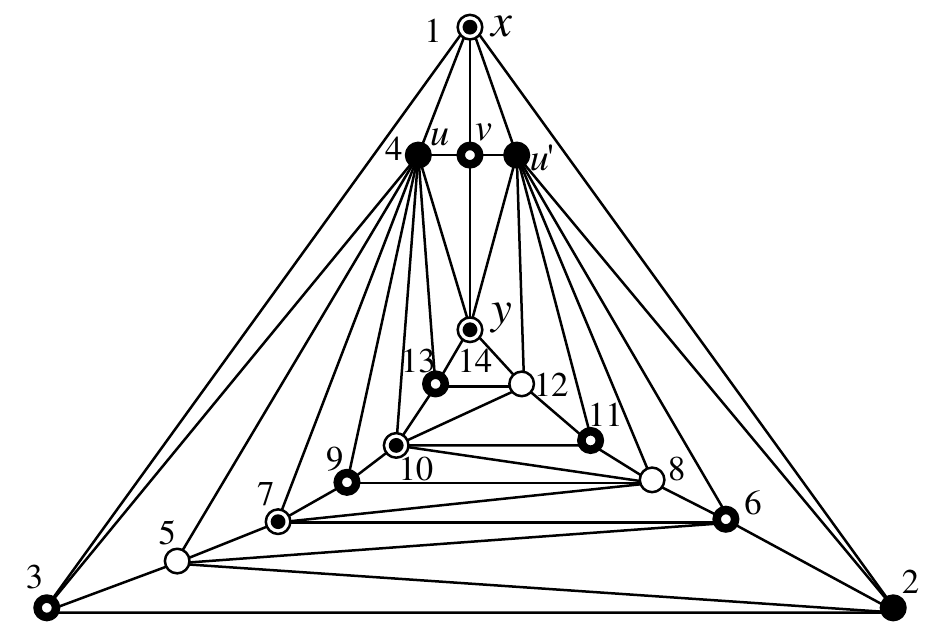}

  \textbf{(d) } The coloring of the induced graph $G*xuy$ based on the step (c) when vertex 2 is recolored with red,
  which generates a pseudo color edge \{$u'$,2\}

  \includegraphics[width=240pt]{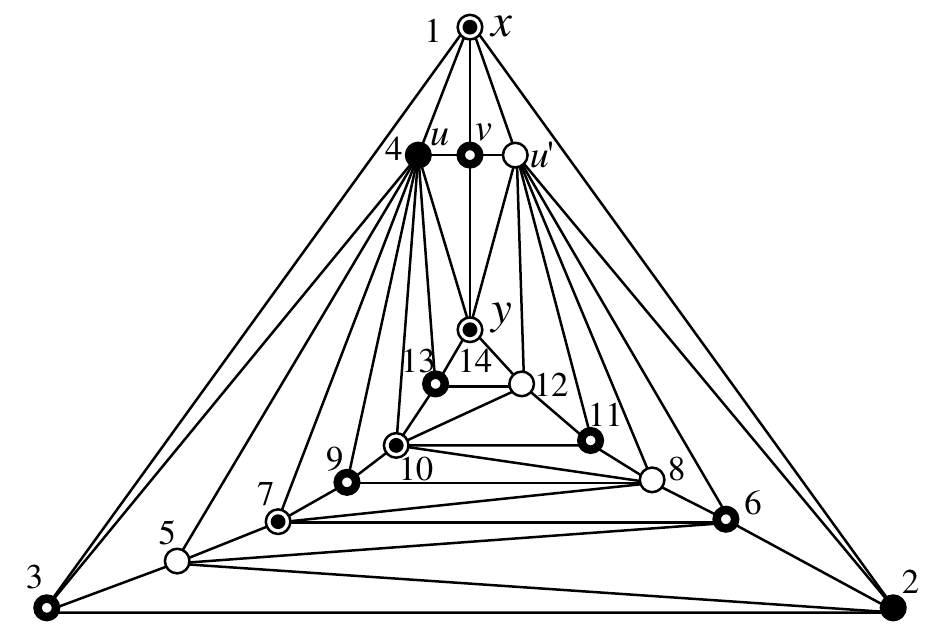}

  \textbf{(e)} The coloring of the induced graph $G*xuy$ based on the step (d) when vertex $u'$ is recolored with yellow,
  which produces several pseudo color edges $\{\{u',u^{\prime \prime}\},u^{\prime \prime}\in C_{4}(u',yellow)\}$

  \includegraphics[width=240pt]{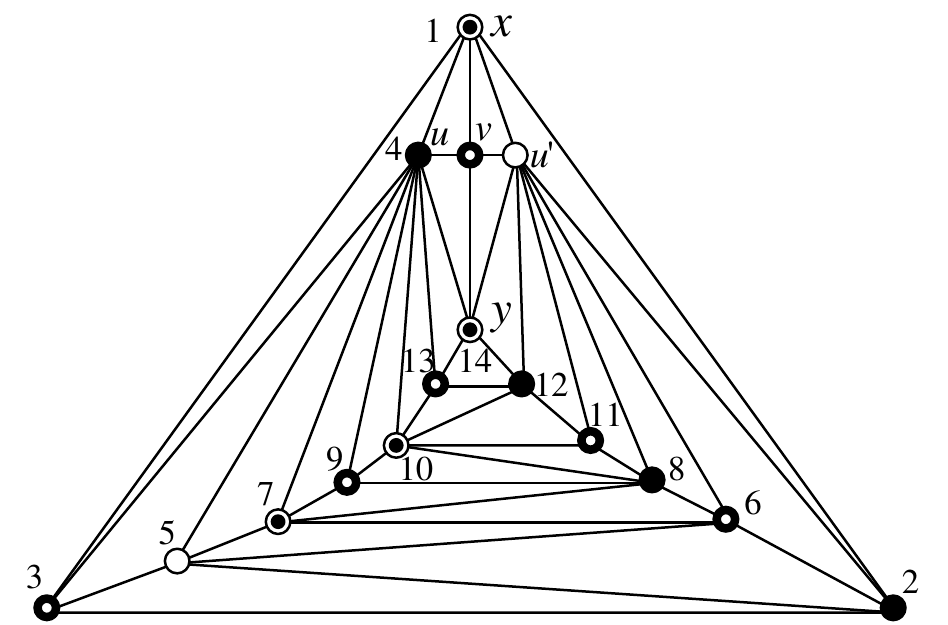}

  \textbf{(f) } A new coloring of the induced graph $G*xuy$ based on the step (e) when the vertices in the set $C_{4}(u',yellow)$
 are colored with red

 \textbf{Figure 4.11.} An illustration for showing that  the induced graph $G*xuy$ by extending 4-wheel operation to a $(2,2)$-FWF graph $G$ belonging to the nonadjacent type of region $II$ is not unique 4-colorable
  \end{center}

\section{Coloring-structure of maximal planar graphs}

  In fact, the essential error, which appeared in the proof of four color conjecture by $Kempe$ in 1879, was that he did not make clear the basic structure of  4-colorings of planar graphs. Although $Heawood$ found this error in 1890, he didn't give a correction to this problem. However, he proved the five color theorem by means of $Kempe$'s idea. The results of this section show that the real reason that later many scholars could not correct $Kempe$'s error might be that the coloring-structure of 4-colorable maximal planar graphs haven't been studied in depth.

    This section aims to make the structure of graph coloring corresponding to  a 4-coloring $f$ of a 4-colorable maximal planar graph $G$ clear. The specific method is: $(1)$ delete the vertices from $G$ that belong to one of the same independent set generated by a 4-coloring $f$ of $G$, thus the 4-coloring problem of a maximal planar graph can be transformed as a 3-coloring problem of a planar graph correspondingly, and the structural problem of six bicolored subgraphs reduces to three bicolored subgraphs' structural problem. So, not only does computation reduce largely, but also the structure becomes simple and easy to study. $(2)$ Furthermore, in the process of researching three bicolored subgraphs, we study the union structure of them and any two of them, respectively. We discover that it is very important to study tree-colorings in $C_4^0(G)$ for attacking JT-conjecture and other problems of graph coloring, so the tree-coloring and cycle-coloring are studied preliminarily in this section.

\subsection{Cycle-colorings and tree-colorings}

\quad\quad Let $G$ be a 4-colorable maximal planar graph and $C(4)$ the set of colors. A \emph{cycle-coloring} $f$ of $G$ is a 4-vertex-coloring such that there exists a cycle $C_{2m}={v_1v_2\cdots v_{2m}v_1}(m\geq 2)$ in $G$ with $|\{f(v_1),f(v_2),\cdots, f(v_{2m})\}|=2$, where $V(C_{2m})=\{v_1,v_2,\cdots, v_{2m}\}$. We refer to $C_{2m}$ as a \emph{bicolored cycle} of $f$, or say $f$ contains a \emph{bicolored cycle}. On the other hand, if $f\in C_4^0(G)$ doesn't contain a bicolored cycle, then $f$ is called a \emph{tree-coloring} of $G$. From the definitions of cycle-coloring and tree-coloring, for any 4-colorable planar graph $G$ and $f\in C_4^0(G)$, $f$ is either a cycle-coloring or a tree-coloring.

For example, there are  eight 4-colorings for the graph shown in Figure $5.1$, and these eight 4-colorings all are cycle-colorings; for the 4-colorings shown in Figure $5.2$, $f_1,f_2$ are cycle-colorings and $f_3,f_4$ are tree-colorings; Figure $5.3$ gives  ten 4-colorings of the icosahedron and they are all tree-colorings. Naturally, we can know a fact that all the 4-colorings of a maximal planar graph maybe contain only tree-colorings, or only cycle-colorings, or both tree-colorings and cycle-colorings. Then, which graphs contain only tree-colorings? Which graphs have only cycle-colorings? Which graphs contain both tree-colorings and cycle-colorings? Obviously, these problems are the basis of studying the coloring properties of 4-colorable maximal planar graphs.

Considering the above three examples, the maximal planar graphs can be divided into three categories according to cycle-coloring and tree-coloring: $\textcircled{1}$ \emph{pure cycle-coloring} graphs, namely these graphs have only cycle-colorings;  $\textcircled{2}$ \emph{pure tree-coloring} graphs, namely such graphs have only tree-colorings; $\textcircled{3}$ \emph{impure coloring} graphs that have both cycle-colorings and tree-colorings.

From now on, we use four different icons shown in Figure 5.1(a), to denote colors 1,2,3 and 4, respectively.
  \begin{center}
        \includegraphics [width=220pt]{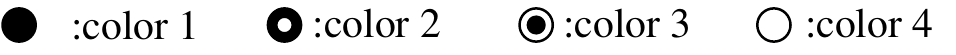}

          \vspace{0.2cm}
        \textbf{Figure 5.1(a).} The check figure between icons and colors
  \end{center}
  \begin{center}
        \includegraphics [width=340pt]{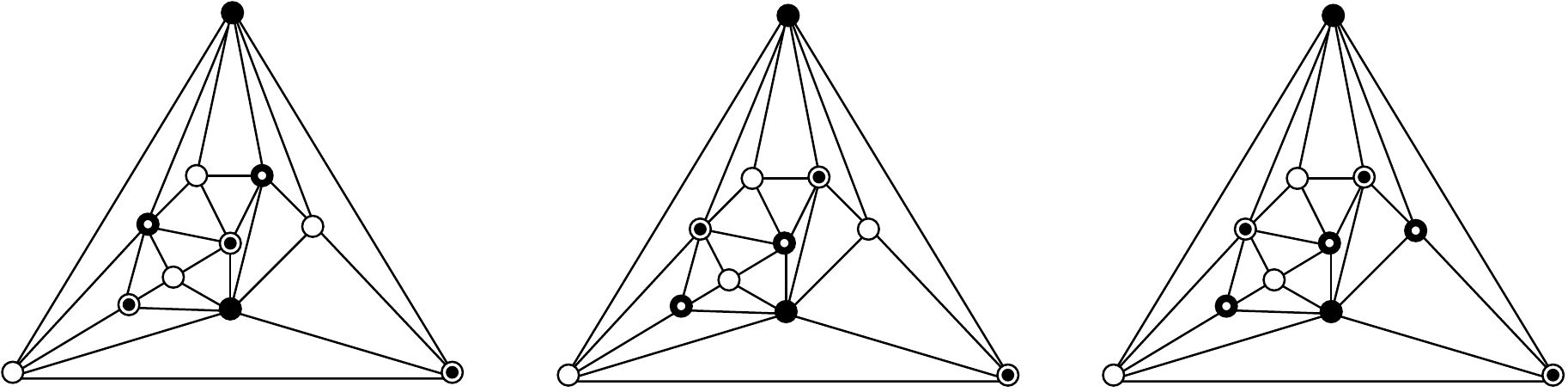}

\small{ \hspace{0.6cm}$f_1$:1-4 cycle \hspace{1.2cm}$f_2$:1-3,1-4 and 3-4 cycles  \hspace{0.4cm}$f_3$:1-3 and 3-4 cycles}

         \vspace{0.2cm}
        \includegraphics [width=340pt]{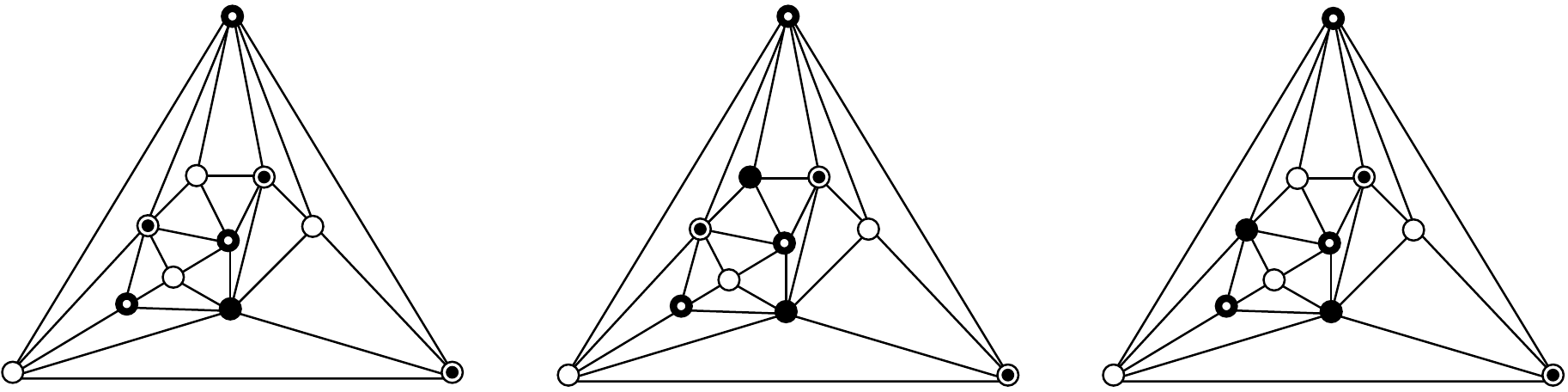}
\small{\hspace{0cm}$f_4$:2-3,2-4 and 3-4 cycles \hspace{0.9cm}$f_5$:2-3 cycle \hspace{1.2cm}$f_6$:1-2,1-4 and 2-4 cycles}

        \vspace{0.2cm}
        \includegraphics [width=260pt]{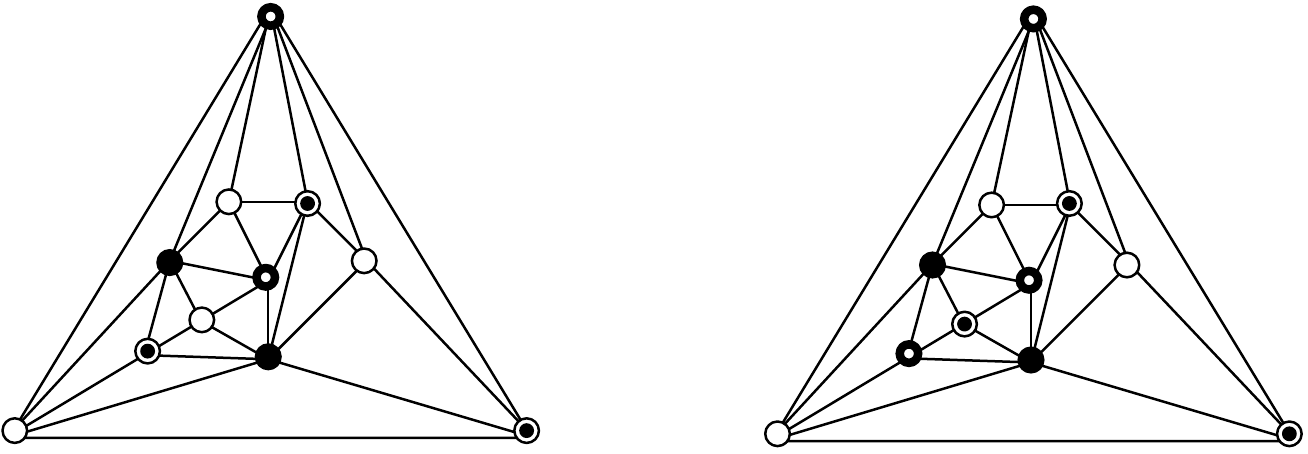}

         \small{\hspace{0cm}$f_7$:1-4 cycle \hspace{3cm}$f_8$:1-2 cycle}\\
         \vspace{0.2cm}
        \textbf{Figure 5.1(b).} All of eight 4-colorings of a  maximal planar graph of order 11
  \end{center}

  \begin{center}

        \includegraphics [width=260pt]{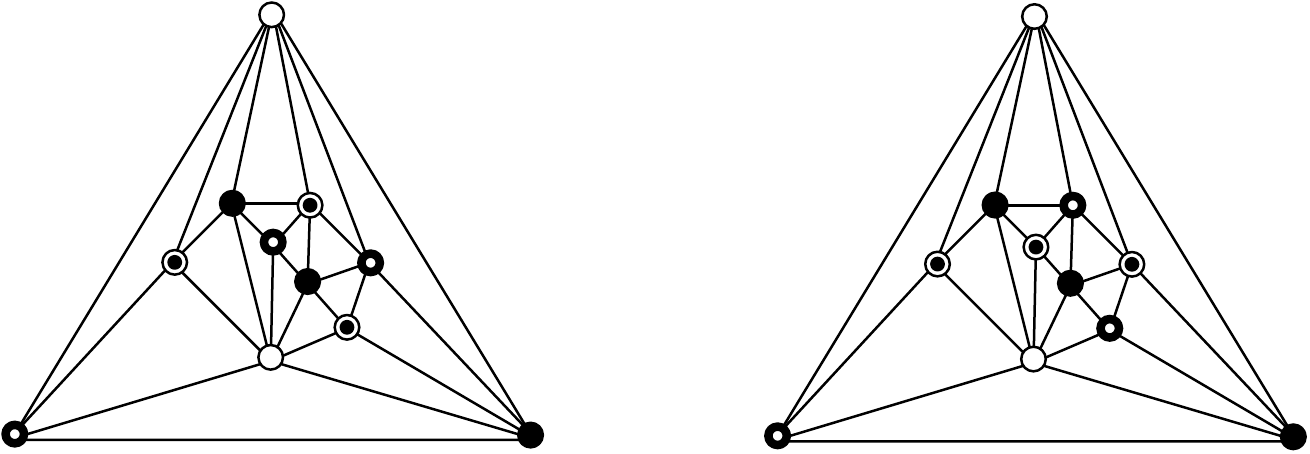}

       \small{\hspace{0cm}$f_1$:1-4 cycle \hspace{3cm}$f_2$:1-4 cycle}

        \vspace{0.2cm}
        \includegraphics [width=260pt]{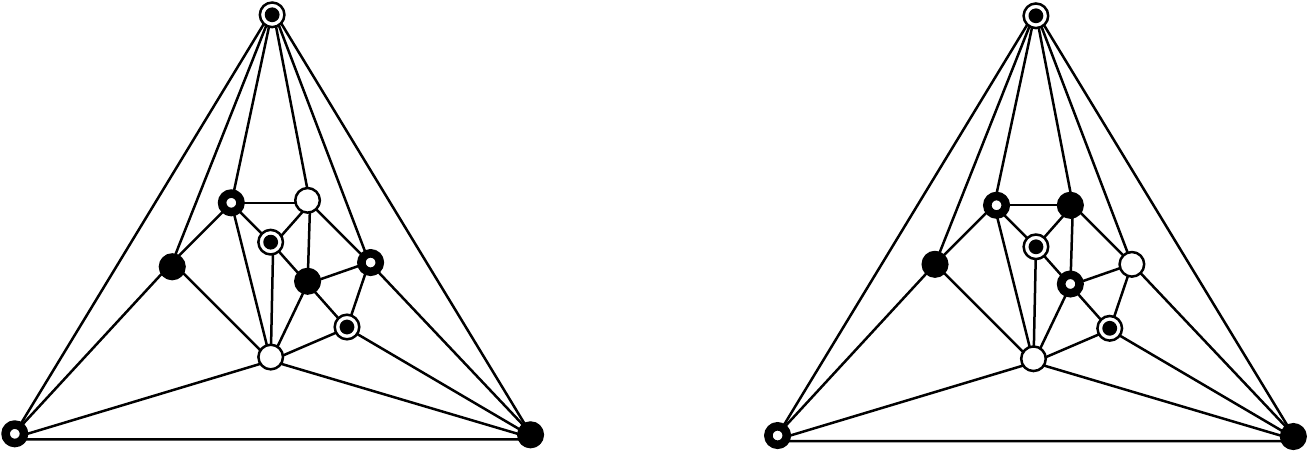}

         \small{\hspace{0cm}$f_3$:no bicolored cycle \hspace{2cm}$f_4$:no bicolored cycle }\\

        \vspace{0.2cm}
        \textbf{Figure 5.2.} All of four 4-colorings of a maximal planar graph with order 11
  \end{center}

  \begin{center}

        \includegraphics [width=320pt]{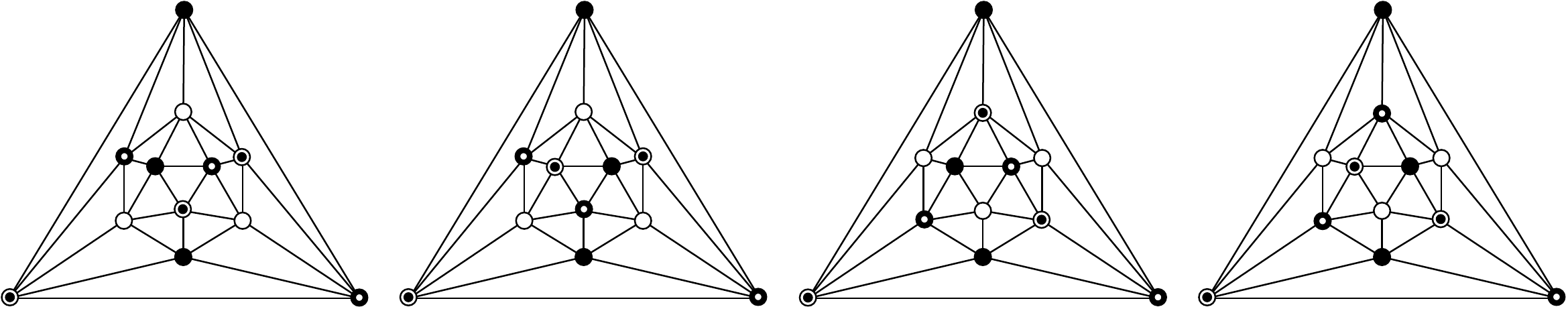}

        \includegraphics [width=320pt]{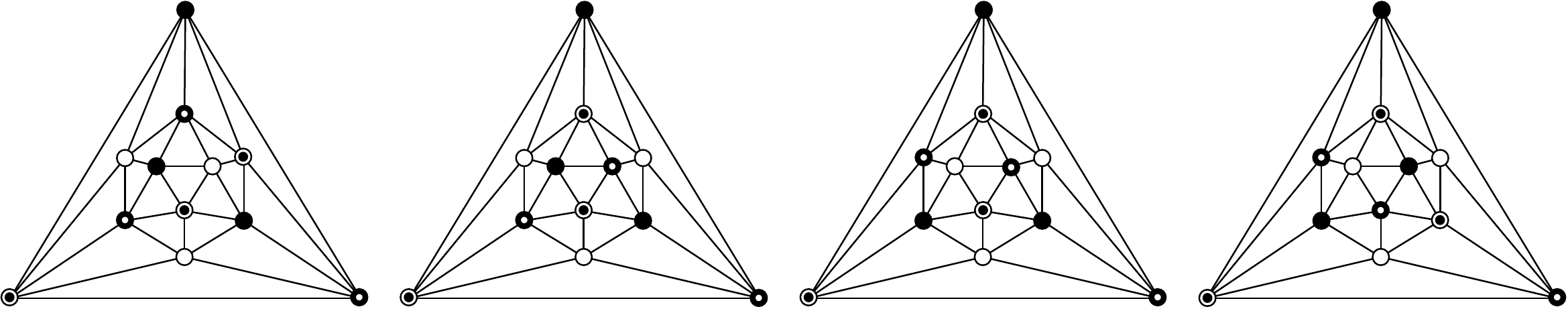}

        \includegraphics [width=320pt]{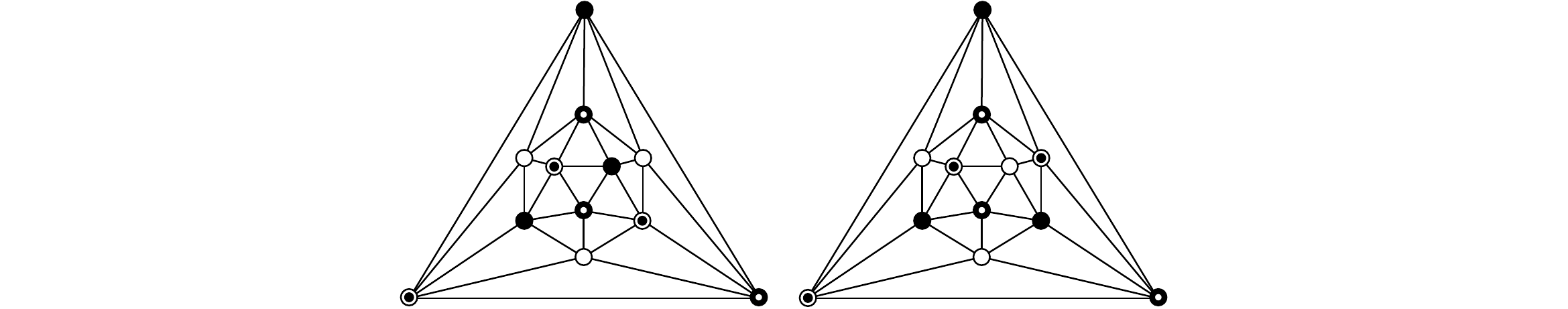}

        \textbf{Figure 5.3.} All of ten 4-colorings of icosahedron
  \end{center}

In terms of the relationship between 4-colorings and the structure of a maximal planar graph, the following results are obvious.

  \begin{theorem2}\label{th5.1} Let $G$ be a double-center wheel graph with $\delta(G)\geq 4$. Then $G$ is a pure cycle-coloring graph.
  \end{theorem2}
  \begin{proof}
   Let $u,v$ be the wheel-center vertices of $G$ and $f$ a 4-coloring of $G$. If $u,v$ are assigned the same color under $f$, it is easy to infer that $f$ contains at least a bicolored cycle of length 4. Otherwise if $u,v$ are assigned different colors, because $G\setminus \{u,v\}$ is a cycle $C$, then the length of $C$ must be even and the vertices of $C$ can be colored only by two colors, so $f$ also contains a bicolored cycle.
  \end{proof}

  \begin{theorem2}\label{th5.2}  For the maximal planar graph $G_1$ and $G_2$  shown in Figures 5.4(a) and (b) , we have: $\textcircled{1}$ when $l$ is even,  $G_1$ has only one tree-coloring; $\textcircled{2}$ when $l$ is odd, $G_2$ has only one tree-coloring.
  \end{theorem2}

The proof of Theorem 5.2 is easy, so omitted here. Now we give some examples of this theorem: the first and fourth graphs shown in Figure 5.5 illustrate the first case of this theorem and the sixth graph shown in Figure 5.5 illustrates the second case.

\begin{center}

        \includegraphics [width=300pt]{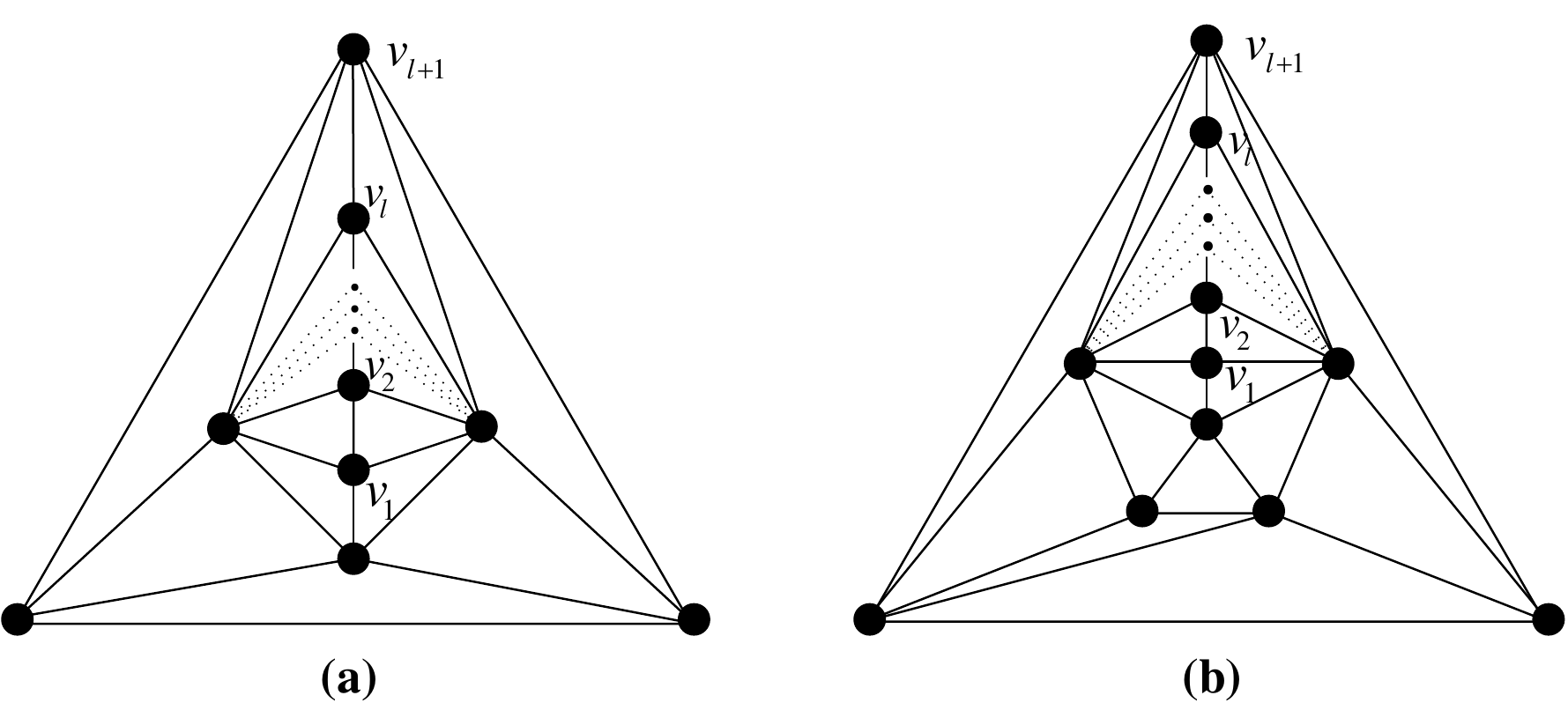}

        \textbf{Figure 5.4.} Two types of maximal planar graph with impure coloring
\end{center}

If a maximal planar graph $G$ contains 3-degree vertices, the
coloring properties of the graph obtained by deleting these 3-degree vertices from $G$
are as the same as the original graph $G$. So, we  need only consider the maximal planar graphs with minimum degree 4 or 5 when we study the coloring properties of them.

Considering the colorings of maximal planar graphs whose orders are from 7 to 11 and whose minimum degrees are not less than 4, we can obtain such a fact that the number of tree-colorings is very small comparing with the number of cycle-colorings. These graphs and their 4-colorings are shown in Appendix $I$, where there is only one graph with order 7 and it is  a double-center wheel graph, so it has no tree-colorings; two graphs with order 8: one is a double-center wheel graph and 3-colorable, the other has two cycle-colorings and one tree-coloring; five graphs with order 9: the first is 3-colorable and divisible, the second has only six cycle-colorings, the third which we refer to as \emph{9-mirror graph} has only two tree-colorings, the fourth is double-center graph, the fifth has only seven cycle-colorings; twelve graphs with order 10: just the fifth and eighth have tree-colorings, and only one respectively, the others do not have tree-colorings; thirty-four graphs with order 11: only the second, fourth, ninth and fourteenth have tree-colorings and each of them at most contain two tree-colorings, the others have no tree-colorings; the seventeenth, nineteenth, twenty-fourth to twenty-eighth, thirtieth and thirty-second are divisible and the twenty-first is 3-colorable.

Table 5.1 \small{The cycle-coloring number and tree-coloring number of the maximal
planar graphs of orders from 7 to 11, and the minimum degrees  are 4 or 5.}

\begin{tabular}{|c|c|c|c|c|c|c|c|c|c|c|c|}  
\hline
\small{GL} &7&$8_1$&$8_2$&$9_1$&$9_2$&$9_3$&$9_4$&$9_5$&$10_1$&$10_2$&$10_3$\\ \hline 
\small{CN} &5&$\ast$&2&$\ast$&6&0&17&9&8&6&$14$\\ \hline          
\small{TN} &0&$\ast$&1&$\ast$&0&2&0&0&0&0&$0$\\ \hline   
\small{GL} &$10_4$&$10_5$&$10_6$&$10_7$&$10_8$&$10_9$&$10_{10}$&$10_{11}$&$10_{12}$&$11_1$&$11_2$\\ \hline
\small{CN} &13&4&$7$&14&10&$\ast$&$\ast$&$\ast$&$\ast$&$8$&6\\ \hline          
\small{TN} &0&1&$0$&0&1&$\ast$&$\ast$&$\ast$&$\ast$&$0$&1\\ \hline         
\small{GL} &$11_3$&$11_4$&$11_5$&$11_6$&$11_7$&$11_8$&$11_9$&$11_{10}$&$11_{11}$&$11_{12}$&$11_{13}$\\ \hline
\small{CN} &$10$&12&16&9&10&11&2&13&22&29&10\\ \hline
\small{TN} &$0$&1&0&0&0&0&2&0&0&0&0\\ \hline
\small{GL} &$11_{14}$&$11_{15}$&$11_{16}$&$11_{17}$&$11_{18}$&$11_{19}$&$11_{20}$&$11_{21}$&$11_{22}$&$11_{23}$&$11_{24}$\\ \hline
\small{CN} &5&8&11&$\ast$&17&$\ast$&14&$\ast$&21&13&$\ast$ \\ \hline
\small{TN} &1&0&0&$\ast$&0&$\ast$&0&$\ast$&0&0&$\ast$ \\ \hline
\small{GL} &$11_{25}$&$11_{26}$&$11_{27}$&$11_{28}$&$11_{29}$&$11_{30}$&$11_{31}$&$11_{32}$&$11_{33}$&$11_{34}$&\\ \hline
\small{CN} &$\ast$&$\ast$&$\ast$&$\ast$&10&$\ast$&41&$\ast$&25&85&\\ \hline
\small{TN} &$\ast$&$\ast$&$\ast$&$\ast$&0&$\ast$&0&$\ast$&0&0&\\ \hline
\end{tabular}
\vspace{0.2cm}

 \small{Here, $GL$ denotes graph label, CN the number of cycle-coloring, TN the number of tree-coloring, $i_j$ the $j$th graph with order $i$ in appendix. And, $\ast$ denotes the corresponding graph is 3-colorable or divisible.}

In Table $5.1$ we can see that there are much more cycle-colorings than tree-colorings. In the total of fifty-four maximal planar graphs with orders from 7 to 11 and  minimum degrees  4 or 5, only one is pure tree-coloring graph and we refer to this graph as \emph{9-mirror graph} (the second graph in Figure 5.5); eight graphs contain at least one tree-coloring (see Figure 5.5). In addition, there are thirty pure cycle-coloring graphs, seven impure coloring graphs and sixteen divisible (or 3-colorable) graphs.  Apart from 3-colorable and divisible graphs in these fifty-four graphs, the number of 4-colorings of the remaining graphs is 533, but the number of tree-colorings is just 10, which shares the proportion of 1.876\%, rarely!

\begin{center}

        \includegraphics [width=320pt]{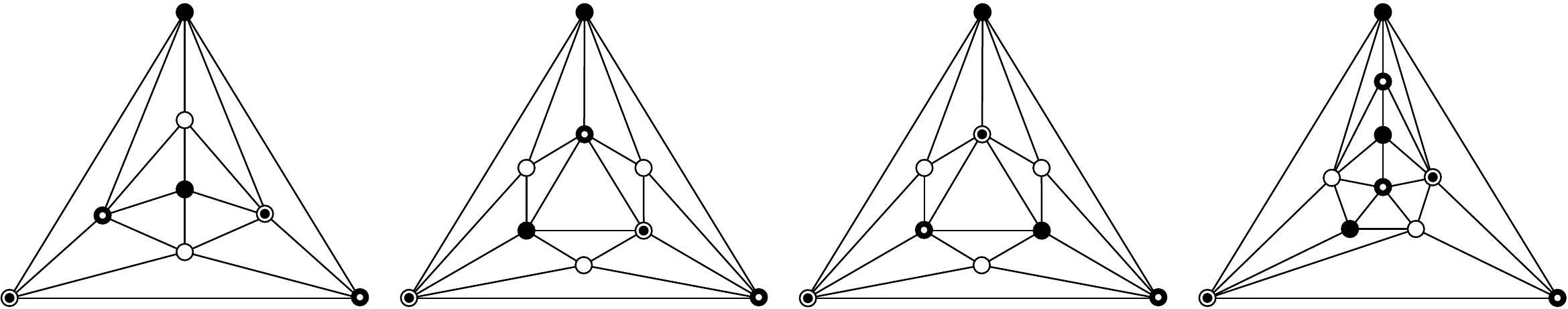}

        \includegraphics [width=320pt]{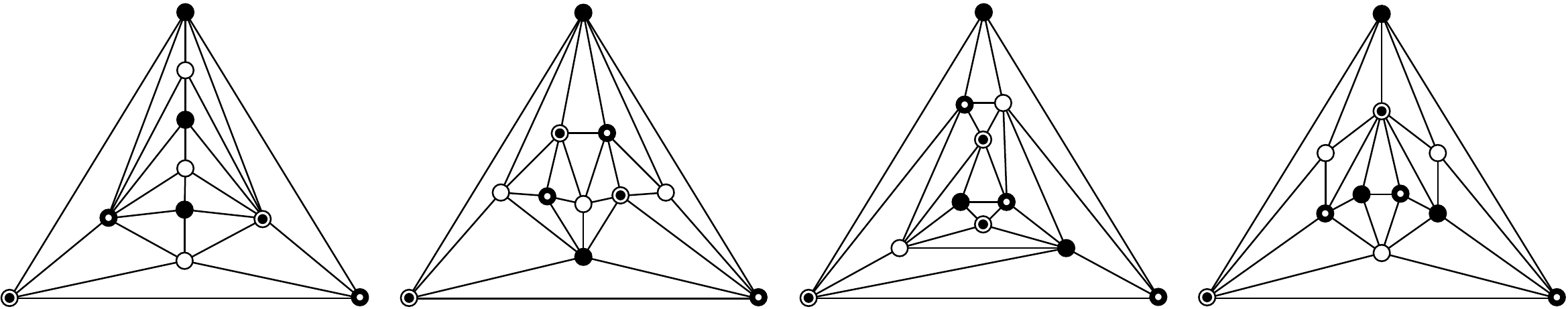}

        \includegraphics [width=320pt]{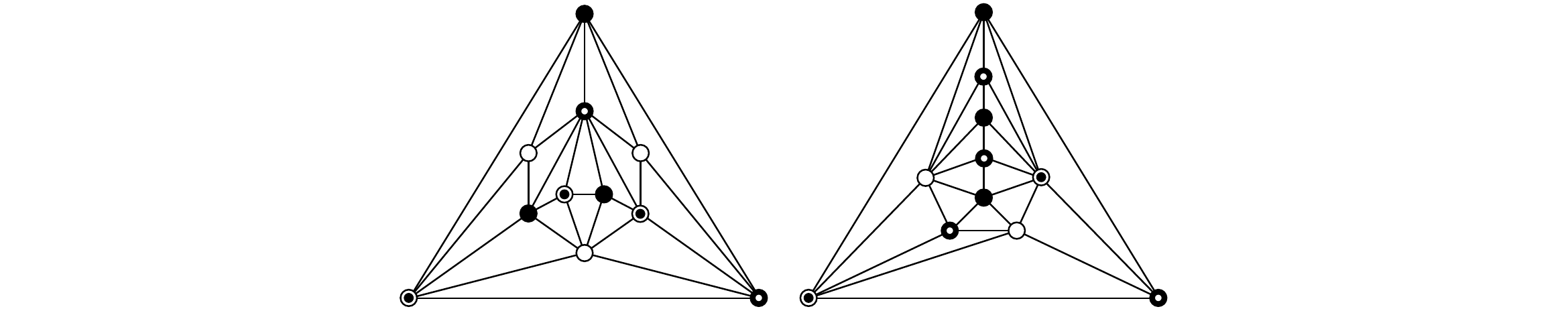}

        \textbf{Figure 5.5.} All tree-colorings of the maximal planar graphs of orders from 7 to 11 and $\delta\geq 4$
\end{center}

For the pure tree-coloring graphs, we have an evident fact as follows:

\begin{Prop} If $G$ is a uniquely 4-colorable maximal planar graph, then $G$ is a pure tree-coloring graph.

\end{Prop}

Hence, all of recursive maximal planar graphs are pure tree-coloring graphs. It is hard to know the basic characterization of pure tree-coloring graphs with minimum degrees at least 4 which are main researching objects in this paper. We have already known that the 9-mirror graph and the icosahedron (see Figure 5.3) are pure tree-coloring graphs. This result is also proved by Xu \cite{Xu2005a,Xu2005b}. What's more, he claimed that there were only two pure tree-coloring maximal planar graphs with orders at most 43 and minimum degrees at least 4. Unfortunately, his claim is incorrect because the maximal planar graph (\emph{13-mirror graph}) shown in Figure 5.6 is also a pure tree-coloring graph.

\begin{center}

        \includegraphics [width=260pt]{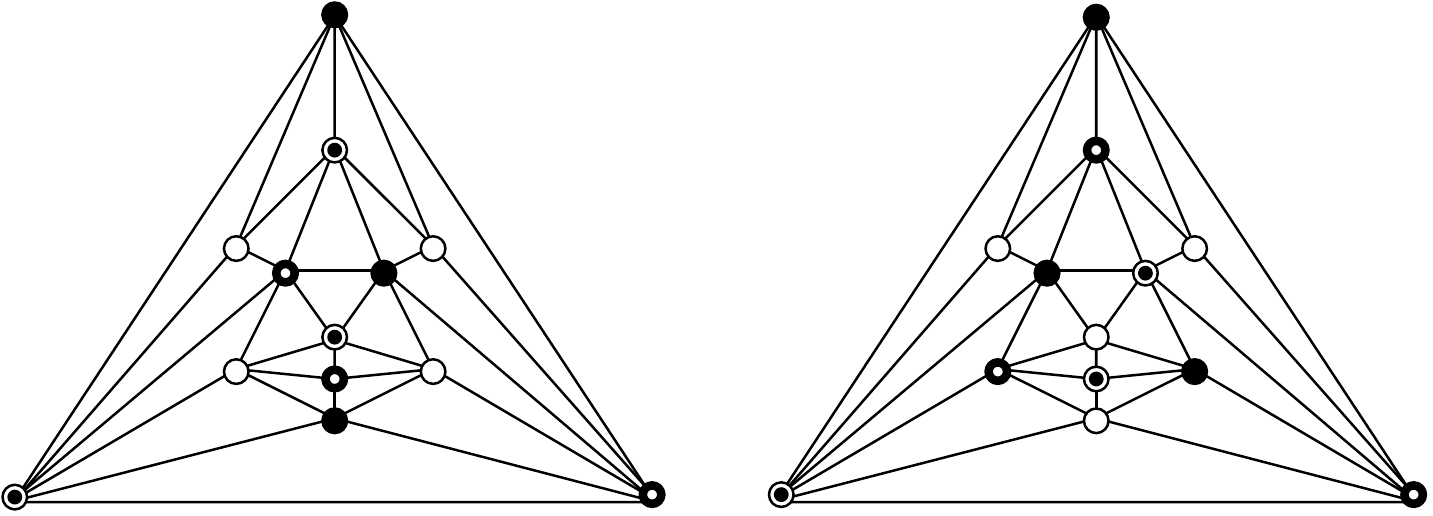}

        \includegraphics [width=260pt]{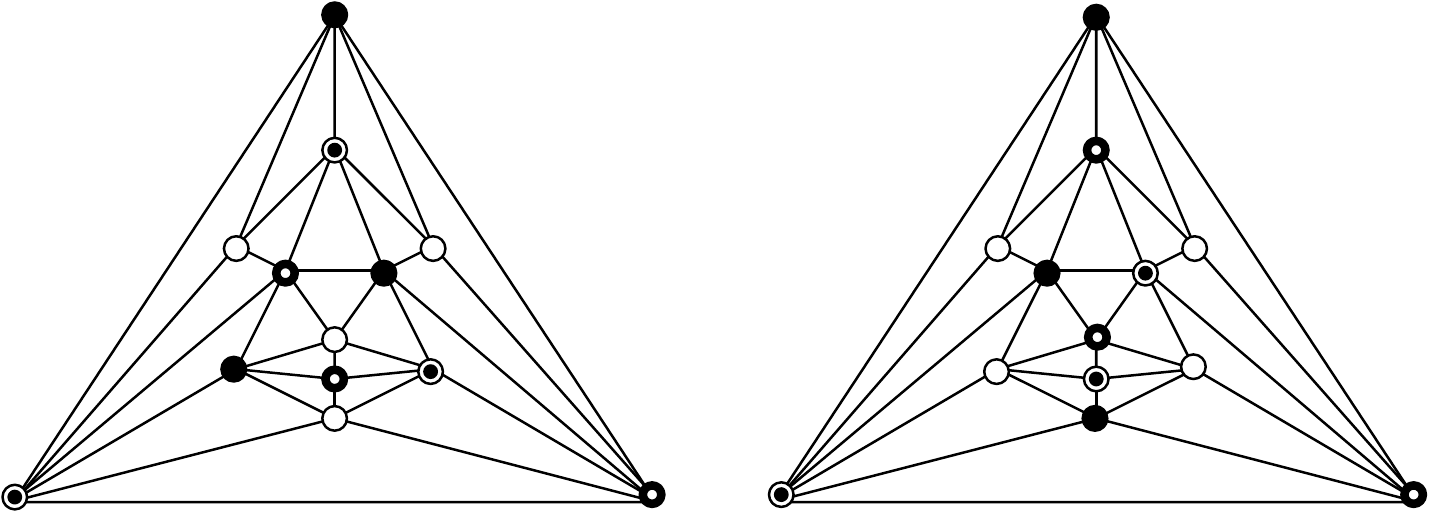}

        \textbf{Figure 5.6.} The third pure tree-coloring graph (13-mirror graph) and its 4-colorings
\end{center}

Naturally, there is an important problem as follows:

\begin{problem} What is the characterization of a pure tree-coloring graph whose minimum degree is not less than $4$? And how many such graphs are there?
\end{problem}

The answer of Problem 5.1 might be: a maximal planar graph $G$ with minimum degree at least 4 is a pure tree-coloring graph if and only if $G$ is icosahedron or $(4k+1)$-mirror graph. The detailed study of this problem will be given in the later sections.

For pure cycle-coloring graphs and impure coloring graphs, we also propose two problems as follows:

\begin{problem}
What are the necessary and sufficient conditions for a maximal planar graph to be a pure cycle-coloring graph?
\end{problem}
\begin{problem}
What are the necessary and sufficient conditions for a maximal planar graph to be an impure coloring graph?
\end{problem}

In the same way, the detailed discussions of these two problems will also be given in the subsequent sections.

\subsection{Equivalency of colorings  between tricolored induced subgraphs and maximal planar graphs}

\quad\quad For a given 4-colorable maximal planar graph $G$, let $f$ be a 4-coloring of $G$ and the color class partition of $f$ be $\{V_1,V_2,V_3,V_4\}$, where $V_i$ denotes the set of vertices assigned color $i$. Obviously, when three color classes of them are determined, the last one is also determined uniquely. So, we only need to make clear the 3-coloring structure of the tricolored induced subgraph that are induced by any three classes partition, such as $V_1,V_2,V_3$.

For the sake of convenience, here we introduce a definition of big-cycle. Let $C$ be a cycle of a  planar graph $G$. If $C$ is a \emph{facial cycle} (the boundary of a face) and has length not less than 4, then we call $C$ a \emph{big-cycle}.

The following gives some examples that illustrate the equivalency of colorings between a maximal planar graph and its tricolored induced graphs keeping  each of  big-cycles  colored with at most three colors. For the first graph shown in Figure 5.7, it has three different 4-colorings totally (see Figures 5.7(a),(b) and (c)). If we denote by $V_4$ the set consisting of vertices received by color 1, then the three 3-colorings of $G-V_4$ corresponding to Figures 5.7(a),(b) and (c) are shown as Figures 5.7(a'),(b') and (c'). For the first graph shown in Figure 5.8(a), we delete the vertices received by color 1, and then obtain its tricolored induced subgraph  and the corresponding 3-coloring (see Figure 5.8(b)). Figures 5.8(c) and (d) exhibit two 4-colorings of this subgraph such that each of its  big-cycles  is colored with three colors.

\begin{center}

        \includegraphics [width=320pt]{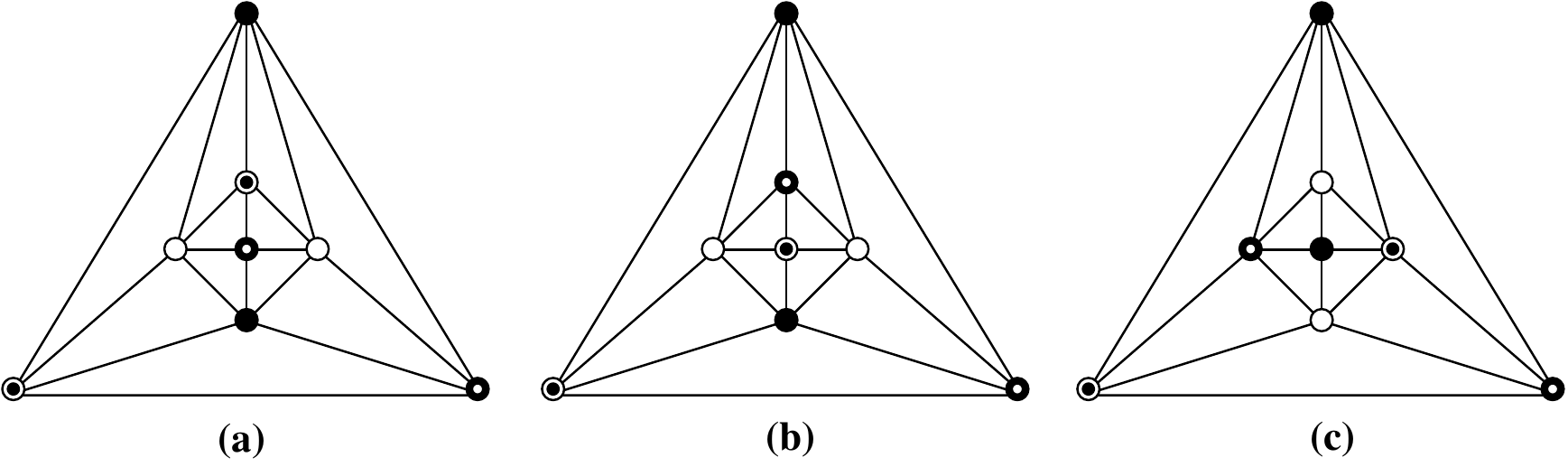}
        \includegraphics [width=320pt]{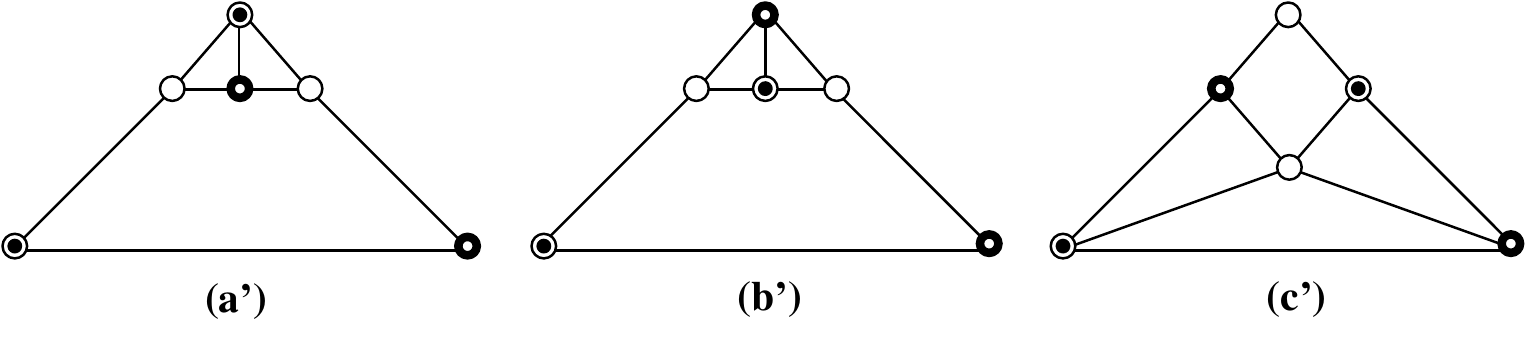}

        \textbf{Figure 5.7.} All the 4-colorings of a maximal planar graph of order 8 and  the corresponding 3-colorings of their tricolored induced subgraphs
\end{center}

\begin{center}

         \includegraphics [width=340pt]{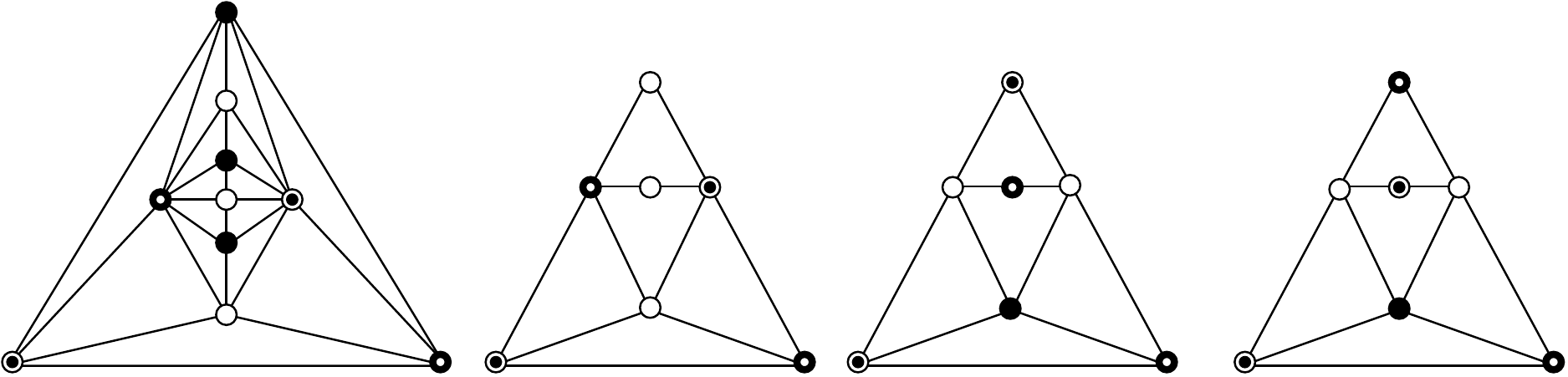}

        \textbf{Figure 5.8.} A 4-coloring of a maximal planar graph with order 10 and the corresponding colorings of its tricolored induced subgraphs
\end{center}

The above two examples show clearly that we have the intuitive and understandable advantage when we study 4-coloring problem on $G-V_4$ rather than on $G$ directly. So, when we study the 4-coloring problem of a maximal planar graph $G$, we need only research the 4-coloring problem of $G-V_4$. However, in the process of studying this problem, we need to pay attention to the following three points:

First, in terms of the choice of $V_4$, we should try to choose such a color class that contains maximal vertices, so that $G-V_4$ can become simple as much as possible.

Second, because the number of triangles in a maximal planar graph of order $n$ is $2n-4$, the number of triangles in $G-V_4$ is
$$
2n-4-\sum\limits_{v\in V_4}d_G(v)
 \eqno{(5.1)}
$$
For example, Figure 5.7 exhibits a maximal planar graph with order $n=8$ and the sum of degrees of the two vertices in $V_4$ is 10, so the number of triangles in $G-V_4$ is $2\times 8-4-10=2$. Observe another maximal planar graph of order 10 shown in Figure 5.8, similarly, we can calculate the number of triangles in $G-V_4$ is 3 by Form (5.1).

Third, the following equality
$$
G[V_1\cup V_2\cup V_3]=G[1,2]\cup G[1,3]\cup G[2,3]
 \eqno{(5.2)}
$$
holds. So, we should study the structure and property of the tricolored induced subgraph $G[V_1\cup V_2\cup V_3]$ from any two bicolored induced subgraphs, such as $G[1,2]\cup G[1,3]$. Studying gradually the structure and property of $G[V_1\cup V_2\cup V_3]$ is the basic idea in this section.

\begin{theorem2}\label{th5.4}
For a 4-colorable maximal planar graph $G$, suppose that $f$ is a 4-coloring of $G$, and the color classes partition of $f$ is $\{V_1,V_2,V_3,V_4\}$. Then

\textcircled{1} $f$ is a tree-coloring of $G$ if and only if $f$ is a tree-coloring  subject to $G[V_1\cup V_2\cup V_3]$;

\textcircled{2} $f$ is a cycle-coloring of $G$ if and only if $f$ is a cycle-coloring (or a disconnected coloring) subject to $G[V_1\cup V_2\cup V_3]$;\\
where``subject to $G[V_1\cup V_2\cup V_3]$" refers to the color class partition of $f$ only presenting on $G[V_1\cup V_2\cup V_3]$; the disconnected coloring means that there are disconnected bicolored induced subgraphs in $G[V_1\cup V_2\cup V_3]$ under the coloring $f$.
\end{theorem2}

Obviously, Theorem 5.4 can bring us some convenience when we judge whether a coloring is a tree-coloring. So, in the later investigation, we will mainly consider the coloring structure of tricolored induced subgraphs.

\subsection{Structure of the union of two bicolored induced subgraphs}
\quad\quad In this section we introduce the concept of \emph{fence}, and prove that the union of any two
bicolored induced subgraphs of 4-colorable maximal planar graph is a fence.
Furthermore, some special fences and the characteristic of maximal planar graphs
corresponding to them are discussed.

\subsubsection{General theory}

\quad\quad The degree of a face is the number of
edges in its boundary, cut edges being counted twice. Let $G$ be a planar graph. If
the degree of every face of $G$ is even and not less than 4, then $G$ is called a \emph{fence}.
The graphs shown in Figures 5.9(a), (c), (d) are fences, but the graph shown in Figure 5.9(b)
is not a fence, because there exists an odd cycle that encloses a face in this graph. For a fence $G$,  it may not have any suspending vertices, of course, it  may also contain suspending vertices. Here the suspending vertices refer to the vertices with degree 0 or 1.
If there exists a suspending vertex $v$ in $G$, then the subtree containing $v$  maybe connect with a cycle  by a common vertex, say $u$ and called a \emph{weld-vertex}. Denote by $t$
the distance between $u$ and $v$, namely there exists a path of length $t$ between them.
Choose a maximum $t$ and refer to $G$ as a \emph{$t$-fence}.  If $G$ has no suspending vertex, we call it a $0$-fence.
The graph shown in
Figure 5.9(a) is a 1-fence; the graph shown in Figure 5.9(c) is a 0-fence; the graph shown in
Figure 5.9(d) is a 2-fence; the graph shown in Figure 5.9(e) is a 3-fence. If there exists no
path between a suspending vertex and any cycle of $G$, that is to say, the graph $G$ is disconnected,
then $G$ is called an $\infty$-fence. The graph shown in Figure 5.9(f) is an $\infty$-fence.

\begin{center}
        \includegraphics [width=260pt]{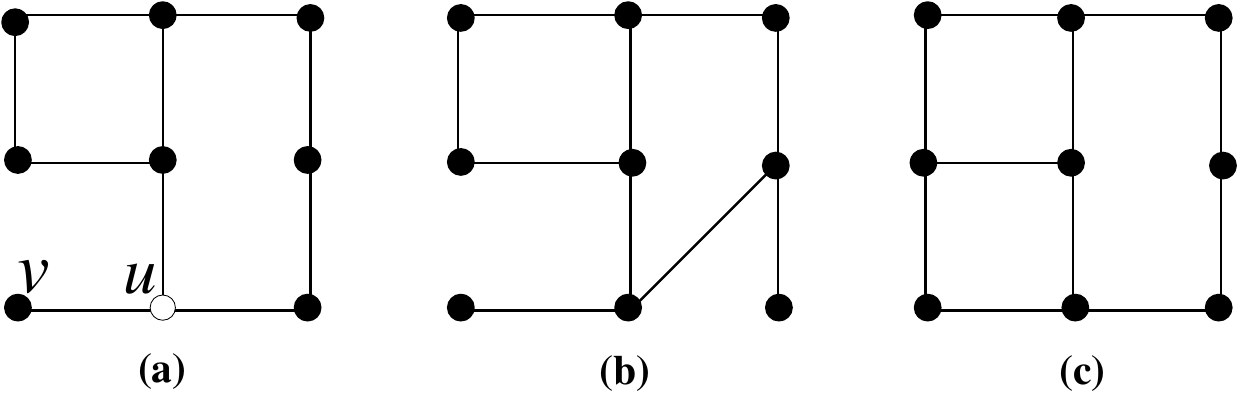}

        \includegraphics [width=260pt]{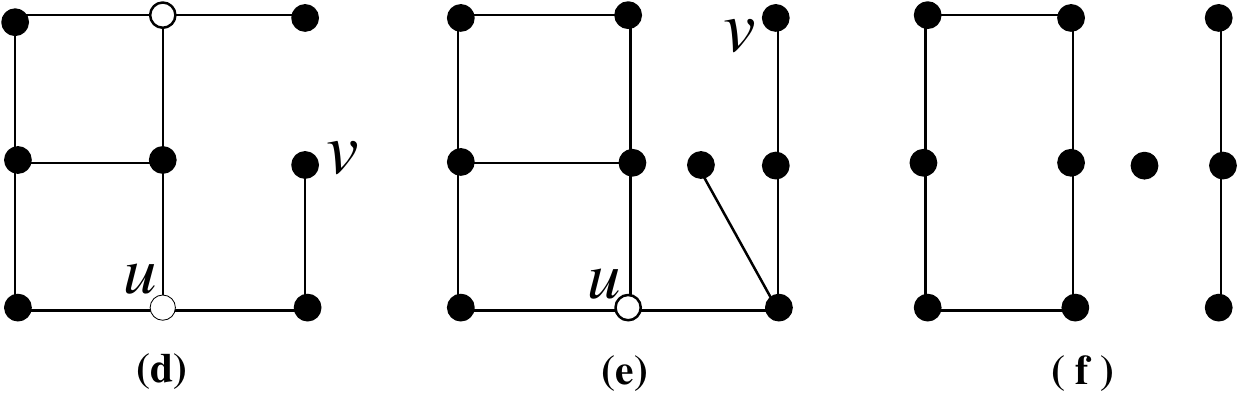}

        \textbf{Figure 5.9.} The illustration of the concept of fence
\end{center}

\begin{theorem2}\label{th5.5}
Let $G$ be a $4$-colorable maximal planar graph, $C(4)=\{1,2,3,4\}$ the color set. For any  $4$-coloring $f$ of $G$ and the union of any two bicolored induced subgraphs with a common color, say $G[1,2]\cup G[1,3]$, we have

 $(1)$ $G[1,2]\cup G[1,3]$ has no odd cycle;

 $(2)$ If the order of both $G[1,2]$ and $G[1,3]$ are at most $3$, then the graph $G[1,2]\cup G[1,3]$ is a cycle under only one case and a tree in other cases;

$(3)$ If $f$ is a tree-coloring and the order of  $G[1,2]$ or $G[1,3]$ is at least $4$, then $G[1,2]\cup G[1,3]$
is a $1$-fence or a $0$-fence;

$(4)$  If $f$ is a tree-coloring, then every suspending vertex of $G[1,2]\cup G[1,3]$ must be adjacent to the vertices
colored by the common color $1$.
\end{theorem2}
 \begin{proof}
 (1) Assume that $G[1,2]\cup G[1,3]$ has an odd cycle $C$. Since the set of colors appeared on the cycle $C$ must contain color 1,
2 and 3, so $C$ contains not only the 1-2 edges and 1-3 edges, but also the 2-3 edges. But $G[1,2]\cup G[1,3]$ has no edges in $G[2,3]$, it is a contradiction, where $i-j$ edge denotes the edge whose two end vertices are colored by color $i$ and $j$ respectively.

 (2) If the order of both $G[1,2]$ and $G[1,3]$ are not more than 3, then all 4 cases are shown in Figure 5.10. It is easy
 to see that only in one case $G[1,2]\cup G[1,3]$ is a cycle, and is a tree in other cases.

  \begin{center}

        \includegraphics [width=260pt]{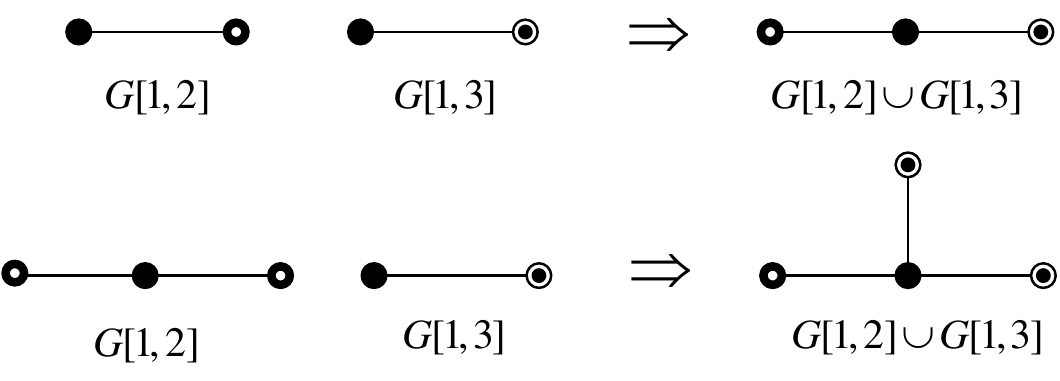}

        \includegraphics [width=260pt]{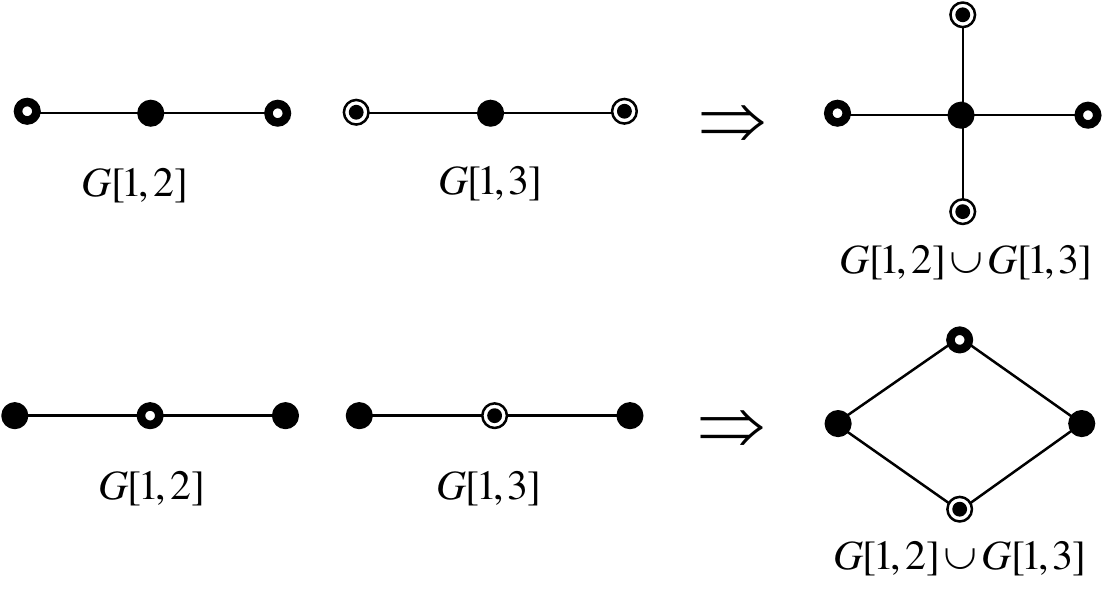}

        \textbf{Figure 5.10.} Four cases of $G[1,2]\cup G[1,3]$ in Theorem 5.5(2)
  \end{center}

(3) Based on the cases (1) and (2), we study the case in which the order of  $G[1,2]$ or $G[1,3]$ is at least 4.
Suppose $v$ is a suspending vertex of $G[1,2]\cup G[1,3]$ and the distance of $v$ to the nearest cycle in $G[1,2]\cup G[1,3]$ is 2. The unique vertex adjacent to $v$ in $G[1,2]\cup G[1,3]$ is denoted by $v^{\prime}$.
 If $f(v)=1$, since $f$ is a tree-coloring, both $G[1,2]$ and $G[1,3]$ are connected, so $v$ is adjacent to some vertices
 colored with color 2 and 3 respectively. Thus, $v$ is adjacent to at least two
 vertices in $G[1,2]\cup G[1,3]$, it is a contradiction to the fact that $v$ is a suspending vertex.   If $f(v)=2$ (or $f(v)=3$),
 then $f(v^{\prime})=1$ and the vertex $w$ ($\neq v$) adjacent to $v^{\prime}$ may be assigned with color 2 or 3. If $f(w)=3$,
 then $G[1,2]$ is disconnected; if $f(w)=2$, then $G[1,3]$ is disconnected. There is a contradiction in these two cases.  So, $G[1,2]\cup G[1,3]$
is either a 1-fence or a 0-fence.

(4) This case is obvious.

The proof of the theorem is complete.
\end{proof}

The results of Theorem \ref{th5.5} show that for a 4-coloring $f$ of a 4-colorable maximal planar graph $G$,  $G[1,2]\cup G[1,3]$ doesn't contain odd cycles, that is to say, it is a tree or contains only even cycles; furthermore, for the structure including the even cycles, we prove that $G[1,2]\cup G[1,3]$ is a 1-fence or a  0-fence when $f$ is a tree-coloring. In addition, for a fence, if it does contain suspending vertices, they must be adjacent to the vertices colored by the common color 1.

In Figure 5.11(b), it is easy to see that the resulted graph by adding a new vertex $v$ to every face of degree at least 4 and connecting
$v$ to all vertices on the cycle of the face including $v$ is a maximal planar graph $G$ with $\delta(G)\geq 4$.  The graph shown in Figure 5.11(a) exhibits the case that the union of two bicolored induced subgraphs of a
maximal planar graph $G$ with $\delta(G)\geq 4$ has no suspending vertex.

\begin{center}

        \includegraphics [width=280pt]{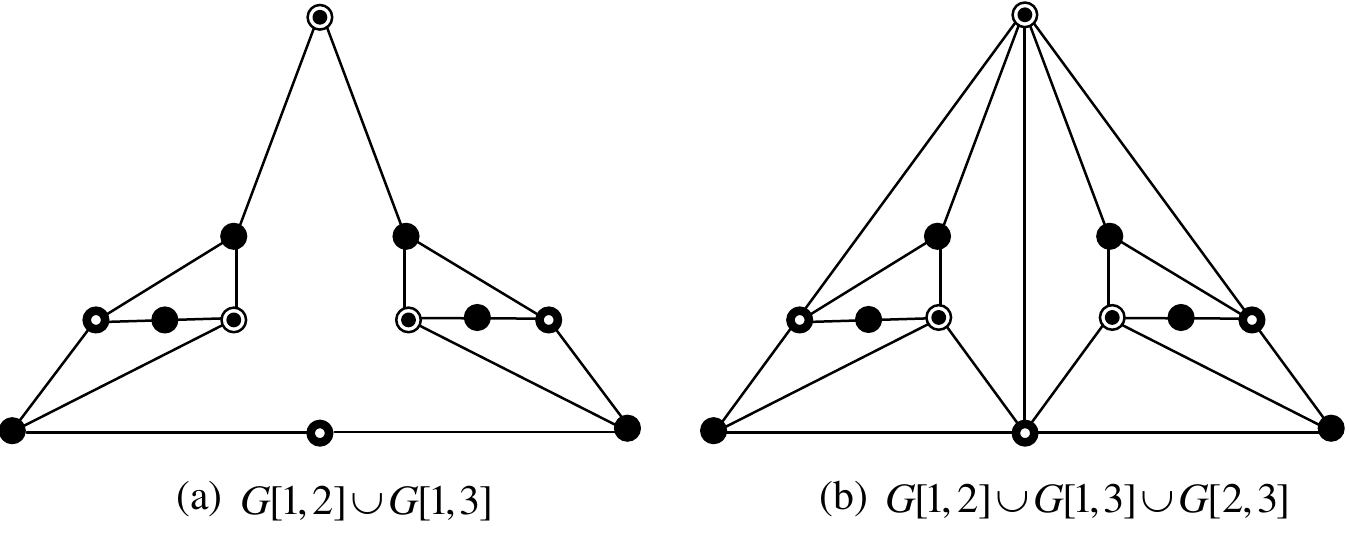}

        \textbf{Figure 5.11.} The illustration for the union of two bicolored induced subgraphs has no suspending vertex
  \end{center}

\subsubsection{Structure of  bicolored induced subgraphs and their union of pure tree-coloring graphs}

The necessary and sufficient condition, for which a maximal planar graph $G$ is a pure
tree-coloring graph is that $G$ is the icosahedron or $(4k+1)$-mirror graph, $k\geq2$. This statement will be discussed in the following chapters.
Now, we analyze the structure of the union of bicolored induced subgraphs of the icosahedron, the 9-mirror graph and the 13-mirror graph.

For the 9-mirror graph (see the second and the third graphs in Figure 5.5), it is easy to prove that all its bicolored induced subgraphs are paths,
and the length of them has only two kinds: one is 4 and another is 3. Concerning the structure of the union of two bicolored induced subgraphs with a common color, there are three cases: $\textcircled{1}$ the union consists of two paths of length 4 (see Figure 5.12(a)): it is a fence without suspending vertices and includes two cycles of length 6 and one cycle of
length 4. $\textcircled{2}$ The union consists of two paths with length 3 and 4 respectively: it is a 1-fence with 3 suspending vertices and
a cycle of length 4 (see Figure 5.12(b)). $\textcircled{3}$ The union consists of two paths of length 3: this union is a 1-fence with 2 suspending vertices and a cycle of length 4 (see Figure 5.12(c)).

For the icosahedron (shown in Figure 5.3),
all of its bicolored induced subgraphs are paths of length 5 (see Figure 5.12(d)). The union of two bicolored induced subgraphs, which consists of two paths of length 5 with a common color (see Figure 5.12(d)), is formed by adding two suspending vertices on the basis of the graph shown in Figure 5.12(a).

\begin{center}

         \includegraphics [width=260pt]{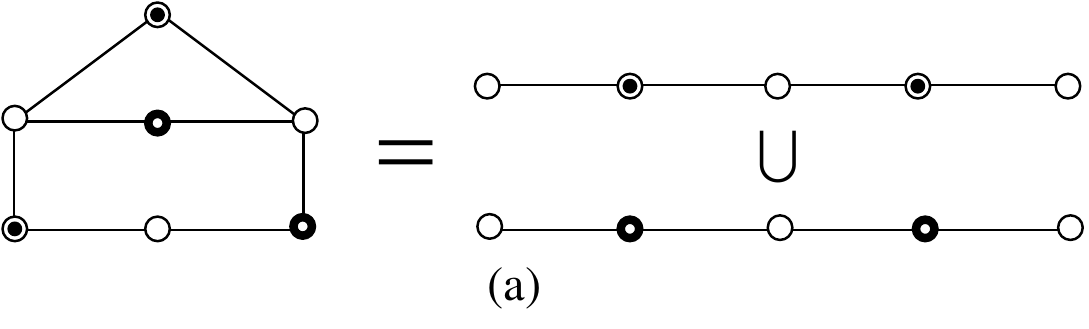}

         \vspace{0.2cm}
         \includegraphics [width=260pt]{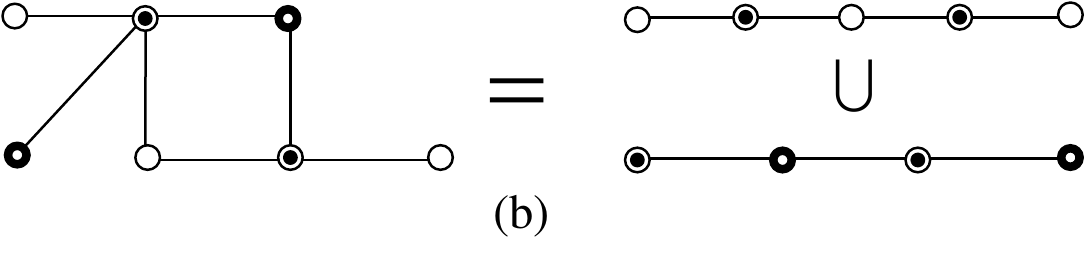}

         \vspace{0.2cm}
         \includegraphics [width=260pt]{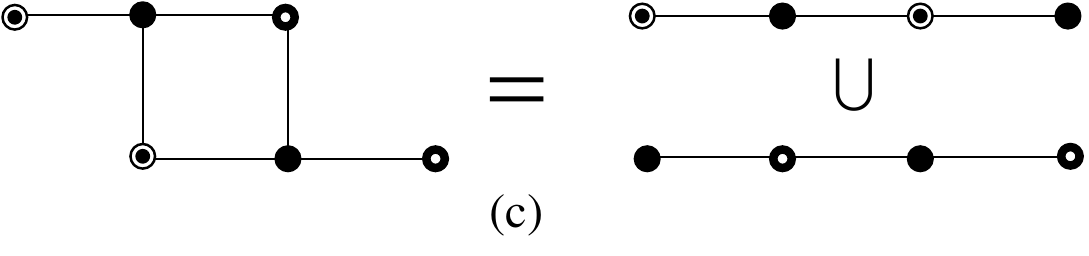}

        \vspace{0.2cm}
         \includegraphics [width=260pt]{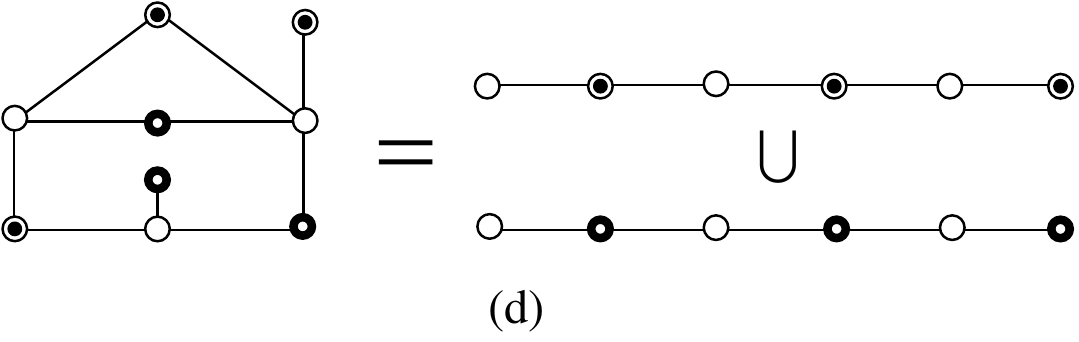}

        \textbf{Figure 5.12.} Structure analysis of all bicolored induced subgraphs of the icosahedron and 9-mirror graph
  \end{center}

 The 13-mirror graph shown in Figure 5.6 has three kinds of bicolored induced subgraphs:  two of them are paths with lengthes  5 and 6, the rest is a tree obtained by adding a vertex with 1-degree to the middle vertex of a path with length 4. Obviously, the maximum degree of this tree is 3. Concerning the structure of the union of two bicolored induced subgraphs with a common color, there are 11 cases shown  in Figures 5.13(a)$\sim$ 5.13(k) respectively.

\begin{center}

         \includegraphics [width=360pt]{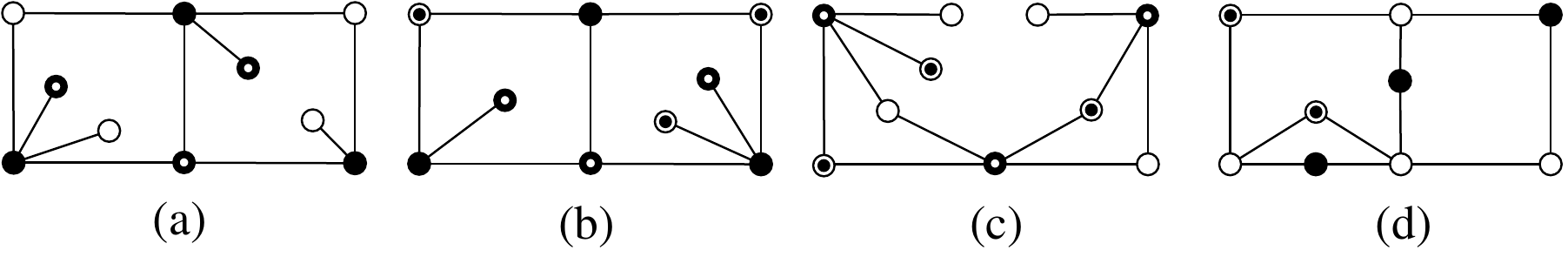}

         \includegraphics [width=360pt]{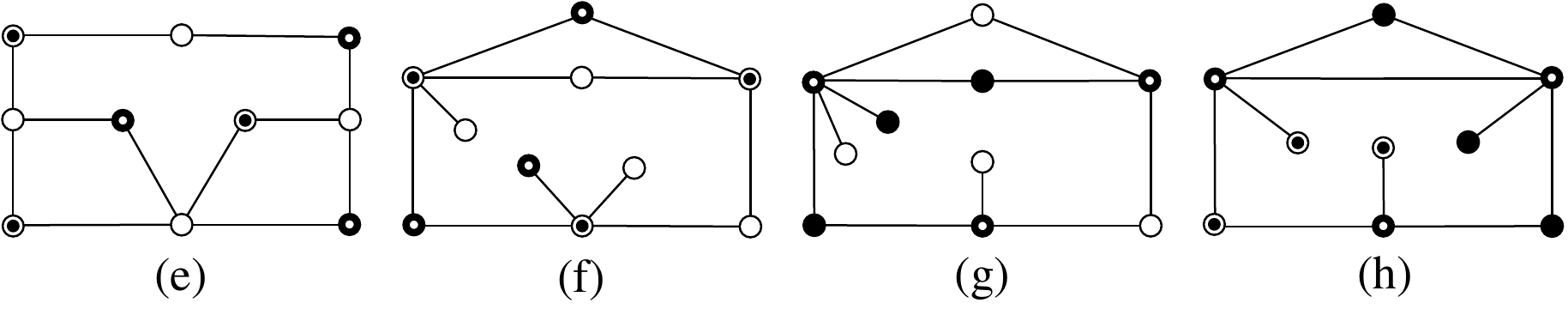}

          \includegraphics [width=360pt]{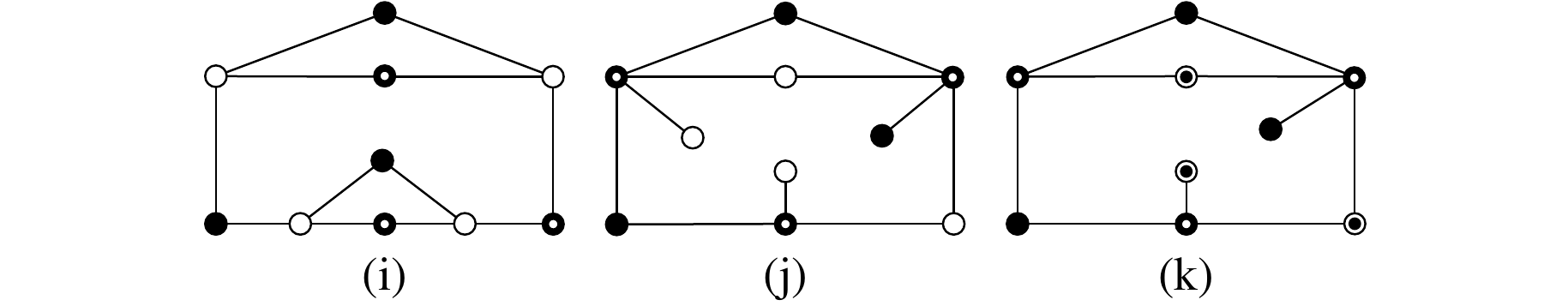}

        \textbf{Figure 5.13.} Structure of all bicolored induced subgraphs of the 13-mirror graph
  \end{center}

Furthermore, the icosahedron, the 9-mirror graph and the 13-mirror graph have several properties as follows:

 $\textcircled{1}$ Each of  bicolored induced subgraphs of the icosahedron is a path whose length is 5; each of bicolored
 induced subgraphs of the 9-mirror graph  is also a path, the length of which is 5, 4 or 3; each of bicolored
 induced subgraphs of the 13-mirror graph is a path with length 6 or 5, or a tree with maximum degree 3.

$\textcircled{2}$ For these three graphs, there are at most two cycles of length 4 in the union of two bicolored
induced subgraphs with a common color and the length of the longest cycle is 8.

 $\textcircled{3}$ These three graphs are symmetrical strongly.

 \subsubsection{Structure of the union of two bicolored induced subgraphs whose common color appears on a cycle}

This subsection is devoted to the structure of a class of specific graphs, in which the union of two bicolored induced subgraphs
whose vertices colored the common color are totally on a cycle.

Let $G$ be a 4-coloring maximal planar graph, $f$ a tree-coloring of $G$, $G[1,2]$ and $G[1,3]$  two bicolored induced subgraphs having a common color. Then
$G[1,2]\cup G[1,3]$ is a 1-fence or a 0-fence.
For the graph $G[1,2]\cup G[1,3]$, suppose the number of vertices assigned  color 1, 2 and 3 are $a, b$ and $c$, respectively.
If all of the vertices assigned  color 1 are in a cycle, then $G[1,2]\cup G[1,3]$ has the following properties:

\textbf{Property 1.} Besides  edges incident with  suspending vertices, each edge is either in one or in two cycles of length 4;

\textbf{Property 2.} There are $a-i$ cycles of length $2i+2$, $i=1,2,\cdots,a-1$, where the cycles of length 4 are adjacent alternately and
each pair of adjacent two cycles have just a common edge;

\textbf{Property 3.}  The number of vertices with degree 1 is $b+c-a$.

It is easy to prove
 \begin{theorem2}\label{th5.6}
Let $G$ be a $4$-colorable maximal planar graph, $C(4)=\{1,2,3,4\}$ the color set and $f$ a $4$-coloring of $G$.
If the order of
$G[1,2]$ or $G[1,3]$ is $3$, say $G[1,2]$, then  $\delta(G)=3$ if there exist  odd suspending vertices in $G[1,2]\cup G[1,3]$, otherwise, $G[1,2]$ contains just two vertices receiving color $1$.
 \end{theorem2}

\subsection{Construction of tricolored induced subgraphs}

Theorem 5.4 tells us that for a 4-coloring $f$ of a maximal planar graph, its properties can be  characterized by its tricolored
induced subgraphs completely.
In fact, this subgraph is the union of three bicolored induced subgraphs. We have made clear the structure on the union of two bicolored induced subgraphs, which is the basis of studying the structure of tricolored induced graphs.

Without loss of generality, we still denote by $G[1,2]\cup G[1,3]$ the union of any two bicolored induced subgraphs here and we will construct the tricolored induced subgraphs based on it by connecting the edges between the vertices receiving color 2 and color 3.
We should remark that if $G[1,2]\cup G[1,3]$ is connected, there is always at least a path $uu_1^{1}u_2^{2-3}u_3^{1}u_4^{2-3}\cdots u_k^{1}u^{\prime}$ for every pair of vertices $u$ and $u'$ that are assigned color 2 and 3 respectively, where
$u_i^{1}$ $(1\leq i \leq k)$  denotes the number of vertices colored 1, $u_i^{2-3} (2\leq i \leq k)$ denotes
the number of vertices colored by 2 or 3. Obviously, there are odd vertices on the path between the vertices $u$ and $u^\prime$. So, we have

\begin{theorem2}\label{th}
Let $G$ be a $4$-colorable maximal planar graph, $C(4)=\{1,2,3,4\}$  the color set, $f\in C_4^0(G)$ and $\{u,u^\prime\}$ a pair of vertices colored by $2$ and $3$ respectively, then any path between
$u$ and $u^\prime$ has odd vertices in the connected subgraph $G[1,2]\cup G[1,3]$.
\end{theorem2}

This theorem shows that: in order to construct $G[V_1\cup V_2\cup V_3]=G[1,2]\cup G[1,3]\cup G[2,3]$ based on
$G[1,2]\cup G[1,3]$, when $G[2,3]$ is a tree with $q$ edges, we need only to connect $q$ edges which can form a tree between the vertices colored by 2 and 3 in $G[1,2]\cup G[1,3]$.

Since $G[1,2]\cup G[1,3]$ has only even cycles and the length of any path between vertices colored by 2 and 3 is an odd number, so every edge of $G[2,3]$ contributes to $G[V_1\cup V_2\cup V_3]$ at least one odd cycle. Now it is discussed in detail
in the following.

 \textbf{Case 1.} Cycle-cycle edge: an edge in $G[2,3]$ and its two ends $u$ and $u'$ are on some cycles in $G[1,2]\cup G[1,3]$. For this case, by Theorem 5.7
 there are at least two different paths with even lengthes. So $G[V_1\cup V_2\cup V_3]$
 is contributed at least two odd cycles when  $u$ and $u'$ are joined in $G[1,2]\cup G[1,3]$.

\textbf{Case 2.} Cycle-suspending edge: an edge in $G[2,3]$  is formed by joining a vertex on a cycle and a suspending vertex in $G[1,2]\cup G[1,3]$.

\textbf{Case 3.} Suspending-suspending edge: an edge in $G[2,3]$ is formed by joining two suspending vertices in $G[1,2]\cup G[1,3]$.

\textbf{Example 5.1.} For the 4-colorable maximal planar graph $G$ shown in Figure 5.14(a), let $f$ be a 4-coloring of $G$. The subgraph shown in Figure 5.14(b) is obtained by deleting  all vertices assigned color 4; the subgraph shown in Figure 5.14(c) is $G[1,2]\cup G[1,3]$, in which there are only two suspending vertices $v_3$ and $v_8$. Thus, the edges $v_4v_8$ and $v_2v_3$ are two so-called
cycle-suspending edges, and the edge $v_3v_8$ is a so-called suspending-suspending edge (see Figure 5.14(b)).

Since there exist two paths $v_4v_7v_8$ and $v_4v_1v_2v_7v_8$ from $v_4$ to $v_8$ in graph $G[1,2]\cup G[1,3]$, so the edge $v_4v_8$
contributes  to $G[V_1\cup V_2\cup V_3]$ two cycles $v_4v_8v_7v_4$ and $v_4v_1v_2v_7v_8v_4$. This also illustrates Theorem 5.7.
Similarly, the edge $v_2v_3$ contributes to $G[V_1\cup V_2\cup V_3]$  two odd cycles $v_2v_1v_3v_2$ and  $v_2v_7v_4v_1v_3v_2$. In addition, because there exist two paths $v_3v_1v_4v_7v_8$ and $v_3v_1v_2v_7v_8$ from $v_3$ to $v_8$ in graph $G[1,2]\cup G[1,3]$,
the edge $v_3v_8$ contributes to $G[V_1\cup V_2\cup V_3]$ two cycles $v_3v_1v_4v_7v_8v_3$ and $v_3v_1v_2v_7v_8v_3$.

In graph $G[1,2]\cup G[1,3]$, since there exist two paths $v_3v_1v_4$ and $v_3v_1v_2v_7v_4$ from $v_3$ to $v_4$, so the two edges
$v_4v_8$ and $v_8v_3$ in $G[2,3]$ contribute to $G[V_1\cup V_2\cup V_3]$ two even cycles $v_3v_1v_4v_8v_3$ and $v_3v_1v_2v_7v_4v_8v_3$. Similarly, because there exist two paths $v_2v_7v_8$ and $v_2v_1v_4v_7v_8$ from $v_2$ to $v_8$,  the two edges
$v_2v_3$ and $v_3v_8$ in $G[2,3]$ contribute to $G[V_1\cup V_2\cup V_3]$ two even cycles $v_2v_7v_8v_3v_2$ and $v_2v_1v_4v_7v_3v_8v_2$.

\begin{center}

         \includegraphics [width=380pt]{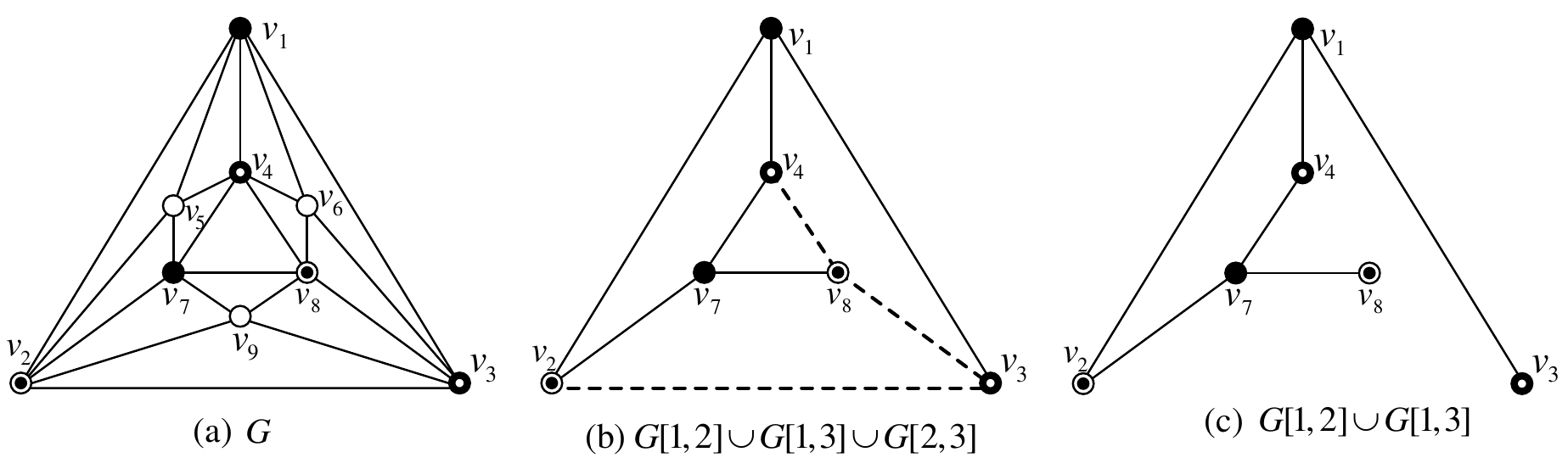}

        \textbf{Figure 5.14.}  An illustration for cycle-cycle edge, cycle-suspending edge and suspending-suspending edge
  \end{center}

Furthermore, we consider  three edges $v_4v_8$, $v_8v_3$ and $v_3v_2$ in $G[2,3]$. Since there exist two paths $v_2v_7v_4$ and $v_2v_1v_4$
from $v_2$ to $v_4$ in  $G[1,2]\cup G[1,3]$,  three edges
$v_4v_8$, $v_8v_3$ and $v_3v_2$ in $G[2,3]$ contribute to $G[V_1\cup V_2\cup V_3]$ two odd cycles $v_2v_7v_4v_8v_3v_2$ and $v_2v_1v_4v_8v_3v_2$.

From the discussion of Example 5.1, we can obtain the following theorem.

\begin{theorem2}\label{th}
Let $G$ be a $4$-colorable maximal planar graph, $C(4)=\{1,2,3,4\}$ the color set and $f\in C_4^0(G)$. Suppose the number of the paths
from a vertex $u$ to another vertex $u^\prime$ is $q$ in $G[1,2]\cup G[1,3]$, for any path $P$ of length $p$ in $G[2,3]$. Then

$\textcircled{1}$ $q$ is an even number;

$\textcircled{2}$ The number of odd (even) cycles including $P$ in $G[1,2]\cup G[1,3]\cup P$ is $q$ more than the number of odd cycles in $G[1,2]\cup G[1,3]$ when $p$ is an odd (even) number.
 \end{theorem2}

 It follows Theorem 5.8 that the coloring structure of a 4-colorable maximal planar graph $G$ under a 4-coloring $f$ has been shown clearly.
Let $V_1,V_2,V_3,V_4$ be four independent sets of $G$ based on $f$, then

$\textcircled{1}$ the coloring structure of $G$ corresponding to $f$ is equivalent to $G[V_1\cup V_2\cup V_3]$;

$\textcircled{2}$ $G[1,2]\cup G[1,3]$ is a fence;

$\textcircled{3}$ $G[V_1\cup V_2\cup V_3]$ is obtained by adding continuously the edges of $G[2,3]$ in $G[1,2]\cup G[1,3]$.
$G[V_1\cup V_2\cup V_3]$ is contributed several odd cycles when an edge of $G[2,3]$ is added to $G[1,2]\cup G[1,3]$. Furthermore, $G[V_1\cup V_2\cup V_3]$ is contributed
 several odd (even) cycles including $P$ when an odd (even) path $P$ of $G[2,3]$ is added to $G[1,2]\cup G[1,3]$.

\section{Black-White coloring, and necessary and sufficient conditions for 2-colorable cycle}

In this Chapter, we will propose a new coloring method, called Black-White coloring.
For a maximal planar graph $G$, a Black-White coloring of $G$ is to assign only two colors black and white, to the vertices of $G$. Beginning with an even cycle, either black or white is colored to the  vertices of the cycle according to a definite rule. Based on this method, if  the vertices of $G$ are colored black or white, then we can deduce that $G$ is 4-colorable if no odd-cycles are contained in the induced subgraphs by the vertices that are colored the same color (black or white). In fact, the Black-White coloring provides a subset of 4-colorings for maximal planar graphs.

Depending on the Black-White coloring, we will show a necessary and sufficient condition for an even cycle in a maximal planar graph to be 2-colorable. For this purpose, we introduce some new concepts, such as closed-maximal planar graphs, opened-maximal planar graphs, semi-maximal planar graphs and 2-colorable cycles. Furthermore, we  research deeply  the characteristics of even-cycles in a maximal planar graph and discuss the enumeration of even-cycles.

Let $G$ be a maximal planar graph and $C$ a cycle in $G$ with length not less than 4. We refer to a subgraph of $G$ that is induced by the vertex set consisting of  vertices of $C$ and the inside component of $C$, as a \emph{semi-maximal planar graph}. In other words, the so-called semi-maximal planar graphs are just a kind of special planar graphs such that each of them contains a unique face bounded by a cycle of length not less than 4 and other faces are triangles. We call every graph $G$ in this sort of graphs a \emph{semi-maximal planar graphs on $C$ }, written $G^C$. Obviously, for $G$ and $C$, there are just two semi-maximal planar graphs on $C$, and we call them the semi-maximal planar graphs on $C$ of $G$.

\subsection{Characteristics and distribution of the even-cycles}

Even-cycle is the basic element of cycle colorings, so considering the characteristics and distribution of the even-cycles in a maximal planar graph is very important to research whether there exist cycle colorings. Suppose that $C=v_1v_2\cdots v_mv_1$ is a cycle with length $m$. A \emph{chord}  on a cycle $C$ is an edge that link two nonadjacent vertices $v_i,v_j$ of $C$. We refer to the cycle containing chords as \emph{chord-cycle}. Both the longest cycles in Figures 6.1(b) and (c) are chord-cycles.

Thus it is clear that for a maximal planar graph and any its connected subgraph $H$, the subgraph induced by the neighbor set of $H$ contains either a cycle or  a chord-cycle except for a tree. For example,
in Figure 6.1(d), the cycle $C=v_1v_2v_3v_5v_{10}v_1$ is a chord-cycle (see Figure 6.1(e)). In addition, there are a special kind of chord-cycles, called \emph{one-side cycle}, which contain no vertices inside (or outside).
For example, the cycle $C=v_1v_2v_3v_4v_1$ in Figure 6.1(d) is a one-side cycle (see Figure 6.1(f)).

 \begin{center}
      \includegraphics [width=180pt]{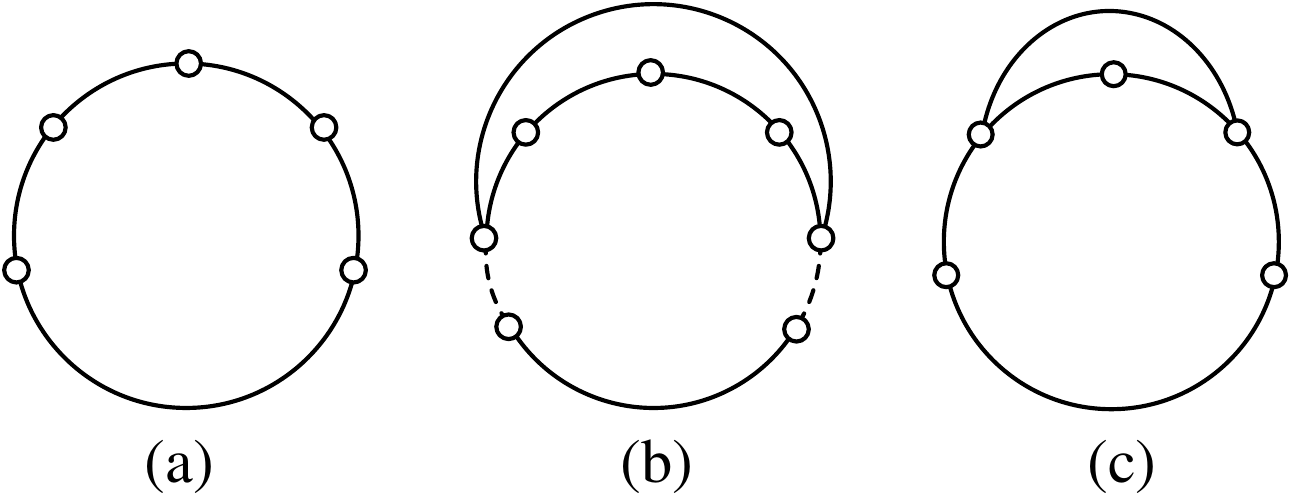}

       \includegraphics [width=320pt]{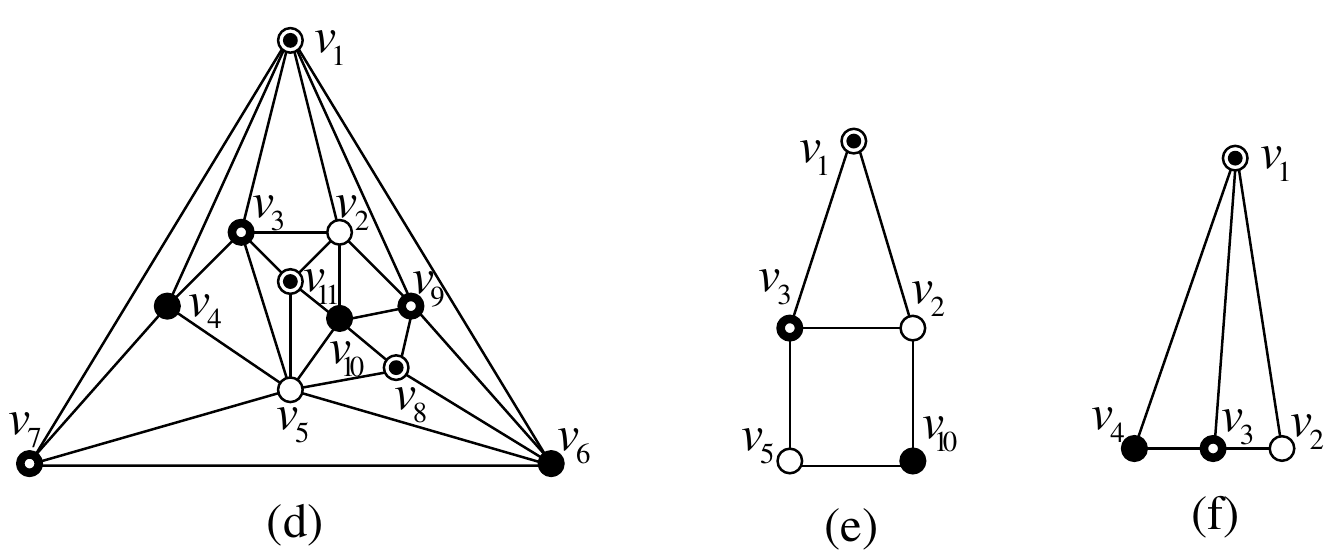}

\textbf{Figure 6.1.} Graphic expression of the definition of cycle, chord-cycle and one-side cycle
  \end{center}

Here, we make an agreement that for $m\geq 4$, a path $P=v_1v_2\cdots v_m$ of a maximal planar graph $G$ is called a \emph{basic path} only if the subgraph induced by $V(P)$ is a path, and a cycle $C=v_1v_2\cdots v_mv_1$ of a maximal planar graph $G$ is called a \emph{basic cycle} only if the subgraph induced by  $V(C)$ is a cycle. Obviously, for a cycle $C$ of a maximal planar graph $G$, $G[V(C)]$ is either a basic cycle or a chord-cycle.

In a maximal planar graph, we are concerned more about the structure, distribution and enumeration of cycles. So, on the basis of the definition defined above, now we are discussing these problems in depth.

\begin{theorem2}\label{th6.1}
Suppose that $G$ is a maximal planar graph with $\delta(G)\geq 4$ and $G$ is not a divisible graph, then the subgraph induced by the set of neighbors of each vertex $v\in V(G)$ is just a cycle with length $d(v)$.
\end{theorem2}

\begin{proof}
Let $\Gamma(v)$ be the neighbor set of $v$. Then there are three possible cases for the induced subgraph $G[\Gamma(v)]$ as follows:

Case 1. $G[\Gamma(v)]$ is a cycle;

Case 2. $G[\Gamma(v)]$ contains triangles;

Case 3. $G[\Gamma(v)]$ contains at least three cycles with length not less than 4.

Case 1 is just the result of this theorem; for the Case 2 and Case 3, it is easy to prove $G$ contains vertices of degree 3 and $G$ is a divisible graph, respectively.
\end{proof}

Denote by $\varsigma(G)$ the set containing all of cycles with lengths at least 4 in $G$, and $\varsigma^{1}(G)$ the set containing all of odd-cycles with lengths at least 5 in $G$, and $\varsigma(G)^{2}$ the set containing all of  even-cycles. Obviously,
$$
\varsigma(G)=\varsigma^{1}(G) \cup \varsigma(G)^{2}.
$$

Actually, Theorem 6.1 says under the condition of  $G$ being a maximal planar graph with $\delta(G)\geq 4$ and a nondivisible graph, each vertex $v$ just contributes a cycle with length $d(v)$ for $\varsigma(G)$. In detail, if $d(v)$ is an odd number, then $v$ just contributes an odd-cycle for $\varsigma^{1}(G)$; if $d(v)$ is an even number, then $v$ contributes an even-cycle for $\varsigma^{2}(G)$.

Naturally, for a maximal planar graph $G$, a correlative problem will be asked about the structure of the subgraph induced by the neighbor set of two adjacent vertices $u,v$ or a connected subgraph having vertices $v_1,v_2,\cdots, v_m$. What are the structures of $G[\Gamma(u,v)]$  and $G[\Gamma(v_1,v_2,\cdots,v_m)]$? Trees? Cycles (odd-cycle or even-cycle)? Chord-cycles or one-side cycles?

Anyway, the length of cycles should be considered largely.

\begin{theorem2}\label{th6.2}
Suppose that $G$ is a maximal planar graph and $P=v_1v_2$ $\cdots v_m$ is a basic path of $G$. If the subgraph induced by the neighbor set of $P$ is a cycle, denoted as $C=G[\Gamma(v_1,v_2,\cdots, v_m)]$, then the length of $C$ is equal to
$$
|C|=\sum\limits_{i=1}^{m}d(v_i)-4(m-1) \eqno{(6.1)}
$$
\end{theorem2}

\begin{proof}
By induction on lengths $m$ of basic paths. When $m=1$, if the subgraph, denoted $H$, which is induced by the neighbors of some vertex $v$ of $G$ is a cycle, then the length of $H$ is equal to $d(v)$ by Theorem 6.1 and the assertion holds. Suppose that it is true for all basic paths of $G$ with fewer than $m$ vertices, where $m\geq 2$, and let $P=v_1v_2\cdots v_{m-1}v_m$ be a basic path of $G$ and the subgraph induced by the neighbor set of $P$ be a cycle, denoted $C$. Choose edge $e=v_{m-1}v_m$ on $P$ and contract $e$ such that two ends $v_{m-1},v_m$ are identified one vertex $v$. Then the resulting graph $G/e$ is a maximal planar graph with a basic path $P^{\prime}=v_1v_2\cdots v_{m-2}v$ and the subgraph of $G/e$ induced by the neighbor set of $P'$ is also the cycle $C$ because the neighbors of $P$ in $G$ and the neighbors of $P^{\prime}$ in $G/e$ are  identical. By the induction hypothesis,
$$
|C|=\sum\limits_{i=1}^{m-2}d_{G/ e}(v_i)+d_{G/ e}(v)-4(m-2) \eqno{(6.2)}
$$
Since
$$
d_{G/e}(v_i)=d_{G}(v_i),i=1,2,\cdots,m-2, d_{G/e}(v)=d_G(v_{m-1})+d_G(v_{m})-4
$$
we obtain
$$
|C|=\sum\limits_{i=1}^{m}d_G(v_i)-4(m-1)
$$
The theorem follows by induction.
\end{proof}
From Formula 6.1, we can see that the parity of the cycle $C$'s length only depends on the number of odd-degree vertices in the path $P$.

\begin{corollary}\label{coro6.3}
Suppose that $G$ is a maximal planar graph and $P=v_1v_2\cdots v_m$ is a basic path of $G$. If the subgraph induced by the neighbor set of $P$ is a cycle $C$, then $C$ is an even-cycle if and only if there are even number of odd-degree vertices among $v_1,v_2,\cdots, v_m$.
\end{corollary}

More generally, we are seeking about such cycles induced by the neighbor set of a connected subgraph of $G$.

\begin{theorem2}\label{th6.4}
Suppose that $G$ is a maximal planar graph and $H=G[\{v_1,$ $v_2,\cdots, v_m\}]$ is a connected subgraph of $G$. If the subgraph induced by the neighbor set of $H$ is a cycle, denoted $C=G[\Gamma(v_1,v_2,\cdots, v_m)]$, then
$$
|C|=\sum\limits_{i=1}^{m}d_G(v_i)-\sum\limits_{i=1}^{m}d_H(v_i)-b(H) \eqno{(6.3)}
$$
where  $b(H)$ denotes the number of edges on the boundary of $H$, in which the cut edges are counted twice.
\end{theorem2}

By a similar proof shown in Theorem 6.2, we can prove Theorem 6.4.

\begin{corollary}\label{coro6.5}
Suppose that $G$ is a maximal planar graph and $H=G[\{v_1,v_2$, $\cdots$, $v_m\}]$ is a connected subgraph of $G$. If the subgraph induced by the neighbor set of $H$ is a cycle, denoted $C=G[\Gamma(v_1,v_2,\cdots, v_m)]$, then $C$ is an even-cycle if and only if $b(H)$ and the number of odd-degree vertices have the same parity among $v_1,v_2,\cdots, v_m$.
\end{corollary}

Let $G$ be a maximal planar graph and $H=G[\{v_1,v_2,\cdots, v_m\}]$ be a connected subgraph of $G$. If the subgraph induced by the neighbor set of $H$ is a cycle or a chord-cycle, then the resulting graph $G-C-H$, denoted by $\bar{H}_{C}$, is called the \emph{complement of $H$ on $C$}. Namely
$$
G-C-H\triangleq \bar{H}_{C} \eqno{(6.4)}
$$

\begin{Prop}\label{pro6.6}
Suppose that $G$ is a maximal planar graph and $H=G[\{v_1,v_2$, $\cdots$, $v_m\}]$ is a connected subgraph of $G$. If the subgraph induced by the neighbor set of $H$ is a cycle or a chord-cycle, then the subgraph induced by the neighbors of $\bar{H}_{C}$ is a cycle when $\bar{H}_{C}$ is connected.
\end{Prop}

Based on Proposition 6.6, the subgraph induced by the neighbors of $\bar{H}_{C}$ is written as $C^{\prime}$, then the following theorem is true.

\begin{theorem2}
Let $G$ be a maximal planar graph and $H=G[v_1,v_2,\cdots, v_m]$ be a connected subgraph of $G$. Suppose that the subgraph induced by the neighbor set of $H$ is not a tree, denoted $C$, we have

\emph{(1)}  If $\bar{H}_{C}$ is unconnected, then $C$ is a chord-cycle;

\emph{(2)} If  $\bar{H}_{C}$ is connected, then either $C^{\prime}=C$ that shows $C$ is a basic cycle of $G$, or $|V(C^{\prime})|<|V(C)|$ and $C$ is a chord-cycle that each chord is contained in a triangle.
\end{theorem2}

\subsection{Enumeration of even-cycles}

In this section, we will find a necessary and sufficient condition for a 2-colorable cycle. Obviously, it is very significant to research how many even-cycles are contained in a maximal planar graph.

Suppose that $G$ is a maximal planar graph with order $n$ and $\delta(G)\geq 4$. Denote by $\pi(G)=(d_1,d_2,\cdots,d_n)$ the degree sequence of $G$.  Now, we consider the number of cycles in $G$. If $G$ is indivisible, then according to Theorem 6.1, $G$ contains $n$ basic cycles that are induced only by the neighbors of a vertex of $G$, and we denote the set of such cycles by $C_1$, obviously $|C_1|=n$. Let $m$ be the number of odd-degree vertices of $G$, then $C_1$ contains $m$ odd-cycles. Thereby, the number of even-cycles in $C_1$ is equal to $n-m$. Especially, when $m=n$, there are no even-cycles in $C_1$.  However, the longest cycle in the subgraph of $G$ induced by any two adjacent vertices is either a basic cycle or a chord-cycle,
so there are such $3n-6$ basic cycles and chord-cycles, in which the length of both basic cycles and chord-cycles are even by Theorem 6.2, denoted by $C_2$ and $|C_2|=3n-6$.

In a maximal planar graph $G$, the induced subgraph by the vertices of the boundary of each semi-maximal planar graph, $H$, of $G$ is either a basic cycle or a chord-cycle. We say $H$ corresponds to a basic cycle or chord-cycle. Conversely, we also say a basic cycle or a chord-cycle corresponds to a semi-maximal planar graph.
In addition, it is easy to see that each basic cycle or chord-cycle in $G$ just corresponds to two semi-maximal planar graph of $G$. Let $Ha(G)$ be the set of all semi-maximal planar graphs of $G$, $Ch(G)$ the set of all chord-cycles of $G$, $Cy(G)$ the set of all basic cycles of length not less than 4 of $G$. Then we have

\begin{theorem2} Suppose that $G$ is a maximal planar graph, then
$$
|Cy(G)|\leq\frac{1}{2}|Ha(G)|-|Ch(G)| \eqno{(6.5)}
$$
\end{theorem2}


\subsection{Black-White coloring operation}

Suppose that $G^{C}$ is a semi-maximal planar graph on an even-cycle $C=v_1v_2\cdots v_{2m}$\\$v_1(m\geq 2)$. Denote by $\Gamma^{*}(C)$ the vertex-set consisting of such vertices of $V(G^{C})-V(C)$ that are adjacent to both the vertices of $C$ with odd-subscript and even-subscript. A so-called Black-White coloring for $G^{C}$, denoted $f_{bw}:V(G)\rightarrow \{b,w\}$, is to divide $V(G^{C})$ into two subsets, $B$ and $W$ that are called black vertex-set and white vertex-set, respectively. Hence $V(G^{C})=B \cup W, B,W\neq \emptyset$, where all of the vertices in $B$ and $W$ are colored by black and white, respectively. Further, for a Black-White coloring $f_{bw}=(B,W)$ of $G^{C}$, if both $G^{C}[B]$ and $G^{C}[W]$ contain no odd-cycles, then we refer to $f_{bw}=(B,W)$ as a \emph{proper Black-White coloring}; otherwise, an \emph{improper Black-White coloring}.

Now, for a semi-maximal planar graph $G^{C}$, we introduce a Black-White coloring operation on $C$, which is closely related to 4-colorings. The following gives the detailed steps.

Step 1. Color vertices of $C$ black;

Step 2. Color vertices of $\Gamma^{*}(C)$ white;

Step 3. Let $\Gamma^{*}(\Gamma^{*}(C))\triangleq \Gamma^{2*}(C)$,  then color vertices of $\Gamma^{2*}(C)$ black;

$\cdots\cdots\cdots\cdots\cdots\cdots\cdots\cdots\cdots\cdots\cdots\cdots\cdots\cdots\cdots\cdots\cdots\cdots\cdots\cdots$

Step $2i$. Color vertices of $\Gamma^{(2i-1)*}(C)$ white;

Step $2i+1$. Color vertices of $\Gamma^{(2i)*}(C)$ black;

$\cdots\cdots\cdots\cdots\cdots\cdots\cdots\cdots\cdots\cdots\cdots\cdots\cdots\cdots\cdots\cdots\cdots\cdots\cdots\cdots$

Until

Step $t+1.$  $\Gamma^{(t)*}(C)=\emptyset$;

Step $t+2.$  Color the remainder vertices of $G^{C}$ that are not colored black or white with grey, and the set containing all of the vertices colored grey is denoted as $A$;

Step $t+3.$ For $\forall v\in V(A)$, when $v$ is colored only black (white), there will be odd-cycles in $G[B]$ (or $G[W]$), then we call $v$ the \emph{fixed-vertex} and color it white (black);

Step $t+4.$ For $\forall u\in V(A)$, no matter which color (black or white) the vertex $u$ is colored, there are always odd-cycles in $G[B]$ (or $G[W]$), then we call the vertex $u$ the \emph{petal-vertex} and color it black or white optionally. Then, if there exist fixed-vertices in $A$, go back to Step $t+3$; otherwise, stop.

Here, if $A\neq \emptyset$ when the operation stops, then we also use $f_{bw}=(B,W,A)$ to denote the Black-White coloring operation. According to this operation, an obvious result can be obtained as follows.

\begin{theorem2}\label{th6.9}
Let $G^{C}$ be a semi-maximal planar graph on an even-cycle $C$, and $f_{bw}=(B,W,A)$ be the Black-White coloring operation on $C$, then $G^{C}[B]$ or $G^{C}[W]$ contains odd-cycles if and only if there are petal-vertices appeared at Step $t+4$ in the procedure of $f_{bw}=(B,W,A)$.
\end{theorem2}

Obviously, a petal-vertex is adjacent to at least two vertices of $B$ and $W$, respectively. Figures 6.2(a) and (b) show a structural characteristic partially of a petal-vertex. And the vertex $u$ shown in Figure 6.2(c) is a petal-vertex.

 \begin{center}
       \includegraphics [width=240pt]{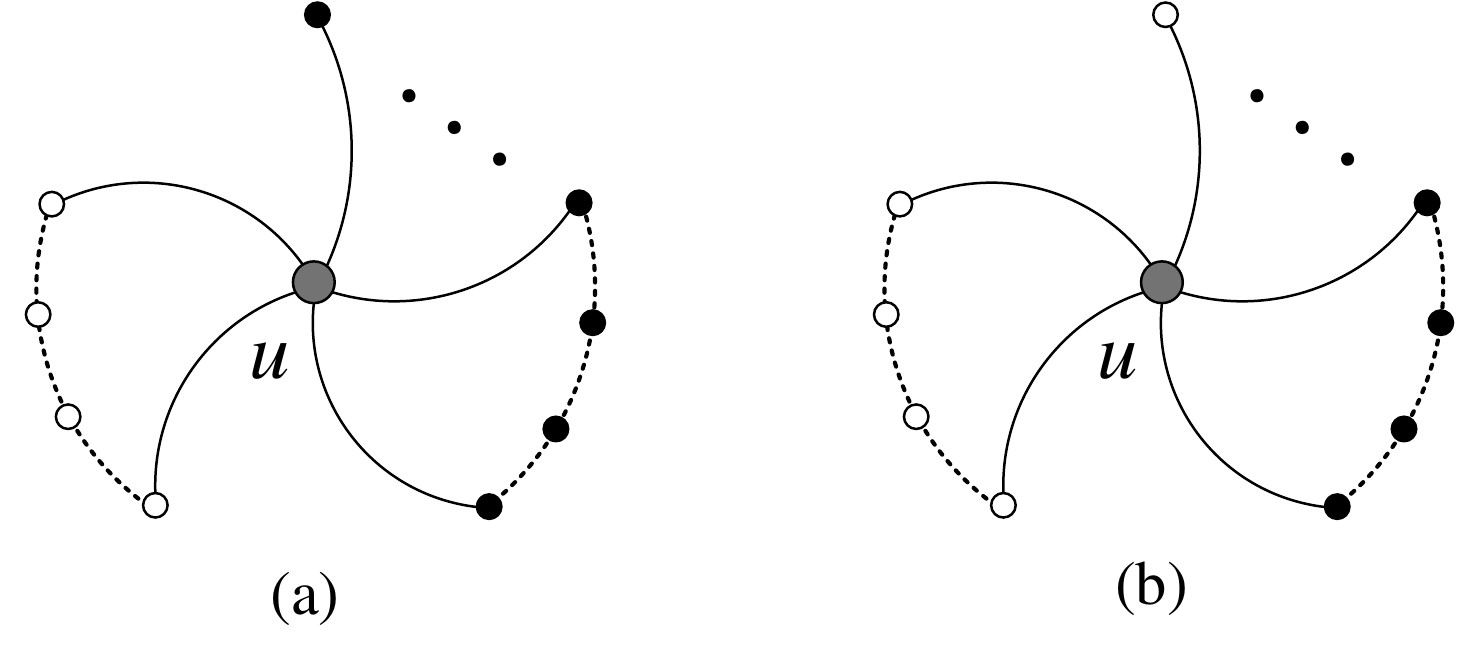}

        \includegraphics [width=320pt]{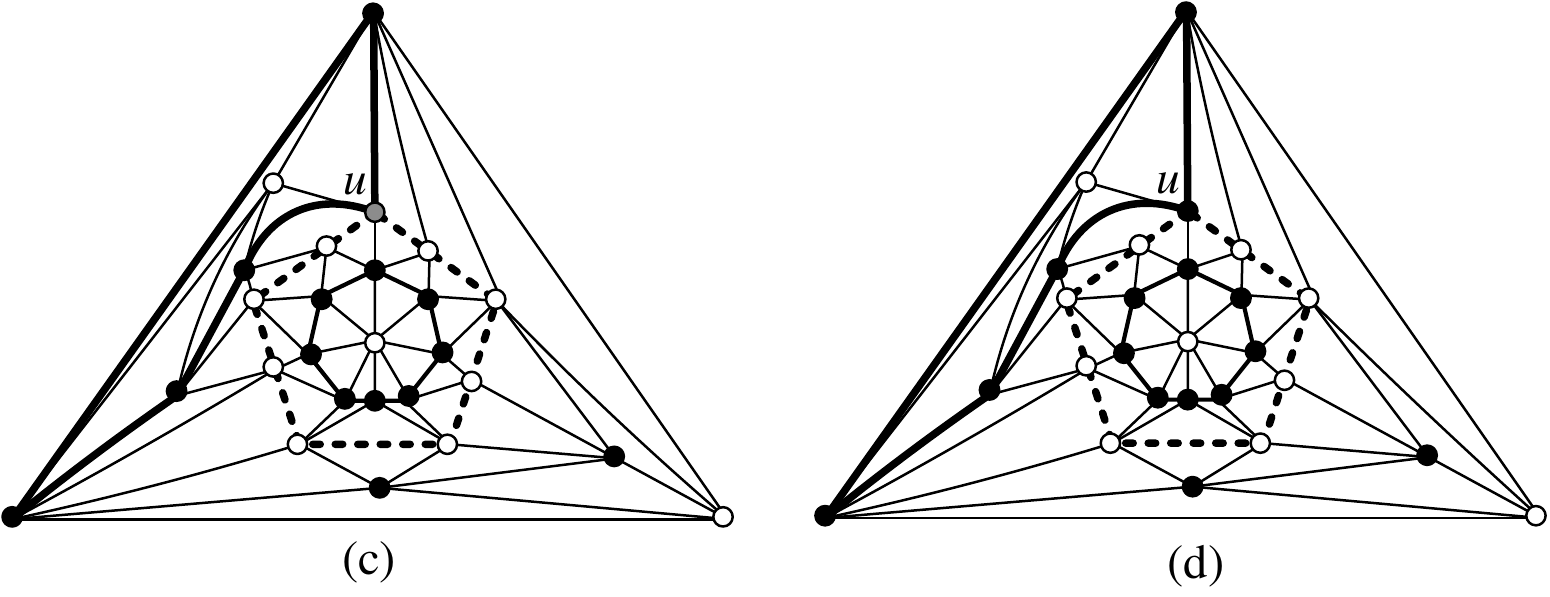}

  \textbf{Figure 6.2.} A diagram of petal-vertices
  \end{center}

 For a maximal planar graph $G$ and a cycle $C$ of $G$, the Black-White coloring operation on $C$ is to conduct a Black-White coloring operation on $C$ for the two semi-maximal planar graphs on $C$ of $G$. For example, Figure 6.3(a) is the resulting graph after conducting a Black-White coloring operation on cycle $v_1v_3v_5v_6v_1$ to the maximal planar graph shown in Figure 6.1; for the graph shown in Figure 6.3(b),
Figure 6.3(b) is the resulting graph after conducting Black-White coloring operation on cycle $v_1v_2v_3v_4v_1$, and Figure 6.3(c) is the resulting graph after conducting Black-White coloring operation on cycle $v_2v_3v_5v_6v_2$. In addition, for the semi-maximal planar graph shown in Figure 6.3(d), implementing a Black-White coloring operation on the inside cycle $C$ produces $\Gamma^{3*}(C)=\emptyset$.

 \begin{center}
       \includegraphics [width=340pt]{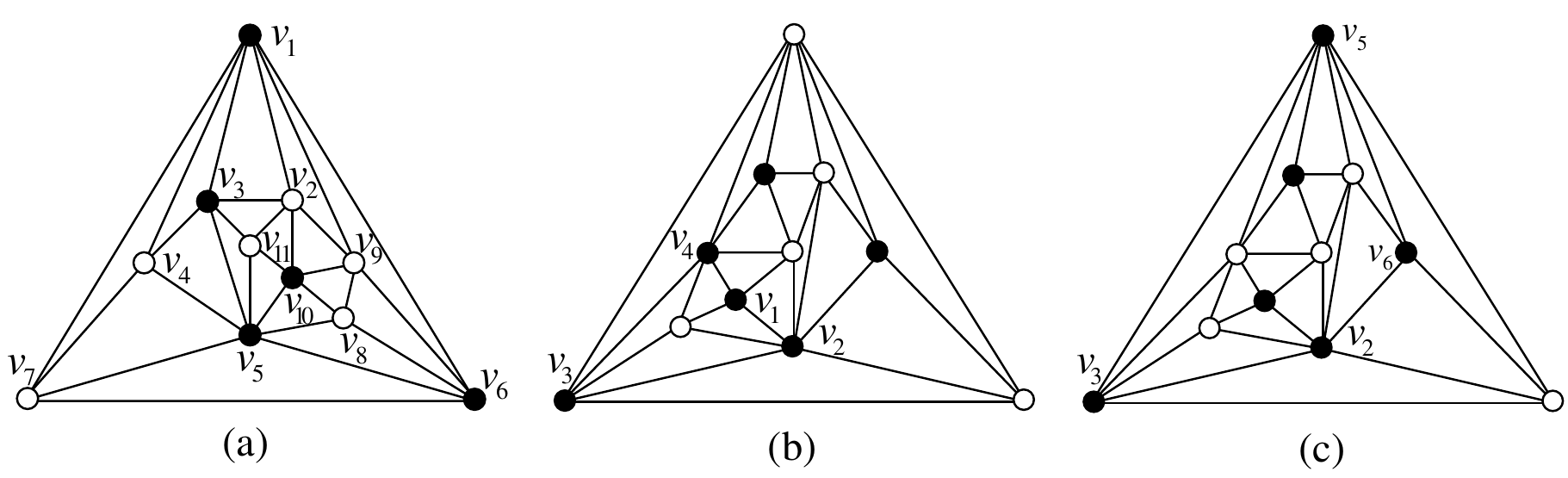}

        \includegraphics [width=200pt]{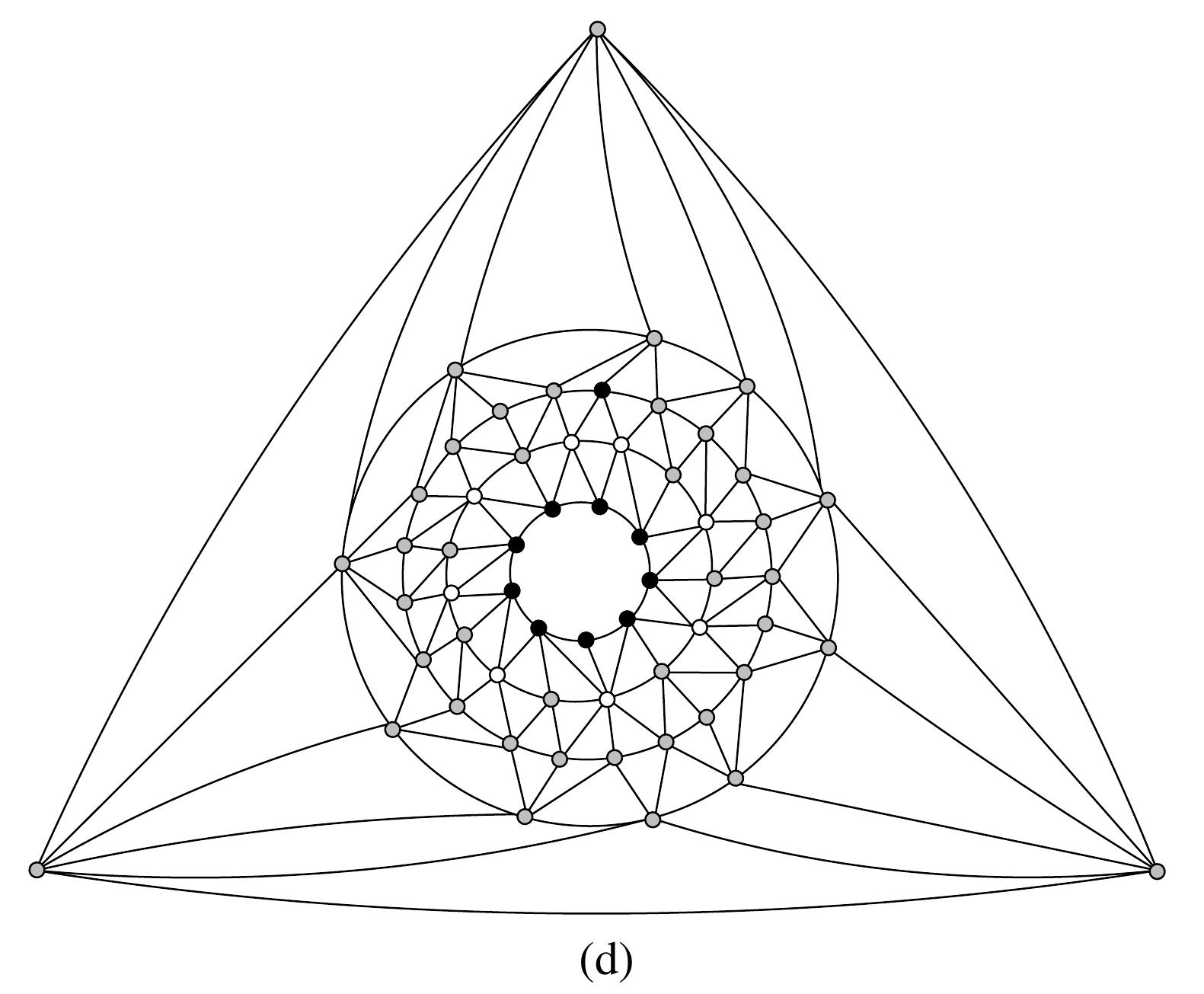}

       \textbf{Figure 6.3.} An example for illustrating a Black-White coloring operation
  \end{center}

 We can know from the above examples that after we conduct Black-White coloring operations to some maximal planar graphs, there may have vertices colored grey, or there may not have such vertices. We refer to the Black-White colorings for the former as \emph{unique Black-White colorings}, and for the latter as \emph{non-unique Black-White colorings}. For example, the colorings in Figures 6.3(a), (b) and (c) are the unique Black-White colorings, but the coloring in Figure 6.3(d) is a non-unique Black-White coloring.

 \subsection{Necessary and sufficient condition of the 2-colorable cycle based on the petal-syndrome}

\quad\quad It is easy to prove the following result:
\begin{theorem2}\label{th6.10} Let $G$ be a maximal planar graph, $C$ be an even-cycle of $G$. Suppose $f_{bw}$ is a Black-White coloring on $C$. If $f_{bw}$ is unique, then $C$ is 2-colorable if and only if $f_{bw}$ is proper.
  \end{theorem2}

According to Theorem 6.10, we can clearly judge that for the maximal planar graph $G$ shown in Figure 6.1(d), the cycle $C=v_1v_3v_5v_6v_1$ is 2-colorable because the Black-White coloring on $C$ of $G$ is unique, and two induced graphs $G[B]$ and $G[W]$ are forest and 1-fence, respectively (see Figure 6.3(a)). Similarly, for the maximal planar graph $G$ shown in Figure 6.1(b), the cycle $C=v_1v_2v_3v_4v_1$ is also a 2-colorable cycle; however, for the icosahedron, the 6-cycles induced by the neighbors of any two adjacent vertices
is not 2-colorable since $G[W]$ contains two triangles (see Figure 6.3(c)).

Let $G^{C}$ be a semi-maximal planar graph on an even-cycle $C$. We conduct the Black-White coloring operation $f_{bw}$ on $C$, which partitions $V(G^{C})$ into three subsets: the black vertex-set $B$, the white vertex-set $W$ and the grey vertex-set $A$. Denote
 $$f_{bw}=(B,W,A) \eqno{(6.6)}$$

Let $f_{bw}=(B,W,A)$ be a Black-White coloring on $C$ of $G^C$, in which $A \neq \emptyset$. If $C$ is an even-cycle and both $G^C[B]$ and $G^C[W]$ contain no odd-cycles, then we recolor any vertex $v \in A$ black or white, and remain the colors of other vertices of $A$ unchanged.  Denote by $f^\prime_{bw}=(B^\prime,W^\prime,A^\prime)$ the new Black-White coloring. Obviously, both $G^C[B^\prime]$ and $G^C[W^\prime]$ still contain no odd-cycles. And then, we call the vertices of $A$ the \emph{free vertices}.

For a Black-White coloring $f_{bw}=(B,W,A)$ on $C$ of $G^C$. Suppose that $|A|\geq 2$ and $u,v \in A$, if $G^C[B\cup \{u,v\}]$ or $G^C[W\cup \{u,v\}]$ contains odd-cycles including $u$ and $v$ when they are recolored black (white), then $\{u,v\}$ is called a \emph{petal-pair}.  In addition, let $S\subseteq A$, if any pair of vertices $u$ and $v$ in $S$ is a petal-pair, then $S$ is called a \emph{petal-set}. For example, both sets $S=\{u_1,u_2,u_3,u_4\}$ in Figure 6.4(a) and $S=\{u_1,u_2,u_3\}$ in Figure 6.4(b) are petal-sets.

We can easily obtain the following theorem by the fact that any planar graph contains no $K_5$ or its subdivision.

\begin{theorem2}\label{th6.11} Let $G^C$ be a semi-maximal planar graph on an even-cycle $C$, and $f_{bw}=(B,W,A)$ be a Black-White coloring on $C$ of $G^C$. If $S\subseteq A$ is a petal-set of $G^C$, then
$$|S|\leq 4 \eqno{(6.7)}$$
\end{theorem2}

An edge $uv$ is called a \emph{petal-edge} if $\{u,v\}$ is a petal-pair. In Figure 6.4(b), three edges $u_1u_2$, $u_2u_3$ and $u_1u_3$ are petal-edges. A path $P$ in $G^C[A]$ is called a \emph{special petal-path} if each edge of $P$ is a petal-edge; analogously, a cycle $C^*$ in $G^C[A]$ is called a \emph{special petal-cycle} if each edge of $C^*$ is a petal-edge. A subgraph $H$ induced by a sequence of $p$ vertices $x_1,x_2,\cdots,x_p$ in $G^C[A]$ is called a \emph{general petal-path} if each pair $\{x_i,x_{i+1}\}$ is a petal-vertex pair for $i=1,2,\cdots,p-1$. Further, if $\{x_1,x_{p}\}$ is a petal-pair, then the subgraph $H$ is called a \emph{general petal-cycle}. Actually, a special petal-path is a special case of general petal-paths; similarly, a special petal-cycle is a special case of general petal-cycles. So a general petal-path and a general petal-cycle are called straightly a \emph{petal-path} and a \emph{petal-cycle}, respectively. For example, the cycle $u_1u_2u_3u_1$ in Figure 6.4(b) is a petal-cycle (special).

\begin{center}
         \includegraphics [width=360pt]{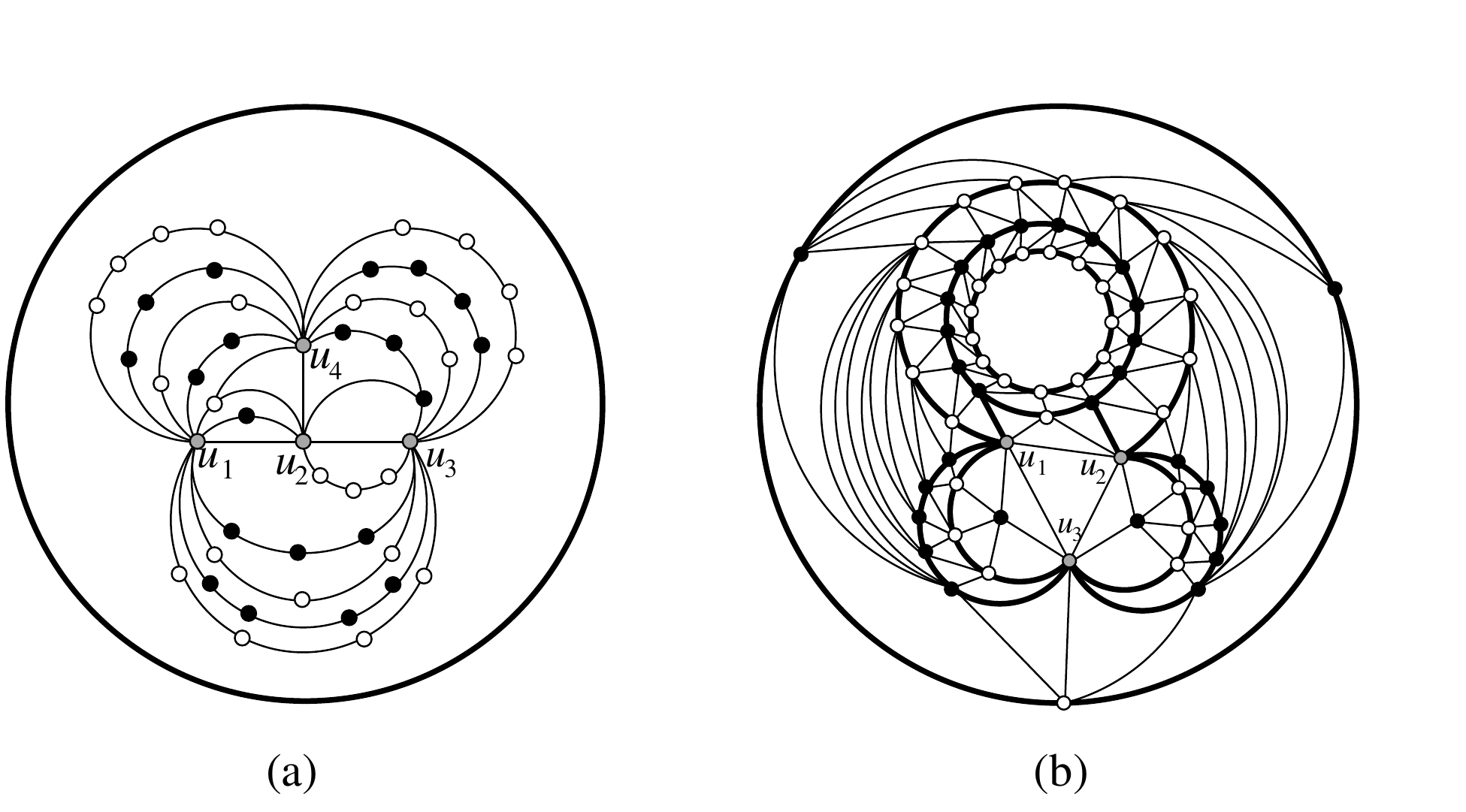}

        \textbf{Figure 6.4.} The illustrations for the concepts of petal-vertex set, petal-edges and petal-cycles
\end{center}

Based on the above arguments, we define a new graph, called a \emph{petal-graph}. For a semi-maximal planar graph $G^C$, we conduct the Black-White coloring operation on $C$, denoted $f_{bw}=(B,W,A)$, $A\neq \emptyset$. The so-called \emph{petal-graph $G_S$ on coloring $f_{bw}$} has its vertex set $\{x_1,x_2,\cdots,x_n\}\subseteq A$ and each its vertex at least belongs to one petal-pair, and two vertices $x_i$ and $x_j$ are adjacent if and only if $\{x_i,x_j\}$ is a petal-pair, where $i,j=1,2,\cdots,n$, $i\neq j$. Obviously, a petal-graph is planar and contains no isolated vertices. For example, if $G^{C}[V(G_S)]$ is a petal-path, then $G_S$ is a path; if $G^{C}[V(G_S)]$ is a petal-cycle, then $G_S$ is a cycle; and if $V(G_S)$ is a petal-set, then $G_S$ is a complete graph.

Let $G^C$ be a semi-maximal planar graph, $f_{bw}=(B,W,A)$ a Black-White coloring on $C$, $A\neq \emptyset$ and $G_S$ the petal-graph on $f_{bw}$. If $G_S$ contains no odd-cycles, namely a bipartite graph with a bipartition $(X,Y)$, then $G_S$ is called an \emph{exclusive petal-graph on coloring $f_{bw}$} if the following two conditions are satisfied:

(1) Both $G^C[B\cup X]$ and $G^C[B\cup Y]$, or both $G^C[W\cup X]$ and $G^C[W\cup X]$ contain odd-cycles;

(2) Both $G^C[B\cup X]$ and $G^C[W\cup X]$, or both $G^C[B\cup Y]$ and $G^C[W\cup X]$ contain odd-cycles.

As shown in Figures 6.5(a) and (b), if there exist a pair of cycles $C_1$ and $C_2$, or a pair of cycles $C_3$ and $C_4$ in $G^C$, then $G^C[X\cup Y]$ is an exclusive petal-graph.

When there exists a path $P=u_1u_2\cdots u_l$ in $G^C[A]$ such that both $G^C[B\cup \{u_{i}u_{i+1}\}]$ and $G^C[W\cup \{u_{i+1}u_{i+2}\}]$ contain odd-cycles, or both $G^C[W\cup \{u_{i}u_{i+1}\}]$ and $G^C[B\cup \{u_{i+1}u_{i+2}\}]$ contain odd-cycles for $i=1,2,\cdots,l-2$, then $P$ is called a \emph{Black-White path}. For the sake of convenience, we always assume that $G^C[B\cup \{u_{i}u_{i+1}\}]$ and $G^C[W\cup \{u_{i+1}u_{i+2}\}]$ contain odd-cycles. Obviously, if a vertex $u$ (not the ends) of the path $P$ is recolored with black or white, then there may appear several fixed-vertices of $A$. Let $D^w_b$ and  $D^w_w$ (contains $u$) be the sets of vertices in $A$ recolored black and white respectively when $u$ is recolored by white, and let $D^b_b$ (contains $u$) and  $D^b_w$ be the sets of vertices in $A$ recolored with black and white respectively when $u$ is recolored with black. A Black-White path $P$ is called an \emph{exclusive Black-White path} if one of the following conditions holds:

 (1) $G^C[ D^b_b ]$ or $G^C[ D^b_w]$, and $G^C[ D^w_b]$ or $G^C[ D^w_w]$ contain simultaneously odd-cycles;

 (2) $G^C[ D^b_b]$ or $G^C[ D^b_w]$ contains two edge-disjoint odd cycles;

 (3) $G^C[ D^w_b]$ or $G^C[D^w_w]$ contains two edge-disjoint odd cycles.

 It is easy to see that if $G^C[A]$ contains an exclusive Black-White path $P$, then there must exist a unicolor odd-cycle  whatever colors (black or white) are assigned to the vertices of $P$. In Figure 6.5(c), when the vertex $u$ is recolored with black,  $D^b_b=\{u,u_2,u_4\}$, $D^b_w=\{u_1,u_3,u_5\}$ and $G^C[D^b_w]$ contains an odd-cycle; when the vertex $u$ is recolored with white, $D^w_b=\{u_6,u_8,u_{10}\}$, $D^w_w=\{u,u_7,u_9,u_{11}\}$, and $G^C[D^w_w]$ contains an odd-cycles, shown in Figure 6.5(d).

\begin{center}
             \includegraphics [width=260pt]{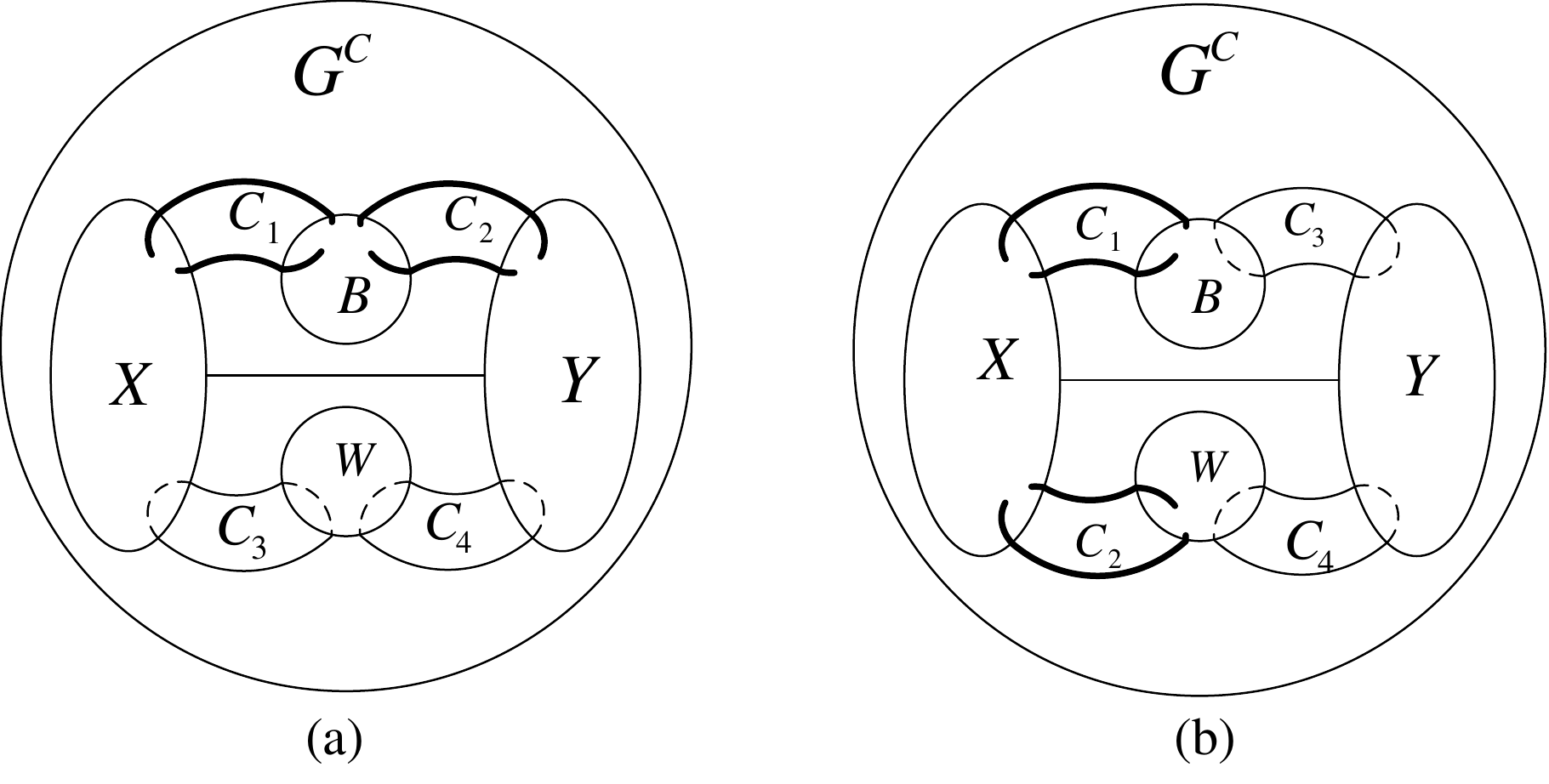}

                            \hspace{1cm}
             \includegraphics [width=320pt]{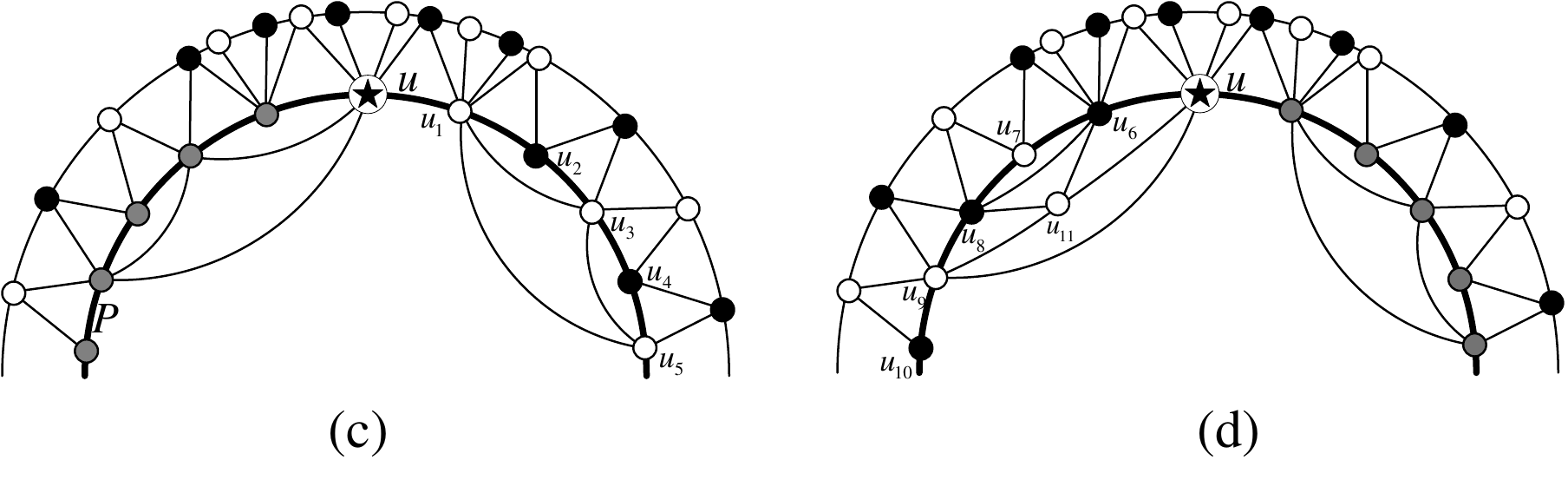}

        \textbf{Figure 6.5.} The illustrations for the concepts of exclusive petal-graphs and exclusive Black-White paths
\end{center}

For a vertex $u$ of $G^C[A]$, when we recolored $u$ with black or white, some fixed-vertices will appear in $\Gamma(u)$. Collect them in $D_1$. Then color them; and another set of fixed-vertices in $\Gamma(D_1)$, $D_2$, will appear correspondingly. Further, when we recolored the vertices of $D_2$, $D_3$ will appear correspondingly. Repeat in this way, until some $D_k$ which has no fixed-vertices appears in $\Gamma(D_k)$. Then we call $D=D_1\cup D_2 \cup \cdots \cup D_k$ the \emph{fixed-set} of $u$.

Suppose that $u$ is a vertex of $G^C[A]$. When $u$ is recolored with white, we denote by $D^w_b$ and $D^w_w$ (contains $u$) the vertex-sets fixed by $u$, in which vertices are recolored with black and white, respectively.
 When $u$ is recolored with black, we denote by  $D^b_b$ (contains $u$) and  $D^b_w$ the vertex-sets fixed by $u$, in which vertices are recolored with black and white, respectively.
 When $u$ is recolored with white,
 let $f^\prime_{bw}=(B^\prime,W^\prime,A^\prime)$ be the extended coloring of $f_{bw}=(B,W,A)$ in the condition that $u$ is recolored with white, in which $B^\prime=B\cup D^w_b$, $W^\prime=W\cup D^w_w$ and $A^\prime=A-(D^w_b\cup D^w_w)$.  When $u$ is recolored with black, let $f^{\prime\prime}_{bw}=(B^{\prime\prime},W^{\prime\prime},A^{\prime\prime})$ be the extended coloring of $f_{bw}=(B,W,A)$ in the condition that $u$ is recolored with black, in which $B^{\prime\prime}=B\cup D^b_b$, $W^{\prime\prime}=W\cup D^b_w$ and $A^{\prime\prime}=A-(D^b_b\cup D^b_w)$. Then, $u$ is called a \emph{general petal-vertex} if one of the following four conditions holds:

(1) $G^C[B^\prime]$ or $G^C[W^\prime]$ contains odd-cycles, and  $G^C[B^{\prime\prime}]$ or $G^C[W^{\prime\prime}]$ contains odd-cycles;

(2) Both petal-graphs on $f^\prime_{bw}$ and $f^{\prime\prime}_{bw}$ of $G^C$ contain odd-cycles;

(3) At least one of $G^C[B^\prime]$ and $G^C[W^\prime]$ contains odd-cycles, and the petal-graph on $f^{\prime\prime}_{bw}$ of $G^C$ contains cycles;

(4) At least one of $G^C[B^{\prime\prime}]$ and $G^C[W^{\prime\prime}]$ contains odd-cycles, and the petal-graph on $f^{\prime}_{bw}$ of $G^C$ contains cycles.

For example, each vertex of an exclusive petal-graph is a general petal-vertex; in Figures 6.5(a) and (b), the vertex $u$ on the exclusive Black-White path is a general petal-vertex; if a petal-graph contains an odd-cycle, then all vertices of the odd-cycle are general petal-vertices. Based on the above arguments, we refer to the phenomenon that $G^C[A]$ contains general petal-vertices as the \emph{petal-syndrome}.

\begin{theorem2}\label{th6.12} Let $G^C$ be a semi-maximal planar graph on an even-cycle $C$, and $f_{bw}=(B,W,A)$ be a Black-White coloring on $C$ of $G^C$. If $A\neq \emptyset$ and $G^C[A]$ contains petal-syndrome, then
$C$ is not a 2-colorable cycle.
  \end{theorem2}

 \begin{theorem2}\label{th6.13}
 Let $G^C$ be a semi-maximal planar graph on an even-cycle $C$, and $f_{bw}=(B,W,A)$ be a Black-White coloring on $C$ of $G^C$. If $A\neq \emptyset$ and promising that no unicolor odd-cycles appear, then for any two adjacent vertices $u,v\in A$:

 \emph{(1)} $u,v$ are recolored only by two different colors if and only if $uv$ is a petal-edge;

 \emph{(2)} $u,v$ are recolored only by the same color if and only if $u,v$ both are on a petal-path and the distance between them is even.
  \end{theorem2}
\begin{proof}
(1) According to the definition of a petal-edge, the sufficient condition holds. Conversely, suppose that $u,v$ are recolored only with two different colors, it shows that both $G^C[B\cup \{u,v\}]$ and $G^C[W\cup \{u,v\}]$ contain  odd-cycles. Therefore, $uv$ is a petal-edge.

(2) Suppose that $u,v$ both are on a petal-path and the distance between them on the path is even, then the colors of $u$ and $v$ must be the some color. Conversely, suppose that $u,v$ are recolored with the same color, namely when one of them, say $u$, is recolored with black (white), then $v$ must be in the vertex-set $D$ fixed by $u$, in which vertices are recolored with black (white). So there exists a petal-path between $u$ and $v$ and the distance between them on the path is even.
\end{proof}

For a semi-maximal planar graph $G^C$ on an even-cycle $C$, in order to deal with the possible problem of the petal-syndrome after conducting  Black-White coloring operation on $C$ for $G^C$, now we give an \emph{improved operation of Black-White coloring on $C$} of $G^C$ as follows.

Step 1. Color vertices of $C$ black;

Step 2. Color vertices of $\Gamma^{*}(C)$ white;

Step 3. Let $\Gamma^{*}(\Gamma^{*}(C))\triangleq \Gamma^{2*}(C)$,  then color  vertices of $\Gamma^{2*}(C)$  black;

$\cdots\cdots\cdots\cdots\cdots\cdots\cdots\cdots\cdots\cdots\cdots\cdots\cdots\cdots\cdots\cdots\cdots\cdots\cdots\cdots$

Step $2i$. Color vertices of $\Gamma^{(2i-1)*}(C)$  white;

Step $2i+1$. Color vertices of $\Gamma^{(2i)*}(C)$ black;

$\cdots\cdots\cdots\cdots\cdots\cdots\cdots\cdots\cdots\cdots\cdots\cdots\cdots\cdots\cdots\cdots\cdots\cdots\cdots\cdots$

Until

Step $t+1$.  $\Gamma^{(t)*}(C)=\emptyset$;

Step $t+2$.  If the subgraph, induced by the set of vertices colored black (or white), contains odd-cycles, then recolor any vertex of the odd-cycles grey and go to stop. Otherwise, if all vertices of $G^C$ are colored black or white, stop; else if there are vertices in $G^{C}$ that are not colored black or white, then color these vertices grey.

For any grey vertex $u$ in $G^{C}$ that is not a general petal-vertex,
if only one color, say white, is assigned to $u$,
the subgraph induced by the set of vertices colored black or white contains odd-cycles, after recoloring the vertices in fixed-set of $u$ black or white properly, then $u$ is called a \emph{restricted-vertex}. Of course, when we recolor $u$ black, no unicolor odd-cycles appear.

Step $t+3$. If there are restricted-vertices in the set of grey vertices, then properly color them and the vertices of their fixed-sets correspondingly;

Step $t+4$. If there exist petal-vertices or petal-syndrome, or the subgraph induced by the set of all black vertices or all white vertices, contains odd-cycles, then stop;

Step $t+5$. If there are restricted-vertices in the set of grey vertices, then go back to Step $t+3$; otherwise, go to next step;

Step $t+6$. Choose a grey vertex $v$, called the \emph{sign-vertex}, which has the most neighbors colored black or white, and color $v$ black \emph{and remark it in sequence}. If there is no grey vertex, stop; otherwise, go to the next step;

Step $t+7$. If there are restricted-vertices in the set of grey vertices, then properly color them and the vertices of their fixed-sets correspondingly;

Step $t+8$. If there are no grey vertices, stop; otherwise, if there exist petal-vertices or petal-syndrome, or the subgraph induced by the set of all black vertices or all white vertices contains odd-cycles, when there are black sign-vertices, choosing the latest black sign-vertex, denoted $w$, then we assign grey to
the vertices, which were colored black or white after $w$ was colored black. At the same time, recolor $w$ white and go back to Step $t+7$; when there are no black sign-vertices, stop; else if there are not petal-vertices or petal-syndrome, or the subgraph, induced by the set of all black vertices or all white vertices, contains no odd-cycles, go to the next step;

Step $t+9$. If there are restricted-vertices in the set of grey vertices, then go back to Step $t+7$; otherwise, go back to Step $t+6$.

Next, we will give an example to illustrate the process of an improved operation of Black-White coloring. Let $G$ be a maximal planar graph and $G^C_1, G^C_2$ be two semi-maximal planar graphs on a cycle $C$, shown in Figure 6.6(a). Now, we run an improved operation of Black-White coloring on $C$ of $G^C_1$.

First, color all vertices of $C$ black, and conduct the operation until Step 3, then $\Gamma^{2*}(C)=\emptyset$. Because not all vertices in $G^C_1$ are colored black or white and the subgraph induced by the set of the vertices colored black (or white) contains no odd-cycles, we color other vertices grey shown in Figure 6.6(b).

Second, conduct Step $t+3$ ($t=2$) of the operation. Since neither grey vertices are restricted-vertices, so conduct Step $t+4$. Because there are no petal-vertices, petal-syndrome, sign-vertices, and the subgraph induced by the set of all black vertices or all white vertices contains no odd-cycle, conduct Step $t+6$ directly. Choose a sign-vertex $v_1$ and color it black; and then conduct Step $t+7$ because $G^C_1$ still contains grey vertices. Here, $v_{11}$ is a restricted-vertex and $\{v_{12}\}$ is the fixed-set of $v_{11}$. Color vertices $v_{11}$ and  $v_{12}$ with the corresponding colors, shown in 6.6(b').

Third, because there are no petal-vertices, petal-syndrome, sign-vertices, and the subgraph induced by the set of all black vertices or all white vertices contains no odd-cycles, go back to Step $t+6$. Choose a sign-vertex $v_2$, and color it black. Here, $v_{21}$ is a restricted-vertex, and $v_{22}$ is the fixed-set of  $v_{21}$. Color $v_{21}$ and  $v_{22}$ with corresponding colors. According to the operation and the above arguments, we choose a sign-vertex $v_3$, then $v_{31}, v_{33}$ both are two restricted-vertices and $\{v_{32},v_{34},v_{35},v_{36},v_{37},v_{38},v_{39},v_{310}\}$ is the fixed-set of them. We color these vertices properly, then the operation stops as if there is no grey vertex in $G^C_1$(see Figure 6.6(b")). So $C$ is 2-colorable in $G^C_1$ since $G^C_1$ contains no unicolor odd-cycle colored.

Analogously, we conduct the improved operation of Black-White coloring on $C$ of $G^C_2$, and the resulting coloring is shown in Figures 6.6(c) and (c'). On the process of the coloring, $v_1,v_2,v_3,v_4$ and $v_5$ are five sign-vertices, and  $\{v_{11},v_{12}\}$, $\{v_{21},v_{22},\cdots,v_{210}\}$, $\{v_{31},v_{32},v_{33}\}$, $\{v_{41},v_{42}\}$, and $\{v_{51},v_{52}\}$ are the fixed-set of them. We can see that $G^C_2$ contains no grey vertices, no unicolor odd-cycle colored  when the operation stops, so $C$ is 2-colorable in $G^C_2$.

Hence, $C$ is a 2-colorable cycle of $G$ and the proper Black-White coloring on $C$ of $G$ is shown in Figures 6.6(b") and (c').

\begin{center}
            \includegraphics [width=240pt]{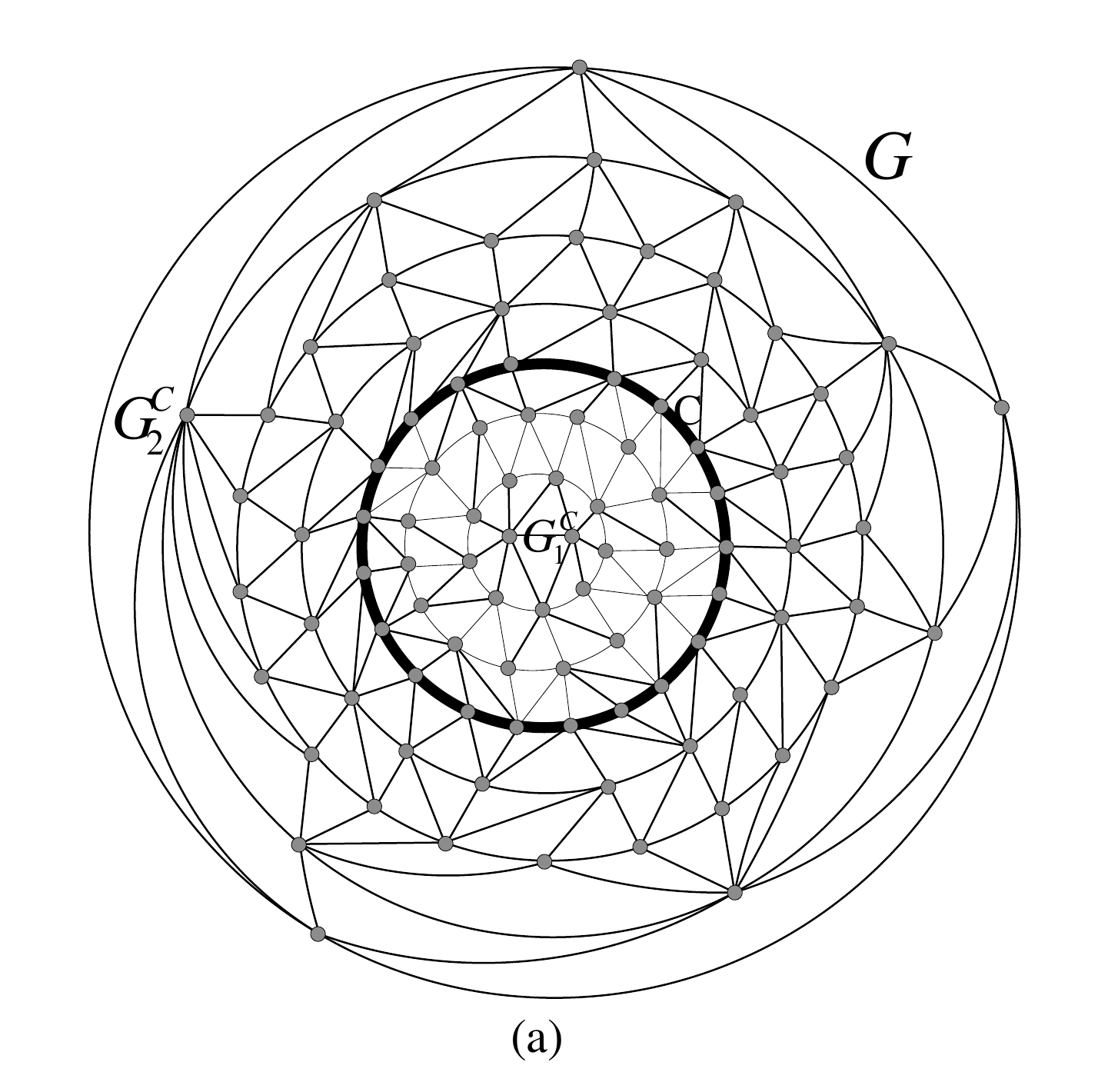}

             \includegraphics [width=380pt]{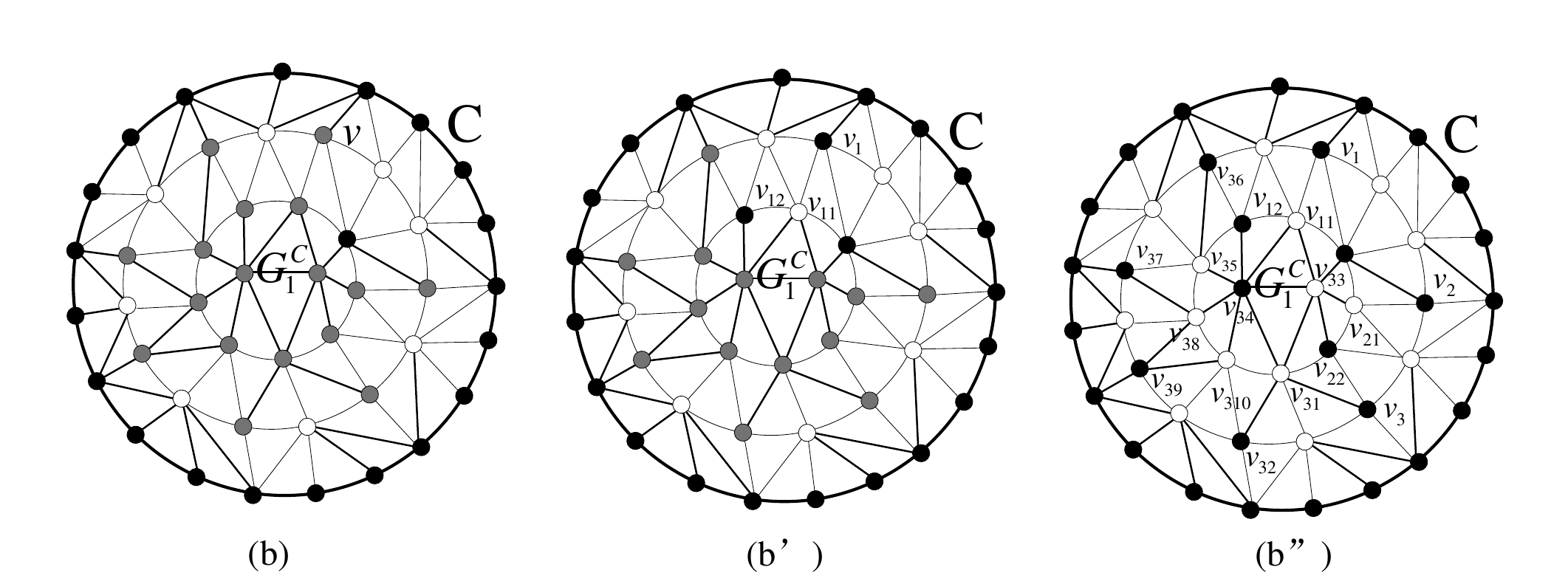}

             \includegraphics [width=360pt]{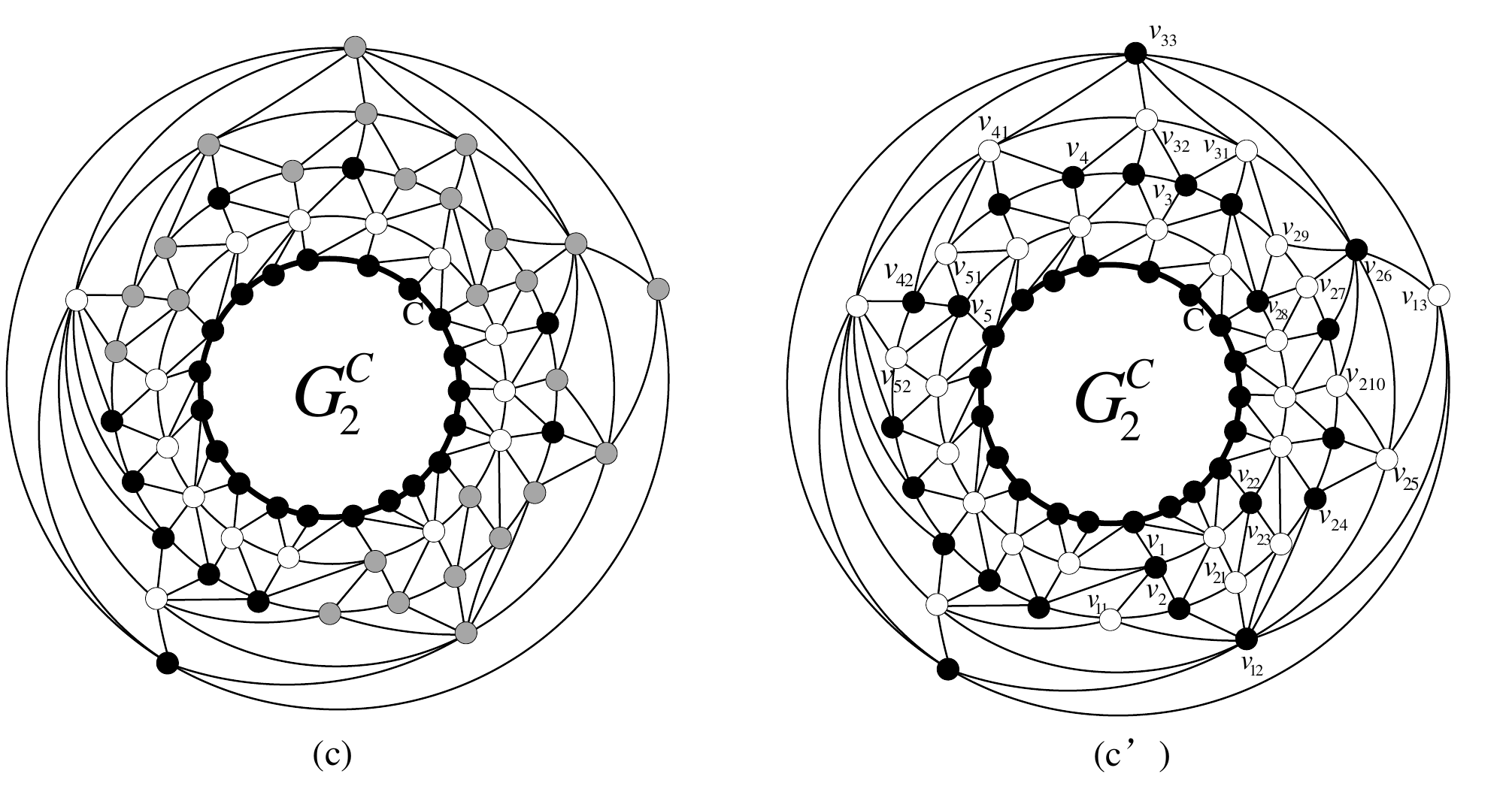}

        \textbf{Figure 6.6.} A maximal planar graph $G$ and two semi-maximal planar graphs $G^C_1, G^C_2$ on $C$
\end{center}

From the improved operation of Black-White coloring, we can obtain the following result.

\begin{theorem2}\label{th6.14} Let $G$ be a maximal planar graph with $\delta(G)\geq 4$, $C$ be an even-cycle of $G$, and $C$ splits $G$ into two semi-maximal planar graphs $G^C_1$, $G^C_2$.  Then $C$ is 2-colorable if and only if neither $G^C_1$ nor $G^C_2$ contains grey vertices after conducting the improved operation of Black-White coloring on $C$ for $G^C_1$ and $G^C_2$, respectively.
  \end{theorem2}
\begin{proof}
  Suppose that $C$ is 2-colorable in $G$. If we conduct the improved operation of Black-White coloring on $C$ for $G^C_1$ and $G^C_2$ respectively, at least one of $G^C_1$ and $G^C_2$, say $G^C_1$, contains grey vertices when the operations stop, then there appears petal-syndrome in the process of conducting the improved operation of Black-White coloring on $C$ for $G^C_1$. That is to say, no matter how to color the vertices of $G^C_1$, there always exist odd-cycles colored with the same color, which contradicts the assumption of $C$ being 2-colorable.

  Conversely, suppose that both $G^C_1$ and $G^C_2$ contain no grey vertices after conducting the improved operation of Black-White coloring on $C$ for $G^C_1$ and $G^C_2$, respectively. Then there appears no odd-cycle in the process of conducting the improved operation of Black-White coloring on $C$ for $G^C_1$ or $G^C_2$. Therefore, $C$ is 2-colorable in $G^C_1$ and $G^C_2$.
\end{proof}

\subsection{Necessary and sufficient conditions of 2-colorable cycles based on graph structure}

On the basis of petal-syndrome, section 6.4 has given a necessary and sufficient condition of 2-colorable cycles. Sometimes, finding or judging petal-syndrome is a very tough task, so this section will study the characteristics of 2-colorable cycles based on the structure of a maximal planar graph.

Suppose that $G^{C}$ is a semi-maximal planar graph on an even-cycle $C$, and $V'=V(G^{C}-V(C)$. Define $G^{C}[V']\triangleq G'$ as a subgraph of $G^{C}$, induced by the inner vertices of $C$.

 For a semi-maximal planar graph $G^{C}$ on an even-cycle $C$, it is easy to prove that $C$ is 2-colorable if and only if $G^{C}$ can be partitioned into two bipartite subgraphs, and $C$ is included in one of them. Suppose that $G_1,G_2$ are two bipartite subgraphs of $G^{C}$ when $C$ is 2-colorable. Without loss of generality, we may assume $C\in G_1$. Thus, we can discuss the structure of $G^C$ through dividing the even-cycle $C$ into two categories according to the relationship between $\Gamma(C)$ and $\Gamma^{*}(C)$.

\textbf{Type 1. Closed-cycles}

If  $\Gamma(C)=\Gamma^{*}(C)$, then we call the cycle $C$ a \emph{closed-cycle} of $G^{C}$. Further, we can divide it into three subcases in detail.

\textbf{Case 1.1. Closed cycle-cycle type}, namely $G^C[\Gamma^{*}(C)]$ is a cycle, denoted $C^*=G^C[\Gamma^{*}(C)]$. We say the cycle $C^*$ encloses $C$. Correspondingly, we say the semi-maximal planar graph $G^{C}$ with the closed-cycle $C$ is a semi-maximal planar graph of \emph{closed cycle-cycle type}. If the Black-White coloring $f_{bw}$ on $C$ for $G^{C}$ is unique, and the subgraphs induced by the fixed-vertices in each step of $f_{bw}$ is a closed-cycle, then we refer to $G^{C}$ as a semi-maximal planar graph on $C$ of \emph{closed type}.

For example, for the semi-maximal planar graph $G^{C_6}$ on a 6-cycle $C_6$ shown in Figure 6.7, $C_6=v_1v_2v_3v_4v_5v_6v_1$ is a closed-cycle of $G^{C_6}$. Because $\Gamma(C_6)=\Gamma^*(C_6)=\{v_1',v_2',v_3',v_4',v_5',v_6'\}$ and $G^{C_6}[\Gamma^*(C_6)]=C_6'$ is also a closed-cycle since $\Gamma(C_6')=\Gamma^*(C_6')=\{x,y,z\}$, then $G^{C_6}$ is a semi-maximal planar graph on the 6-cycle $C_6$ of closed type.

\begin{center}
         \includegraphics [width=240pt]{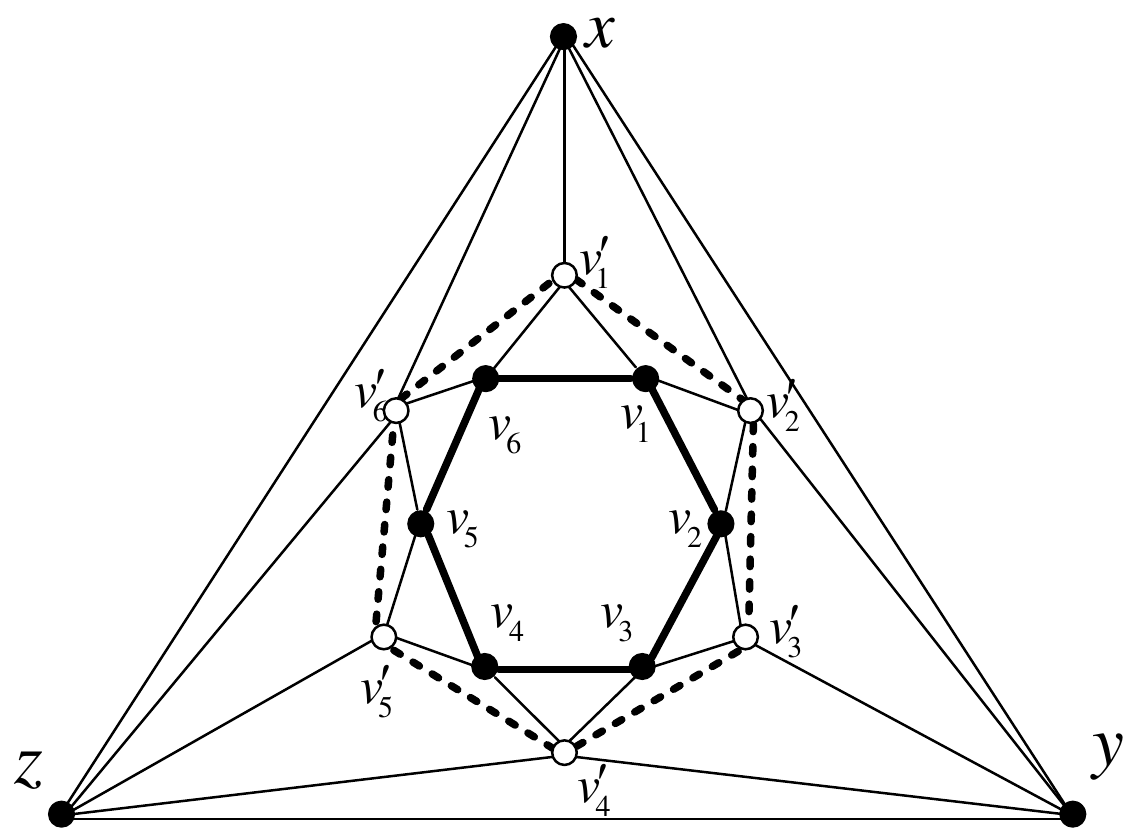}

        \textbf{Figure 6.7.} A diagram for illustrating a closed-cycle and a semi-maximal planar graph of closed type
\end{center}

For semi-maximal planar graphs of closed type, we have an obvious result as follows.

\begin{theorem2}\label{th6.15}
Suppose that $G^{C}$ is a semi-maximal planar graph on even-cycle $C$, and $\Gamma(C)=\Gamma^{*}(C)$. If $G^C[\Gamma^{*}(C)]$ is a cycle, denoted $C^{*}$, then $C$ is 2-colorable in $G^{C}$ if and only if $C^{*}$ is 2-colorable in $G^{C}-V(C)$.
\end{theorem2}

\begin{proof}
If $C^{*}$ is 2-colorable in $G^{C}-C$, then $C$ is also 2-colorable in $G^{C}$, obviously. Conversely, if $C$ is a 2-colorable cycle in $G^{C}$, then there exists a coloring $f$ of $G^{C}$ satisfying $|f(C)|=2$, without loss of generality, $f(C)=\{1,2\}$. By the definition of $C^{*}$, we have $f(v)\neq 1,2$ for any $v\in V(C^*)$. So $f(C^*)=\{3,4\}$, namely $C^{*}$ is 2-colorable in $G^{C}$. Thus, $C^{*}$ is also 2-colorable in $G^{C}-V(C)$.
\end{proof}

For example, in Figure 6.7, because cycle $xyzx$ is not a 2-colorable cycle,  then $C_6$ is not a 2-colorable cycle.

\textbf{Case 1.2. Closed cycle-tree type}, namely $G^C[\Gamma^{*}(C)]=G'$ is a tree. In this case, $G_1=C$, $G_2$ is a tree, and we say  $G^C$ is a semi-maximal planar graph of \emph{closed cycle-tree type}. Obviously, $C$ is a 2-colorable cycle.

\textbf{Case 1.3. Closed cycle-fence type}, namely $G^C[\Gamma^{*}(C)]\neq G'$ and $\Gamma^{*}(C)$ is a fence. Then, we say that $G^C$ is a semi-maximal planar graph of \emph{closed cycle-fence type}. In fact, if  $\Gamma^{*}(C)$ isn't a fence, then $\Gamma^{*}(C)$ contains odd-cycles, which means that $C$ isn't a 2-colorable cycle obviously according to Theorem 6.15. For a semi-maximal planar graph $G^C$ of closed cycle-fence type, $\Gamma^{*}(C)$ is a fence that contains one or more even-cycles. Figure 6.8 gives two examples, in which $\Gamma^{*}(C)$ contains one even-cycle. In this case, obviously, $C$ is 2-colorable if and only if all of the even-cycles of $\Gamma^{*}(C)$ are 2-colorable.

\begin{center}
          \includegraphics [width=260pt]{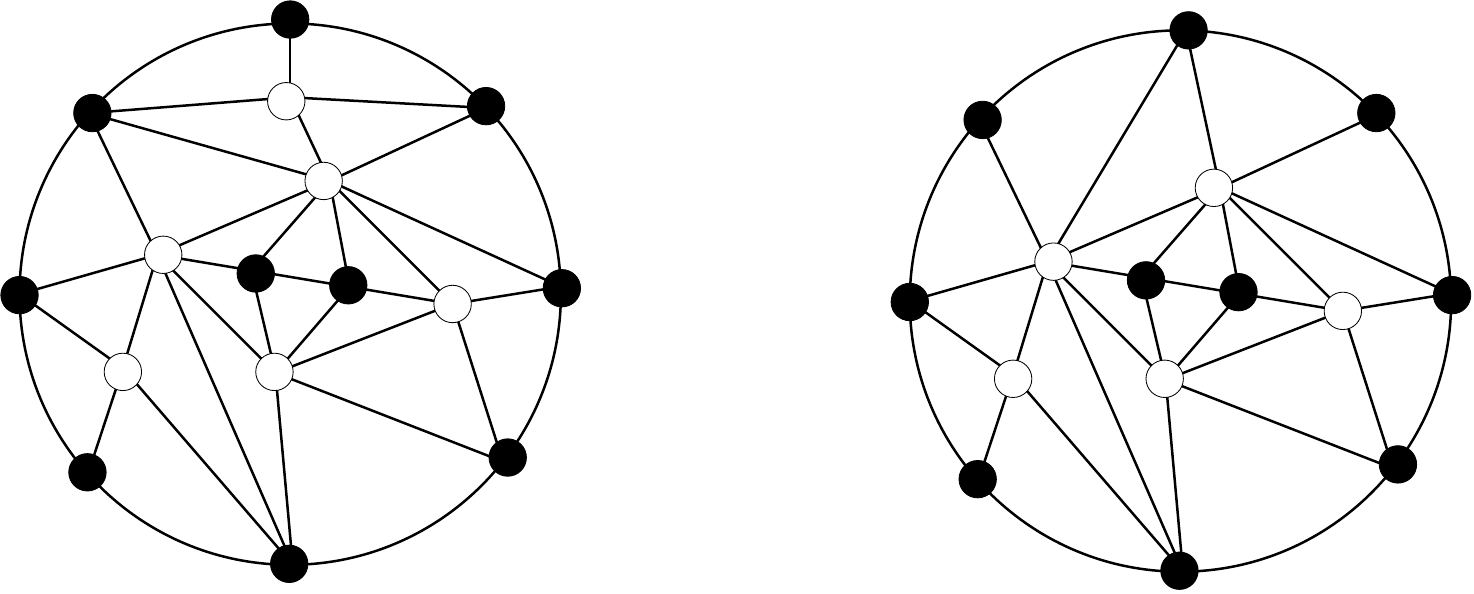}

        \textbf{Figure 6.8.} Two semi-maximal planar graphs of closed cycle-fence type
\end{center}

\textbf{Type 2. Opened-cycles}

If  $\Gamma(C) \neq \Gamma^{*}(C)$, namely $\Gamma^{*}(C)\subset \Gamma(C)$,  then we call cycle $C$ an \emph{opened-cycle} of $G^{C}$, and $O(C)=\Gamma(C)-\Gamma^{*}(C)$ an \emph{opened-vertex set} of $C$,  in which the vertices are called the \emph{opened-vertices}. In this case, we refer to $G^{C}$ as a semi-maximal planar graph of \emph{opened type}. For example, in Figure 6.1(d), consider the 4-cycle $C=v_1v_3v_5v_6v_1$ and its inner components, which is a semi-maximal planar graph on $C$. Obviously, $C$ is an opened-cycle, and we have an opened-vertex set $O(C)=\{v_{10}\}$; in Figure 6.3(b), for the semi-maximal planar graph including the 4-cycle $C=v_1v_2v_3v_4v_1$ and its outer components together, $C$ is also an opened-cycle.

 If $C$ is an opened-cycle, then $G_1$ can only contain the opened-vertices of $\Gamma(C)$. Similar with the condition of closed-cycles, we can also partition all opened type semi-maximal planar graphs into six subtypes as follows:

\textbf{Case 2.1. Fence-tree type.} $G_1$ is a fence that contains only an even-cycle $C$, and $G_2=G^{C}-G_1$ is a tree. In this case, we call $G^{C}$ a semi-maximal planar graph of \emph{fence-tree type}. Consider a special situation that $G_1$ is a 0-fence and $G_2$ is a tree, which is similar as Case 1.2.

\textbf{Case 2.2. Fence-cycle type.} The connected components of $G_1$ including $C$ is a fence $G_F$ that contains only an even-cycle $C$, and $G_2=G^{C}-G_1$ contains even-cycles that are 2-colorable. Then,  we call $G^{C}$ a semi-maximal planar graph of \emph{fence-cycle type} (see Figure 6.9(a)).

\textbf{Case 2.3. Scycles-forest type.} Besides $C$, $G_1$ also contains other even-cycles, which either connect mutually through a path, or connect to $C$ through a path containing opened-vertices, and $G_2=G^{C}-G_1$ is a forest. In this case, we call $G^{C}$ a semi-maximal planar graph of \emph{scycles-forest type} (see Figure 6.9(b)).

\textbf{Case 2.4. Scycles-Scycles type.} Besides $C$, $G_1$ also contains other even-cycles. Let $G_
s$ be the connected components including $C_1$ in $G_1$, then $G_2=G^{C}-G_1$ also contains some even-cycles which are 2-colorable in $G^{C}-G_S$. In this case, we call $G^{C}$ a semi-maximal planar graph of \emph{scycles-scycles type}.

\begin{center}
         \includegraphics [width=260pt]{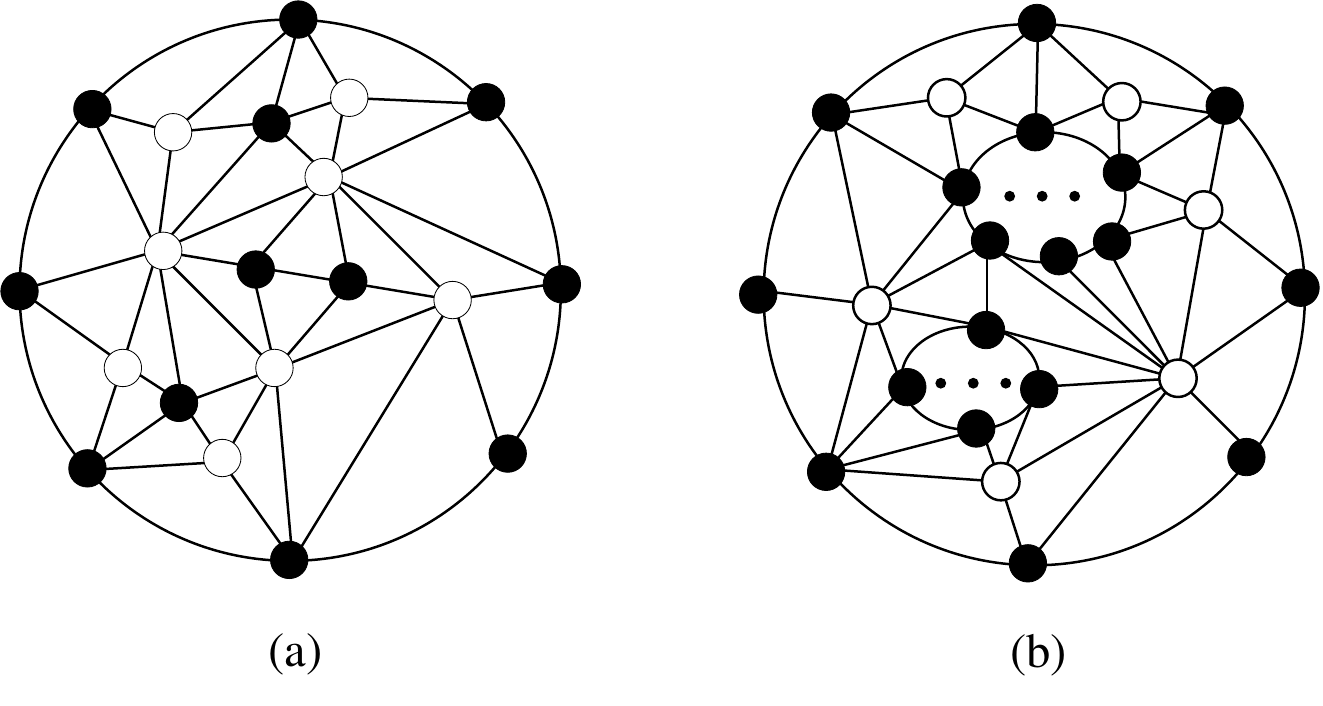}

        \textbf{Figure 6.9.} Semi-maximal planar graphs of fence-cycle, cycles-forest, scycles-scycles type
\end{center}

\textbf{Case 2.5. Intersected cycles-forest type.} $G_1$ consists of some even-cycles $C,C_1,\cdots, C_m$, each of which has at least two vertices of $C$. Namely, the vertices of $C$ are partitioned into $m$ subsets so that each subset is included in at least one even-cycle. Simultaneously, $G_2=G^{C}-G_1$ is a forest. In this case, we call $G^{C}$ a semi-maximal planar graph of \emph{intersected cycles-forest type} (see Figure 6.10).

\textbf{Case 2.6. Intersected cycles-cycles type.} $G_1$ contains $q$ cycles $C,C_1,\cdots, C_q$. Among them there are $m$ even-cycles $C_1,\cdots, C_m$, each of which has at least two vertices of $C$.
Namely, the vertices of $C$ are partitioned into $m$ subsets so that each subset is included in at least one even-cycle. $G_2=G^{C}-G_1$ also contains even-cycles that are 2-colorable. In this case, we call $G^{C}$ a semi-maximal planar graph of \emph{intersected cycles-cycles type}.

\begin{center}
             \includegraphics [width=260pt]{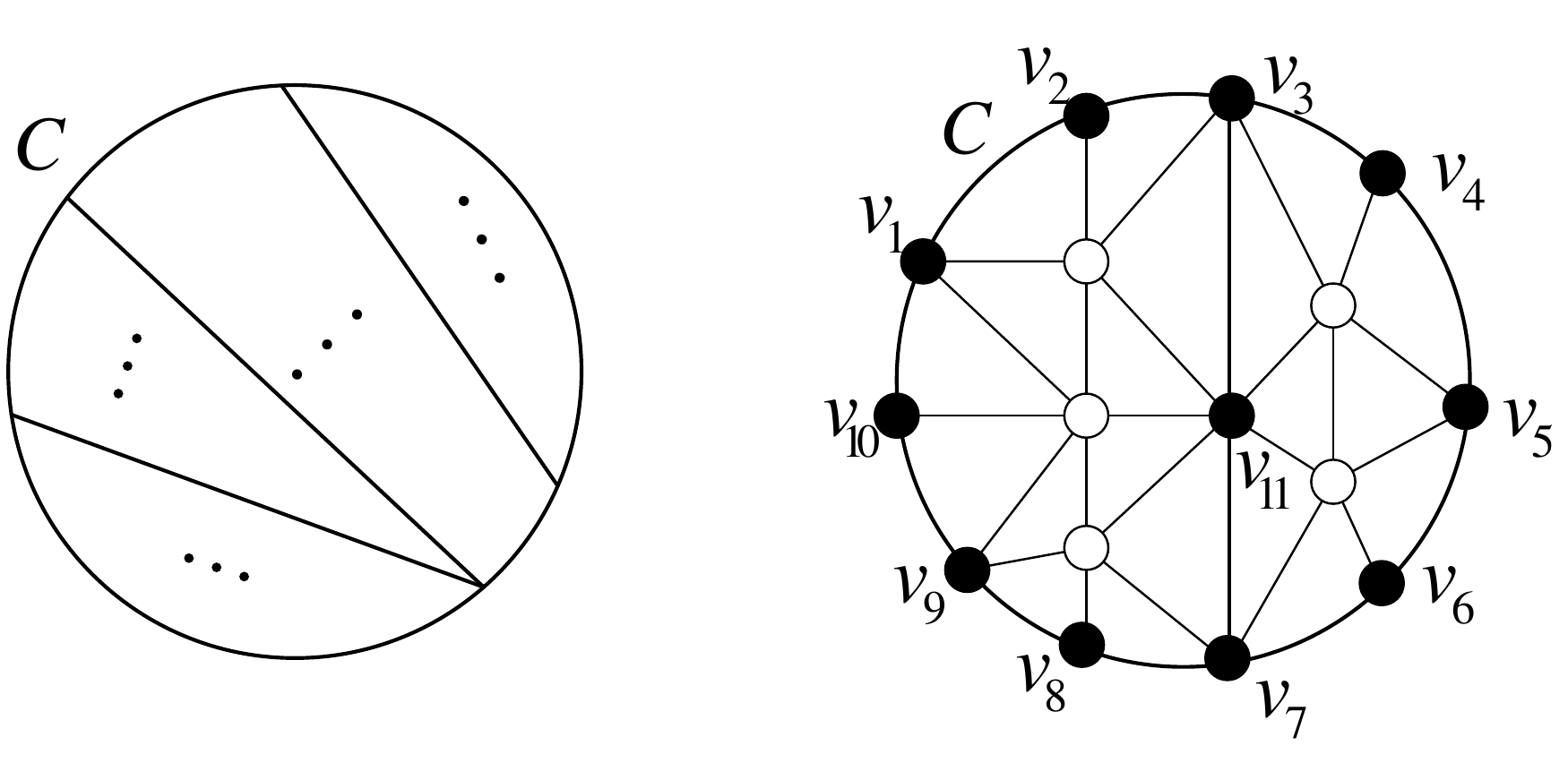}

        \textbf{Figure 6.10.} Semi-maximal planar graphs of intersected cycles-cycles type
\end{center}

The above arguments of the nine cases on closed-cycles and opened-cycles, in fact, give a necessary and sufficient condition of $C$ being 2-colorable. In addition, these nine cases are the categories when we give a classification for semi-maximal planar graphs according to closed-cycle and opened-cycle. So, we can obtain the following result.

\begin{theorem2}\label{th6.16}
Suppose that $G^{C}$ is a semi-maximal planar graph on an even-cycle $C$. Then $C$ is $2$-colorable if and only if $G^{C}$ belongs to one of the following:

\emph{(1)} closed cycle-cycle type, in which $G^{C}[\Gamma^{*}(C)]$ is 2-colorable in $G^{C}-C$;

\emph{(2)} closed cycle-tree type;

\emph{(3)} closed cycle-fence type, in which the even-cycles in the fence are 2-colorable;

\emph{(4)} fence-tree type;

\emph{(5)} fence-cycle type;

\emph{(6)} scycles-forest type;

\emph{(7)} scycles-Scycles type;

\emph{(8)} intersected cycles-forest type;

\emph{(9)} intersected cycles-cycles type.
$$
\eqno{\Box}
$$
\end{theorem2}

\subsection{Construction of semi-maximal planar graphs having 2-colorable cycles}

In a semi-maximal planar graph $G^{C}$ on $C$, when we implement a Black-White coloring operation on $C$, if $C$ is a bicolored cycle, and there is no other bicolored cycles in the component inside $C$, then we call $C$ a \emph{basic bicolored cycle} of $G^{C}$ and say $G^{C}$ is to be a semi-maximal planar graph of \emph{basic type}; otherwise, a \emph{compound bicolored cycle} and  say $G^{C}$ belongs to the \emph{compound type}. In Theorem 6.16, only in the second and the fourth cases, \emph{closed cycle-tree and fence-tree type},  $C$ is a basic bicolored cycle. Obviously, compound bicolored cycles can be gained from basic bicolored cycles through some given operations. Then, what are these operations? This section will  reply this question. In fact, just three operations involved: one is the cycle-spliced operation, and another two are the bicolored path-split operation and its inverse operation, the bicolored cycle-contracted operation.  At the end of this section, we study the characteristics and structure of a semi-maximal planar graph $G^{C}$ with a basic bicolored cycle $C$. For this goal, we introduce two new operations: the \emph{folded operation on even-cycles} and its inverse operation, the \emph{unfolded operation on even-cycles}.  Further, we give the structural characteristics of a semi-maximal planar graph that belongs to  fence-tree type.

\subsubsection{Generating operation system of semi-maximal planar graphs having compound bicolored cycles}

Denote by $G_B^C$ and $G_C^C$ the semi-maximal planar graphs having basic bicolored cycles and compound bicolored cycles, respectively. This subsection will give a generating operation system of semi-maximal planar graphs having compound bicolored cycles, denoted $\zeta{(G_C^C)}$, in the following:
\begin{center}
    \textbf{$\zeta{(G_C^C)}$=($G_B^C$, $S$)}
\end{center}
in which $S$ contains three operations: the cycle-spliced operation, the bicolored path-split operation and the bicolored cycle-contracted operation. Now, we will describe these three operations at length.

\textbf{I. Cycle-spliced operation}

Suppose that $G^{C_1}$ and $G^{C_2}$ are two semi-maximal planar graphs on two 2-colorable cycle $C_1$ and $C_2$, respectively. $C_i$ contains a path $P_i$ of length $x$ for $i=1,2$ satisfying
$$
1\leq x\leq \left\{\begin{array}{c}
              \min\{|C_1|,|C_2|\}-1, |C_1|\neq |C_2|; \\
             \hspace{0.8cm} |C_1|-2, \hspace{1cm} |C_1|=|C_2|.
            \end{array}\right.
$$
Then, we merge $P_1$ and $P_2$ into one path $P$ of length $x$ such that $G^{C_1}$ and $G^{C_2}$ form a new semi-maximal planar graph $G^{C}$ on a cycle $C$, whose length is
$$
|C|=|C_1|+|C_2|-2x.
$$
An example for illustrating this process is shown in Figure 6.11. Easily, it follows that $|C|\geq 4$, and $C$ is a 2-colorable cycle of  $G^{C}$.

\begin{center}
          \includegraphics [width=160pt]{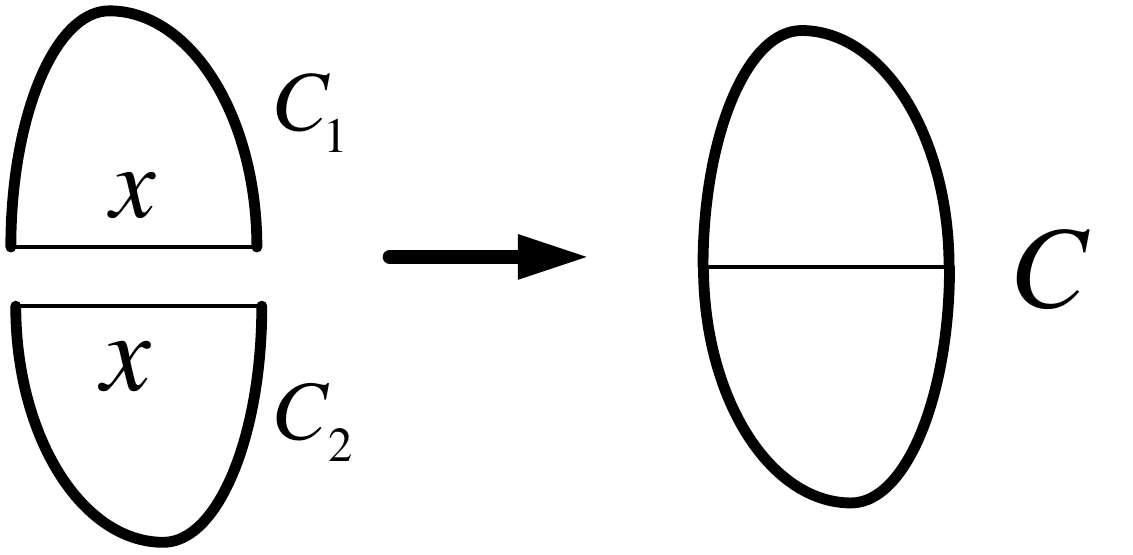}

        \textbf{Figure 6.11.} A diagram for illustrating a cycle-spliced operation
\end{center}

\textbf{II. Bicolored path-splitting operation}

Let $G^C$ be a 4-colorable semi-maximal planar graph on an even-cycle $C$, which is 2-colorable.  Suppose that $f\in C_4^0(G^C)$, and the vertices of $C$ are colored 1,2 alternately. Under $f$, let $P$ be a bicolored path with length $l\geq 2$, the vertices of which are colored 1,2 (or 3,4) alternately (see Figure 6.12(a)). Conducting bicolored path-splitting operation on $P$ in $G^C$, we can gain a new semi-maximal planar graph $G_1^C$ on $C$, and $C$ is also a 2-colorable cycle in $G_1^C$. Following, we will give the detailed description of a bicolored path-splitting operation.

(1) Similarly with the extending $4$-wheel operation, along the direction from one end to another end of $P$, cut a crack to inner vertices and edges of $P$ in accordance with edge-vertex-$\cdots$-vertex-edge order. That is to say except the initial and terminus of $P$, other vertices and edges are cut a crack from their inner side. In this way, each vertex $v$ of $P$ (except the ends) reproduces a new vertex inheriting its color, namely a copy of $v$; and each edge of $P$ reproduces a new edge correspondingly (see Figure 6.12(b)).

(2) Extending path $P$ from the crack produces a bicolored cycle $C^{\prime}$ with length $2(l-1)$ (see Figure 6.12(c)).

(3) Finally, add a tree $T$ to inside of $C^{\prime}$ properly, and connect the vertices of $T$ with the vertices of $C^{\prime}$ properly so that all the faces inside $C^{\prime}$ are triangle, and no vertices have degrees less than 4.

\textbf{III. Bicolored cycle-contracted operation}

Let that $G^C$ be a 4-colorable semi-maximal planar graph on an even-cycle $C$, which is 2-colorable.  Suppose that $f\in C_4^0(G^C)$, and the vertices of $C$ are colored 1,2 alternately.

(1) Under $f$, suppose that $C_p$ is a bicolored cycle that differs from $C$, and the vertices of $C$ are colored 1,2 (or 3,4) alternately. Choose a vertex of $C_p$ arbitrarily, say $v_1$, and then choose another vertex of $C_p$, say $v_2$, which has the longest distance to $v_1$ (see Figure 6.12(c)).

(2) Delete all of vertices inside $C_p$ (see Figure 6.12(c)).

(3) Starting with $v_1$, we identify any pair of vertices that have the same distance to $v_1$ (see Figures 6.12(a) and (b)). Thus, a new 4-colorable semi-maximal planar graph, $G^{C}_{C}$ is obtained, and $C$ is also a 2-colorable in $G^{C}_{C}$.

\begin{center}
        \includegraphics [width=260pt]{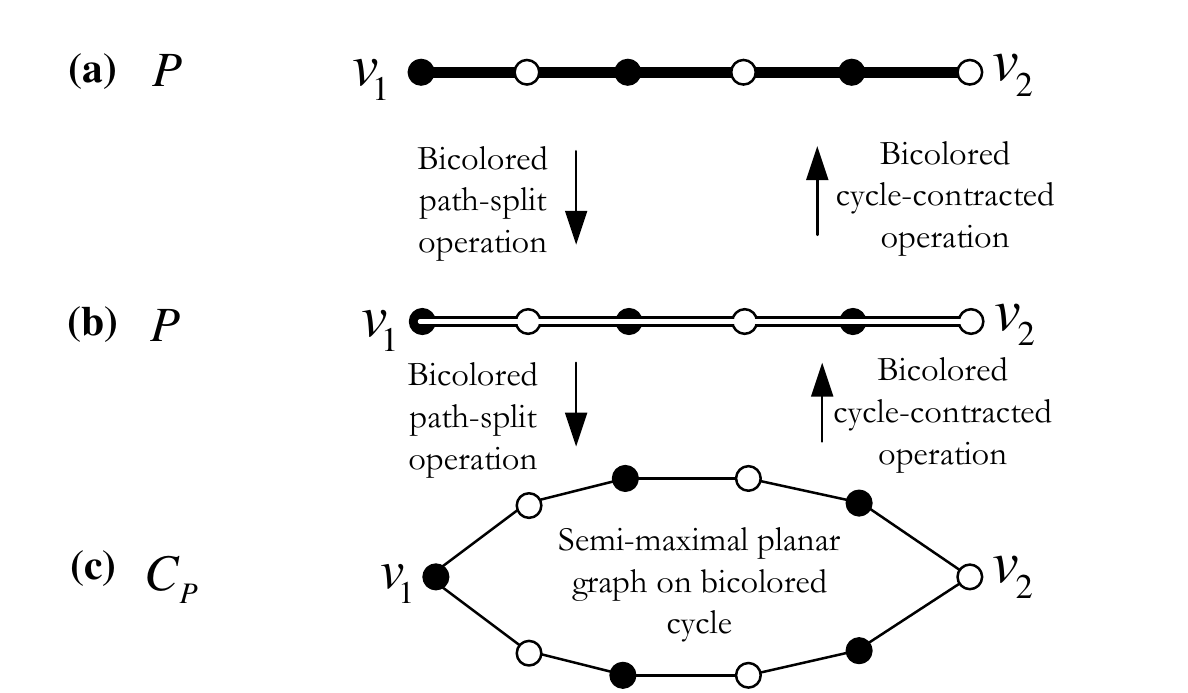}

        \textbf{Figure 6.12.} A Bicolored path-splitting operation and a bicolored cycle-contracted operation
\end{center}

On the basis of the above three operations, now we give a method to construct a semi-maximal planar graph having 2-colorable cycles as follows.

Step 1. Choose some semi-maximal planar graphs $G^{C_i}$ on basic bicolored cycles $C_i$ for $i=1,2,\cdots,m$, $m\geq 2$. Namely, $G^{C_1},G^{C_2},\cdots,G^{C_m}$ belong to the closed-tree type or the fence-tree type;

Step 2. Conducting the cycle-spliced operation to $G^{C_1},G^{C_2},\cdots,G^{C_m}$, we can produce some new semi-maximal planar graphs having 2-colorable cycles, which belong to the intersected-cycle type.

Step 3. The semi-maximal planar graphs having 2-colorable cycles, which belong to the closed-cycle type and closed-fence type, can be constructed by conducting some bicolored path-splitting operation on the semi-maximal planar graphs having 2-colorable cycles, which belong to the closed-path type and the closed-tree type, respectively;

Step 4. The semi-maximal planar graphs having 2-colorable cycles, which belong to the fence-cycle type, the scycle-forest type and the scycle-scycle type, can be constructed by conducting some bicolored path-splitting operation on the semi-maximal planar graphs having 2-colorable cycles, which belong to the fence-path type;

Step 5. For any semi-maximal planar graphs with 2-colorable cycles belonging to one of the nine types in Theorem 6.16, we can always produce various semi-maximal planar graphs having 2-colorable cycles by conducting the bicolored path-splitting operation and the cycle-spliced operation, simultaneously.

Step 6. For any semi-maximal planar graph $G_C^{C}$ on compound bicolored cycles, when we conduct some bicolored cycle-contracted operation to $G_C^{C}$, a semi-maximal planar graph $G_C^{C}$ based on a basic bicolored cycle will be produced.

\subsubsection{Folded operation, unfolded operation and the characteristics of 2-colorable cycles}

We have known that a semi-maximal planar graph $G^{C}$ on the basic bicolored cycle $C$ belongs to either the closed-tree type or the fence-tree type. For the former, its structure is very clear that the subgraph $G^{\prime}$ induced by the vertices inside cycle $C$ is a tree; but for the latter, $G^{C}$ can be viewed as the union of a fence $G_F$ on cycle $C$ and a tree $T$. Namely, $G^{C}=G_F\cup T$ and $V(G_F)\cap V(T)=\emptyset$. So, $G^{\prime}$ contains not only the vertices in $T$, but also the vertices in $G_F-C$. Then, what are the characteristics of $G^{\prime}$? The following will discuss this problem, for which we need a pair of operation on even-cycles: the \emph{folded operation} and the \emph{unfolded operation}.

\textbf{I. Folded operation}

Actually, the so-called folded operation is to conduct the action of identifying some succussive pairs of vertices. Suppose that $G^{C}$ is a semi-maximal planar graph on a cycle $C$, which belongs to closed-tree type, and $f\in C_4^0(G^C)$. When $|C|\geq 6$, choose a vertex of $C$ arbitrarily, say $u$, and identify the two vertices adjacent to $u$ on $C$, say $v_1,v_1'$. Define that the new 4-colorable semi-maximal planar graph $ G^{C}\circ\{v_1,v_1'\}\triangleq G_1^{C_1}$,  $\{v_1,v_1'\}\triangleq v_1$ and $C_1$ is the new outer cycle (see Figure 6.13).
If $|C_1|\geq 6$, then identify the two vertices adjacent to $v_1$ on $C_1$, say $v_2,v_2'$, and we will produce another new semi-maximal planar graph $G_2^{C_2}$ on a cycle $C_2$. This procedure can be continued until some new outer cycle $C_i$ satisfying \emph{$|C_i|=4$}.

\textbf{II. Unfolded operation}

In fact, the so-called unfolded operation is the inverse operation of a folded operation. Differing from the bicolored path-splitting operation, for a bicolored path $P=v_1v_2\cdots v_l$, we can split one of its ends $v_l$, the terminal end,  but it is not permitted in bicolored path-splitting operation. The detailed description is as follows:

Step 1. For a semi-maximal planar graph $G^{C_1}$ on a cycle $C_1$, which belongs to the fence-tree type, and $f\in C_4^0(G^{C_1})$.  Suppose $G_F$ and $T$ are the fence and the tree, respectively. Namely, $G_F\cup T=G^{C_1}$, and $C\subseteq G_F$. We choose a bicolored $t$-path in $G_{F}$, denoted $P=u-v_1$, which only contains one vertex $v_1$ of $C_1$, namely $v_1$ is one end of $P$;

Step 2.  Starting with the vertex $u$, along the direction from $u$ to $v_1$ on the path $P$, cut a crack the inner vertices and edges of $P$ in accordance with edge-vertex-$\cdots$-edge -vertex order. That is to say except the end $u$, all vertices and edges of $P$ are cut a crack from their inner side. In this way, each vertex $v$ of $P$ (except $u$) reproduces a new vertex inheriting its color, namely a copy of $v$ and each edge of $P$ reproduces a new edge correspondingly;

Step 3. Extending the path $P$ from the crack and a new semi-maximal planar graph will be obtained.

Obviously, when we conduct some unfolded operations for all the paths of $G_F$, which only contain one vertex of $C_1$, every resulting graph is going to be a semi-maximal planar graph belonging to the closed-tree type.
For example, Figure 6.13(b) is a semi-maximal planar graph belonging to the fence-tree type. Of course, it is opened and $u$ is an opened-vertex. Conducting an unfolded operation on path $P=uv_1$, we can obtain the graph shown in Figure 6.13(a).

\begin{center}
         \includegraphics [width=280pt]{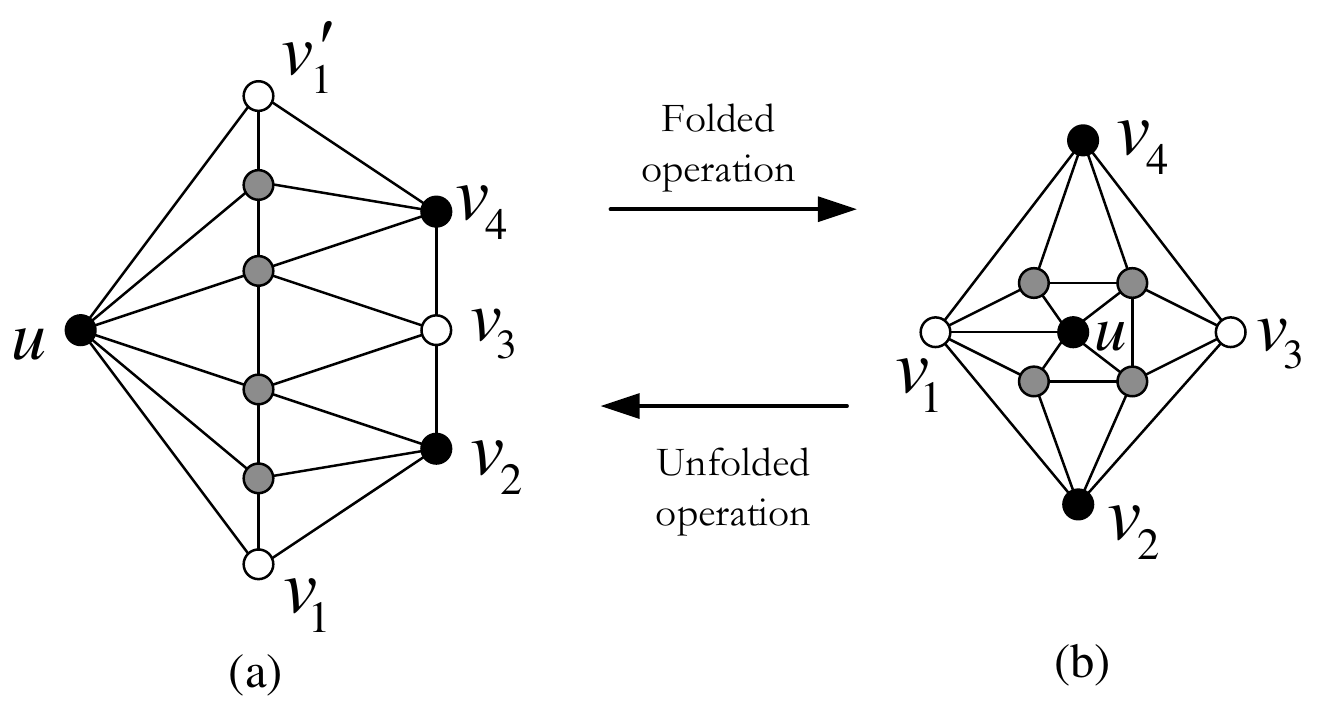}

        \textbf{Figure 6.13.} Schematic diagram for illustrating folded operation and unfolded operation
\end{center}

Considering the case that in $G_F$, some suspending vertices connect with the vertex $v_1$ of $C_1$ through a tree. That is to say, in $G_F -(V(C_1)\backslash \{v_1\})$, the connected branch containing $u$ is a tree, denoted $T'$, not a path. In this case, first, we should choose a path $P$ that starts with a vertex $v_1$ of $C_1$, which is adjacent to opened-vertices, to a vertex $u'$ of $T'$, which has degree not less than 3. Conducting the unfolded operation on $P$, we can obtain a new semi-maximal planar graph $G^{C_2}$, which  belongs to the fence-tree type. Obviously, if we remain the colors appearing in the vertices of $G^{C_1}$ unchanged, including the copies of the vertices of $P$, then the 4-coloring of $G^{C_2}$ is also proper.
Next, similarly with the above process, conducting the  unfolded operation to $G^{C_2}$ enables us to obtain another 4-colorable semi-maximal planar graph $G^{C_3}$. Continue this procedure until to some semi-maximal planar graph $G^{C_i}$, which belongs to the closed-tree type. Figure 6.14 gives an illustration of this procedure, in which Figure 6.14(a) gives a semi-maximal planar graph belonging to the fence-tree type, Figure 6.14(b) gives the new semi-maximal planar graph after conducting the unfolded operation to a path $P=uv$, Figure 6.14(c) gives another semi-maximal planar graph after conducting the unfolded operation to a path $P=uu_1$, and Figure 6.14(d) gives another semi-maximal planar graph after conducting unfolded operation on path $P=uu_2$.

\begin{center}
        \includegraphics [width=280pt]{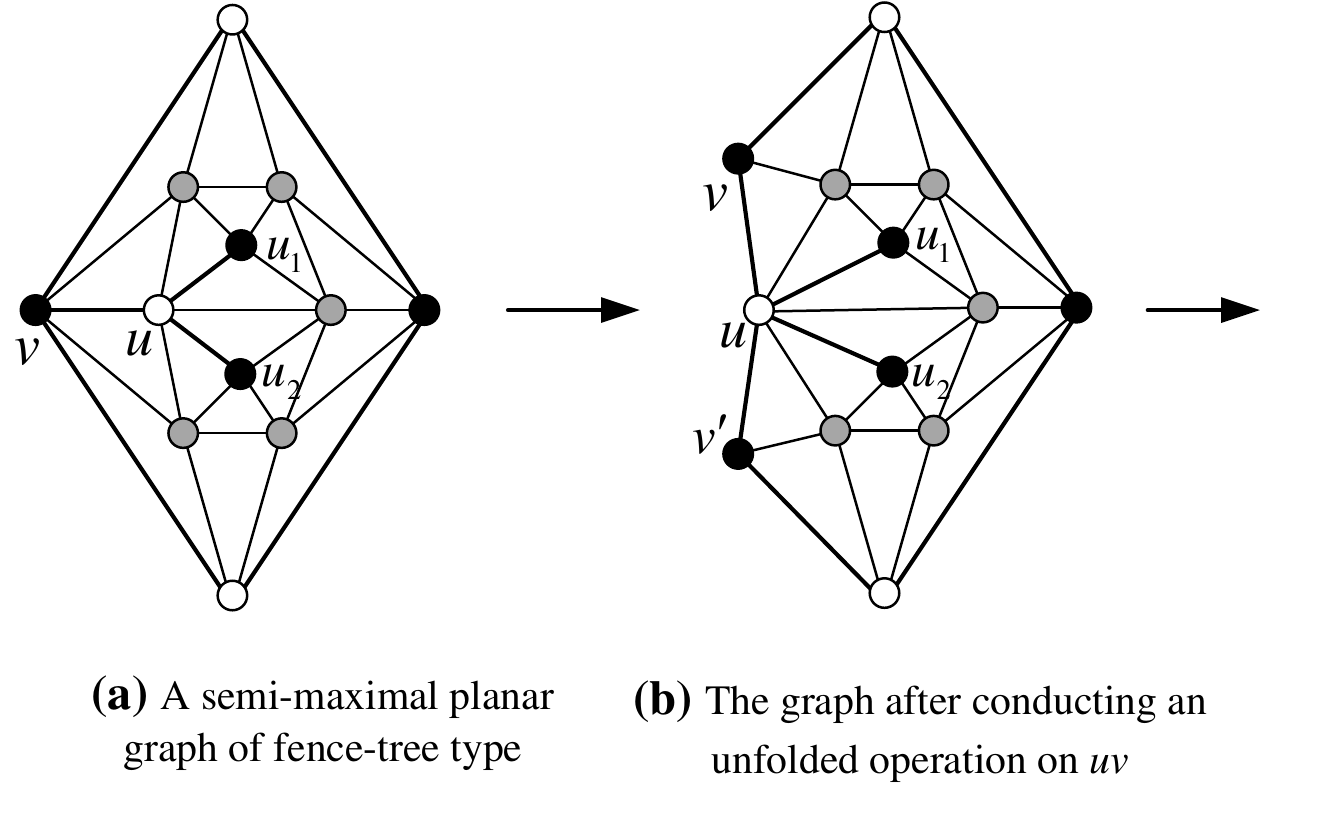}

         \includegraphics [width=300pt]{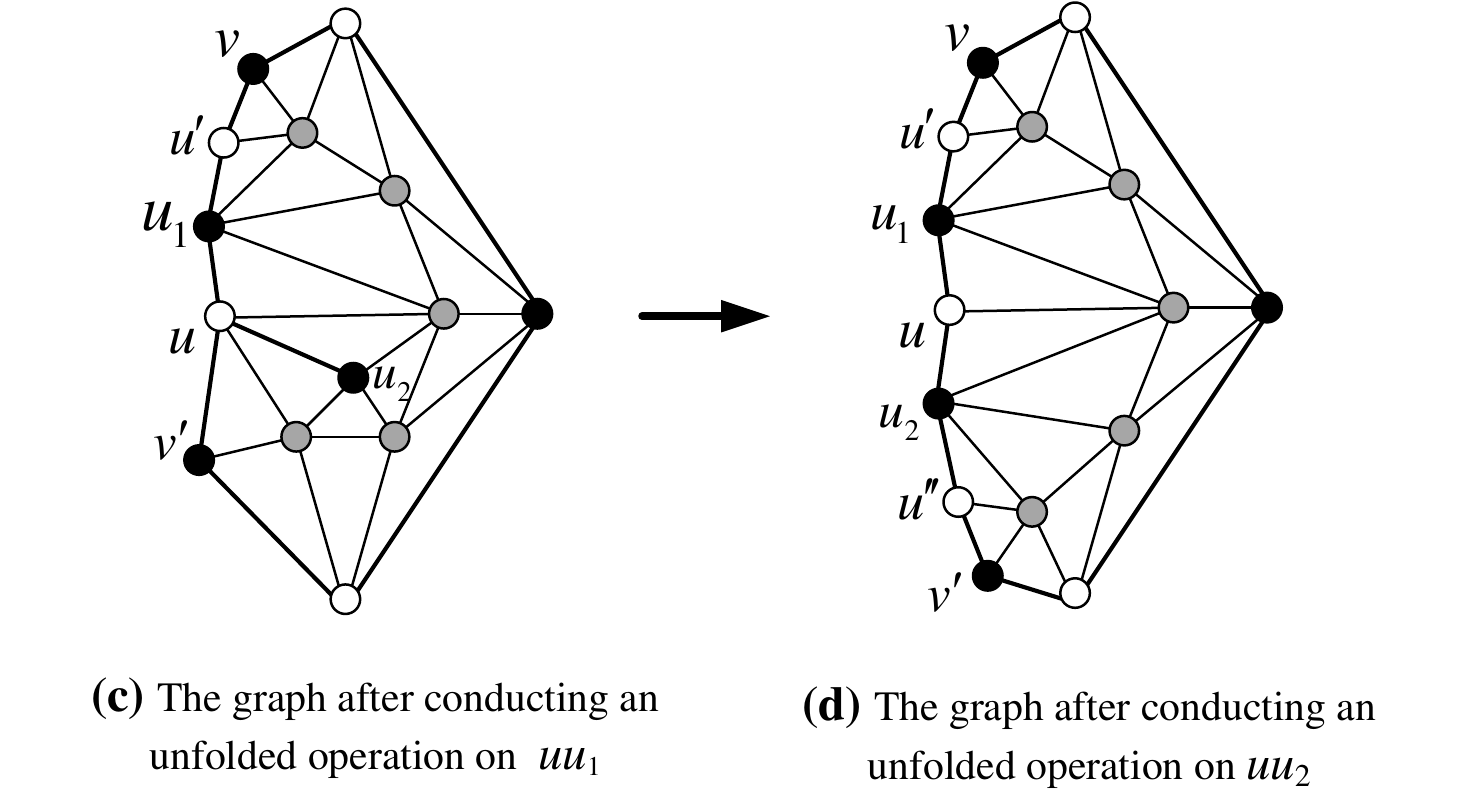}

        \textbf{Figure 6.14.} A diagram of unfolded operations in special case
\end{center}

It follows the following fact from the above two operations.

\begin{theorem2}\label{th6.17}
Suppose that $G^{C}$ is a semi-maximal planar graph belonging to the closed-tree type with $|C|\geq 6$, $V'=V(G^{C})-V(C)$, and $G[V']\triangleq G'$ is a tree. When conducting a folded operation to $C$, we can obtain a new semi-maximal planar graph $G^{C_1}$, which belongs to the fence-tree type. Denote by $G_F$ and $T$ the fence and tree of $G^{C_1}$ respectively, then $T=G'$. That is to say, when we conduct the folded operation to a semi-maximal planar graph belonging to the closed-tree type, the structure of its fence gets changed, but the tree remains unchanged. Conversely, for a semi-maximal planar graph $G^{C_1}$ belonging to the fence-tree type, when we conduct unfolded operation to $G^{C_1}$ repeatedly, finally, we can obtain a semi-maximal planar graph belonging to the closed-tree type, and the structure of the tree $T$ will not be affected in the process of unfolded operation.

\end{theorem2}

This theorem actually tells us that fence-tree structures can be obtained from closed-tree structures by conducting some folded operations to them, repeatedly. Because the tree structure remained unchanged in the process of folded operations, we need only to study the change of edges on the cycle after conducting this operation.

Here, we give an example for constructing a semi-maximal planar graph on two paths. Suppose $P,P'$ are two paths with length not less than 3. Now, we will construct a semi-maximal planar graph $G^{C}$ from $P$ and $P'$ (see Figure 6.15), and the resulting graph is called the semi-maximal planar graph of \emph{cycle-path type}.

Step 1. Connect one end $u$ of $P'$ to the two ends of $P$, respectively, then connect another end of $P'$ to at least three vertices of $P$;

Step 2. For the inner vertices of $P'$, connect them to some vertices of $P$ such that each connected edge is in a triangle, and the degrees of vertices on $P$, except $u$, have to increase to 4 at least.

\begin{center}
         \includegraphics [width=360pt]{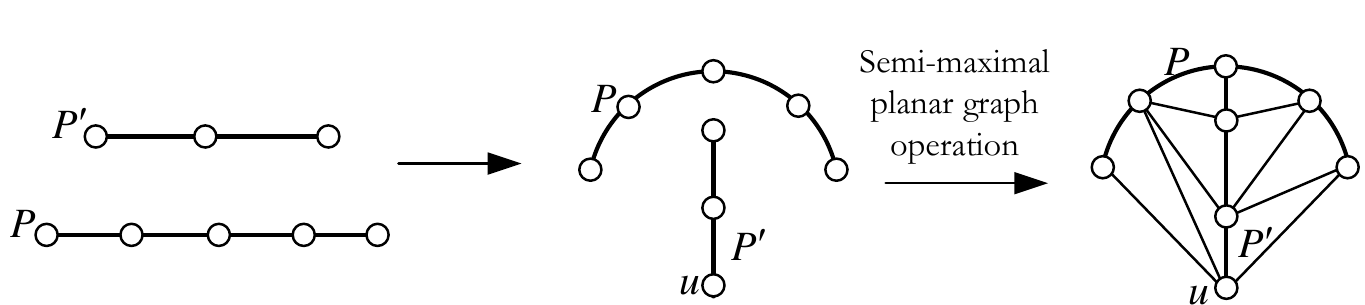}

        \textbf{Figure 6.15.} A diagram of constructing a semi-maximal planar graph belonging to the closed-tree type on two paths
\end{center}

In the process of constructing semi-maximal planar graphs on two paths, we can deem that the path $P'$ is obtained by identifying pairs of vertices with the same color in a semi-maximal planar graph $G^{C}$ belonging to the closed-tree type repeatedly. That is to say, $P'$ generated in the process of conducting the folded operation to $G^{C}$, is a path in a semi-maximal planar graph $G^{C_1}$ belonging to fence-tree type. If $G^{C_1}$ can be obtained only by conducting folded operation once in $G^{C}$, then the induced graph $G^{C_1}[V_1' \cup\{u\}]$ is the join of a semi-maximal planar graph belonging to the cycle-path type
and some possible trees, where $V'=V(G^{C_1})-V(C_1)$ and $u$ is a weld-vertex in the fence; if $G^{C_1}$ is obtained by conducting the folded operation many times to $G^{C}$, then this process can be seen turned back as conducting cycle-spliced operation among semi-maximal planar graphs belonging to the closed-path type, and then connecting some possible trees on the outer cycle. We call these graphs \emph{barrette-structure} graphs. Thus, we have proved that if $G^{C_1}$ is a semi-maximal planar graph belonging to fence-tree type, then the subgraphs induced by the inside vertices of $C_1$ and the vertices on $C_1$ joined the trees in the fence is a barrette-structure graph.

\begin{theorem2}\label{th6.18}
Suppose that $G^{C_1}$ is a semi-maximal planar graph belonging to the fence-tree type, $u_1,u_2,\cdots,u_k$ are all weld-vertices in the fence, and $V'=V(G^{C_1})-V(C_1)$, then $G^{C_1}[V_1' \cup \{u_1,u_2,\cdots,u_k\}]$ is a barrette-structure graph.
\end{theorem2}

\subsection{Summary}

The contents of this chapter mainly cover the following:

\textbf{First}, we point out that there were two categories of cycles in maximal planar graphs, called the basic cycles and the chord cycles, and study their distribution and enumeration, which paves the way to the latter research.

\textbf{Second}, we propose a new method, the Black-White coloring, to study maximal planar graphs. The advantages of this method are that its process can be realized simply, and for a maximal planar graph, each of its Black-White coloring consists of a subset of its 4-colorings set. Especially, it is a powerful technique to study the 2-colorable cycles.

\textbf{Third}, we set up the petal-syndrome, on which we find a necessary and sufficient condition that an even-cycle is 2-colorable in a maximal planar graph.

\textbf{Fourth}, we find a necessary and sufficient condition of a 2-colorable cycle based on structure. Namely, independently satisfy two basic types and seven compound types (see Theorem 6.16).

\textbf{Fifth}, we prove that each compound type can be obtained from basic type by conducting the cycle-spliced operation, the bicolored path-splitting operation and the bicolored cycle-contracted operation.

\textbf{Sixth}, we make clear the relationship between two basic types through introducing the folded operation and the unfolded operation, and depict the inner structure of the semi-maximal planar graphs belonging to fence-tree type: a barrette-structure graph. Thereby, the structure of the semi-maximal planar graphs with 2-colorable cycles are described deeply.

However, for the compound types, it is still a tough problem to judge which type they belong to the seven cases in Theorem 6.16. In order to deal with this problem, we need to argue it combined with the open-vertices properly.

The more in-depth study on this problem will be given in later articles.


    \vspace{5mm}
    \begin{center}
    \textbf{Acknowledgements}
    \end{center}
    \vspace{5mm}

This series of articles is completed based on the modification and improvement of the paper "Mathematical proofs of two conjectures: the four color problem and the uniquely 4-colorable planar graph", which was finished in May, 2010. This is our first article. I am indebted to many friends and colleagues for their help with this paper, especially my two Ph.D students, Enqiang Zhu and Zepeng Li, who deserve a special word of thanks. They spent much time and energy reading through the entire manuscript, and discussed meaningfully with me the problems arisen in the manuscript. In particular, almost all of the figures drawing, languages checking, and translation from Chinese to English were accomplished by them together. In addition, I would like to thank my graduate student Yang Yang, who also drew some figures that appear in the Appendix, and post-doctor Yan Qu, who participated in the final checking of the English version.

As soon as the manuscript was completed, many experts in this field carried on detailed scrutiny. They provided numerous unfailingly pertinent comments, corrected some technical errors and linguistic infelicities, and made valuable suggestions.  Professors Bing Yao and Xiang'en Chen in Northwest Normal University organized a seminar to discuss the manuscripts weekly, where students who participated in the seminar as a speaker include Yuping Gao, Xiangqian Zhou, Jiajing Wei, Jiajuan Zhang, Weihua Lu, Fanghong Zhang and Zhitao Hu. In addition, there were twelve more students who took part in the discussion class, including Zhiqiang Wang, Jingxia Guo, Chunhu Sun, Wenjuan Liu, Xiaohui Liu, Wangfa Liu, Chunyan Ma, Fang Yang, Yuanyuan Liu, Yan Gou, Hongyu Wang, Chao Wang. I thank them all warmly for their various contributions. I am grateful also to professor Suixiang Gao in Chinese Academy of Sciences and his student Wei Zhang for over half a year's reviewing, and Professor Zhixiang Yin in Anhui University of Science and Technology for  nearly a year's reviewing, and professor Daoheng Yu in Peking University for more than two months' reading. I am most fortunate to benefit from their excellent knowledge and taste. So, this paper also embodied their energies and painstaking effort. Here, I express my sincere gratitude to them once again. In particular, Prof. Bing Yao not only reviewed the Chinese version for more than two years, but also checked the English writing at length. Here, I express deep acknowledgment to him.

Finally, the accomplishment of this paper was closely related to the facilities and comfortable environment provided kindly by Peking University, and the encouragement and support provided by colleagues, especially Academicians Fuqing Yang, Xingui He, and Hong Mei, and Professors Wanling Qu and Hanpin Wang.  I would like to express my sincere thanks to them also.

\section{Appendix}
    This appendix gives all 4-colorings of the maximal
 planar graphs whose orders are from 6 to 11 and
 $\delta(G)\geq 4$ excepting 3-colorable and divisible graphs.

 1. There is only one maximal planar graph of order 6 whose minimal
 degree is 4. Its degree sequence is 444444, and it is a 3-colorable graphs.

 The following figures show the drawing of this graph and its unique
 3-coloring:
\begin{center}
        \includegraphics [width=380pt]{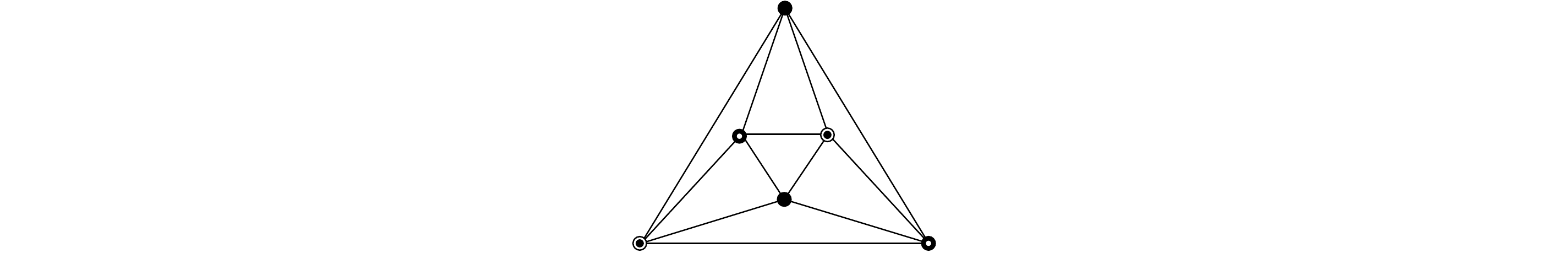}

  \end{center}

 2. There is only one maximal planar graph of order 7 whose
 degree sequence is 4444455. It has 5 different 4-colorings.
   \begin{center}
        \includegraphics [width=380pt]{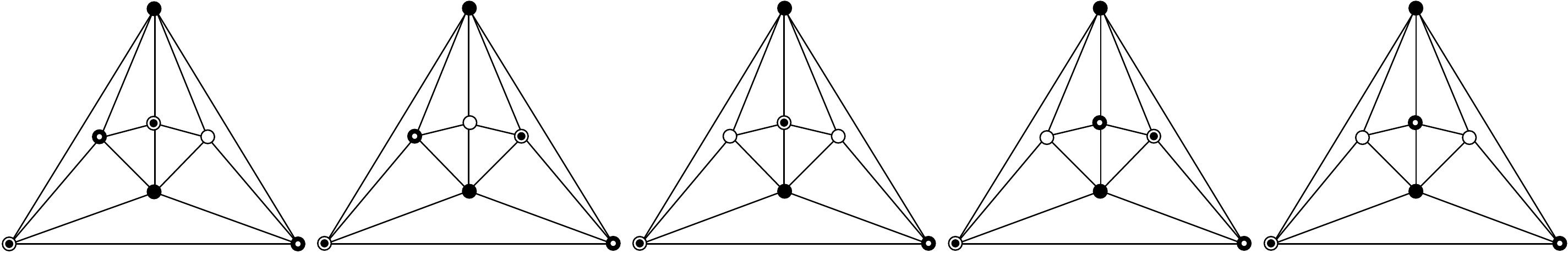}

   \end{center}

 3. There are two  maximal planar graphs of order 8 whose minimal
 degree are 4.

 3.1 Degree sequence is 44444466, and it is uniquely 3-colorable.

  \begin{center}
        \includegraphics [width=380pt]{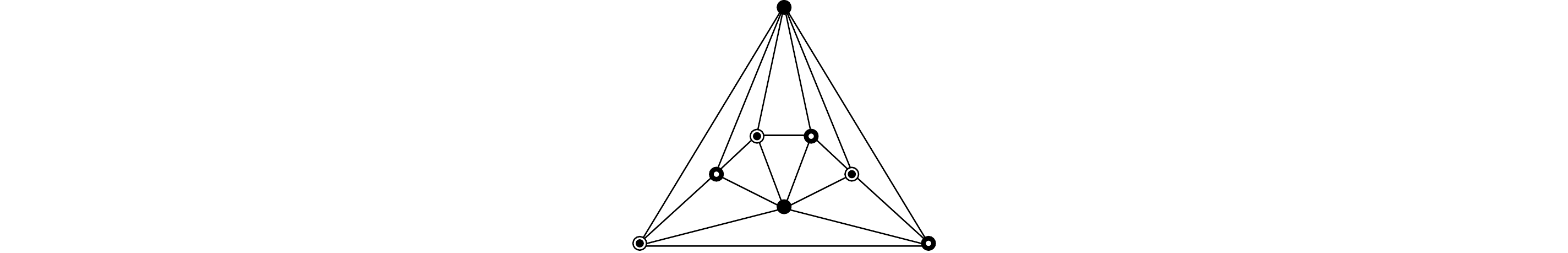}
  \end{center}

 3.2 Degree sequence is 44445555, and it has 3 kinds of different 4-colorings.

  \begin{center}
        \includegraphics [width=380pt]{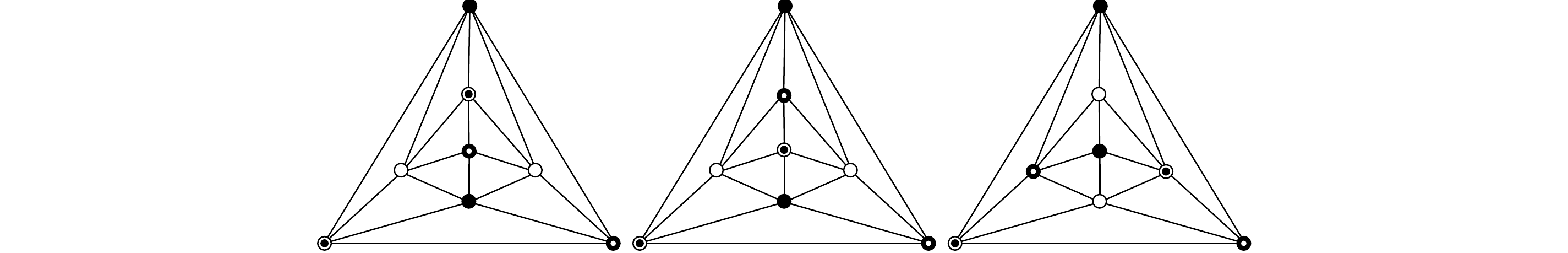}
  \end{center}

 4. There are five  maximal planar graphs of order 9 whose minimal
 degree is 4.

 4.1 Degree sequence is 444444666, and it is 3-colorable.
  \begin{center}
        \includegraphics [width=380pt]{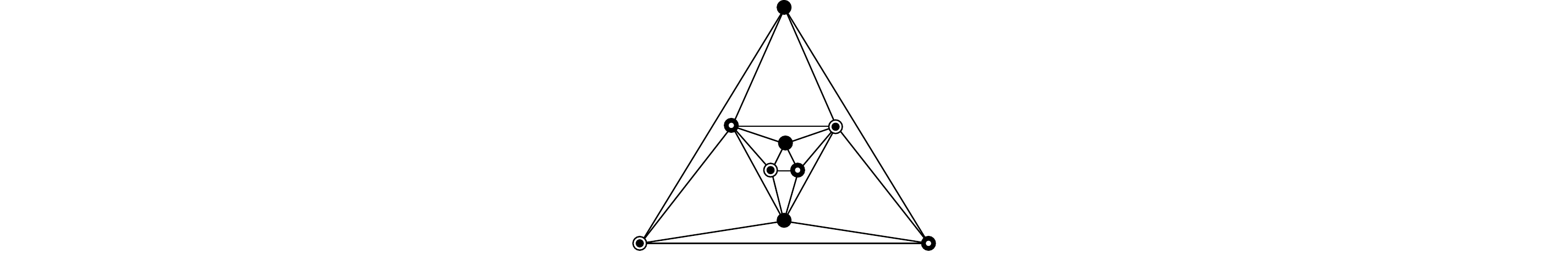}
  \end{center}

 4.2 Degree sequence is 444455556, and it has 6 kinds of different colorings.
 \begin{center}
        \includegraphics [width=380pt]{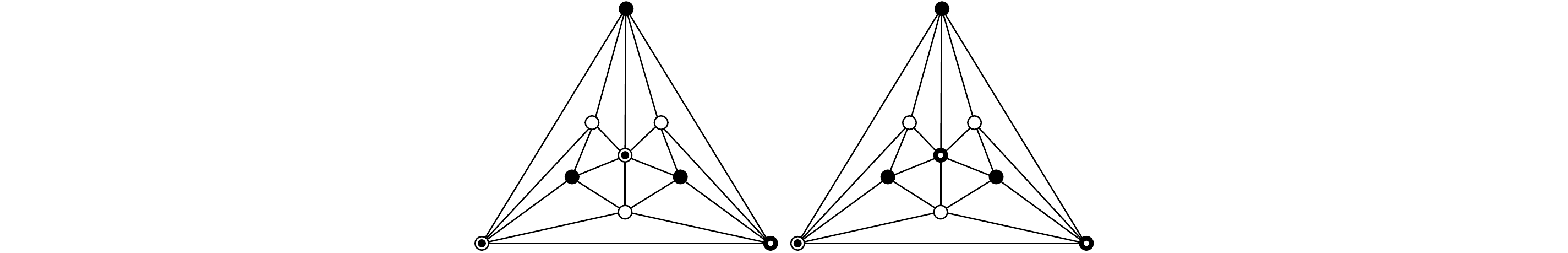}

        \vspace{2mm}
        \includegraphics [width=380pt]{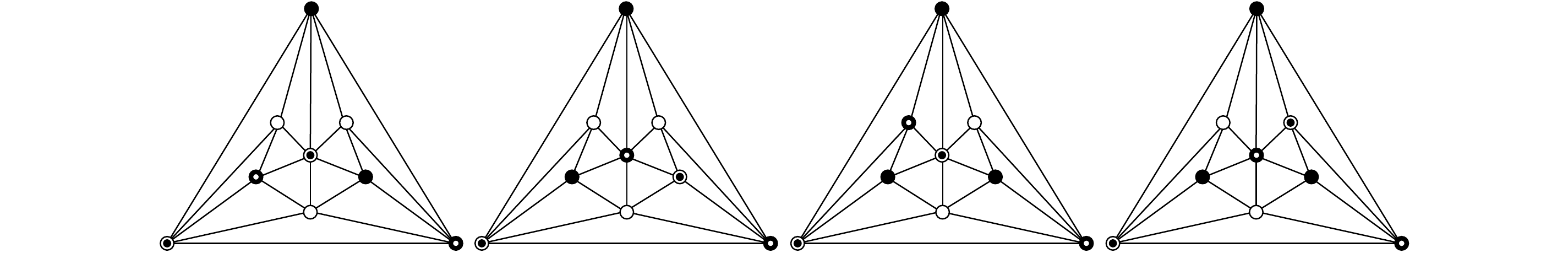}
  \end{center}

 4.3 Degree sequence is 444555555, and it has 2 kinds of different colorings.
 \begin{center}
        \includegraphics [width=380pt]{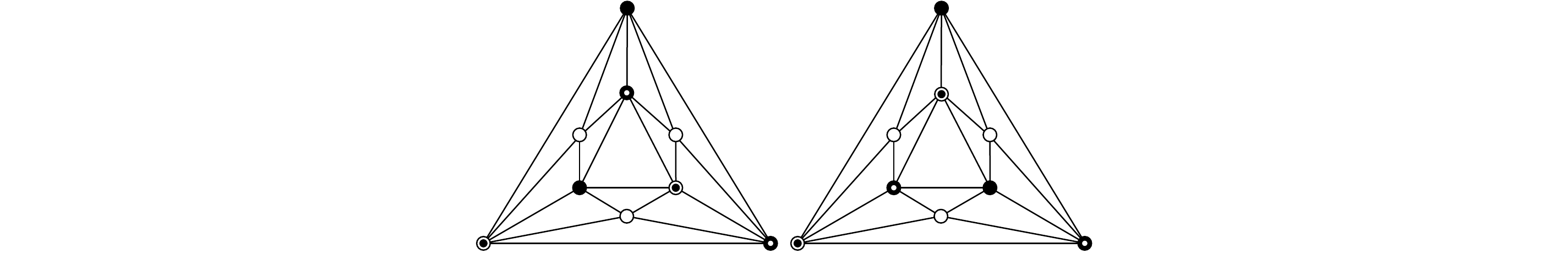}
  \end{center}

 4.4 Degree sequence is 4444444477, and it has 17 kinds of different colorings.
  \begin{center}
        \includegraphics [width=380pt]{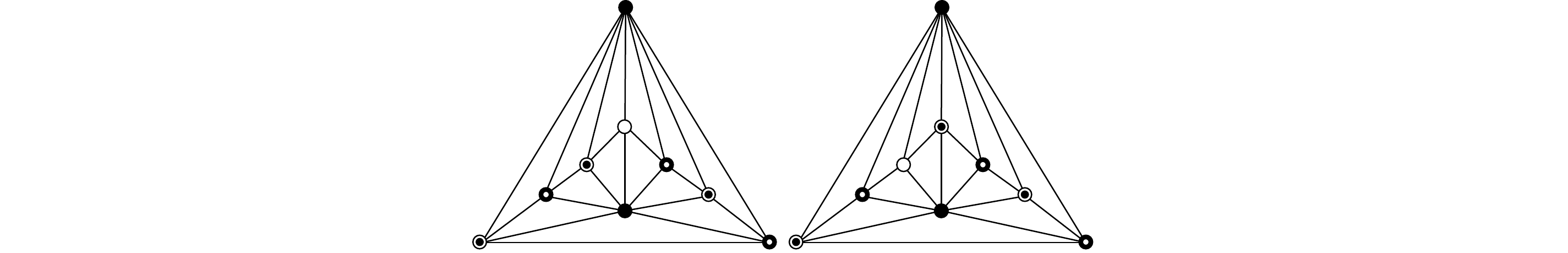}

        \vspace{2mm}
        \includegraphics [width=380pt]{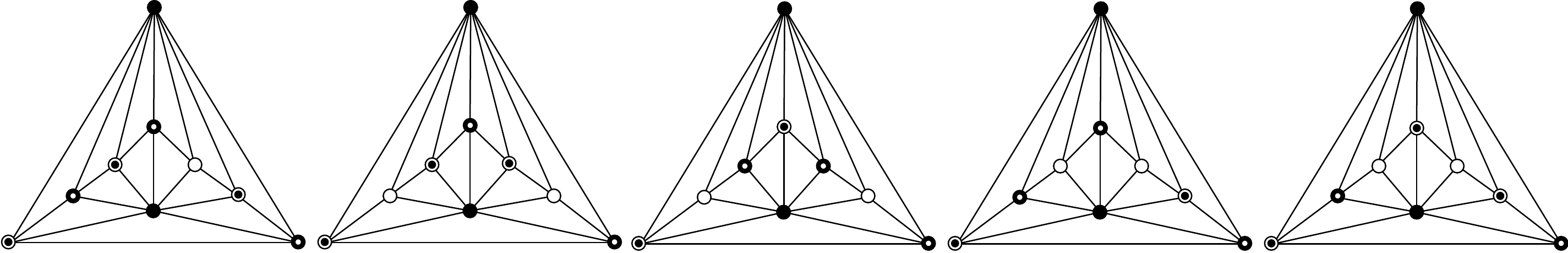}

        \vspace{2mm}
        \includegraphics [width=380pt]{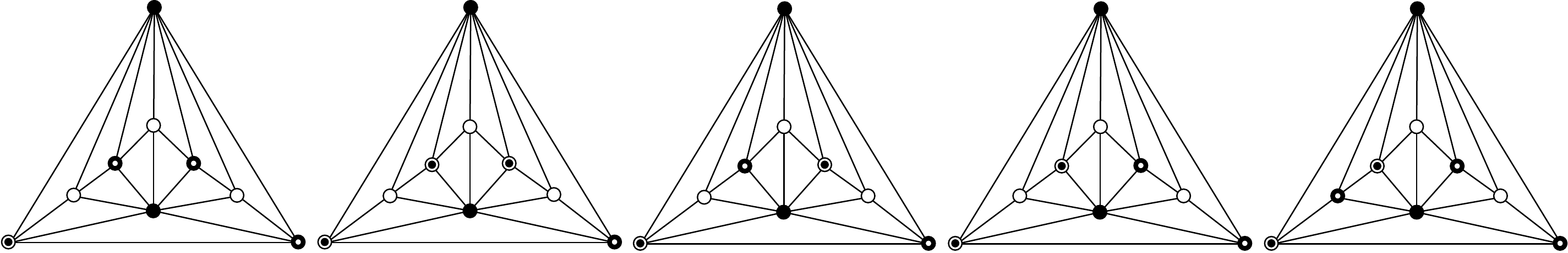}

        \vspace{2mm}
        \includegraphics [width=380pt]{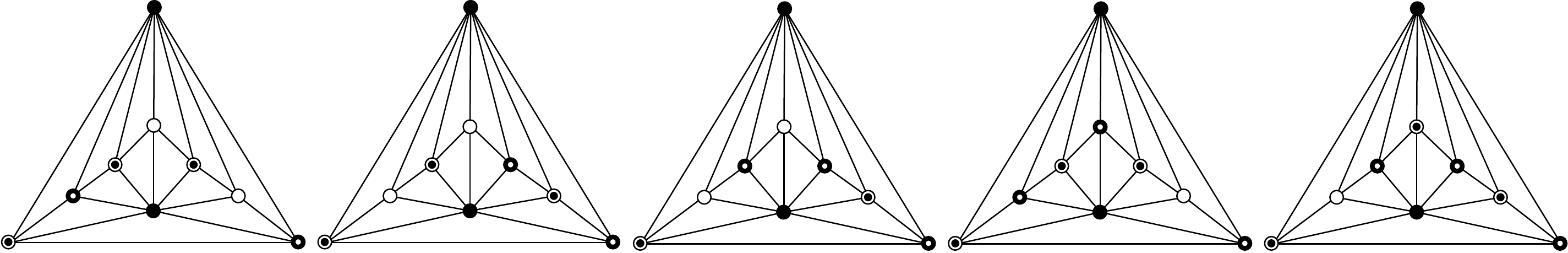}
  \end{center}

 4.5 Degree sequence is 4444445566, and it has 9 kinds of different colorings.
  \begin{center}
        \includegraphics [width=380pt]{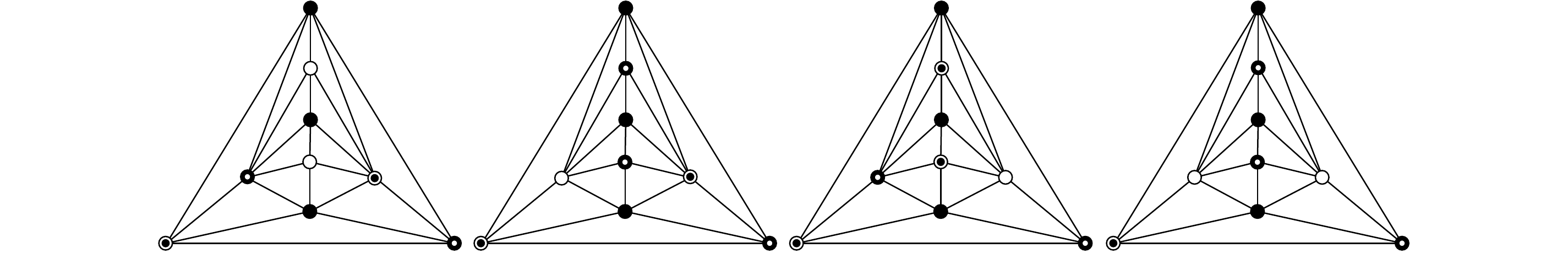}

        \vspace{2mm}
        \includegraphics [width=380pt]{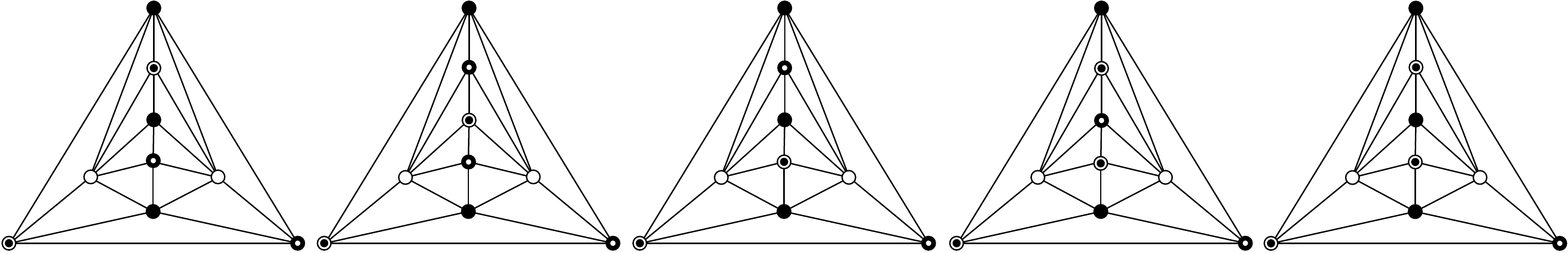}
  \end{center}

 5. There are 12 maximal planar graphs of order 10 whose minimal
 degree is 4.

 5.1 Degree sequence is 4455555555, and it has 8 kinds of different colorings.
    \begin{center}
        \includegraphics [width=380pt]{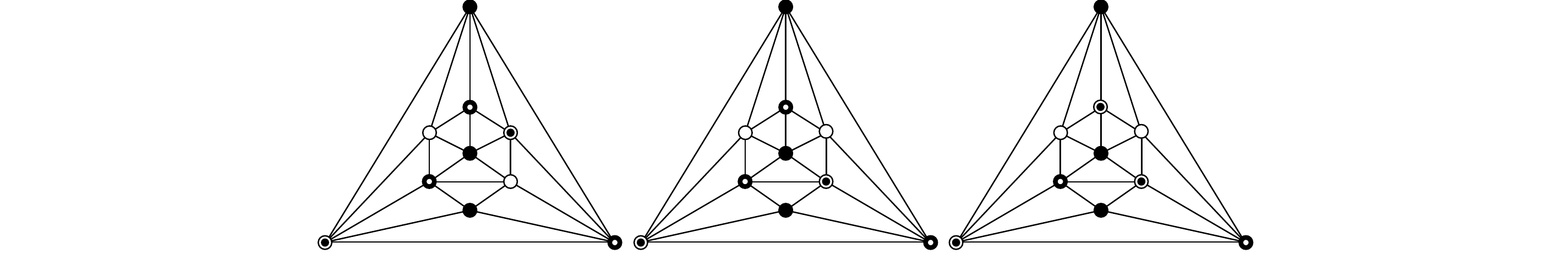}

        \vspace{2mm}
        \includegraphics [width=380pt]{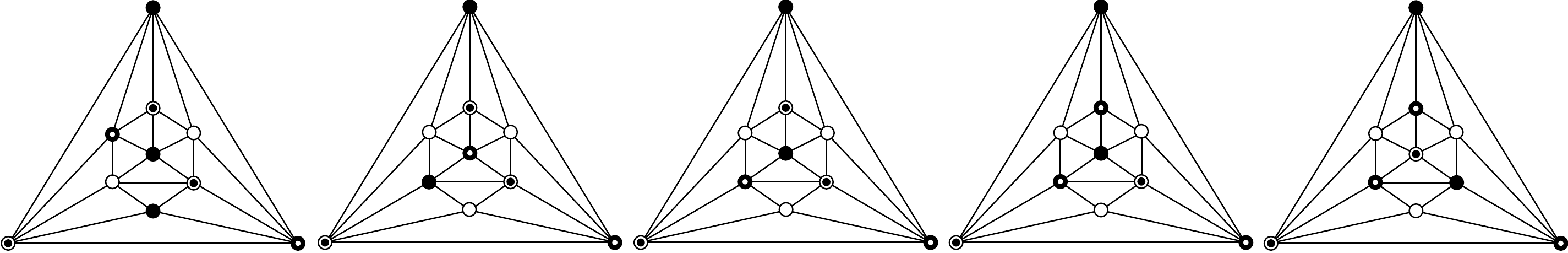}
   \end{center}

 5.2 Degree sequence is 4445555556, and it has 6 kinds of different colorings.
    \begin{center}
        \includegraphics [width=380pt]{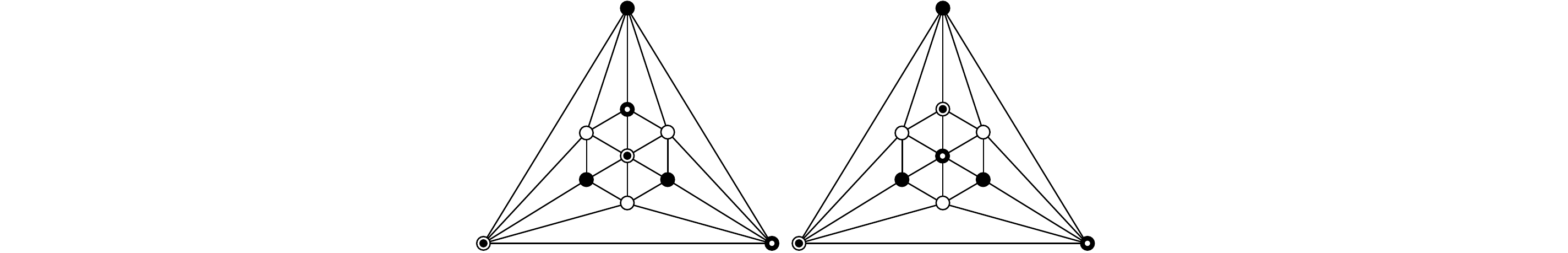}

        \vspace{2mm}
        \includegraphics [width=380pt]{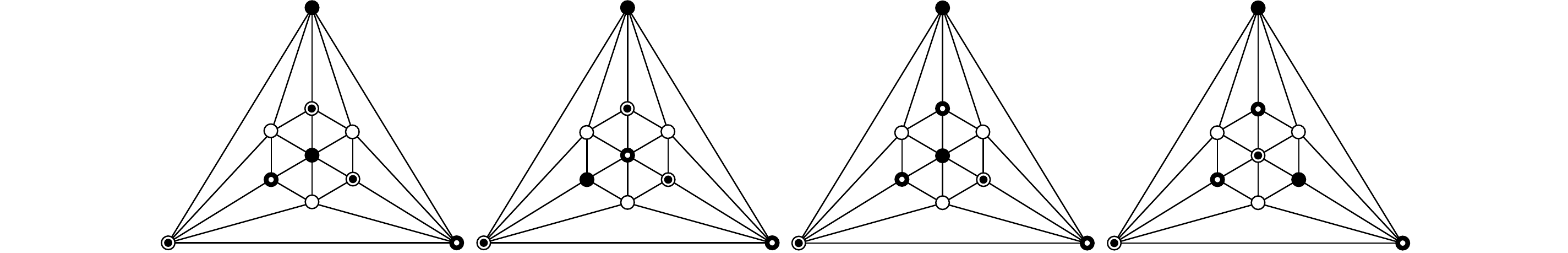}
   \end{center}

 5.3 Degree sequence is 4444555566, and it has 14 kinds of different colorings.
      \begin{center}
        \includegraphics [width=380pt]{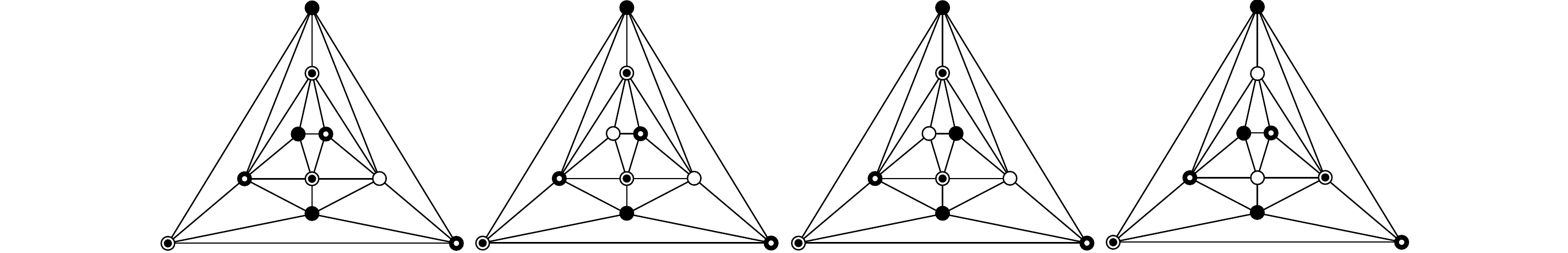}

        \vspace{2mm}
        \includegraphics [width=380pt]{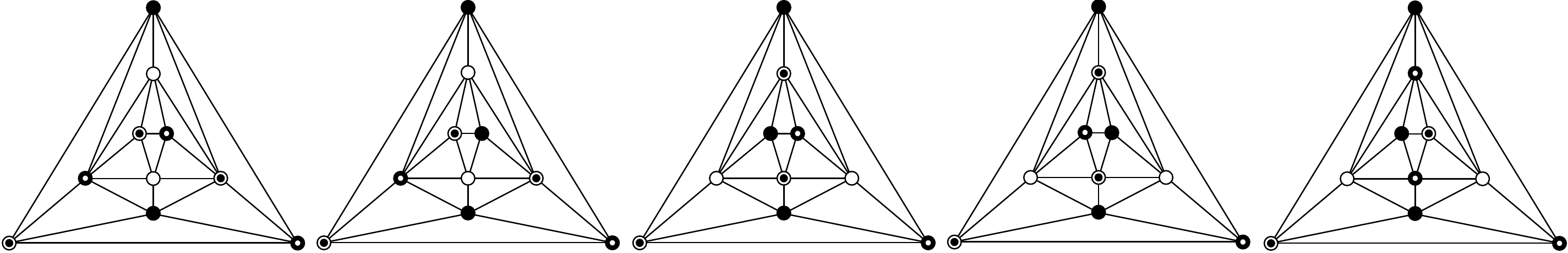}

        \vspace{2mm}
      \includegraphics [width=380pt]{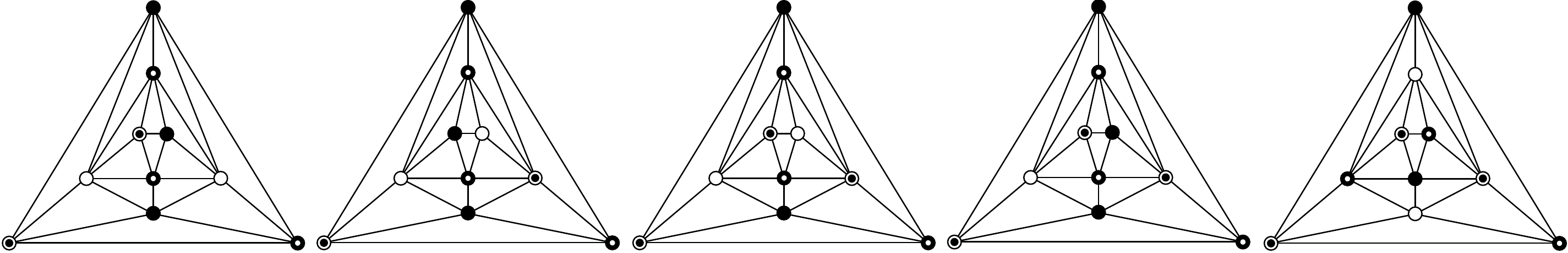}
    \end{center}
 5.4 Degree sequence is 4444555566, and it has 13 kinds of different colorings.
  \begin{center}
        \includegraphics [width=380pt]{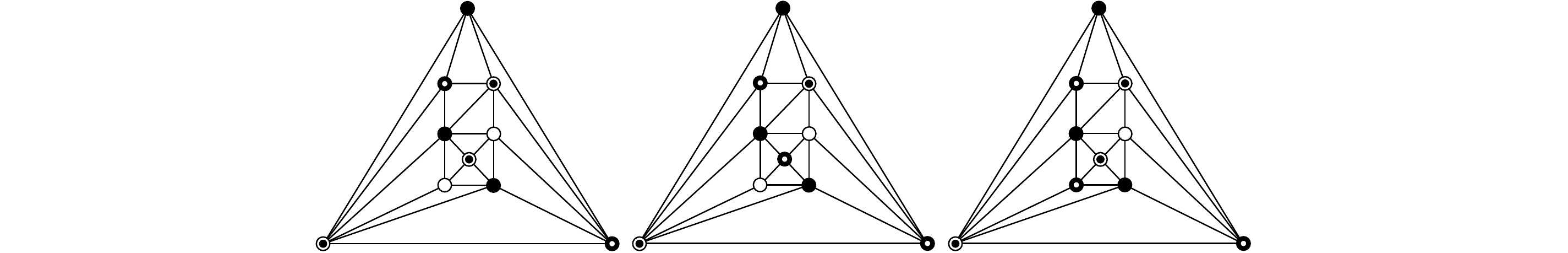}

        \vspace{2mm}
        \includegraphics [width=380pt]{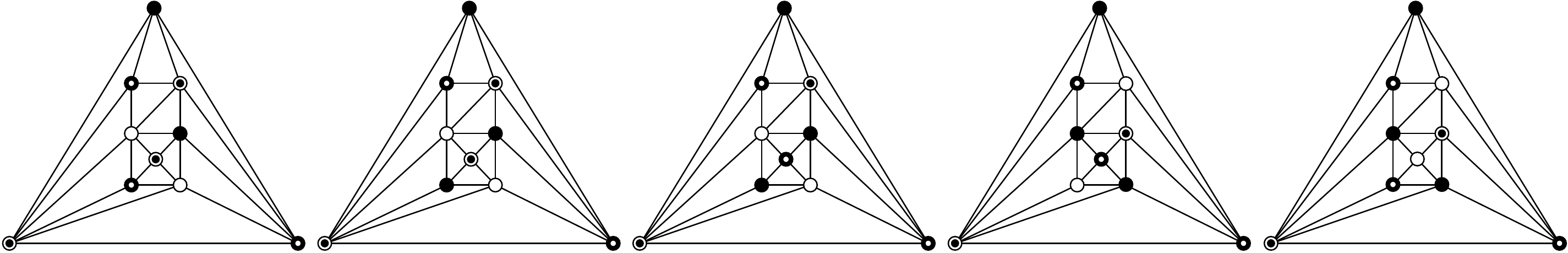}

        \vspace{2mm}
      \includegraphics [width=380pt]{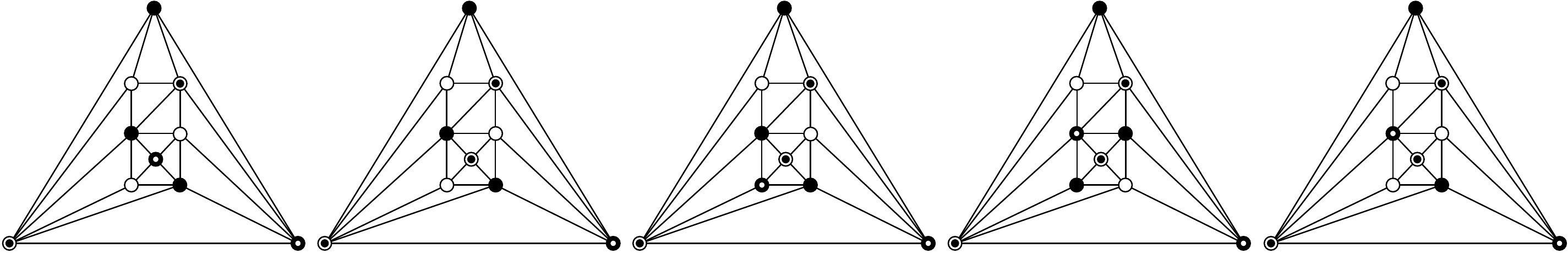}
    \end{center}

 5.5 Degree sequence is 4444555566, and it has 5 kinds of different colorings.
  \begin{center}
        \includegraphics [width=380pt]{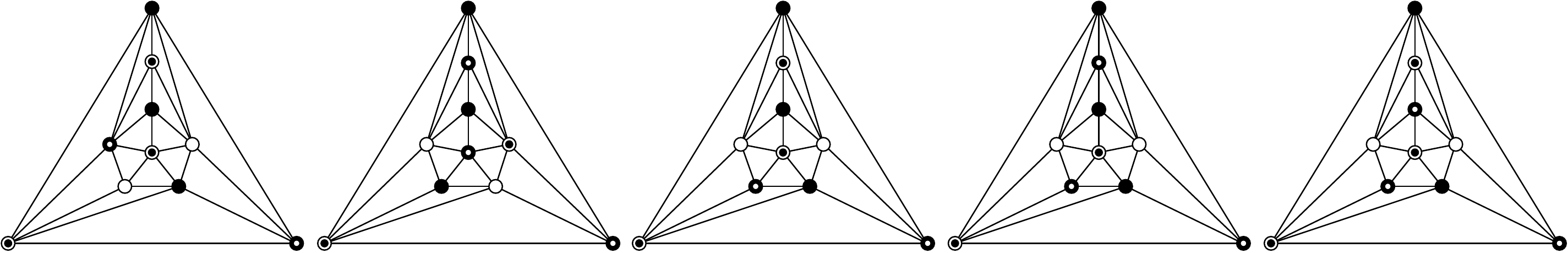}

    \end{center}

 5.6 Degree sequence is 4444455567, and it has 7 kinds of different colorings.

   \begin{center}
        \includegraphics [width=380pt]{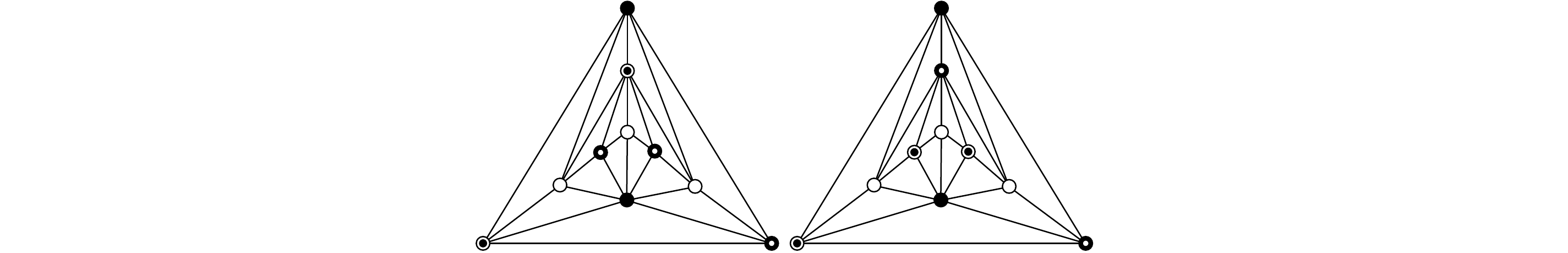}

          \vspace{2mm}
        \includegraphics [width=380pt]{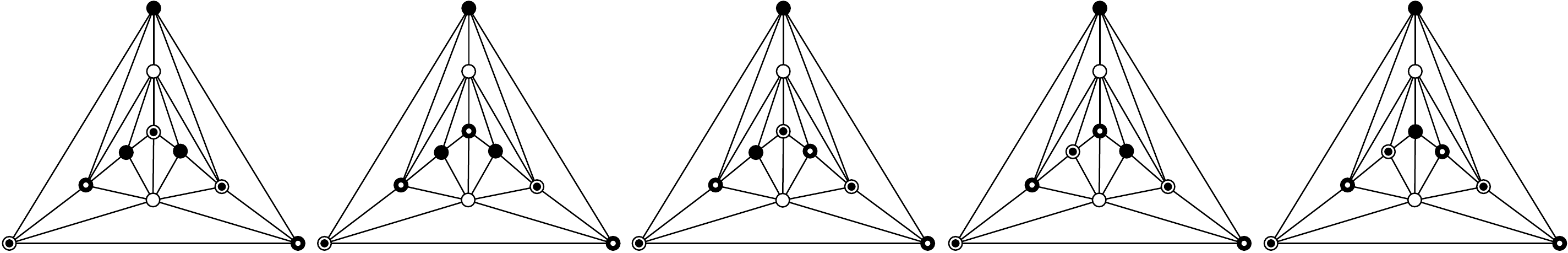}

    \end{center}

 5.7 Degree sequence is 4444455666, and it has 14 kinds of different colorings.
   \begin{center}
        \includegraphics [width=380pt]{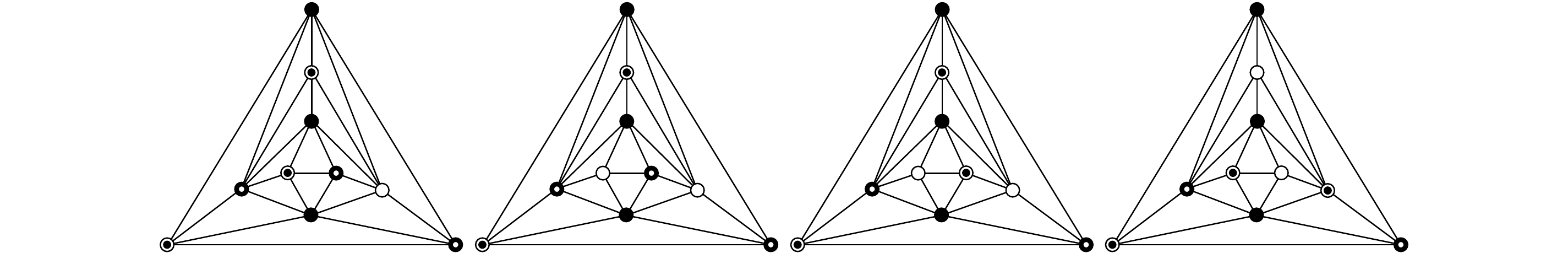}

         \vspace{2mm}
        \includegraphics [width=380pt]{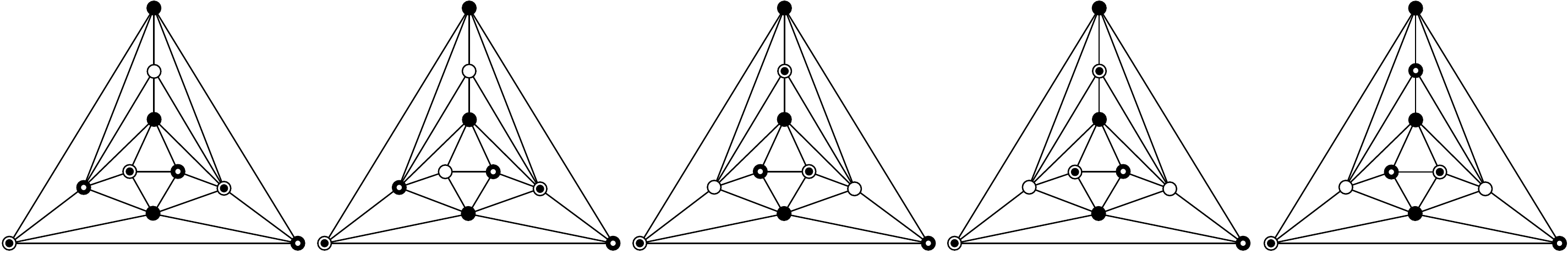}

          \vspace{2mm}
      \includegraphics [width=380pt]{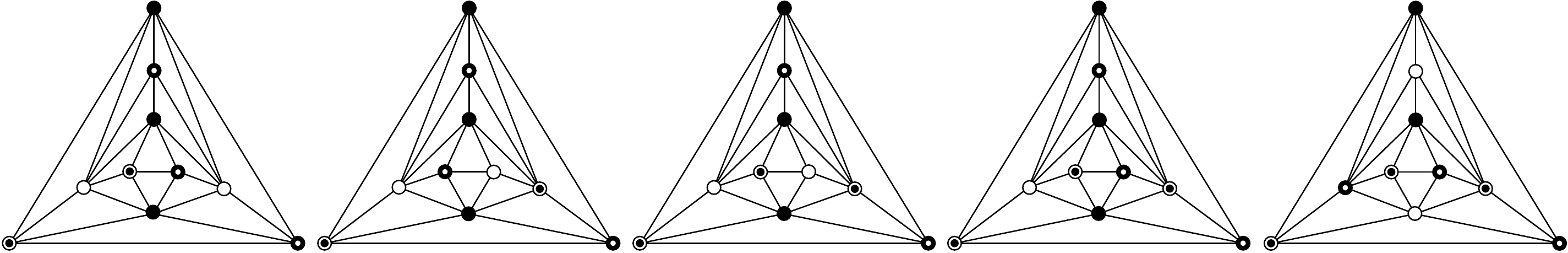}
    \end{center}

 5.8 Degree sequence is 4444445577, and it has 11 kinds of different colorings.

    \begin{center}
        \includegraphics [width=380pt]{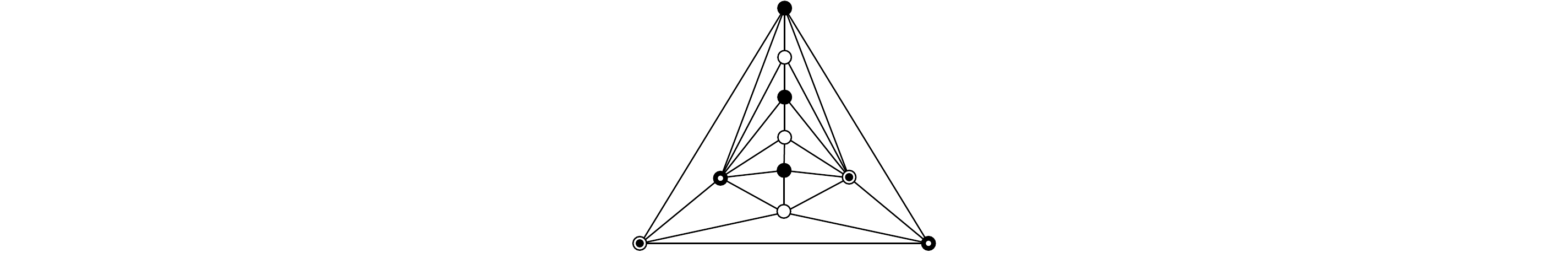}

          \vspace{2mm}
        \includegraphics [width=380pt]{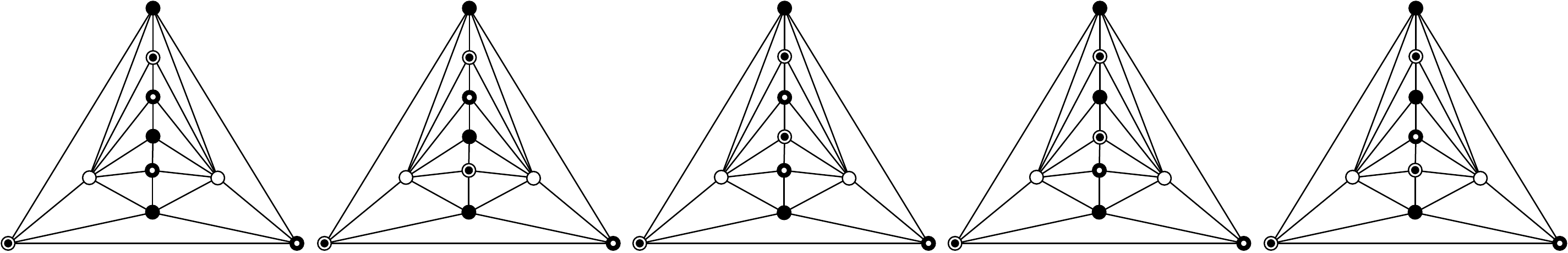}

         \vspace{2mm}
      \includegraphics [width=380pt]{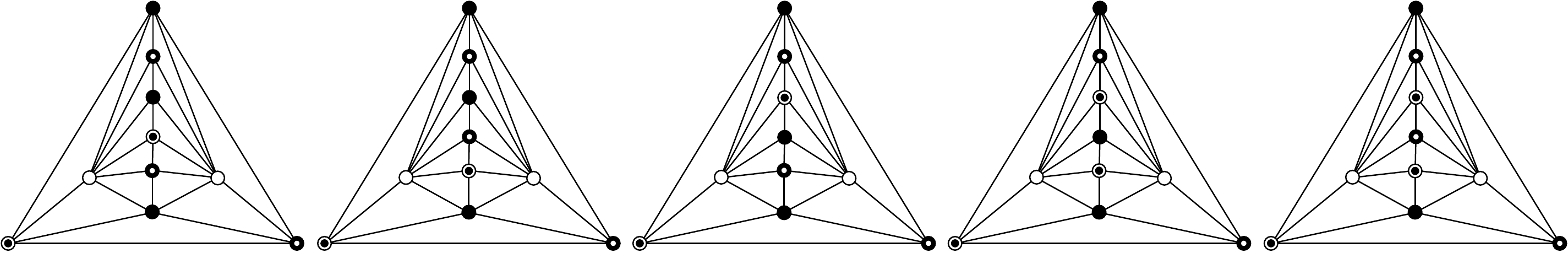}
    \end{center}

 5.9 Degree sequence is 4444445577, and it is a divisible graph.
     \begin{center}
        \includegraphics [width=380pt]{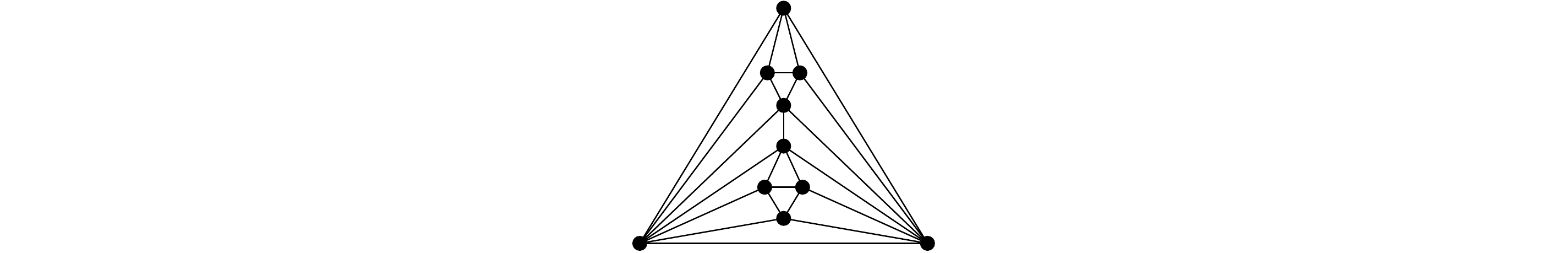}

    \end{center}

 5.10 Degree sequence is  4444446666, and it is uniquely 3-colorable.
     \begin{center}
        \includegraphics [width=380pt]{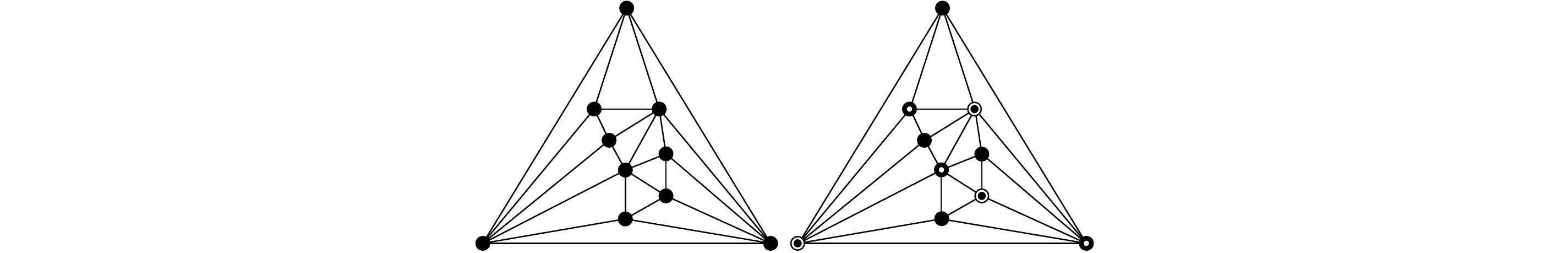}

    \end{center}

 5.11 Degree sequence is 4444445667, and it is divisible.
    \begin{center}
        \includegraphics [width=380pt]{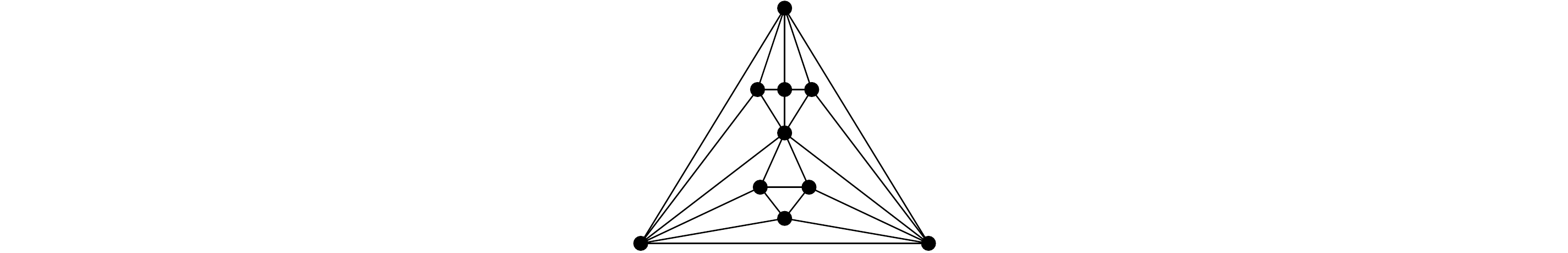}

    \end{center}

 5.12 Degree sequence is 4444444488, and it is 3-colorable.
    \begin{center}
        \includegraphics [width=380pt]{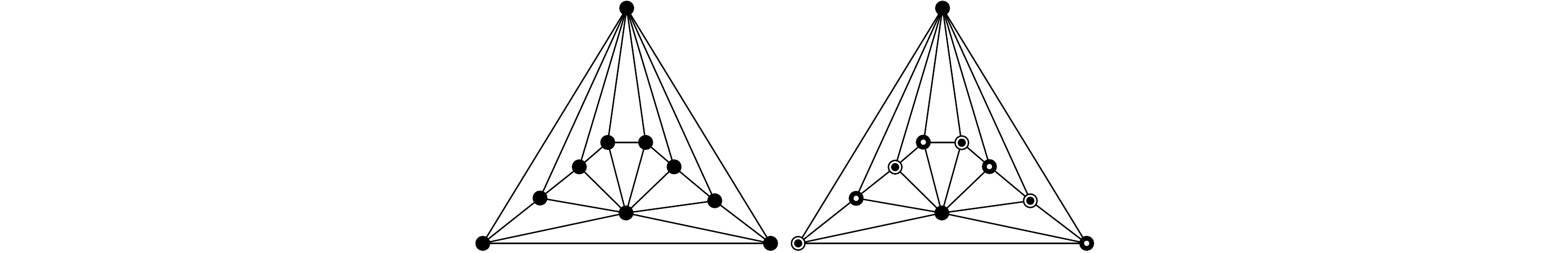}

    \end{center}

 6. There are 34 maximal planar graphs of order 11 whose minimal
 degree is 4.

 6.1 Degree sequence is 44555555556, and it has 8 kinds of different colorings.
    \begin{center}
        \includegraphics [width=380pt]{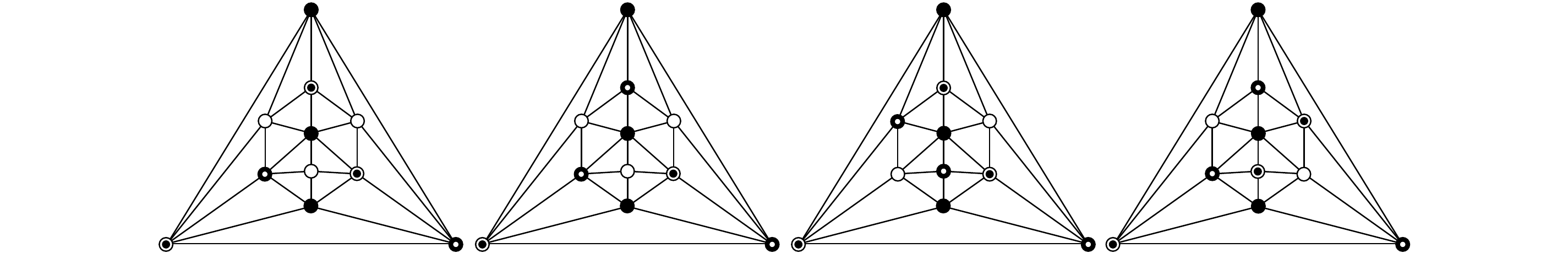}

        \vspace{2mm}
        \includegraphics [width=380pt]{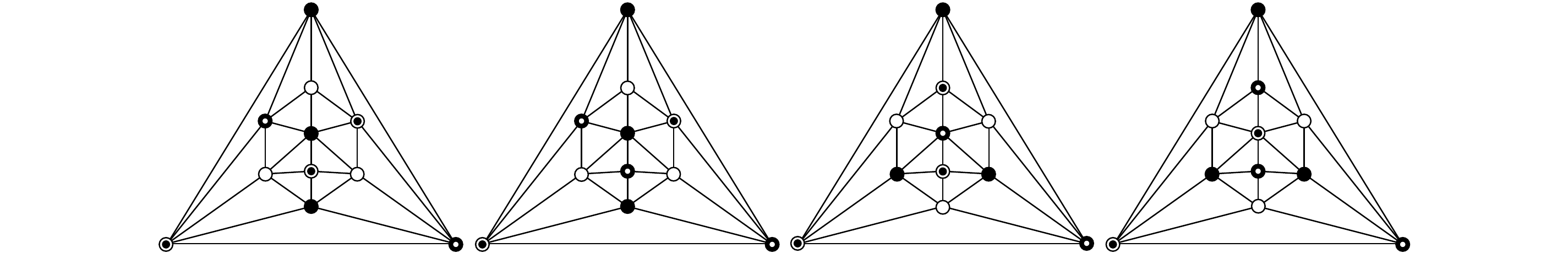}
   \end{center}

 6.2 Degree sequence is 44455555566, and it has 7 kinds of different colorings.
    \begin{center}
        \includegraphics [width=380pt]{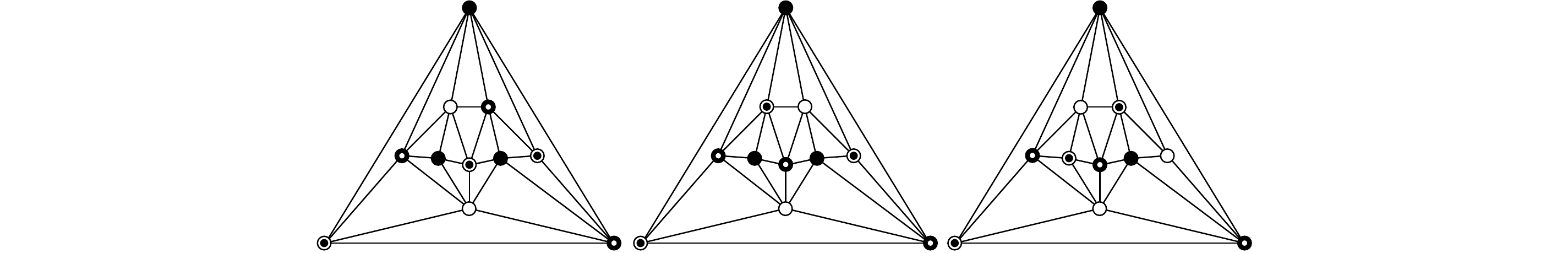}

        \vspace{2mm}
        \includegraphics [width=380pt]{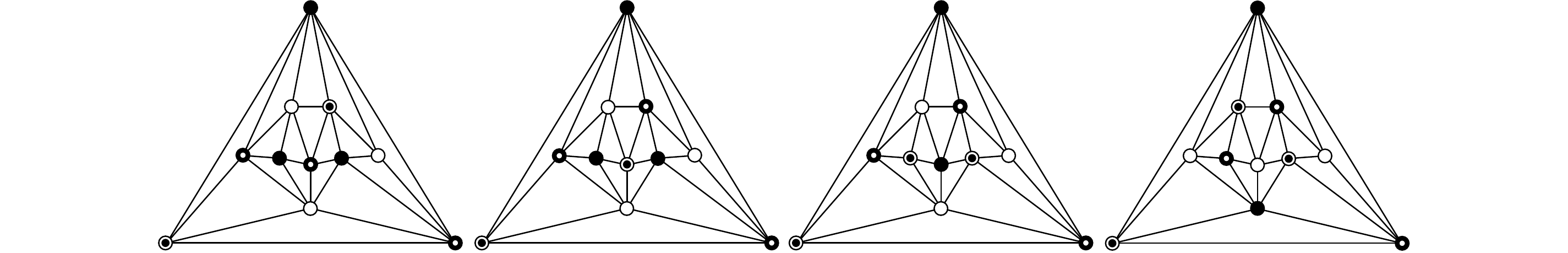}
   \end{center}

6.3 Degree sequence is 44455555566, and it has 10 kinds of different colorings.
\begin{center}
        \includegraphics [width=380pt]{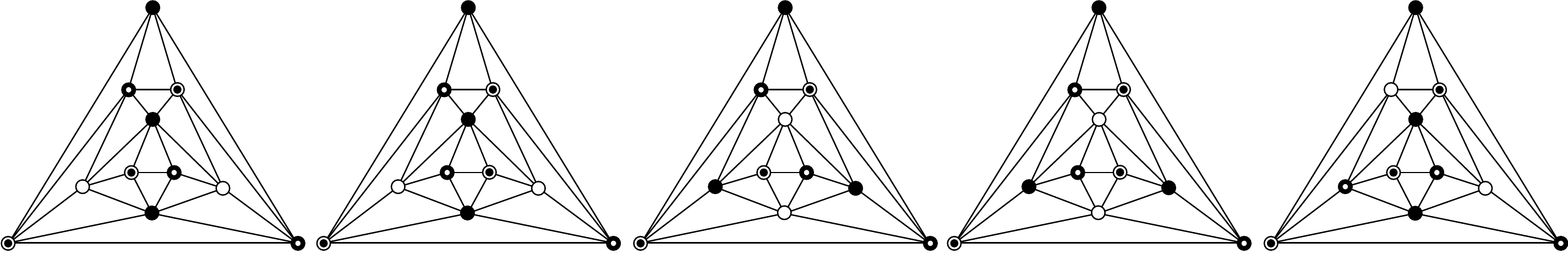}

        \vspace{2mm}
        \includegraphics [width=380pt]{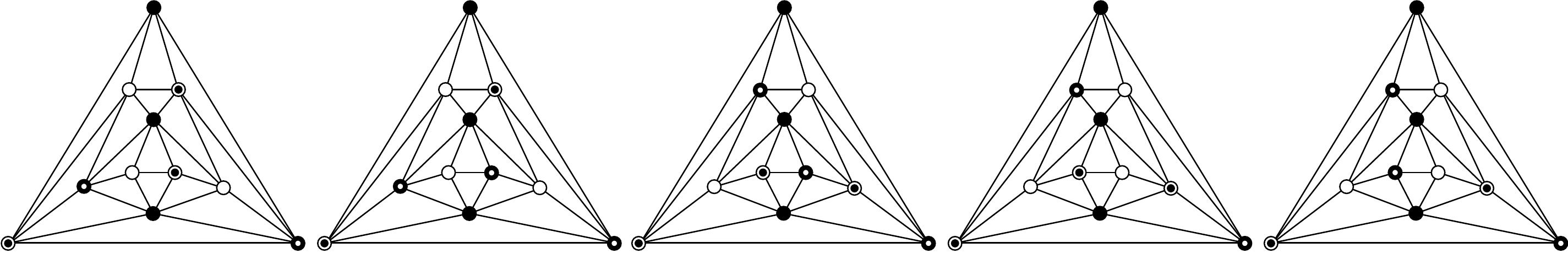}
   \end{center}

6.4 Degree sequence is 44445555666, and it has 13 kinds of different colorings.
\begin{center}
        \includegraphics [width=380pt]{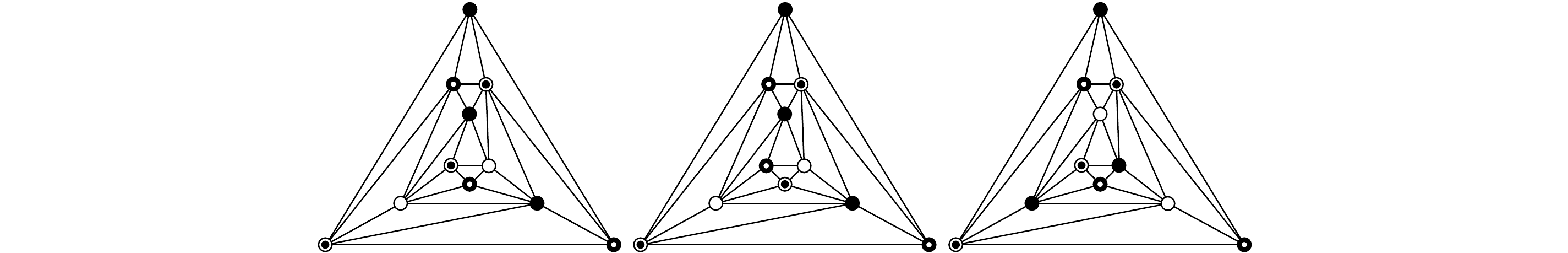}

        \vspace{2mm}
        \includegraphics [width=380pt]{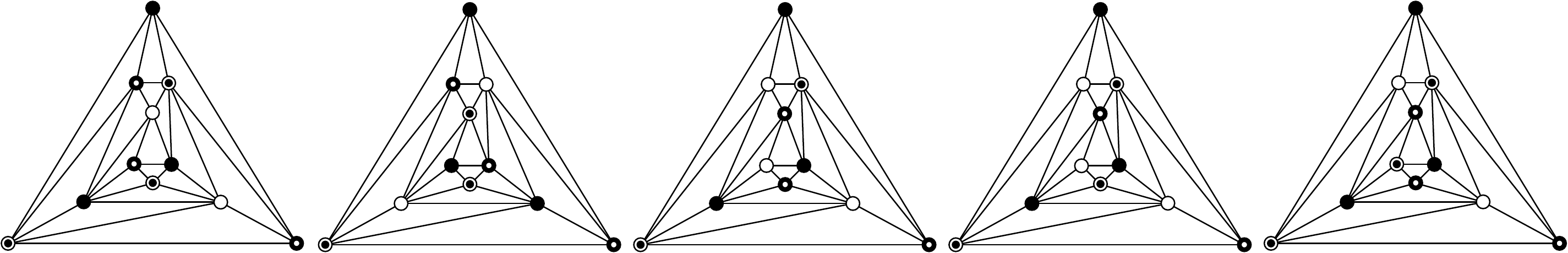}

        \vspace{2mm}
        \includegraphics [width=380pt]{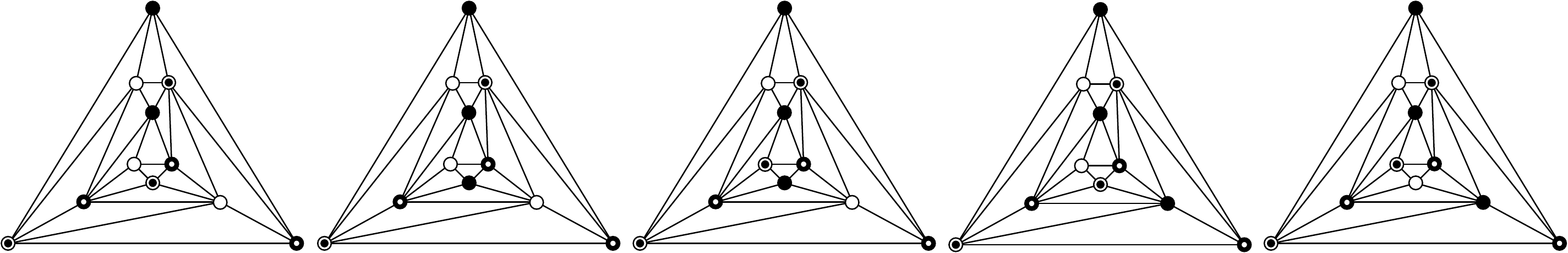}
   \end{center}

6.5 Degree sequence is 44445555666, and it has 16 kinds of different colorings.

\begin{center}
        \includegraphics [width=380pt]{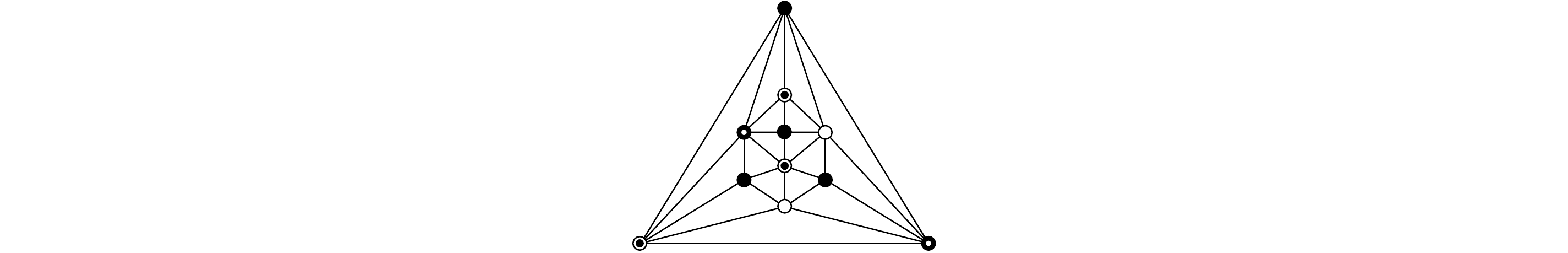}

        \vspace{2mm}
        \includegraphics [width=380pt]{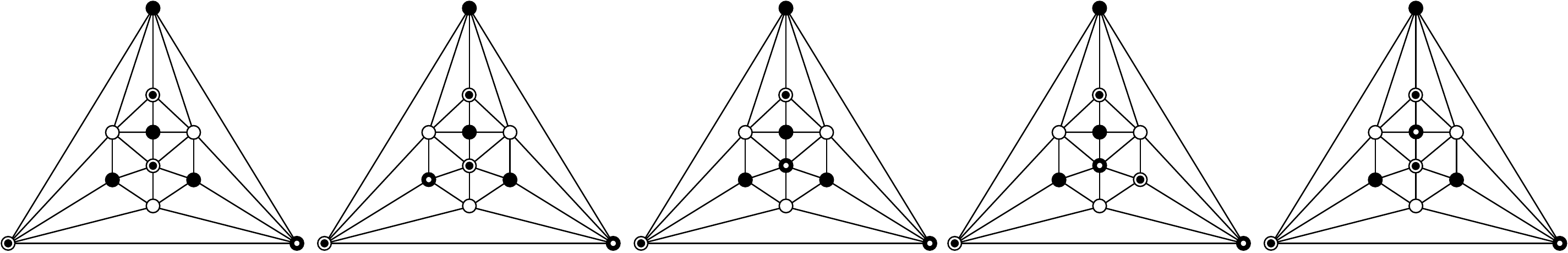}

        \vspace{2mm}
        \includegraphics [width=380pt]{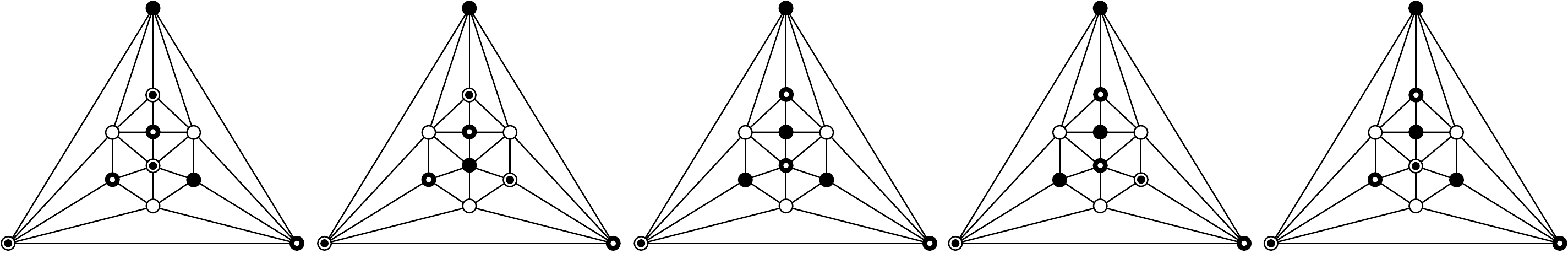}

         \vspace{2mm}
        \includegraphics [width=380pt]{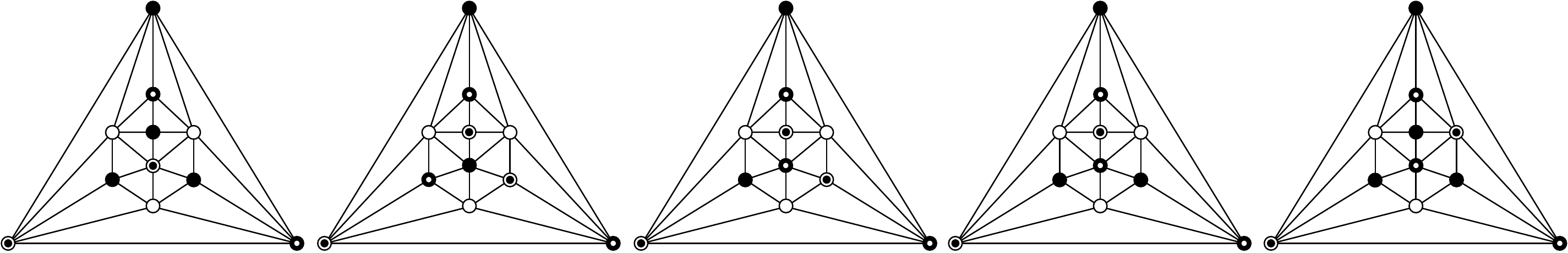}
   \end{center}

6.6 Degree sequence is 44445555666, and it has 9 kinds of different colorings.

\begin{center}
        \includegraphics [width=380pt]{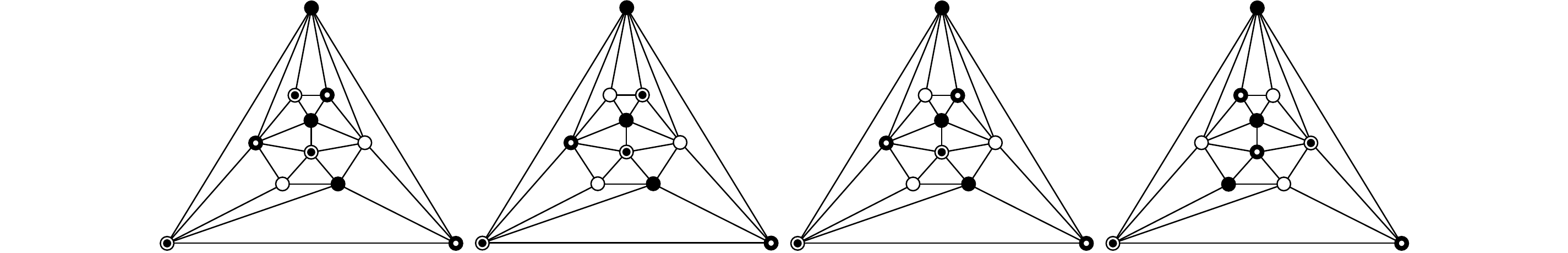}

        \vspace{2mm}
        \includegraphics [width=380pt]{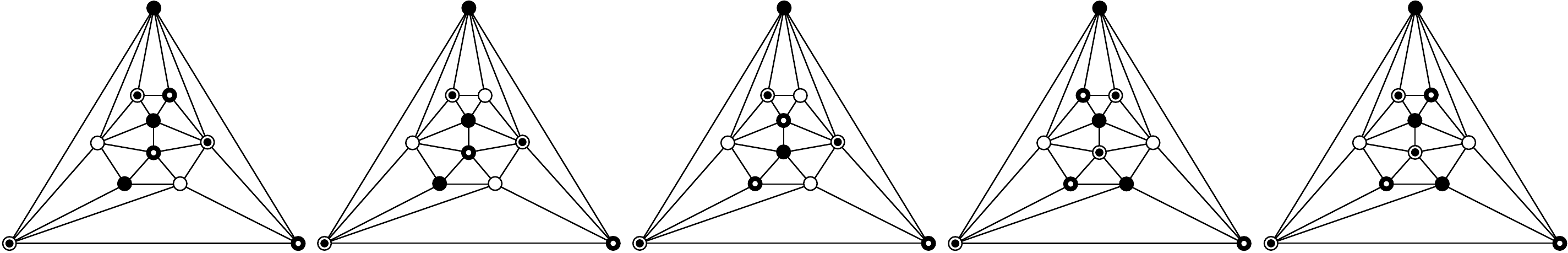}
   \end{center}

6.7 Degree sequence is 44445555666, and it has 10 kinds of different colorings.

\begin{center}
        \includegraphics [width=380pt]{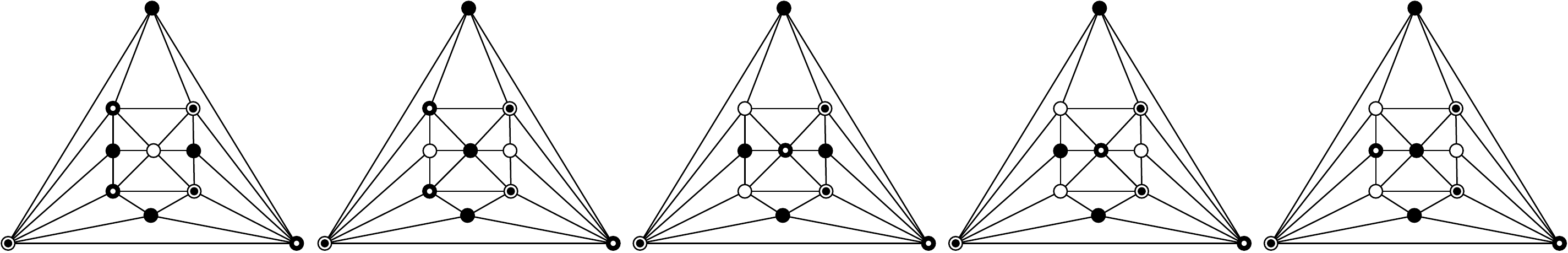}

        \vspace{2mm}
        \includegraphics [width=380pt]{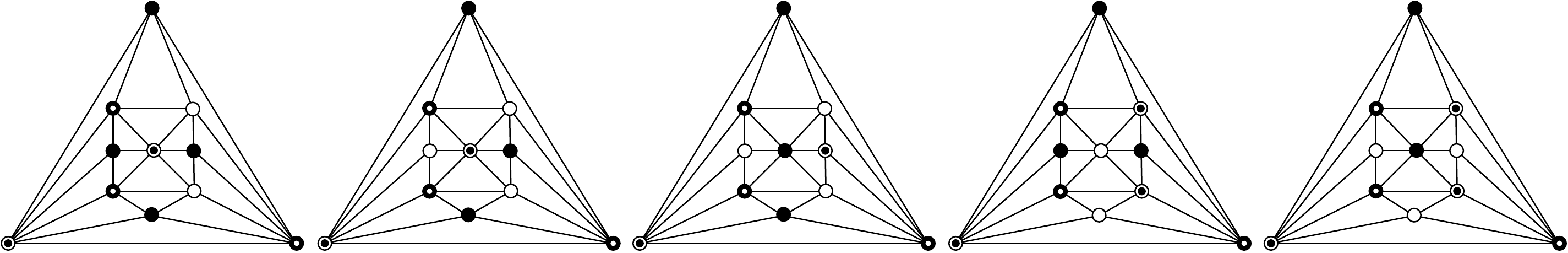}
   \end{center}

6.8 Degree sequence is 44445555567, and it has 11 kinds of different colorings.

\begin{center}
        \includegraphics [width=380pt]{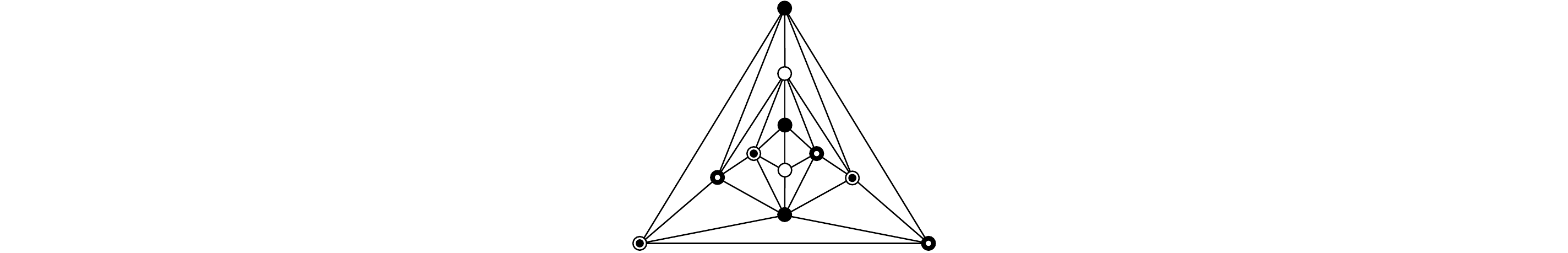}

        \vspace{2mm}
        \includegraphics [width=380pt]{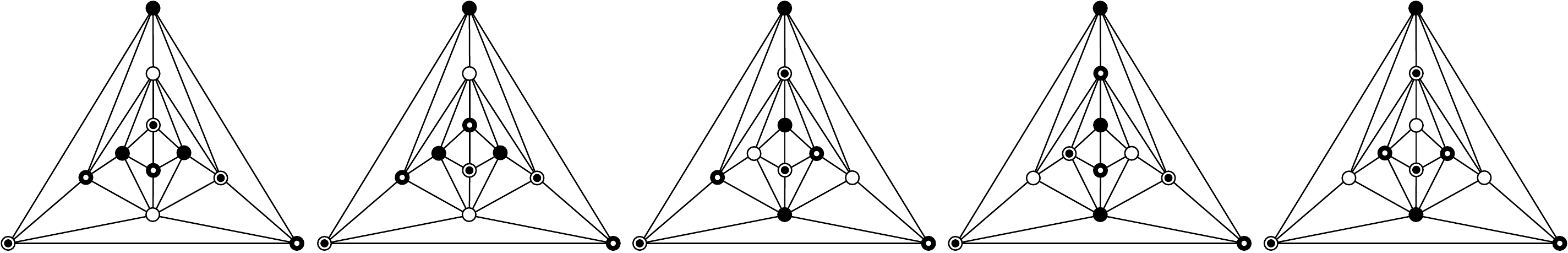}

        \vspace{2mm}
        \includegraphics [width=380pt]{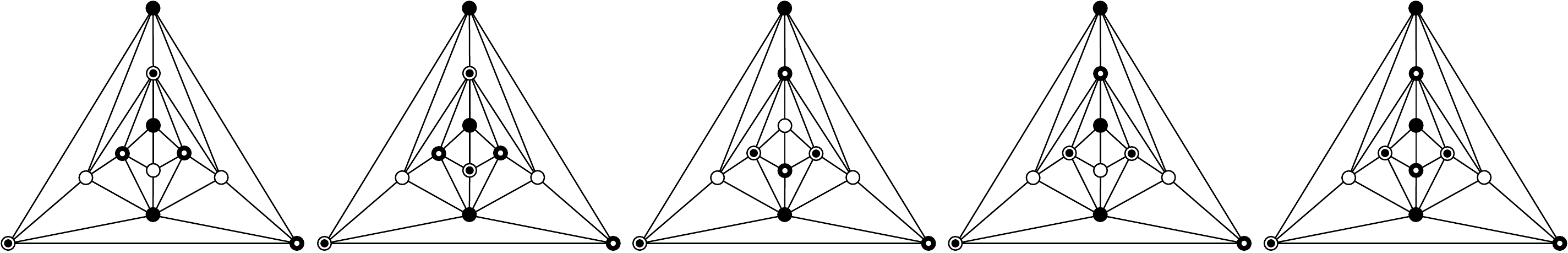}
   \end{center}

6.9 Degree sequence is 44445555567, and it has 4 kinds of different colorings.
\begin{center}
        \includegraphics [width=380pt]{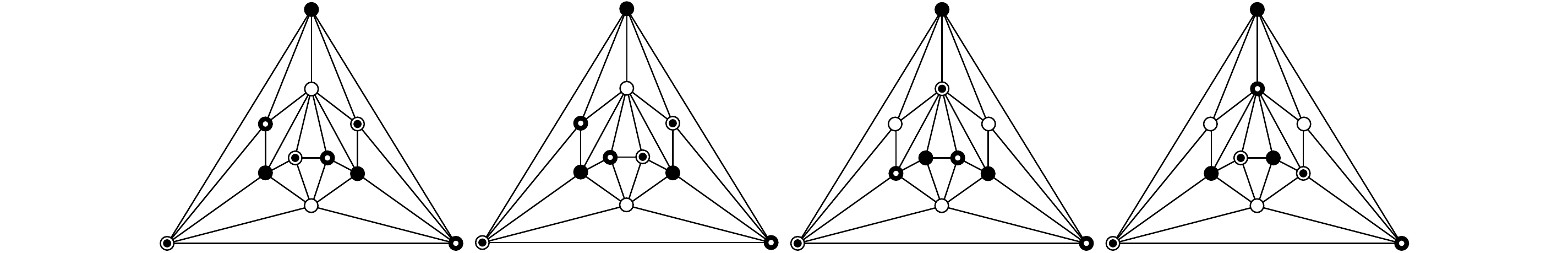}

   \end{center}

6.10 Degree sequence is 44445555567, and it has 13 kinds of different colorings.

\begin{center}
        \includegraphics [width=380pt]{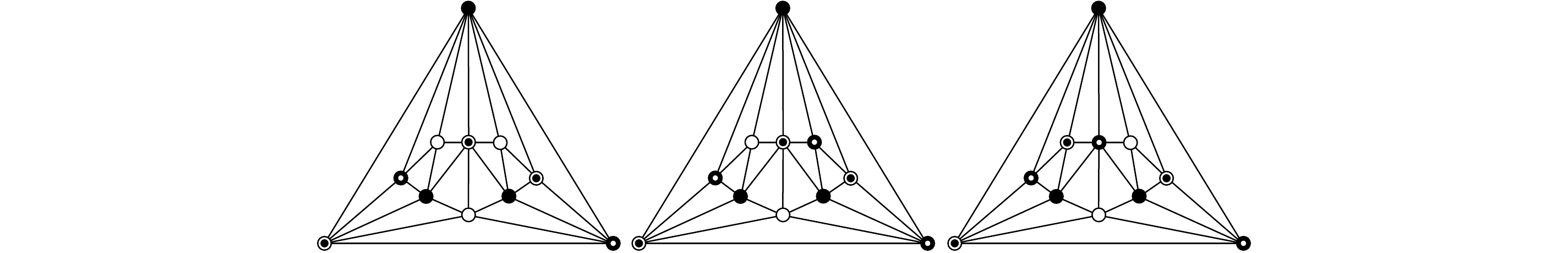}

        \vspace{2mm}
        \includegraphics [width=380pt]{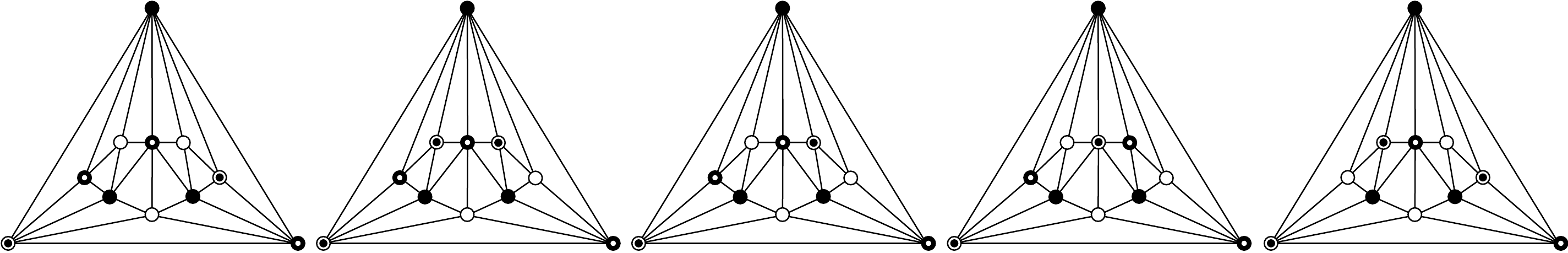}

        \vspace{2mm}
        \includegraphics [width=380pt]{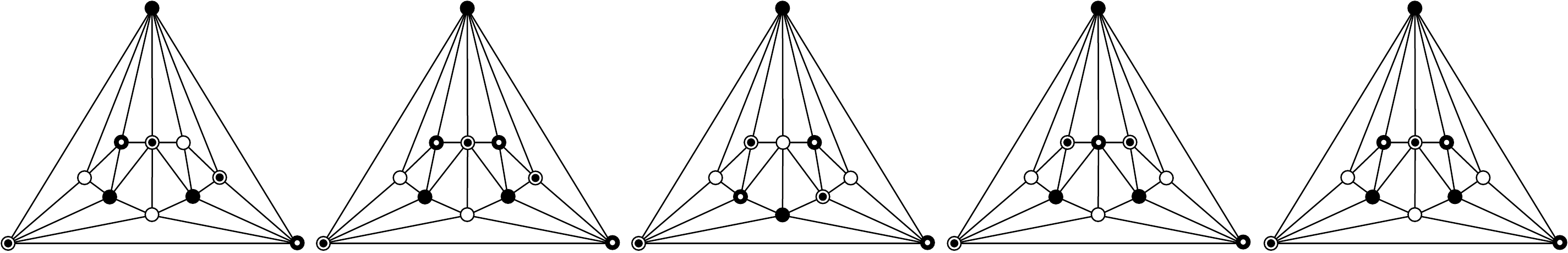}
   \end{center}

6.11 Degree sequence is 44444555667, and it has 22 kinds of different colorings.

\begin{center}
        \includegraphics [width=380pt]{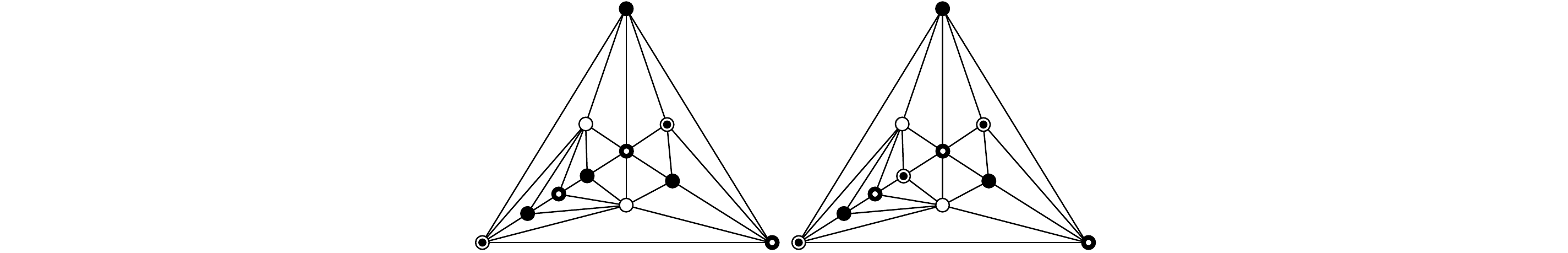}

        \vspace{2mm}
        \includegraphics [width=380pt]{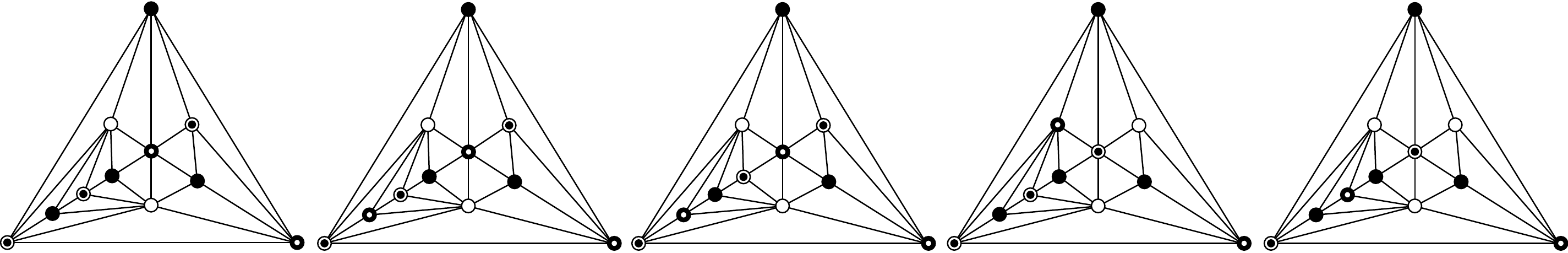}

        \vspace{2mm}
        \includegraphics [width=380pt]{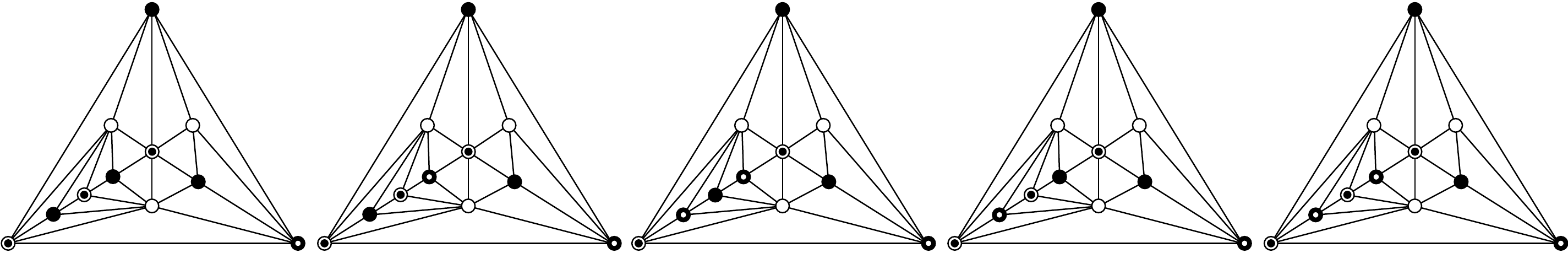}

         \vspace{2mm}
        \includegraphics [width=380pt]{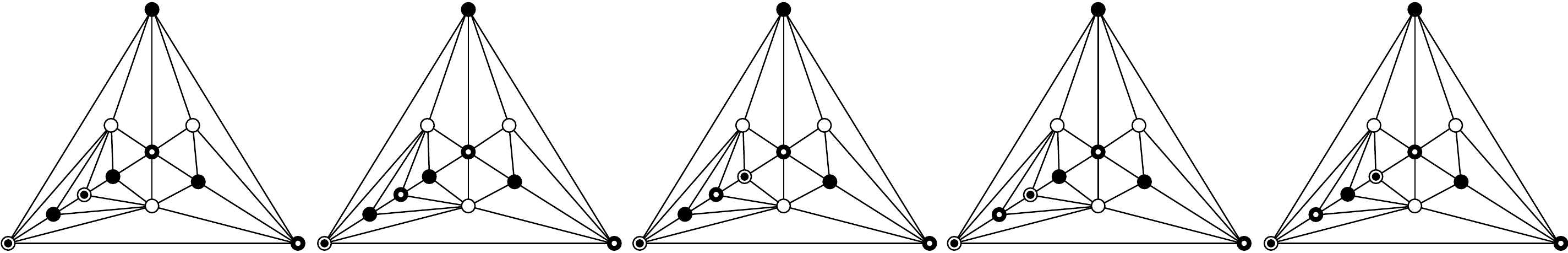}

         \vspace{2mm}
        \includegraphics [width=380pt]{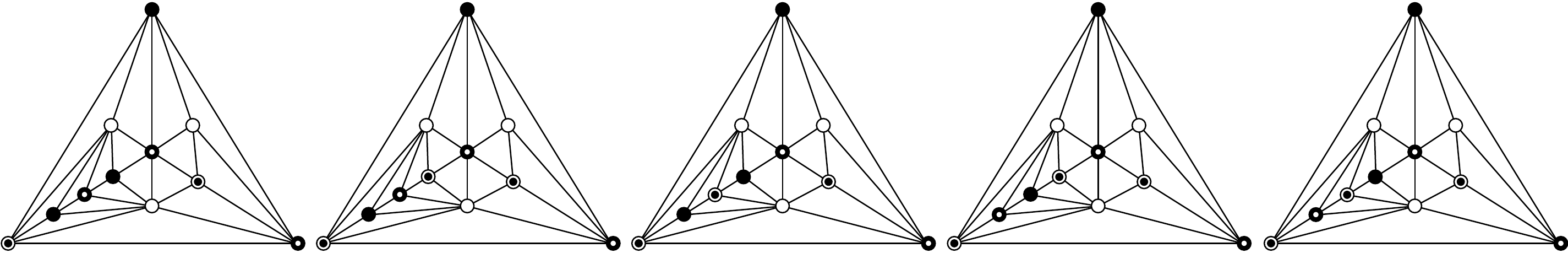}
   \end{center}

6.12 Degree sequence is  44444556666, and it has 29 kinds of different colorings.

\begin{center}
        \includegraphics [width=380pt]{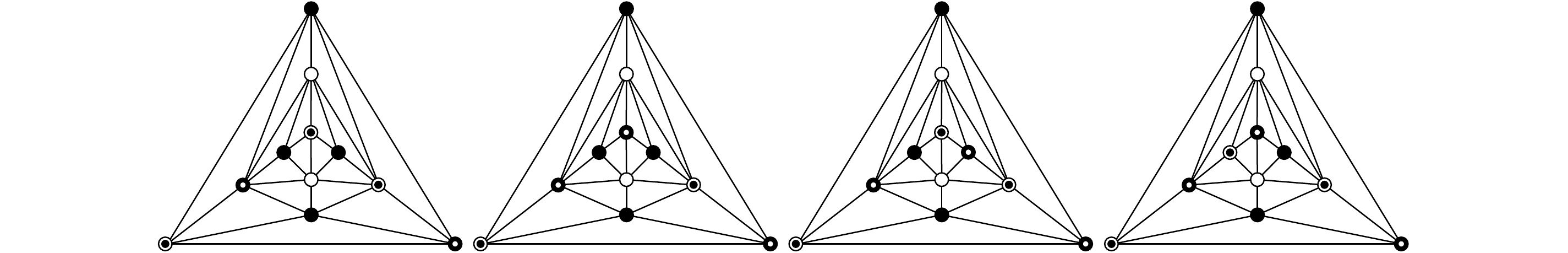}

        \vspace{2mm}
        \includegraphics [width=380pt]{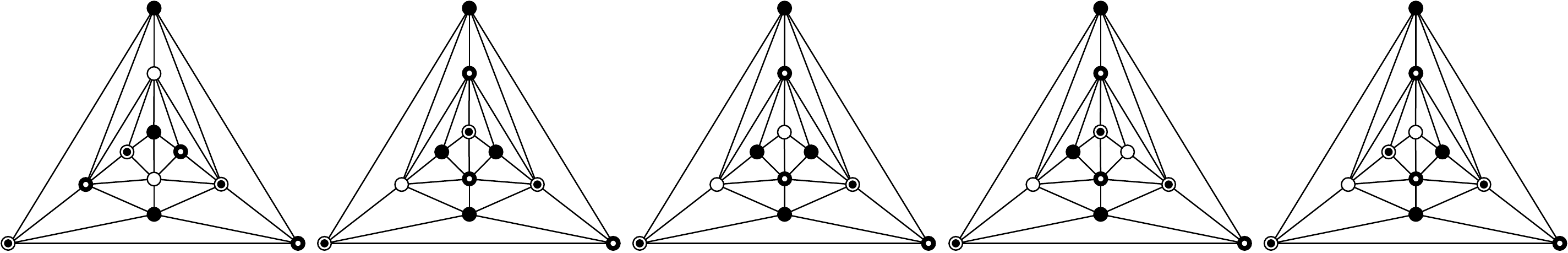}

        \vspace{2mm}
        \includegraphics [width=380pt]{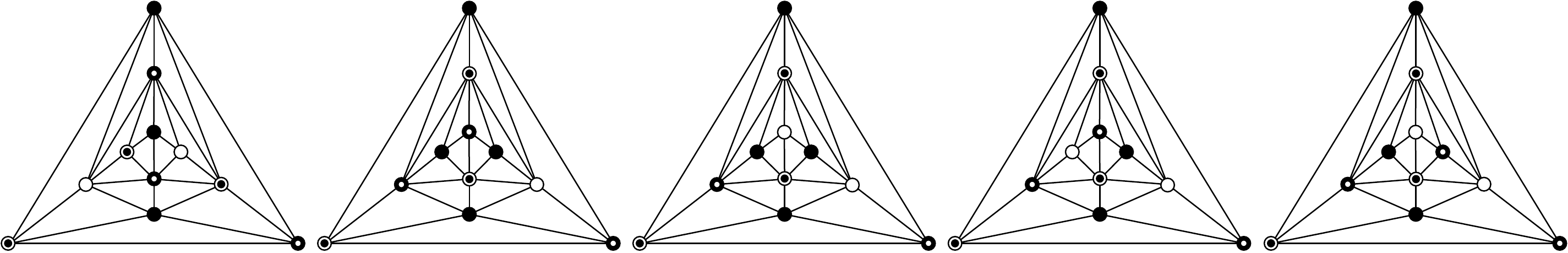}

         \vspace{2mm}
        \includegraphics [width=380pt]{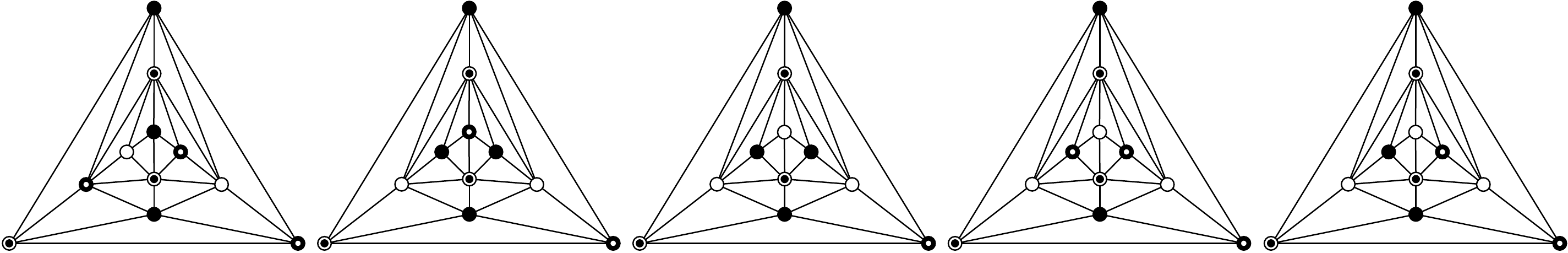}

         \vspace{2mm}
        \includegraphics [width=380pt]{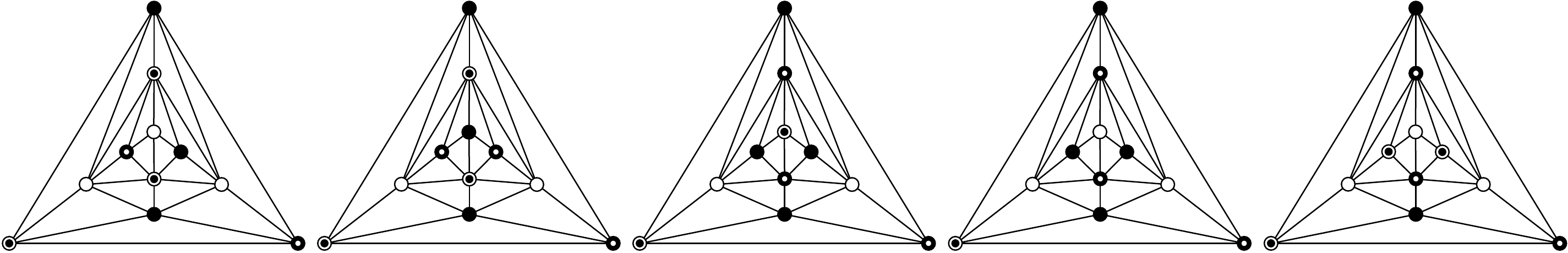}

         \vspace{2mm}
        \includegraphics [width=380pt]{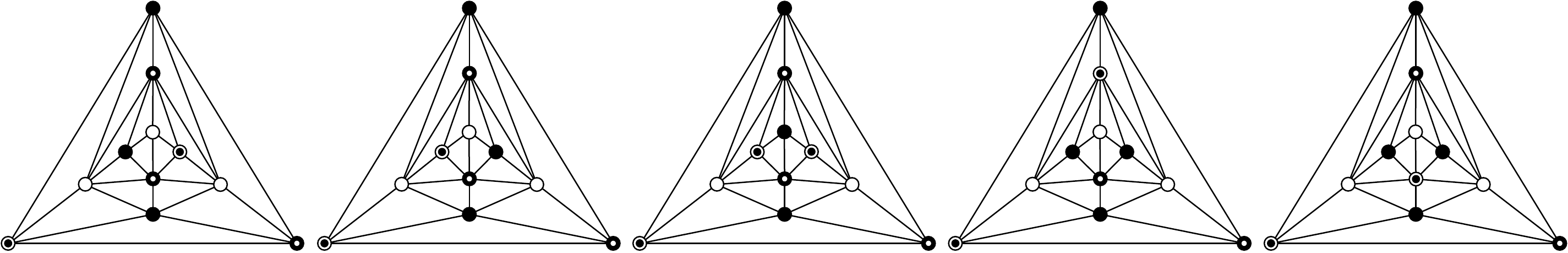}
   \end{center}

6.13 Degree sequence is 44444555577, and it has 10 kinds of different colorings.

\begin{center}
        \includegraphics [width=380pt]{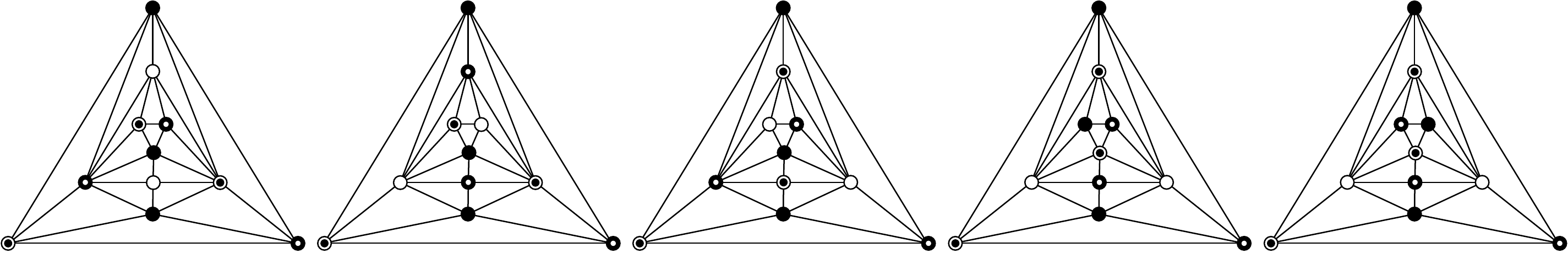}

        \vspace{2mm}
        \includegraphics [width=380pt]{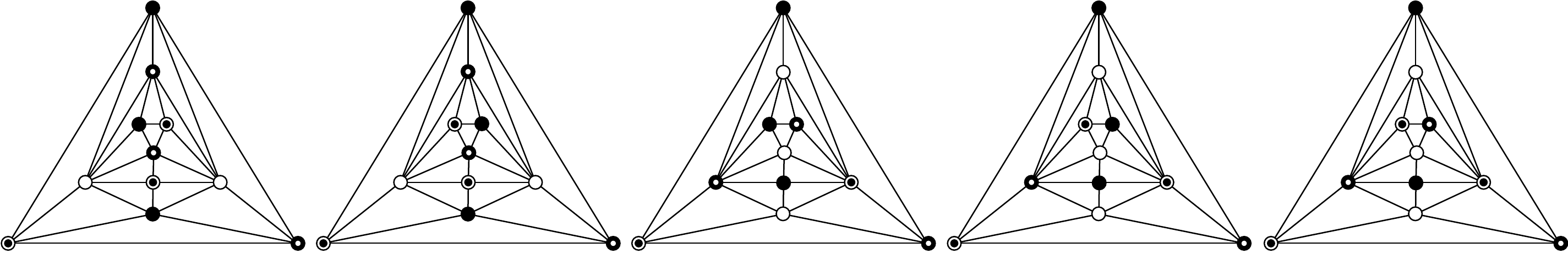}
   \end{center}

6.14 Degree sequence is 44444555577, and it has 6 kinds of different colorings.
\begin{center}
        \includegraphics [width=380pt]{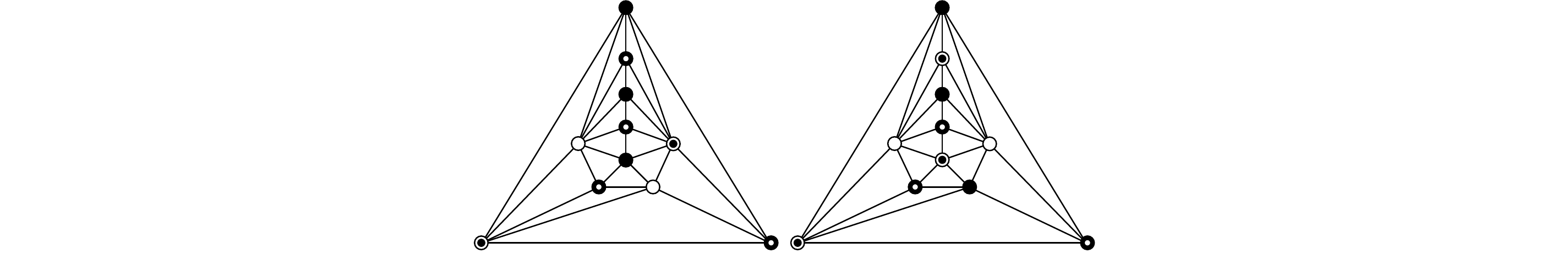}

        \vspace{2mm}
        \includegraphics [width=380pt]{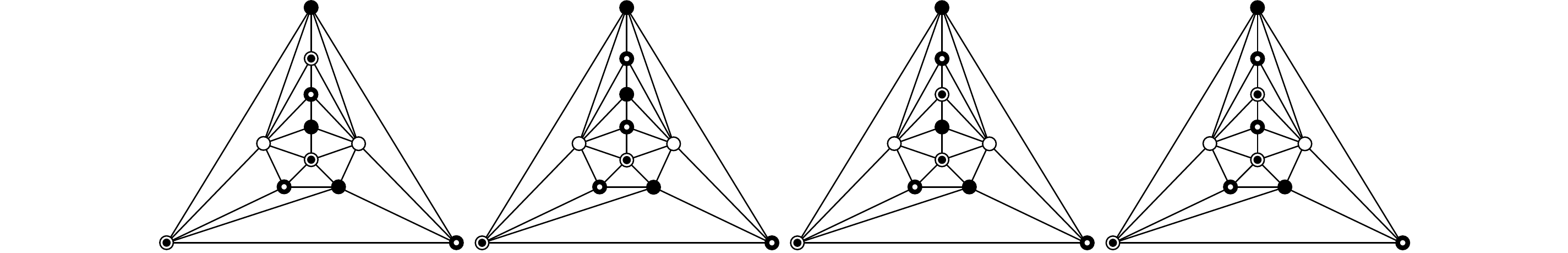}
   \end{center}

6.15 Degree sequence is 44444555667, and it has 8 kinds of different colorings.
\begin{center}
        \includegraphics [width=380pt]{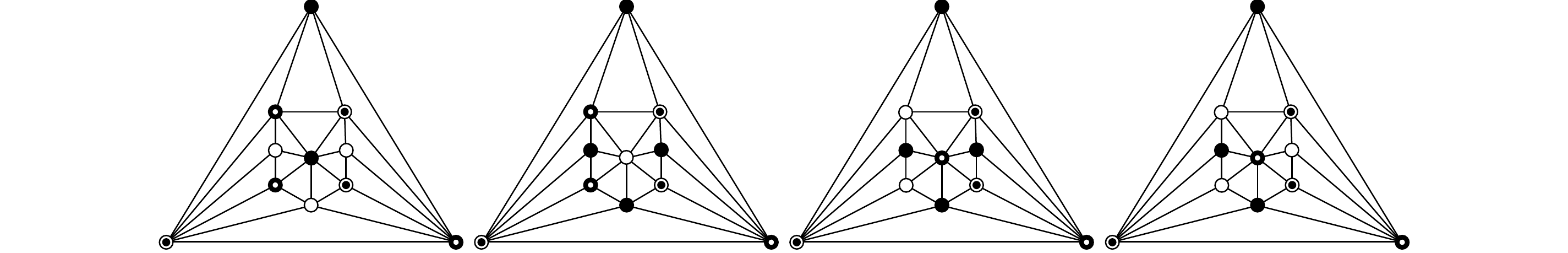}

        \vspace{2mm}
        \includegraphics [width=380pt]{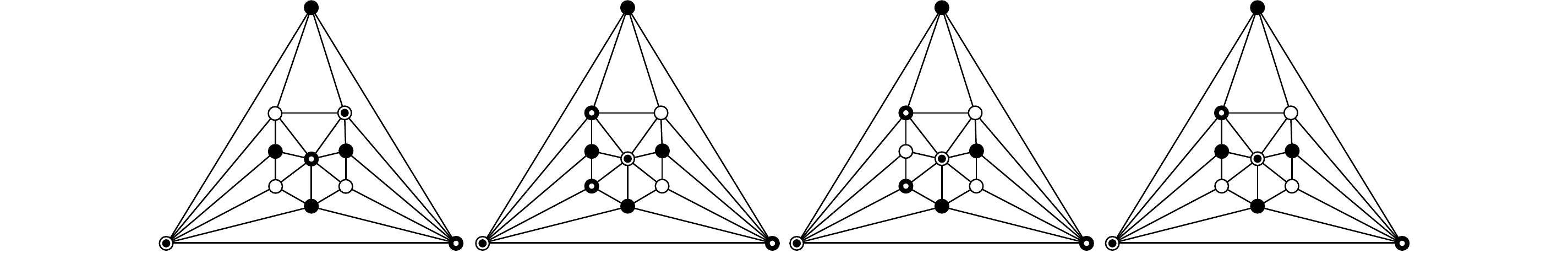}
   \end{center}

6.16 Degree sequence is 44444555667, and it has 11 kinds of different colorings.

\begin{center}
        \includegraphics [width=380pt]{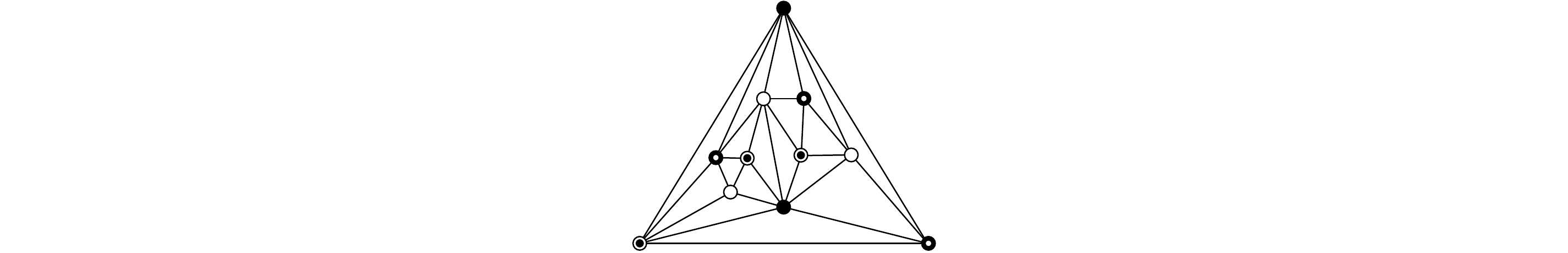}

        \vspace{2mm}
        \includegraphics [width=380pt]{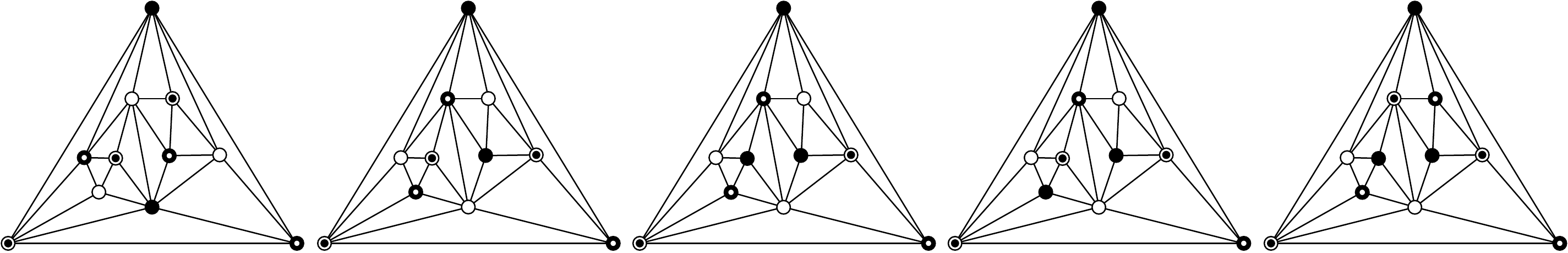}

        \vspace{2mm}
        \includegraphics [width=380pt]{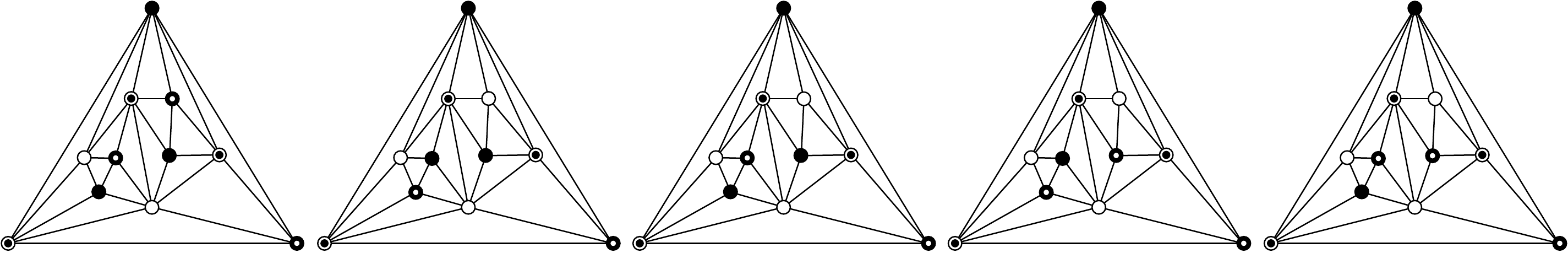}
   \end{center}

6.17 Degree sequence is 44444555667, and it is divisible.

 \begin{center}
        \includegraphics [width=380pt]{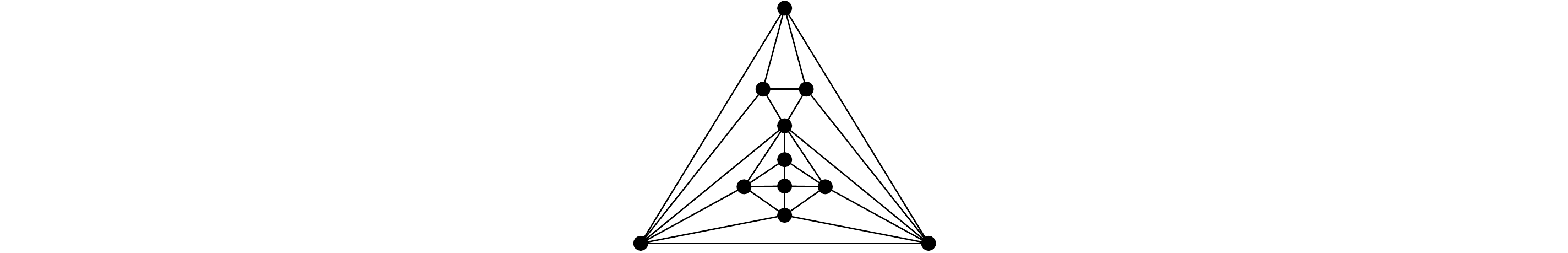}

   \end{center}

6.18 Degree sequence is 44444556666, and it has 17 kinds of different colorings.

  \begin{center}
        \includegraphics [width=380pt]{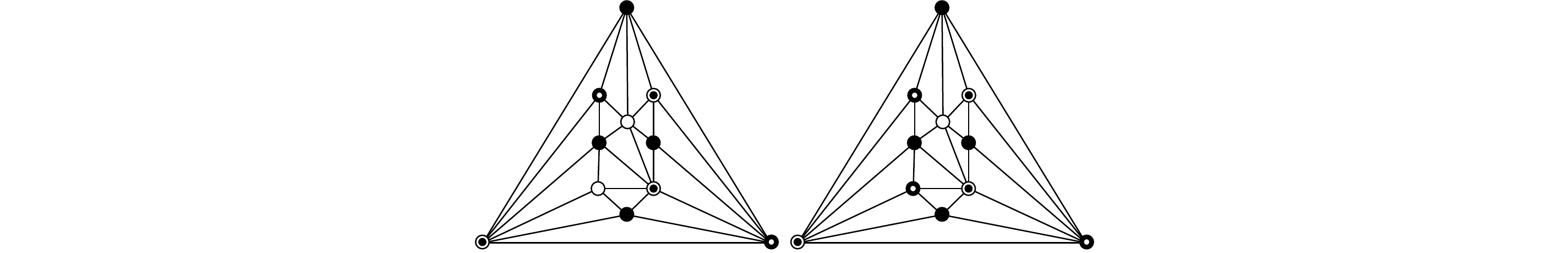}

        \vspace{2mm}
        \includegraphics [width=380pt]{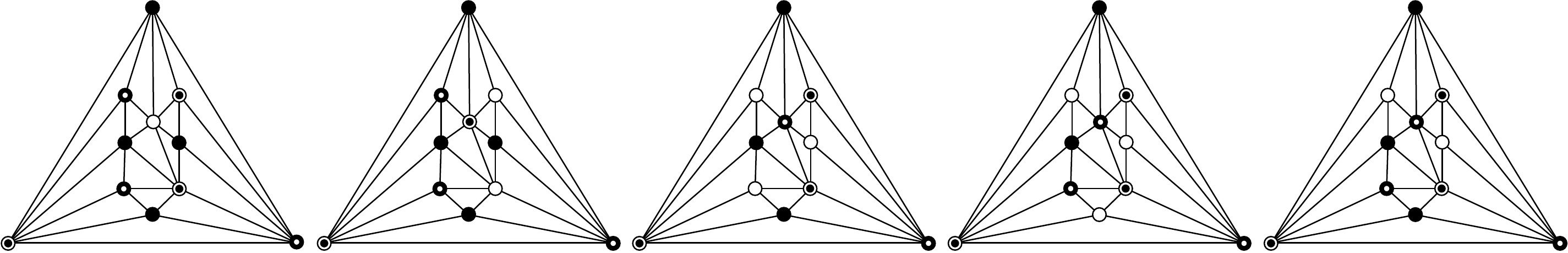}

        \vspace{2mm}
        \includegraphics [width=380pt]{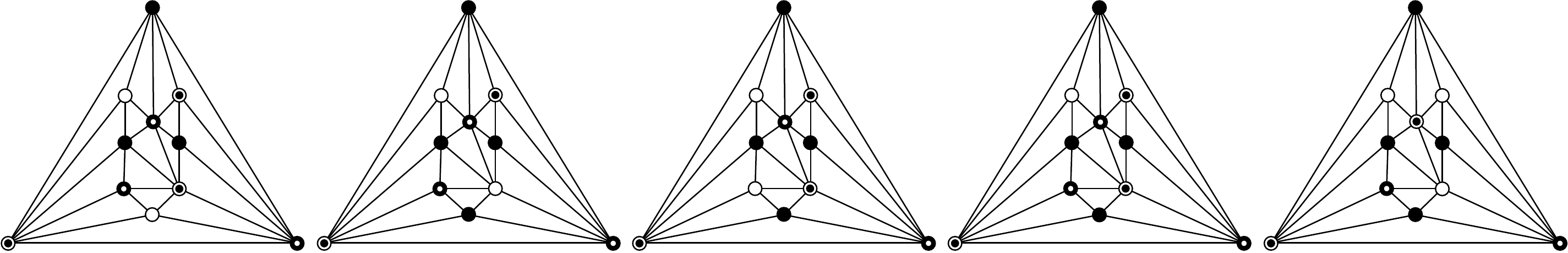}

        \vspace{2mm}
        \includegraphics [width=380pt]{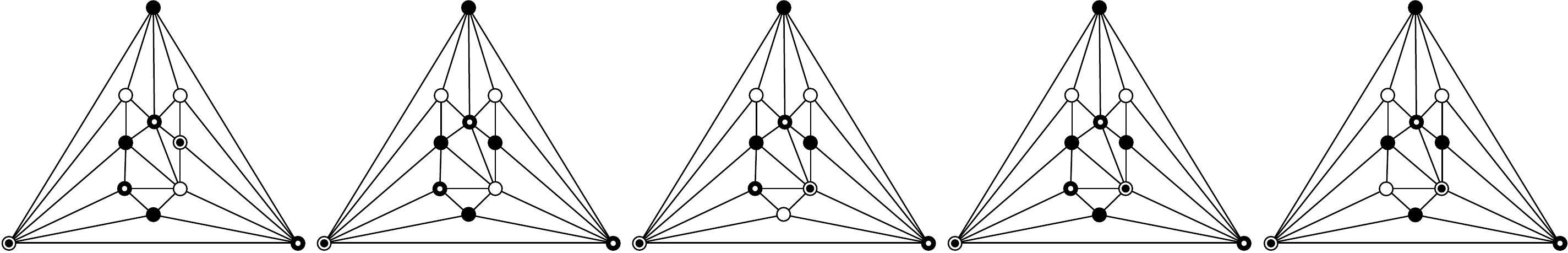}
   \end{center}

6.19 Degree sequence is 44444455677, and it is divisible.

\begin{center}
        \includegraphics [width=380pt]{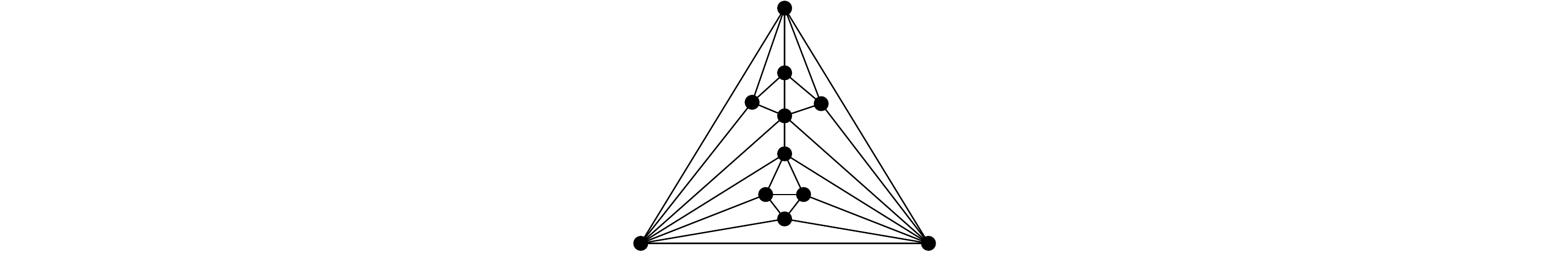}

   \end{center}

6.20 Degree sequence is 44444455578, and it has 14 kinds of different colorings.

    \begin{center}
        \includegraphics [width=380pt]{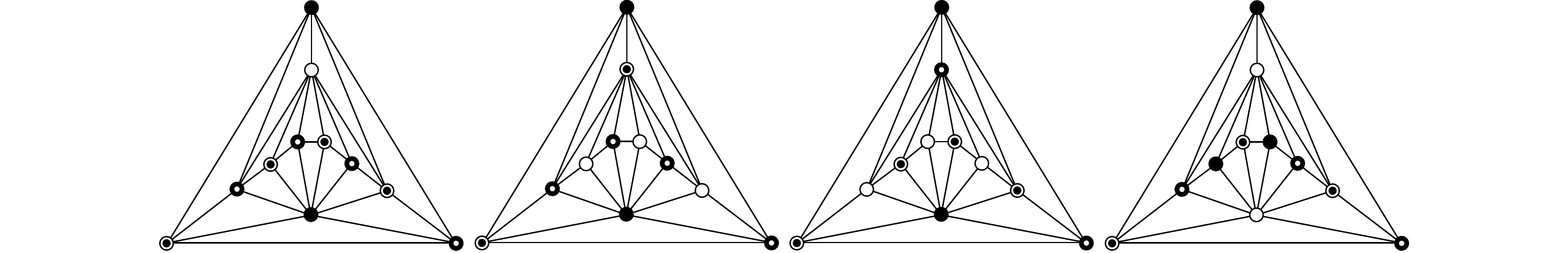}

        \vspace{2mm}
        \includegraphics [width=380pt]{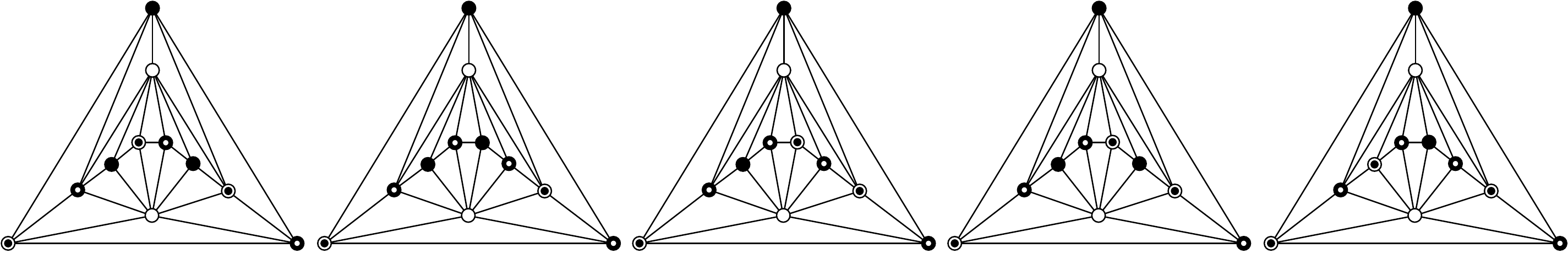}

        \vspace{2mm}
        \includegraphics [width=380pt]{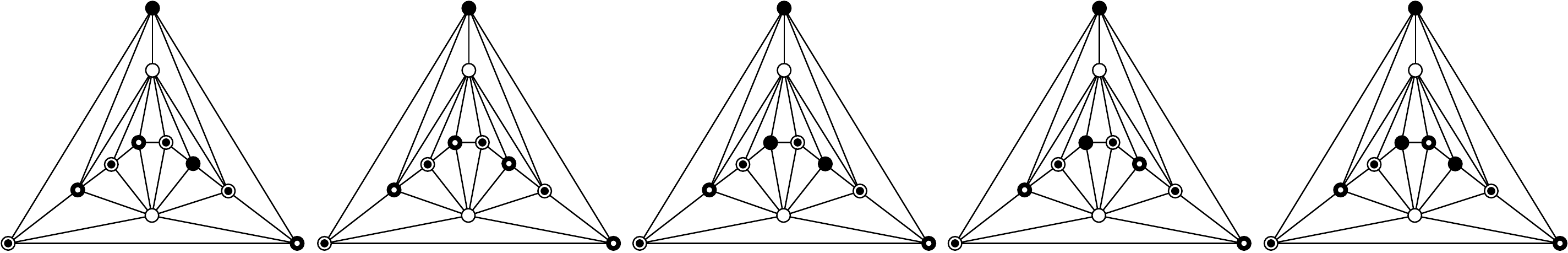}
   \end{center}

6.21 Degree sequence is 44444466666, and it is uniquely 3-colorable.

   \begin{center}
        \includegraphics [width=380pt]{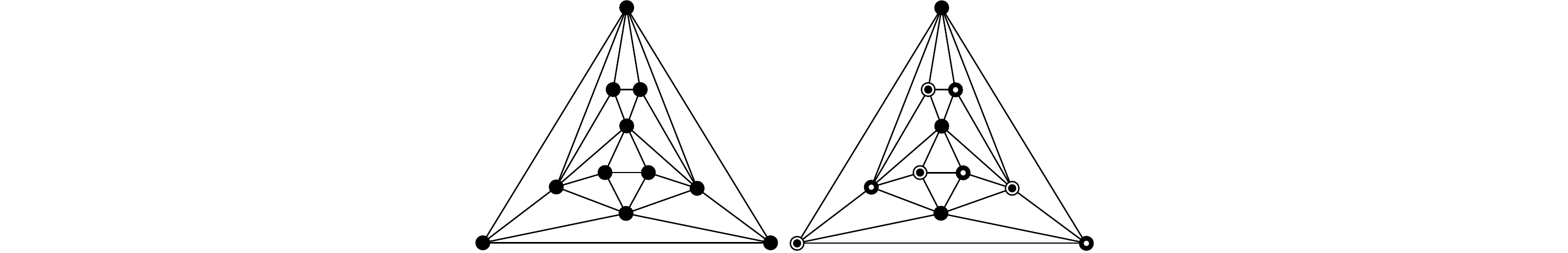}

   \end{center}

6.22 Degree sequence is  44444456667, and it has 21 kinds of different colorings.

\begin{center}
        \includegraphics [width=380pt]{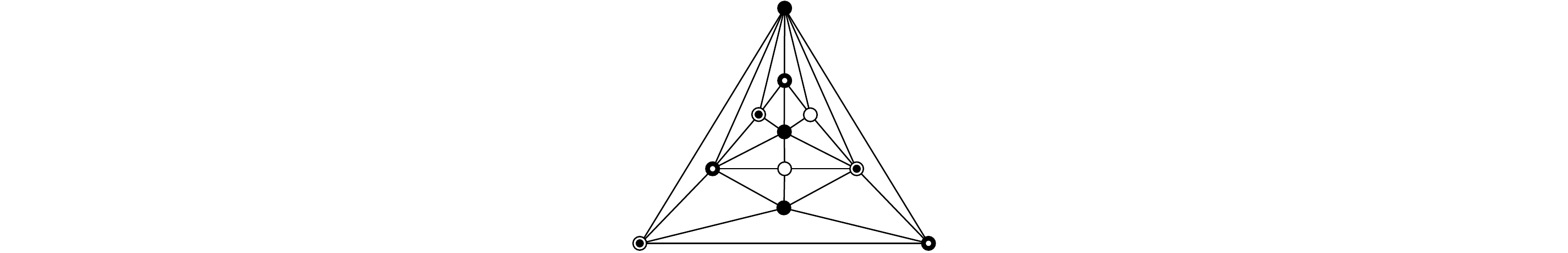}

        \vspace{2mm}
        \includegraphics [width=380pt]{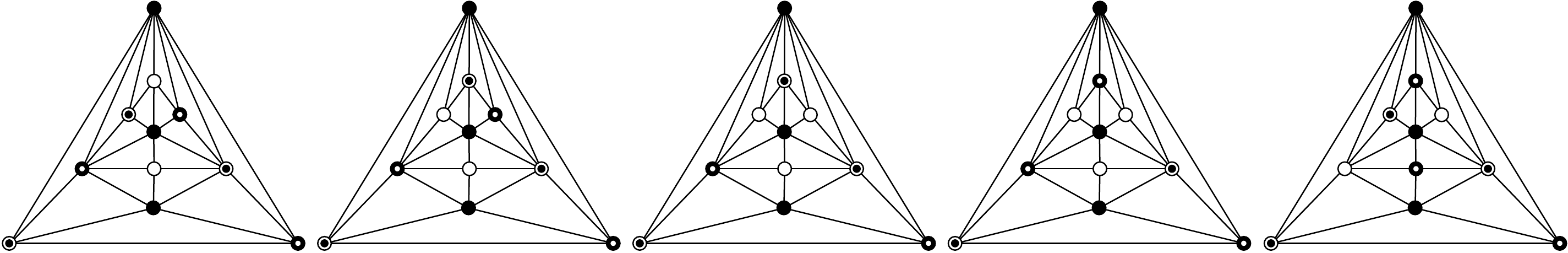}

        \vspace{2mm}
        \includegraphics [width=380pt]{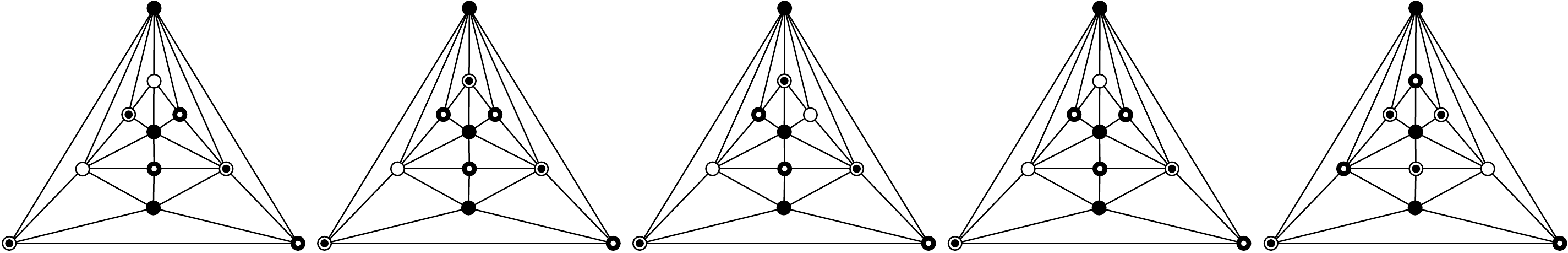}

        \vspace{2mm}
        \includegraphics [width=380pt]{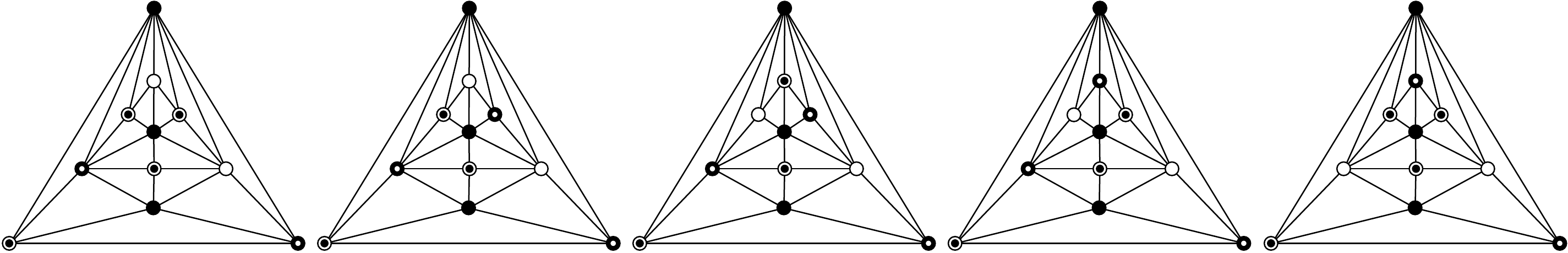}

        \vspace{2mm}
        \includegraphics [width=380pt]{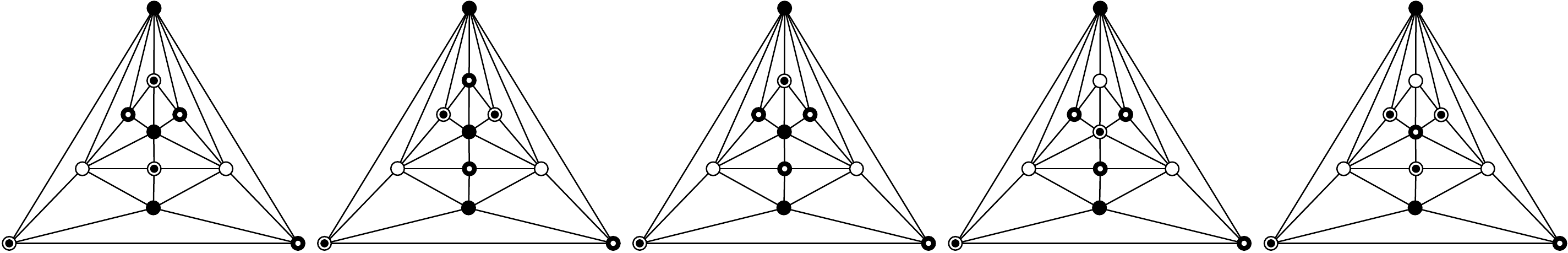}
   \end{center}

6.23 Degree sequence is  44444455668, and it has 13 kinds of different colorings.

\begin{center}
        \includegraphics [width=380pt]{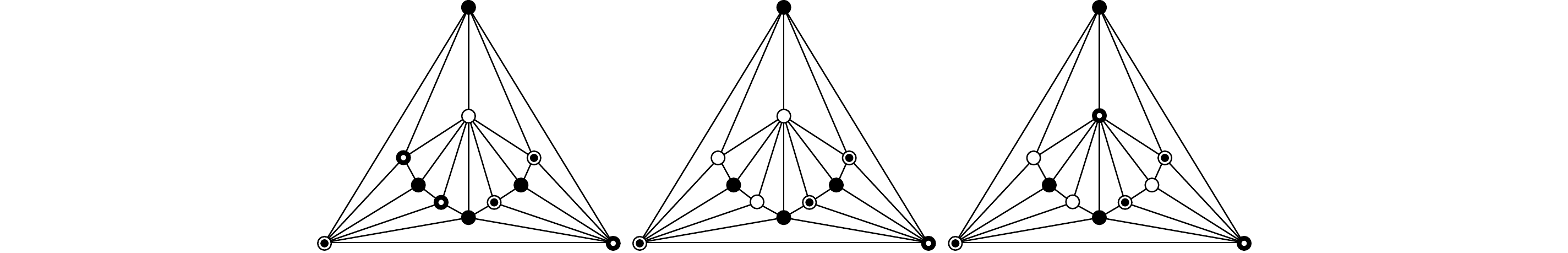}

        \vspace{2mm}
        \includegraphics [width=380pt]{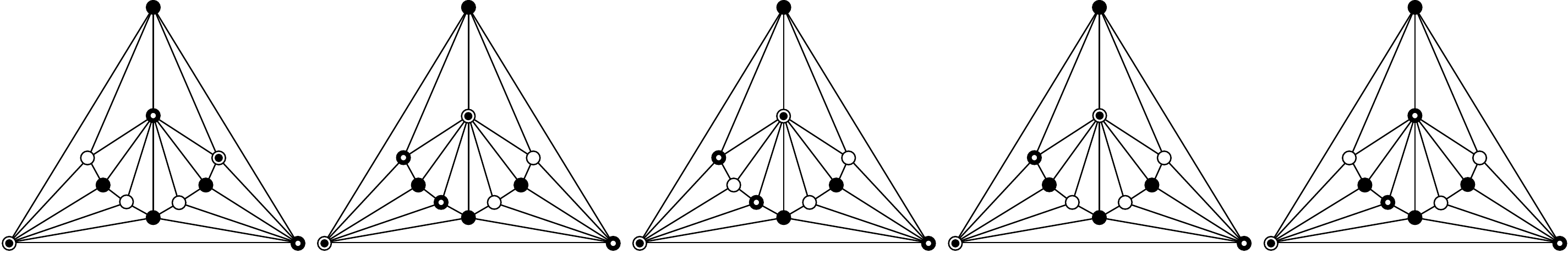}

        \vspace{2mm}
        \includegraphics [width=380pt]{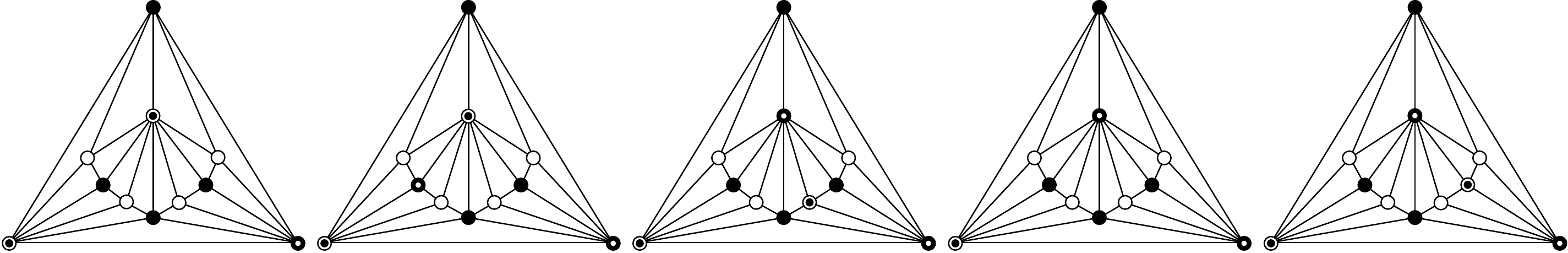}
   \end{center}

6.24 Degree sequence is 44444455578, and it is divisible.

  \begin{center}
        \includegraphics [width=380pt]{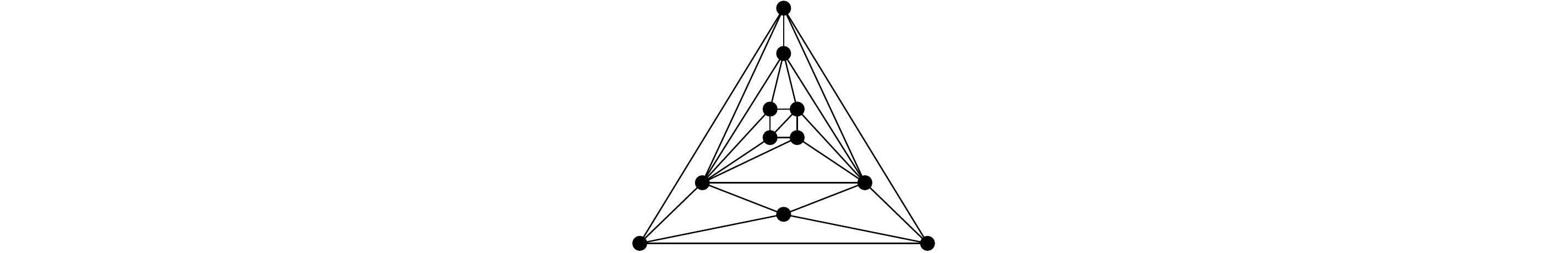}

   \end{center}

6.25 Degree sequence is 44444446668, and it is divisible and uniquely 3-colorable.

\begin{center}
        \includegraphics [width=380pt]{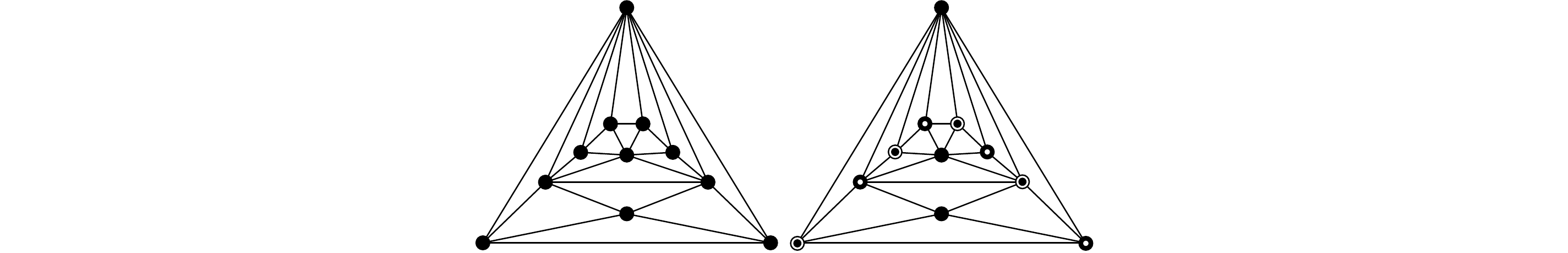}

   \end{center}

6.26 Degree sequence is  44444455668, and it is divisible.

\begin{center}
        \includegraphics [width=380pt]{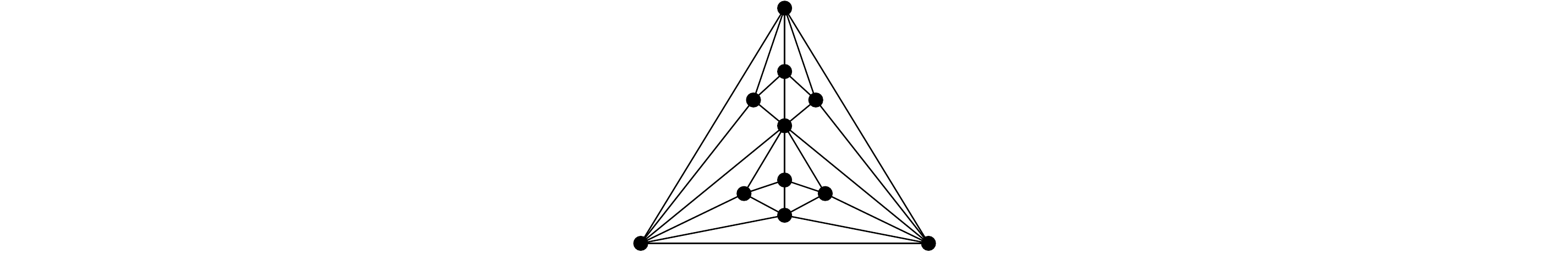}

   \end{center}

6.27 Degree sequence is  44444455668, and it is divisible.

\begin{center}
        \includegraphics [width=380pt]{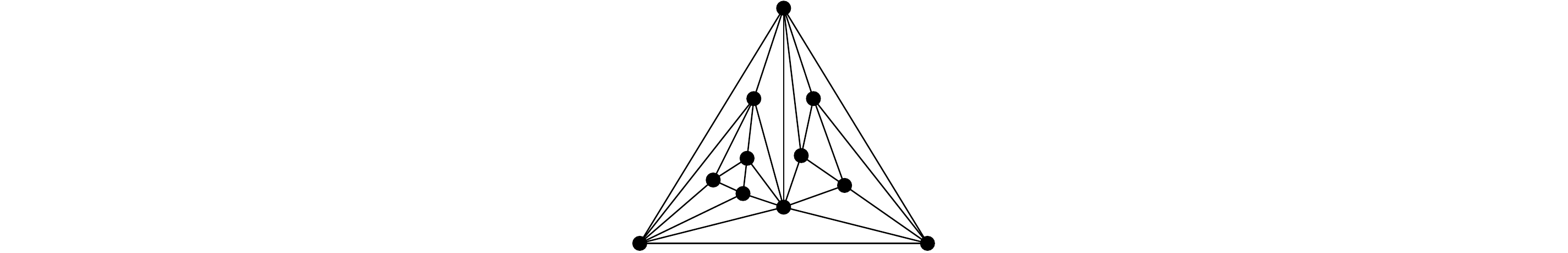}

   \end{center}

6.28 Degree sequence is 44444455677, and it is divisible.

\begin{center}
        \includegraphics [width=380pt]{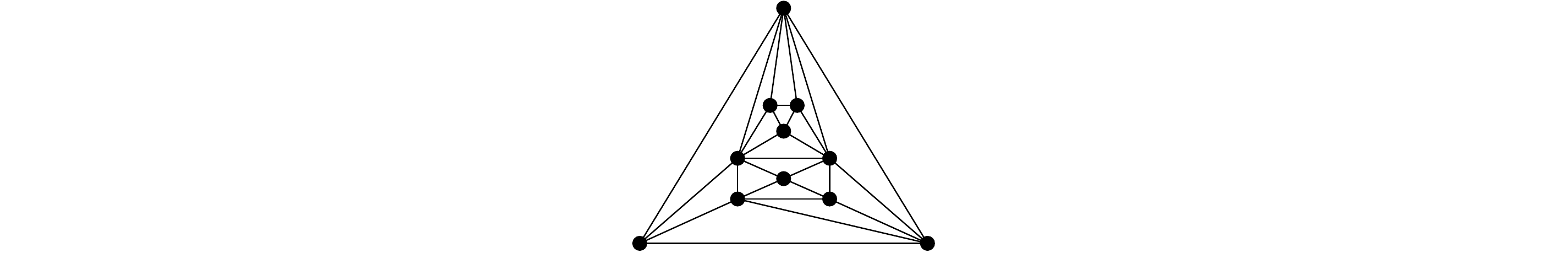}

   \end{center}

6.29 Degree sequence is 44444455677, and it has 10 kinds of different colorings.

\begin{center}
        \includegraphics [width=380pt]{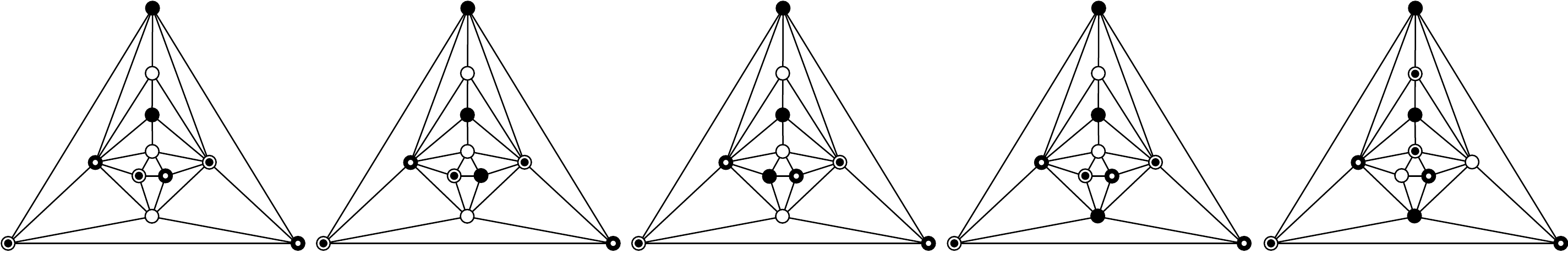}

        \vspace{2mm}
        \includegraphics [width=380pt]{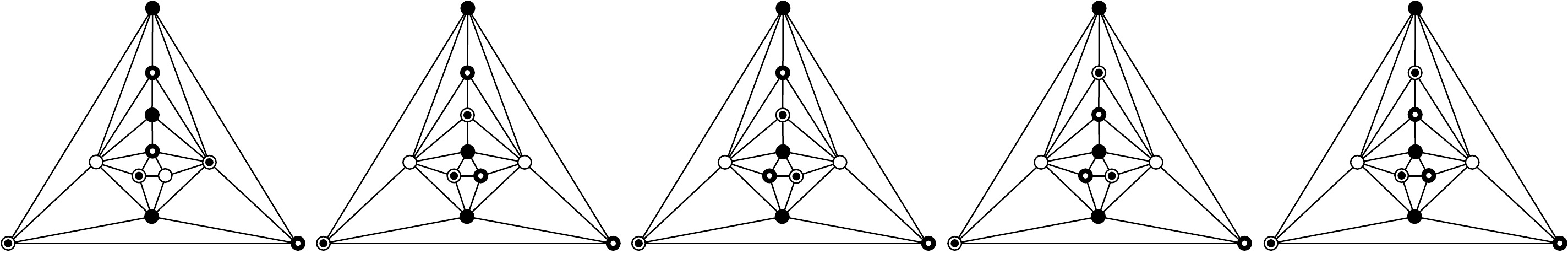}

   \end{center}

6.30 Degree sequence is 44444455677, and it is divisible.

\begin{center}
        \includegraphics [width=380pt]{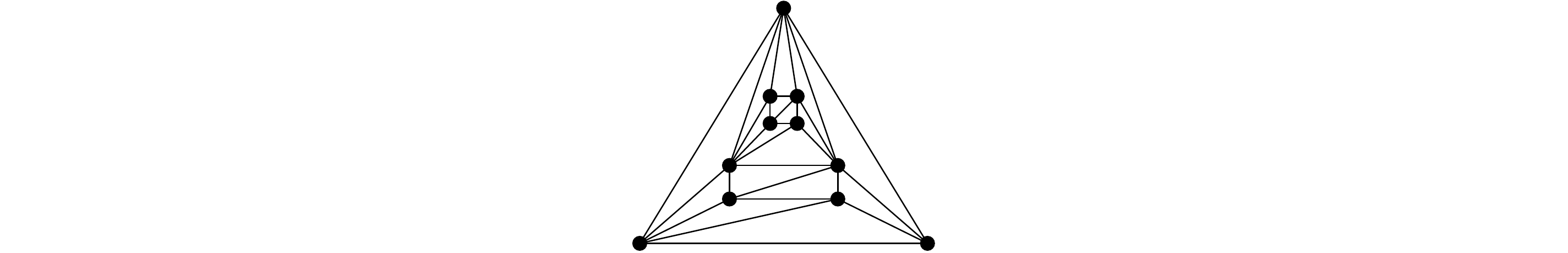}

   \end{center}

6.31 Degree sequence is 44444446677, and it has 41 kinds of different colorings.

\begin{center}
        \includegraphics [width=380pt]{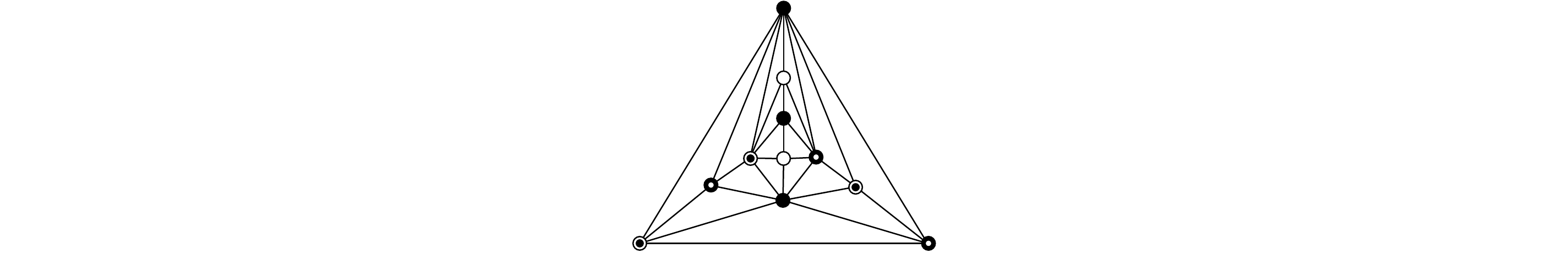}

        \vspace{2mm}
        \includegraphics [width=380pt]{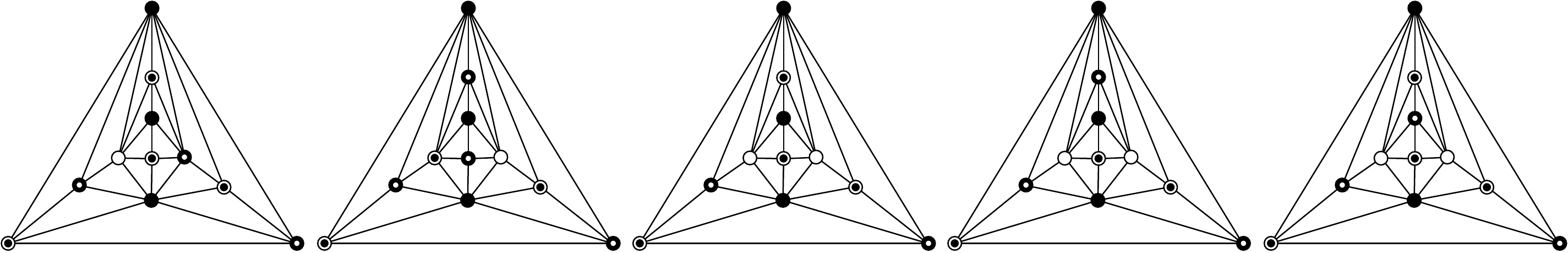}

        \vspace{2mm}
        \includegraphics [width=380pt]{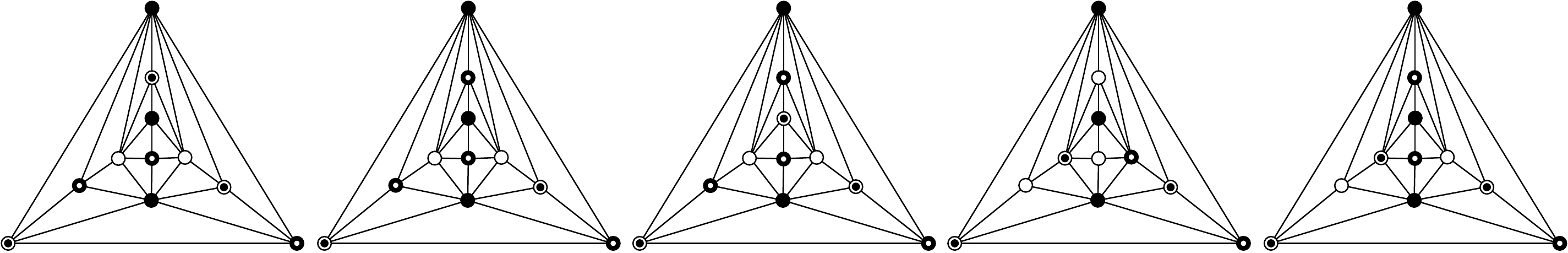}

        \vspace{2mm}
        \includegraphics [width=380pt]{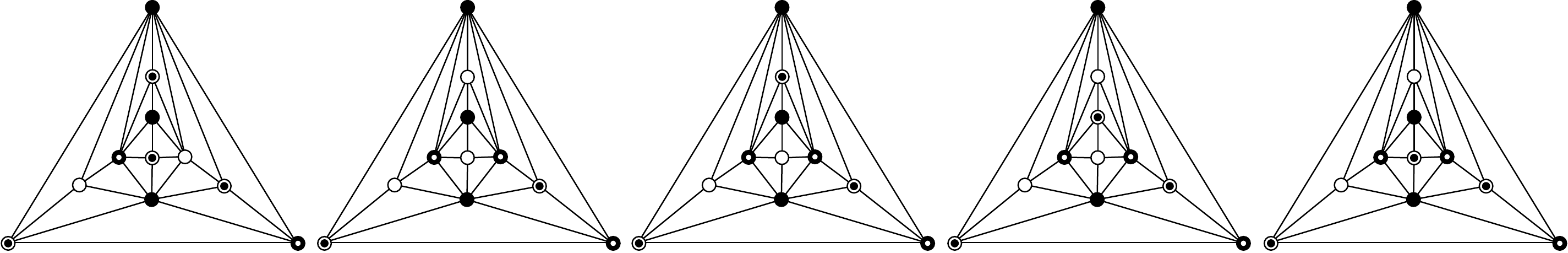}

        \vspace{2mm}
        \includegraphics [width=380pt]{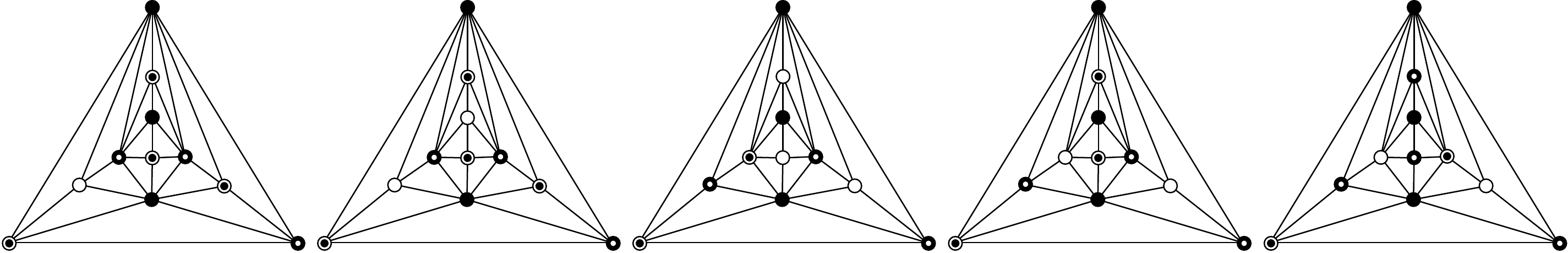}

        \vspace{2mm}
        \includegraphics [width=380pt]{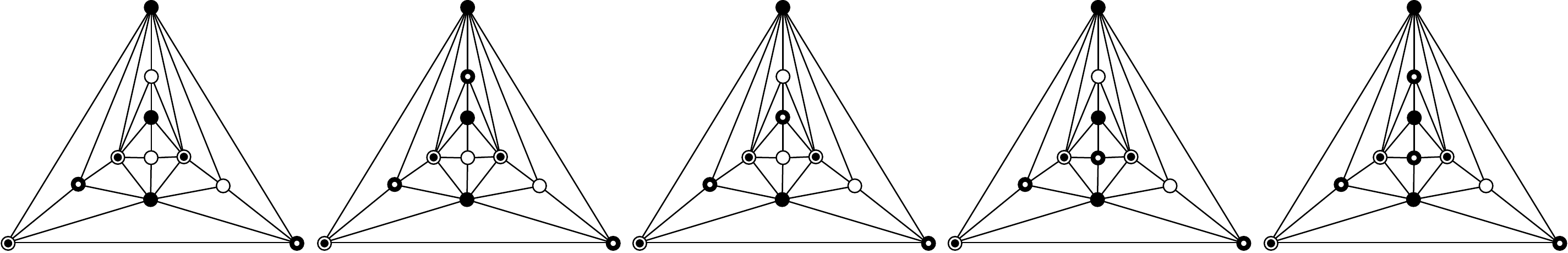}

        \vspace{2mm}
        \includegraphics [width=380pt]{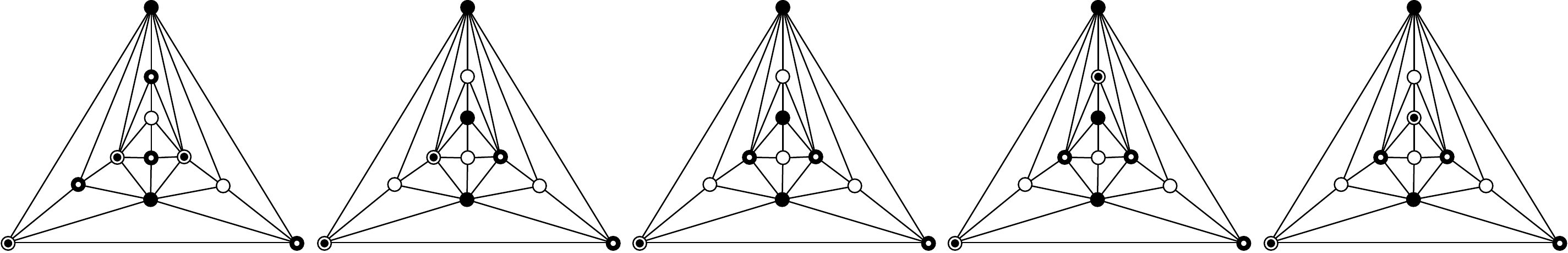}

        \vspace{2mm}
        \includegraphics [width=380pt]{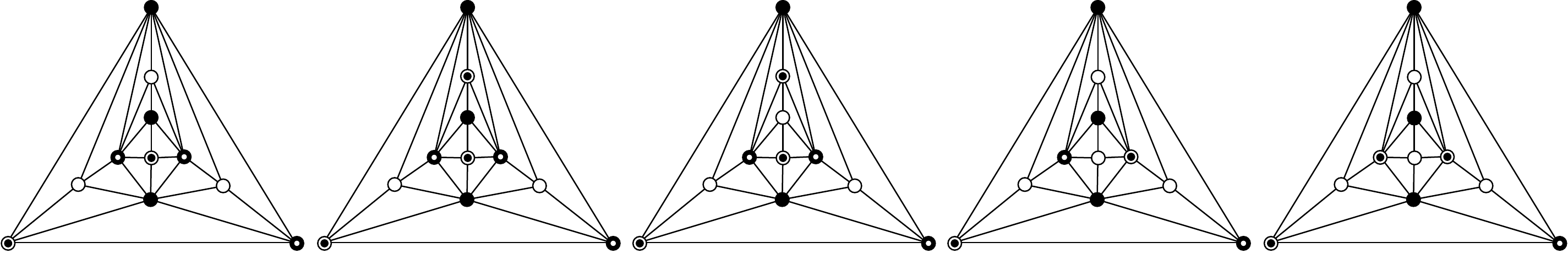}

        \vspace{2mm}
        \includegraphics [width=380pt]{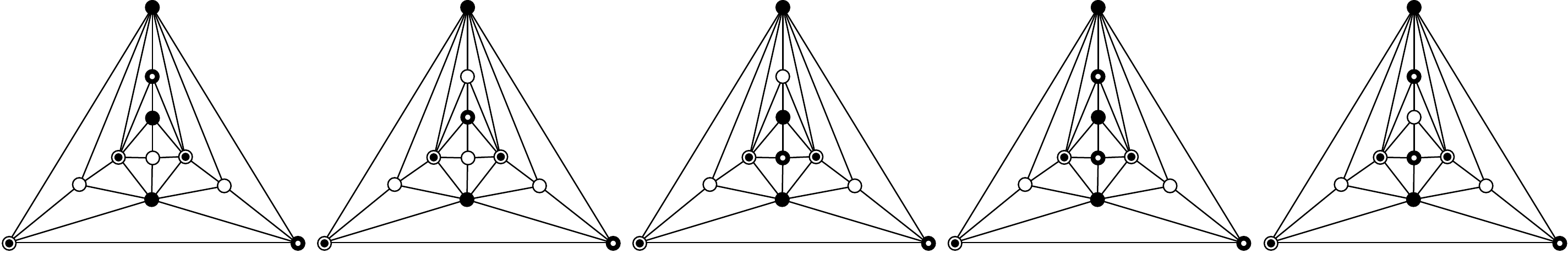}

   \end{center}

6.32 Degree sequence is 44444445588, and it is divisible.

\begin{center}
        \includegraphics [width=380pt]{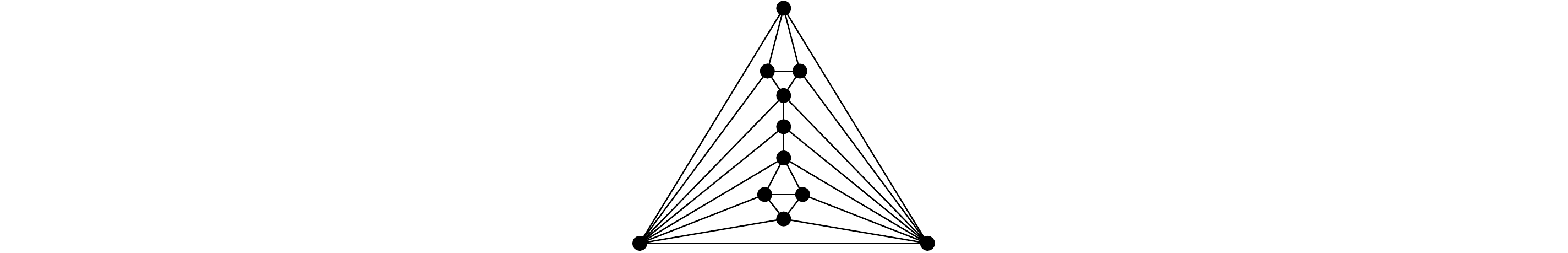}

   \end{center}

6.33 Degree sequence is 44444445588, and it has 25 kinds of different colorings.

\begin{center}
        \includegraphics [width=380pt]{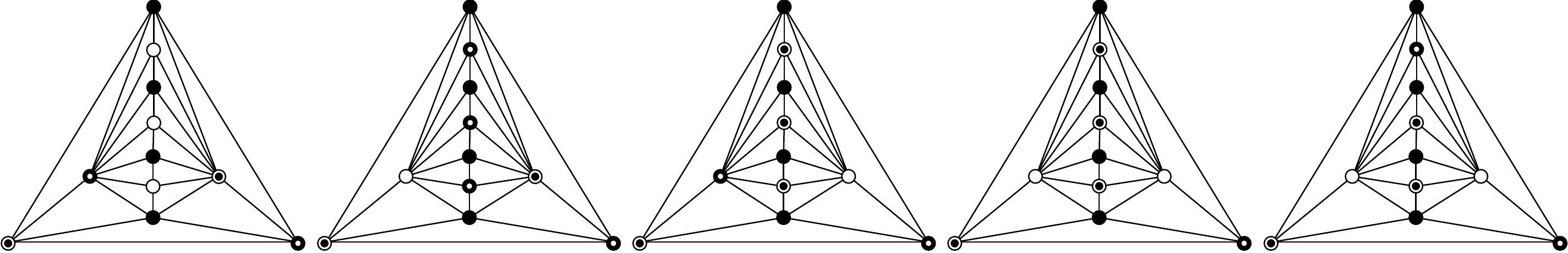}

        \vspace{2mm}
        \includegraphics [width=380pt]{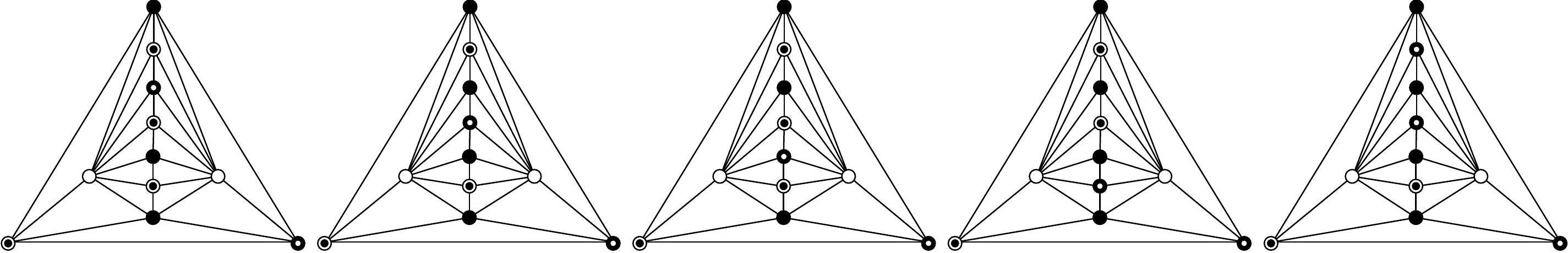}

        \vspace{2mm}
        \includegraphics [width=380pt]{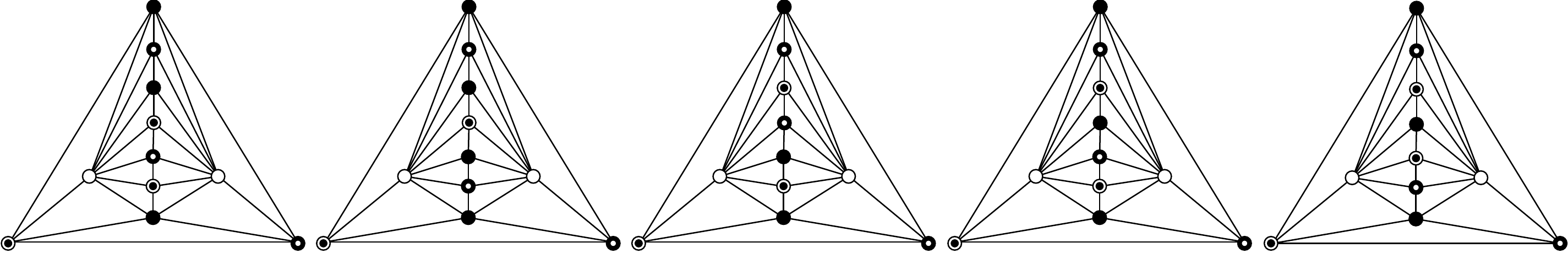}

        \vspace{2mm}
        \includegraphics [width=380pt]{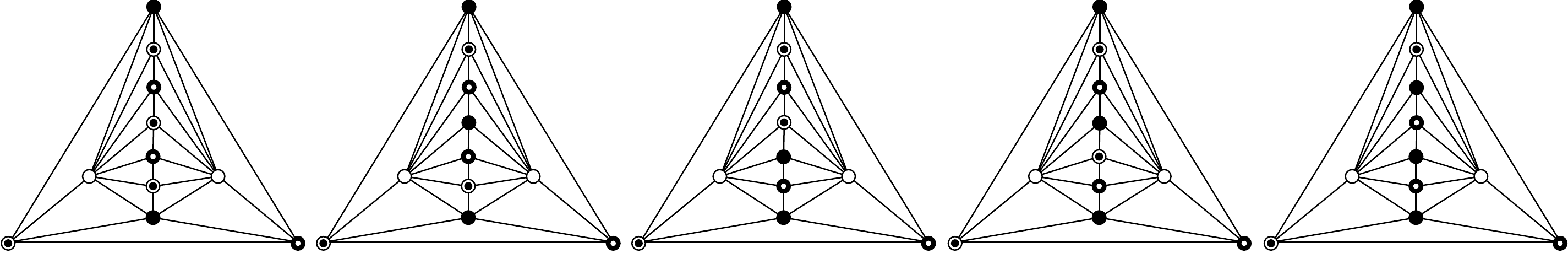}

        \vspace{2mm}
        \includegraphics [width=380pt]{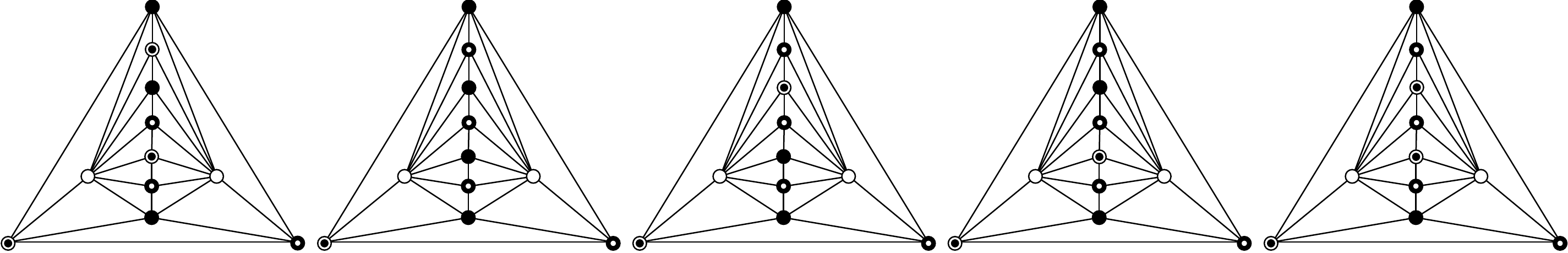}

   \end{center}

6.34 Degree sequence is 44444444499, and it has 85 kinds of different colorings.

\begin{center}
        \includegraphics [width=380pt]{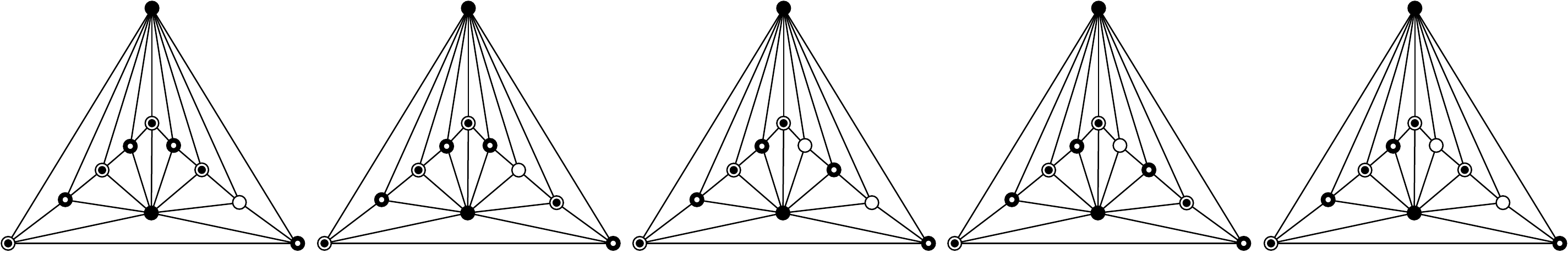}

        \vspace{2mm}
        \includegraphics [width=380pt]{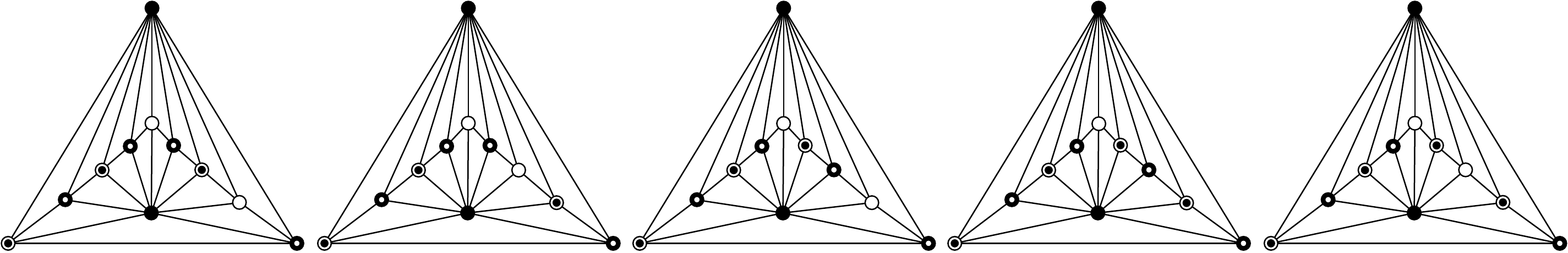}

        \vspace{2mm}
        \includegraphics [width=380pt]{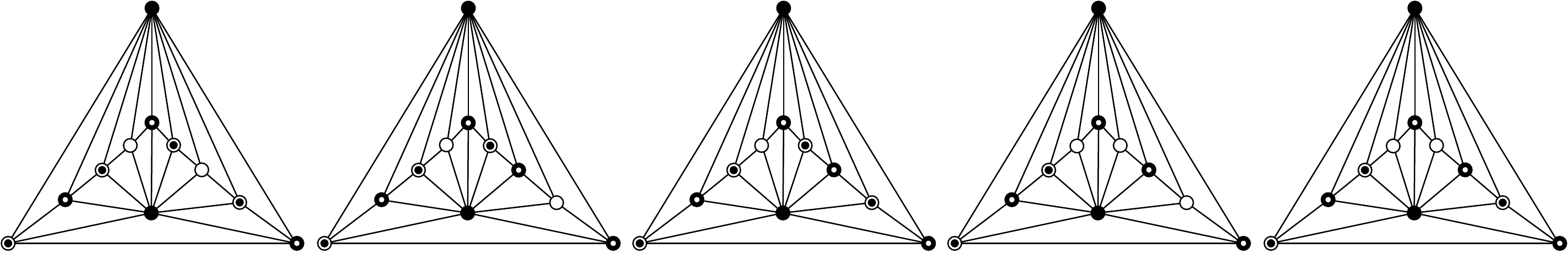}

        \vspace{2mm}
        \includegraphics [width=380pt]{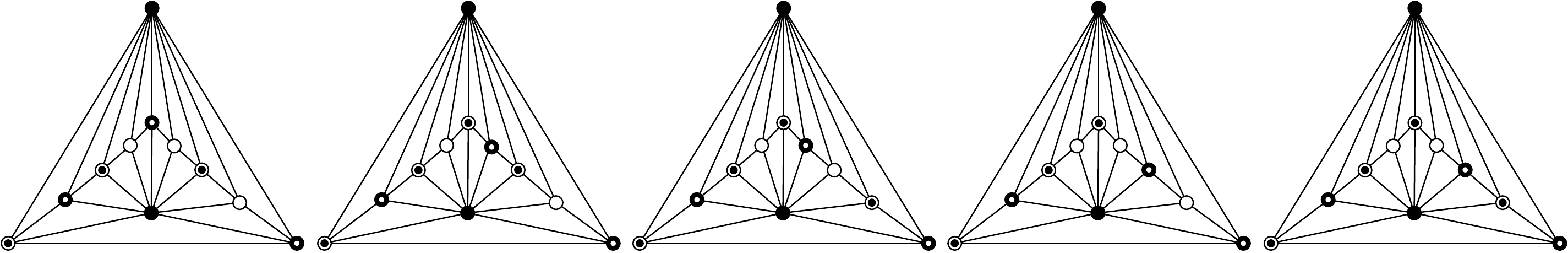}

        \vspace{2mm}
        \includegraphics [width=380pt]{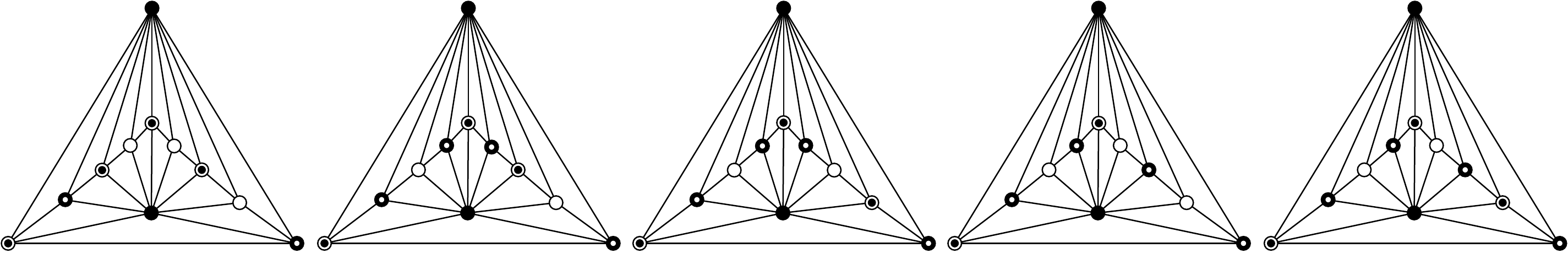}

        \vspace{2mm}
        \includegraphics [width=380pt]{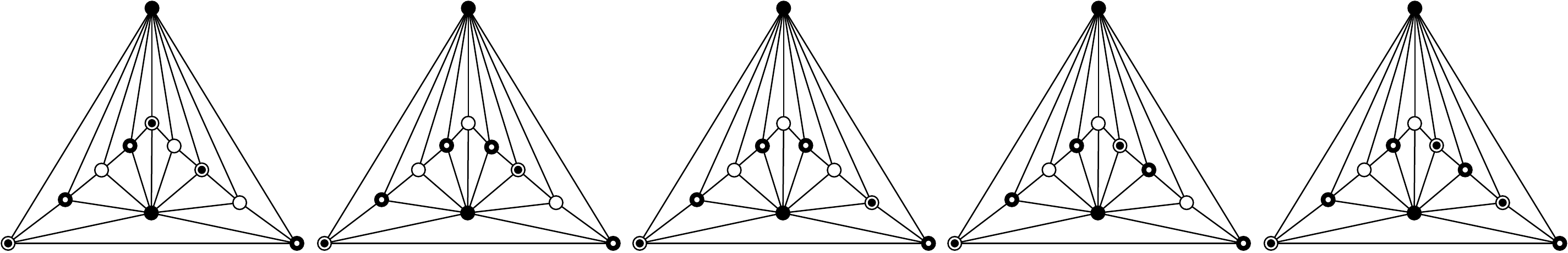}

        \vspace{2mm}
        \includegraphics [width=380pt]{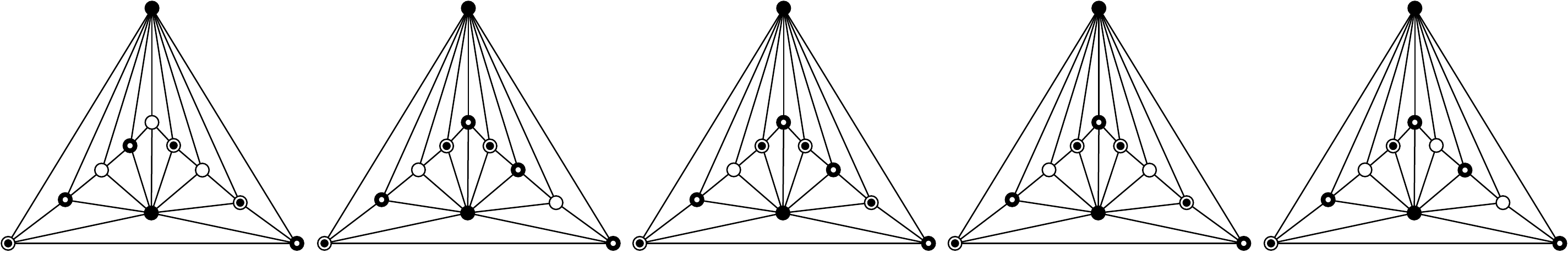}

        \vspace{2mm}
        \includegraphics [width=380pt]{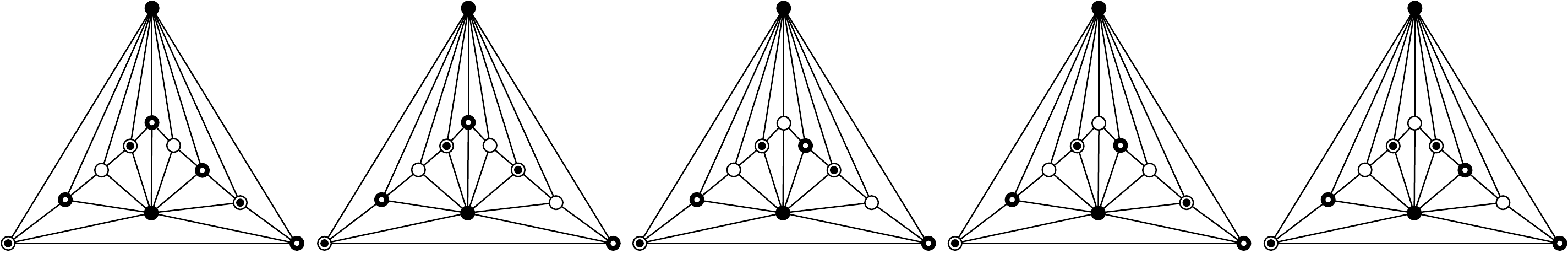}

        \vspace{2mm}
        \includegraphics [width=380pt]{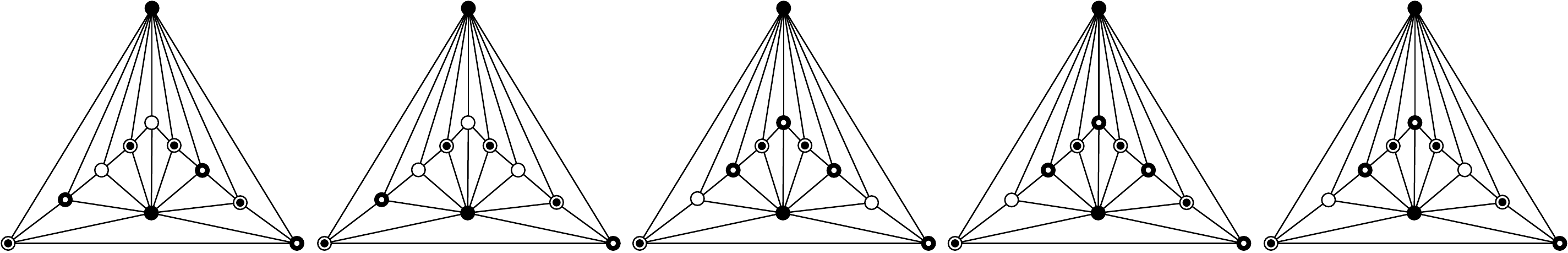}

        \vspace{2mm}
        \includegraphics [width=380pt]{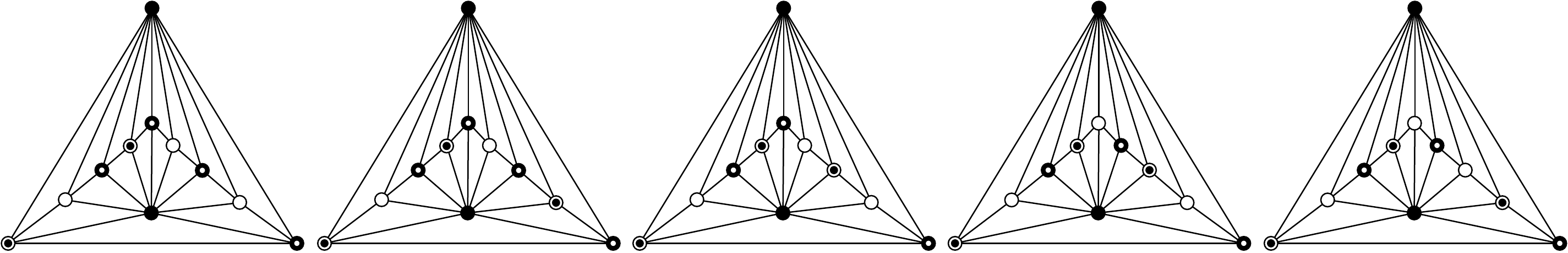}

        \vspace{2mm}
        \includegraphics [width=380pt]{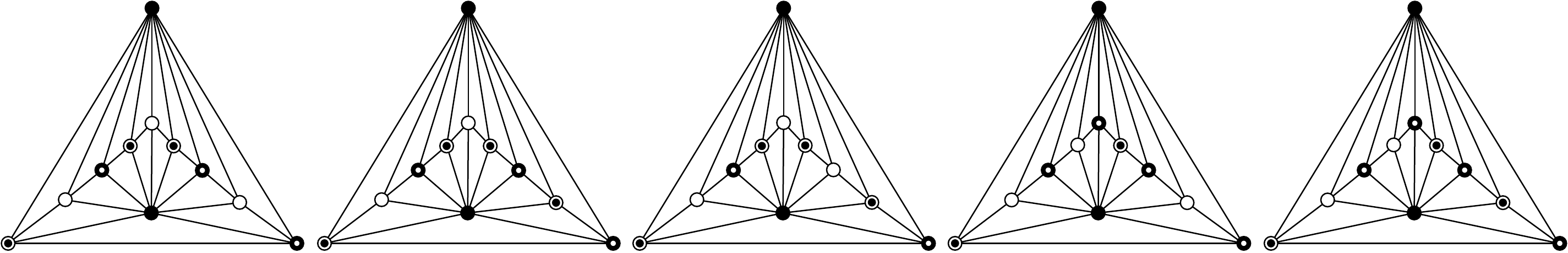}

        \vspace{2mm}
        \includegraphics [width=380pt]{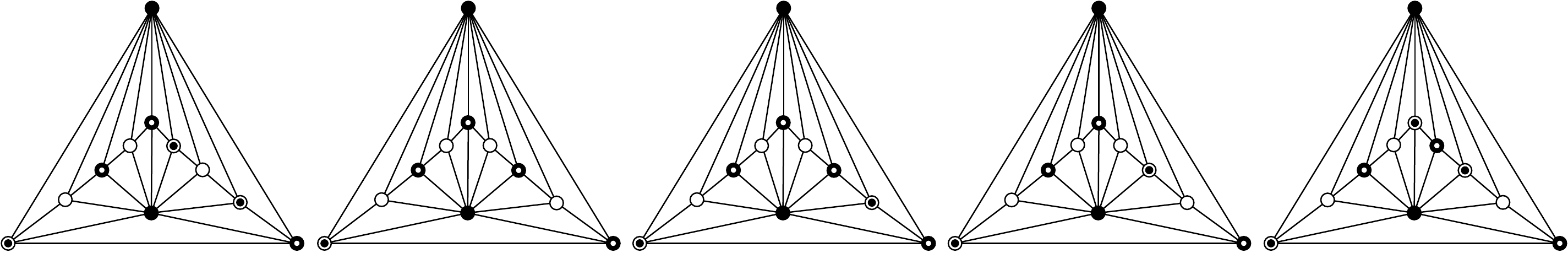}

        \vspace{2mm}
        \includegraphics [width=380pt]{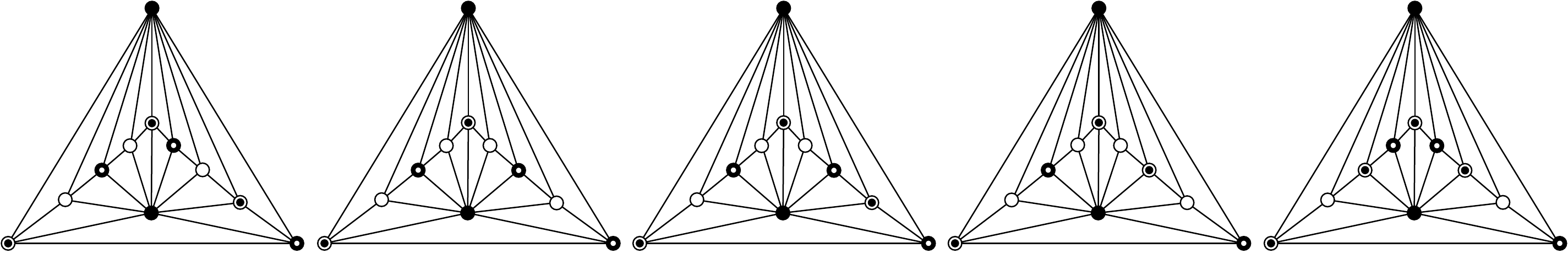}

        \vspace{2mm}
        \includegraphics [width=380pt]{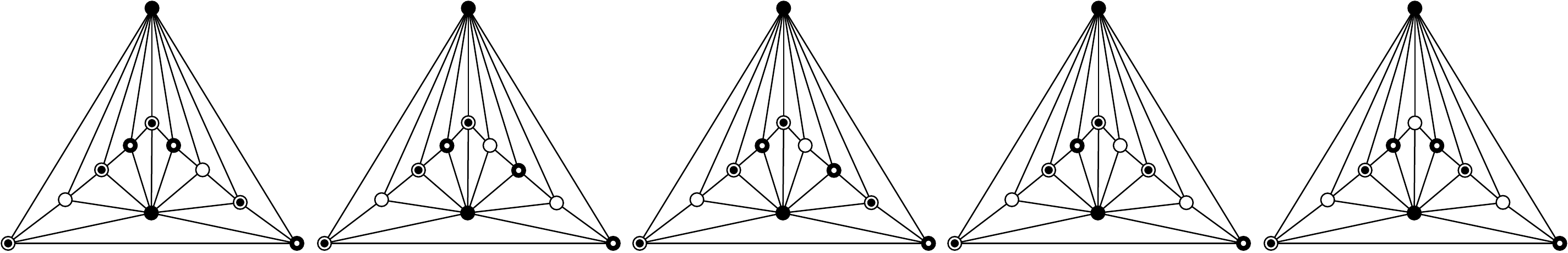}

        \vspace{2mm}
        \includegraphics [width=380pt]{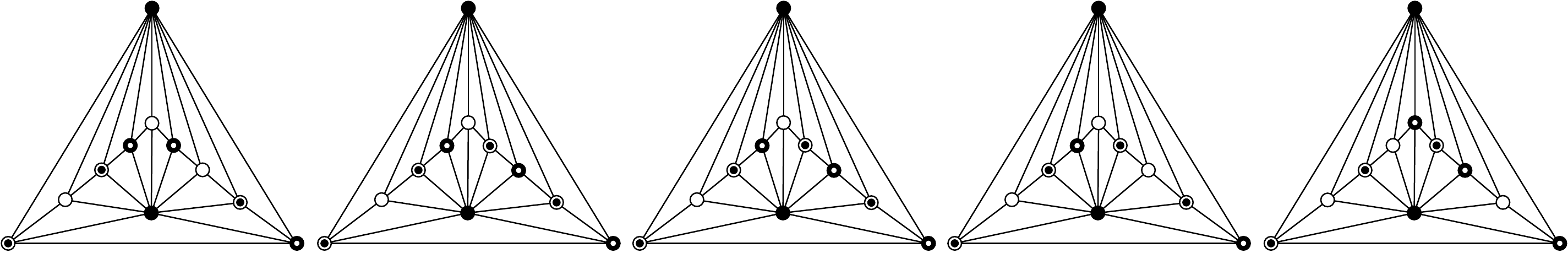}

        \vspace{2mm}
        \includegraphics [width=380pt]{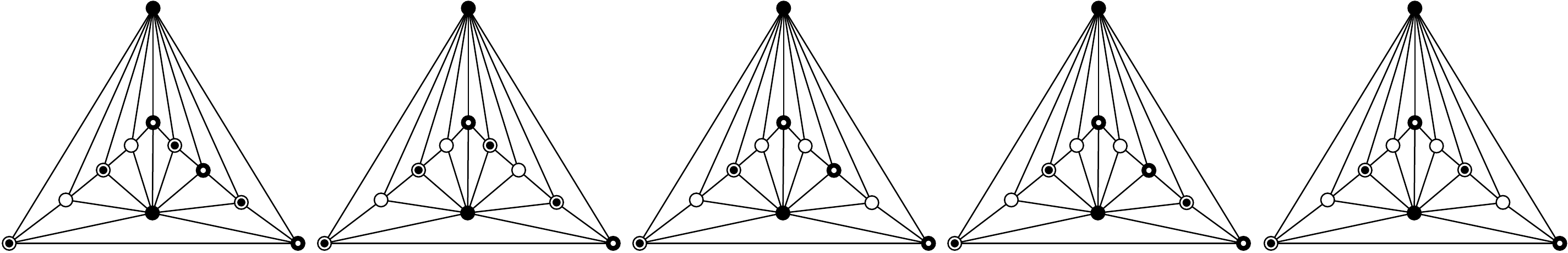}

        \vspace{2mm}
        \includegraphics [width=380pt]{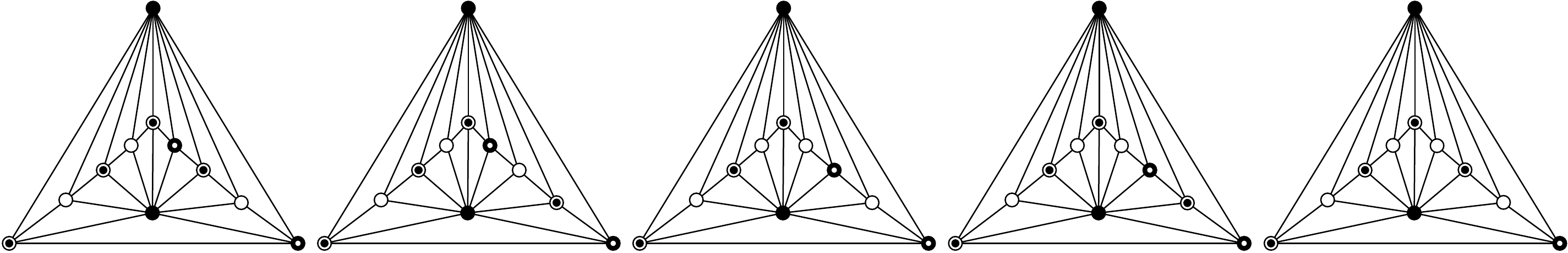}

   \end{center}

\end{document}